\documentclass[12pt,psamsfonts,leqno,oneside,letterpaper]{amsart}
\usepackage[dvips,text={6.5truein,9truein},left=1truein,top=1truein]{geometry}
\usepackage{amssymb,amsmath,amscd,enumerate,
}
\usepackage[pdftex]{graphicx}
\usepackage{url}

\usepackage[colorlinks,linkcolor=blue,citecolor=blue,pdfstartview=FitH]{hyperref}
\input xy
\xyoption{all}
\SelectTips{cm}{12}

\usepackage{color}
\usepackage{mathrsfs}

\newcommand\nonessentialproofreadingnote[1]{\empty}
\newcommand\tinymissingref[1]{\empty}
\newcommand\abstractcomment[1]{\empty}

\usepackage{amsthm}
\theoremstyle{definition}
\newtheorem{para}{}[section]

\newtheorem{remark}[para]{Remark}
\newtheorem{remarks}[para]{Remarks}
\newtheorem{notation}[para]{Notation}

\newtheorem{convention}[para]{Convention}
\newtheorem{definition}[para]{Definition}
\newtheorem{definitions}[para]{Definitions}
\newtheorem{definitionnotation}[para]{Definition and Notation}
\newtheorem{remarksnotation}[para]{Remarks and Notation}

\newtheorem{remarksdefinitions}[para]{Remarks and Definitions}
\newtheorem{remarknotation}[para]{Remark and Notation}
\newtheorem{remarksdefinition}[para]{Remarks and Definition}
\newtheorem{notationremark}[para]{Notation and Remark}
\newtheorem{notationremarks}[para]{Notation and Remarks}
\newtheorem{definitionsnotation}[para]{Definitions and Notation}

\newtheorem{reviewdefinition}[para]{Review and Definition}
\newtheorem{definitionremark}[para]{Definition and Remark}

\newtheorem{definitionsremarks}[para]{Definitions and Remarks}
\newcommand\Alternatives{\begin{enumerate}[(i)]}
\newcommand\EndAlternatives{\end{enumerate}}
\newcommand\Conditions{\begin{enumerate}[(1)]}
\newcommand\EndConditions{\end{enumerate}}

\theoremstyle{plain}
\newtheorem{theorem}[para]{Theorem}
\newtheorem{lemma}[para]{Lemma}
\newtheorem{proposition}[para]{Proposition}
\newtheorem*{propA}{Proposition A}
\newtheorem*{propB}{Proposition B}
\newtheorem*{propC}{Proposition C}

\newtheorem{corollary}[para]{Corollary}
\newtheorem{conjecture}[para]{Conjecture}

\newtheorem{claim}[equation]{}
\numberwithin{equation}{para}
\numberwithin{figure}{section}

\newcommand\Number{\begin{para}}
\newcommand\EndNumber{\end{para}}
\newcommand\Definition{\begin{definition}}
\newcommand\EndDefinition{\end{definition}}
\newcommand\Definitions{\begin{definitions}}
\newcommand\DefinitionsNotation{\begin{definitionsnotation}}
\newcommand\EndDefinitionsRemarks{\end{definitionsremarks}}
\newcommand\DefinitionsRemarks{\begin{definitionsremarks}}
\newcommand\RemarksNotation{\begin{remarksnotation}}
\newcommand\EndRemarksNotation{\end{remarksnotation}}
\newcommand\RemarkNotation{\begin{remarknotation}}
\newcommand\EndRemarkNotation{\end{remarknotation}}
\newcommand\NotationRemark{\begin{notationremark}}
\newcommand\NotationRemarks{\begin{notationremarks}}
\newcommand\EndNotationRemark{\end{notationremark}}
\newcommand\EndNotationRemarks{\end{notationremarks}}
\newcommand\EndDefinitionsNotation{\end{definitionsnotation}}
\newcommand\DefinitionNotation{\begin{definitionnotation}}
\newcommand\EndDefinitionNotation{\end{definitionnotation}}
\newcommand\ReviewDefinition{\begin{reviewdefinition}}
\newcommand\EndReviewDefinition{\end{reviewdefinition}}
\newcommand\EndDefinitions{\end{definitions}}
\newcommand\Theorem{\begin{theorem}}
\newcommand\EndTheorem{\end{theorem}}
\newcommand\Conjecture{\begin{conjecture}}
\newcommand\EndConjecture{\end{conjecture}}
\newcommand\Remark{\begin{remark}}
\newcommand\EndRemark{\end{remark}}
\newcommand\Remarks{\begin{remarks}}
\newcommand\EndRemarks{\end{remarks}}
\newcommand\Convention{\begin{convention}}
\newcommand\EndConvention{\end{convention}}
\newcommand\Notation{\begin{notation}}
\newcommand\EndNotation{\end{notation}}
\newcommand\Lemma{\begin{lemma}}
\newcommand\EndLemma{\end{lemma}}
\newcommand\Proposition{\begin{proposition}}
\newcommand\EndProposition{\end{proposition}}
\newcommand\Corollary{\begin{corollary}}
\newcommand\EndCorollary{\end{corollary}}
\newcommand\Claim{\begin{claim}}
\newcommand\EndClaim{\end{claim}}
\newcommand\Proof{\begin{proof}}
\newcommand\EndProof{\end{proof}}
\newcommand\Equation{\begin{equation}}
\newcommand\EndEquation{\end{equation}}

\newcommand\Bullets{\begin{itemize}}
\newcommand\EndBullets{\end{itemize}}

\newcommand\pl{_{\rm PL}}
\newcommand\PL{^{\rm PL}}
\newcommand\smooth{_{\rm smooth}}
\newcommand\oldOmega{{\mathfrak M}} 

\newcommand\oldGamma{{\mathfrak G}}
\newcommand\oldDelta{{\mathfrak D}}
\newcommand\oldTheta{{\mathfrak T}}
\newcommand\oldLambda{{\mathfrak L}}   
\newcommand\otheroldLambda{{\mathfrak O}}   
\newcommand\oldXi{\mathfrak X} 
\newcommand\oldPi{\mathfrak P} 
\newcommand\oldSigma{\mathfrak S} 
\newcommand\oldUpsilon{\mathfrak W} 
\newcommand\oldPhi{\mathfrak F} 
\newcommand\oldPsi{\mathfrak N}
\newcommand\frakA{\mathfrak A} 
 
\newcommand\frakB{\mathfrak B}
\newcommand\frakD{\mathfrak D}
\newcommand\frakC{\mathfrak C}
\newcommand\frakE{\mathfrak E}
\newcommand\frakI{\mathfrak I}
\newcommand\frakJ{\mathfrak J}
\newcommand\frakK{\mathfrak K}

\newcommand\frakZ{\mathfrak Z}
\newcommand\frakU{\mathfrak U}
\newcommand\tfrakU{{\widetilde{\mathfrak U}}}
\newcommand\tfrakV{{\widetilde{\mathfrak V}}}
\newcommand\tfrakW{{\widetilde{\mathfrak W}}}

\newcommand\tfrakA{{\widetilde{\mathfrak A}}}
\newcommand\tfrakE{{\widetilde{\mathfrak E}}}

\newcommand\frakW{{{\mathfrak W}}}
\newcommand\tfrakZ{{\widetilde{\mathfrak Z}}}
\newcommand\frakV{\mathfrak V}
\newcommand\frakP{\mathfrak P}
\newcommand\frakR{\mathfrak R}
\newcommand\tfrakB{\widetilde{\mathfrak B}}

\newcommand\toldOmega{\widetilde{\oldOmega}} 

\newcommand\toldGamma{\widetilde{\mathfrak G}}

\newcommand\toldTheta{\widetilde{\mathfrak T}}
\newcommand\toldLambda{\widetilde{\mathfrak L}}   
\newcommand\toldXi{\widetilde{\mathfrak X}}
\newcommand\toldPi{\widetilde{\mathfrak P}}

\newcommand\toldPsi{\widetilde{\oldPsi}}

\newcommand\compnum{{\frak c}}
\renewcommand\epsilon{\varepsilon}

\newcommand\spair{S-pair}
\newcommand\Ssuborbifold{S-suborbifold}
\newcommand\Asuborbifold{A-suborbifold}
\newcommand\pagelike{page-like}
\newcommand\bindinglike{binding-like}
\newcommand\untwisted{untwisted}
\newcommand\twisted{twisted}

\newcommand\simple{atoral}

\newcommand\torifold{solid toric orbifold}

\newcommand\book{\mathop{\rm book}}

\newcommand\kish{\mathop{\rm kish}}

\newcommand\discup{\mathbin{\rotatebox[origin=c]{90}{$\vDash$}}}
\newcommand\RR{{\bf R}}
\newcommand\fraks{{\Sigma}}

\newcommand\hatC{{\widehat C}}
\newcommand\hatX{{\widehat X}}

\newcommand\plusX{{ X^+}}

\newcommand\plusN{{ N^+}}

\newcommand\chibar{\overline\chi}

\newcommand\inter{\mathop{\rm int}}

\newcommand\tM{\widetilde M}

\newcommand\tT{{\widetilde T}}

\newcommand\tW{{\widetilde W}}
\newcommand\tR{{\widetilde R}}

\newcommand\tq{\widetilde q}

\newcommand\tcals{\widetilde \cals}

\newcommand\overp{\overline{p}}

\newcommand\barH{\overline{H}}

\newcommand\tx{{\widetilde x}}

\newcommand\calc{{\mathcal C}}

\newcommand\cala{{\mathscr A}}

\newcommand\calp{{\mathcal P}}

\newcommand\ch{\calc_{\rm h}}
\newcommand\cald{{\mathcal D}}
\newcommand\frakH{{\mathfrak H}}
\newcommand\frakf{{\mathfrak f}}
\newcommand\frakh{{\mathfrak h}}

\newcommand\frakY{{\mathfrak Y}}

\newcommand\calf{{\mathcal F}}

\newcommand\calK{{\mathcal K}}

\newcommand\frakQ{{\mathfrak Q}}

\newcommand\calq{{\mathcal Q}}
\newcommand\calt{{\mathcal T}}
\newcommand\calb{{\mathcal B}}

\newcommand\cale{{\mathcal E}}

\newcommand\calx{{\mathcal X}}
\newcommand\caly{{\mathcal Y}}

\newcommand\calv{{\mathcal V}}

\newcommand\grock{\tau}

\newcommand\smock{\zeta}

\newcommand\area{\mathop{{\rm area}}}

\newcommand\cut{^\dagger_}

\newcommand\ZZ{{\mathbb Z}}

\newcommand\NN{{\mathbb N}}
\newcommand\CC{{\mathbb C}}
\newcommand\QQ{{\mathbb Q }}
\newcommand\HH{{\mathbb H}}

\newcommand\calk{{\mathcal K}}
\newcommand\cals{{\mathscr S}}

\newcommand\calr{{\mathscr R}}

\newcommand\wt{\mathop{{\rm wt}}}

\newcommand\obd{\mathop{\omega}}

\newcommand\silv{\mathop{\rm silv}}

\newcommand\vol{\mathop{\rm vol}}
\newcommand\volG{\mathop{\rm vol_G}}
\newcommand\volorb{\mathop{\rm vol_{ASTA}}}
\newcommand\voct{V_{\rm oct}}
\newcommand\vtet{V_{\rm tet}}

\newcommand\Fix{\mathop{{\rm Fix}}}
\newcommand\card{\mathop{{\rm card}}}
\newcommand\rank{\mathop{{\rm rank}}}
\newcommand\image{\mathop{{\rm Im}}}

\newcommand\Fr{\mathop{\rm Fr}}

\DeclareFontFamily{U}{rcjhbltx}{}
\DeclareFontShape{U}{rcjhbltx}{m}{n}{<->rcjhbltx}{}
\DeclareSymbolFont{hebrewletters}{U}{rcjhbltx}{m}{n}

\let\aleph\relax\let\beth\relax
\let\gimel\relax\let\daleth\relax

\DeclareMathSymbol{\aleph}{\mathord}{hebrewletters}{39}
\DeclareMathSymbol{\beth}{\mathord}{hebrewletters}{98}
\DeclareMathSymbol{\gimel}{\mathord}{hebrewletters}{103}
\DeclareMathSymbol{\daleth}{\mathord}{hebrewletters}{100}

\DeclareMathSymbol{\lamed}{\mathord}{hebrewletters}{108}
\DeclareMathSymbol{\mem}{\mathord}{hebrewletters}{109}
\DeclareMathSymbol{\ayin}{\mathord}{hebrewletters}{96}
\DeclareMathSymbol{\tsadi}{\mathord}{hebrewletters}{118}
\DeclareMathSymbol{\qof}{\mathord}{hebrewletters}{114}
\DeclareMathSymbol{\shin}{\mathord}{hebrewletters}{152}

\begin{document}

\author{Peter B. Shalen}
\address{Department of Mathematics, Statistics, and Computer Science (M/C 249)\\  University of Illinois at Chicago\\
 851 S. Morgan St.\\
  Chicago, IL 60607-7045} \email{shalen@math.uic.edu}
\thanks{Partially supported by NSF grant DMS-1207720}


\abstractcomment{
The things below the stars are most likely non-problems, but I'm not quite ready to remove them.

$$\qquad*\qquad*\qquad*\qquad*\qquad*\qquad*\qquad*\qquad*$$

It has just occurred to me (7/14/17) that with my current def. of strong \simple ity, I could probably prove that a covering of a strongly \simple\ orbifold is strongly \simple\ by using Meeks-Simon-Yau for the manifold case. Something tp think about...

I want to check that I have not assumed anywhere that when $\oldXi$ is a $2$-orbifold and $x$ is a point of a $1$-dim. component of $\fraks_\oldXi$, the group associated to $x$ is generated by a reflection alone. This may once have been an issue in  the proof of Lemma \ref{affect}, which in the earlier draft kplus.tex contained a note about an issue raised by Shawn, but the second assertion, the part affected by this, is gone now. So...

I'm having doubts again about whether I've argued correctly about
$I$-fibrations. The vertical boundary of an $i$-fibered orbifold can
have components that are annular orbifolds with underlying surface a
disk and two singular points of order $2$. Is everything I've said
consistent with that example? The issue has come up again in
connection with Proposition \ref{when vertical}.The issue is understanding the case where
the annular orbifolds in the hypothesis of that prop. have singularities.

Here's a more specific point: if the base of an $I$-fibration is a
$2$-manifold, have I assumed that the total space is a $3$-manifold? I
don't think it's true. But I may have argued that one could reduce
certain facts to the manifold case by removing singularities from the
base. When one does this, the base becomes a $2$-manifold, but it's a
mistake to assume that the total space becomes a $3$-manifold.

What I believe at the moment is that if the base is a manifold and the total space is orientable, then the total space is a manifold. This illustrates how confused I get. In general, I somehow ought to
make sure I have never overlooked the case where the total space is an
orientable orbifold but the base space is a non-orientable orbifold,
which is a major source of complication. I've certainly been careful
about it in some places, but I wonder whether I have ignored it in others.

Every time I look at a missingref, I find mistakes in the text that
are independent of the missingref. Stupid things like using a letter
that denotes one object when I mean to refer to another object. I'm
not sure how I can catch all these things.

 I have been concerned about confusion
  between the calligraphic and fraktur fonts, specifically in the proof of Lemma \ref{pre-modification}. I tried
  mathscr in place of mathcal, and it's not at all clear that it's an improvement.

}

\title{Volume and Homology for Hyperbolic $3$-Orbifolds, I}

\begin{abstract}
 Let $\oldOmega$ be a
closed, orientable, hyperbolic 3-orbifold whose singular set is a
link, and such that $\pi_1(\oldOmega)$
contains no hyperbolic triangle group. We show that if the underlying manifold $|\oldOmega|$
is irreducible, and $|\toldOmega|$ is irreducible for every
two-sheeted (orbifold) covering $\toldOmega|$ of $\oldOmega$, and if $\vol \oldOmega\le1.72$, then
$\dim H_1(\oldOmega;\ZZ_2)\le
15$. Furthermore, if
$\vol \oldOmega\le1.22$ then
$\dim H_1(\oldOmega;\ZZ_2)\le
11$, and if $\vol \oldOmega\le0.61$ then
$\dim H_1(\oldOmega;\ZZ_2)\le
7$.  
The proof is an application of results that will be used in the sequel
to this paper to obtain qualitatively similar results without the assumption of
irreducibility of $|\oldOmega|$ and $|\toldOmega|$.
\end{abstract}

\maketitle


\section{Introduction}

It is a standard consequence of the Margulis Lemma \cite[Chapter D]{bp} that an upper bound on the volume of a closed hyperbolic $3$-orbifold $\oldOmega$ imposes an upper bound on the rank of $\pi_1(\oldOmega)$, and hence on the dimension of $H_1(\oldOmega;\ZZ_p)$ for any prime $p$. In \cite{rankfour}, \cite {last}, \cite{lastplusone}, and  \cite{fourfree}, for the case of a closed, orientable hyerbolic $3$-manifold $M$, relatively small upper bounds for $\vol M$ are shown to imply explicit bounds for $\dim _{\ZZ_2}H_1(M;\ZZ_2)$ that are close to being sharp. The purpose of this paper and its sequel \cite{second} is to provide explicit bounds
for hyperbolic $3$-orbifolds that are not manifolds and are subject to relatively small upper bounds on volume; like the bounds given in \cite{rankfour}, \cite {last}, \cite{lastplusone}, and  \cite{fourfree}, the bounds given here, when they apply, will improve the naive bounds by several orders of magnitude.

The approach used here to bounding the dimension of the $\ZZ_p$-vector space $H_1(\oldOmega;\ZZ_2)$, where $\oldOmega$ is a closed, orientable hyperbolic $3$-orbifold whose volume is subject to a given bound, is to begin by finding bounds for  $\dim_{\ZZ_2} H_1(|\oldOmega|;\ZZ_2)$ and $\dim_{\ZZ_2} H_1(|\toldOmega|;\ZZ_2)$, where $\toldOmega$ is an arbitrary two-sheeted covering orbifold of $\oldOmega$, and absolute value signs denote the underlying space of an orbifold. (If $\oldOmega$ is an orientable $3$-orbifold then $|\oldOmega|$ is a $3$-manifold.) In the case where the singular set of $\oldOmega$ is a link, these bounds can be parlayed into bounds for $\dim H_1(|\oldOmega|;\ZZ_2)$ by means of the following result, which is proved in the body of this paper as Proposition \ref{boogie-woogie bugle boy}; the proof is an application of the Smith inequalities.

\begin{propA}
Let $\oldOmega$ be a 
closed, orientable, hyperbolic 3-orbifold whose singular set
is a link. Then $\oldOmega$ is covered with degree at most $2$ by some
orbifold $\toldOmega$ such that
$$\dim_{\ZZ_2}H_1(\oldOmega;\ZZ_2)\le1+ \dim_{\ZZ_2}H_1(|{\toldOmega}|;\ZZ_2)+\dim_{\ZZ_2}H_1(|{\oldOmega}|;\ZZ_2).$$
\end{propA}

In the present paper we establish the following result, which provides bounds for $\dim H_1(M;\ZZ_2)$ in certain situations, and is proved in the body of the paper as Proposition \ref{lost corollary}: 

\begin{propB}
Let $\oldOmega$ be a closed,
orientable, hyperbolic $3$-orbifold containing no embedded negative
turnovers (see \ref{wuzza turnover}). 
Suppose that $M:=|\oldOmega|$ is irreducible.
Then the following assertions are true.
\begin{itemize}
\item If the singular set of $\oldOmega$ is a link and $\vol(\oldOmega)\le3.44$ then $\dim H_1(M;\ZZ_2)\le7$. 
\item If 
the singular set of $\oldOmega$ is a link and $\vol(\oldOmega)\le1.22$, then $\dim H_1(M;\ZZ_2)\le3$.
\item If 
the singular set of $\oldOmega$ is a link and  $\vol(\oldOmega)<1.83$, then $\dim H_1(M;\ZZ_2)\le6$.
\item If $\vol(\oldOmega)<0.915$, then $\dim H_1(M;\ZZ_2)\le3$. 
\end{itemize}
\end{propB}

The hypothesis in Proposition B that $M$ is irreducible is not a particularly natural one, and as we shall explain below, the main importance of Proposition B is that it (or the main result used in proving it) gives information about orbifolds whose underlying manifolds may be reducible. However, as a by-product,  Propositions A and B yield the following result, which relates $\vol \oldOmega$ and
$\dim H_1(\oldOmega;\ZZ_2)$ for certain hyperbolic $3$-orbifolds $\oldOmega$, and is proved in the body of the paper as Proposition \ref{orbifirst}:

\begin{propC} Let $\oldOmega$ be a 
closed, orientable, hyperbolic 3-orbifold whose singular set 
is a link, and such that $\pi_1(\oldOmega)$
contains no hyperbolic triangle group. Suppose that $|\oldOmega|$ is irreducible, and that $|\toldOmega|$ is irreducible for every two-sheeted (orbifold) covering $\toldOmega|$ of $\oldOmega$. If $\vol \oldOmega\le1.72$ then
$\dim H_1(\oldOmega;\ZZ_2)\le
15$. Furthermore, if
$\vol \oldOmega\le1.22$ then
$\dim H_1(\oldOmega;\ZZ_2)\le
11$, and if $\vol \oldOmega\le0.61$ then
$\dim H_1(\oldOmega;\ZZ_2)\le
7$. 
\end{propC}

In \cite{second}, we will obtain analogues of Propositions B and C that do not involve the hypothesis that the underlying manifold of $\oldOmega$ or of its two-sheeted coverings is irreducible. These results will involve stronger upper bounds on volume than those assumed in Proposition B or C, and will yield weaker upper bounds on the rank of homology. For example, the analogue of Proposition B will imply that if a closed orientable hyperbolic $3$-orbifold $\oldOmega$ has volume strictly less than $1.83$, has a link as singular set and contains no embedded negative turnover, then $\dim H_1(|\oldOmega|;\ZZ_2)\le19$. Like Proposition B, its analogue in \cite{second} also provides stronger upper bounds on 
 $\dim H_1(|\oldOmega|;\ZZ_2|)$ in the presence of
stronger upper bounds on $\vol \oldOmega$. Also in analogy with Proposition B, the conclusion becomes weaker if the assumption that the singular set is a link is removed. The analogue of Proposition C proved in \cite{second} will imply that if $\oldOmega$ has volume strictly less than $0.915$, has a link as singular set and contains no embedded negative turnover, then $\dim H_1(\oldOmega;\ZZ_2)\le30$. As in Proposition C, stronger upper bounds on $\vol \oldOmega$ give stronger upper bounds on $\dim H_1(\oldOmega;\ZZ_2)$. Like Proposition C, the analogous result in \cite{second} is deduced from the result giving a bound on $\dim H_1(\oldOmega;\ZZ_2)$ by applying Proposition A, and therefore does not give information when the singular set is not a link.

It should be pointed out that these quantitative results are different from the versions tentatively stated in the expository article \cite{arithmetic}. This confirms the statement made in \cite{arithmetic} that ``the exact bound may `be slightly different when the paper is finished.''

The author's personal interest in the problems addressed in this paper and in \cite{second} arises in large part from their connection with the theory of arithmetic groups. It was established in \cite{borel}, and is explained in \cite{arithmetic}, that the most difficult step in listing the arithmetic subgroups of at most a given covolume in ${\rm PGL}(2,\CC)$ is to bound the dimension of the first homology of such groups with mod $2$ coefficients; this is a special case of the problem of bounding $\dim H_1(\oldOmega;\ZZ_2)$, where $\oldOmega$ is an orientable hyperbolic $3$-orbifold. For this application, it turns out to be possible to restrict attention to the case where $\oldOmega$ is closed and $\pi_1(\oldOmega)$ contains no triangle groups. However, the assumption in Proposition C above that the singular set is a link is too restrictive for this application. The author hopes to obtain results in the future that do not depend on this hypothesis.



The starting point for the proof of Proposition B is the well-known fact (essentially contained in Proposition \ref{hepcat} of this paper) that if the  irreducible manifold $M=|\oldOmega|$ itself admits a hyperbolic structure, then  $\vol M\le\vol\oldOmega$. In this case, the conclusion of Proposition B can be deduced from the results of 
\cite{rankfour}, \cite {last}, and  \cite{fourfree}.
If  $M$ itself admits a hyperbolic structure, then it follows from Perelman's geometrization theorem \cite{bbmbp}, \cite{Cao-Zhu}, \cite{kleiner-lott}, \cite{Morgan-Tian} that $M$ is either a small Seifert fibered space, in which case $\dim H_1(M)$ can be shown to be at most $3$, or $M$ contains an incompressible torus $T$. One can choose $T$ within its isotopy class so that $T=|\oldTheta|$, where $\oldTheta$ is some incompressible $2$-suborbifold of $\oldOmega$. The challenge then becomes to prove that if $\oldOmega$ contains an incompressible $2$-orbifold whose underlying surface is a torus, then certain upper bounds for $\vol\oldOmega$ imply certain upper bounds for $\dim H_1(M;\ZZ_2)$; or, contrapositively, that certain lower bounds for $\dim H_1(M;\ZZ_2)$, together with the existence of such a $2$-suborbifold, imply certain upper bounds for  $\vol\oldOmega$.

For the case of a hyperbolic $3$-manifold, it was shown in \cite{ast} that certain extrinsic topological invariants of an incompressible surface give lower bounds on the hyperbolic volume of the $3$-manifold. This was exploited in \cite{hodad}, \cite{after-hodad}, \cite{last}, and \cite{lastplusone}, to relate volume to the dimension of the mod-$2$ homology. To meet the challenge described above, it is necessary to adapt the results of \cite{ast} to the context of orbifolds. This was first done in \cite{atkinson}. In this paper we need a more systematic version of the orbifold analogue of the results of \cite{ast}. The latter results are stated in terms of the relative characteristic submanifold of the manifold obtained by splitting a hyperbolic $3$-manifold along an incompressible surface. We therefore need a systematic version of the theory of the characteristic suborbifold of a $3$-orbifold. After preliminary material in Sections \ref{prelim section} and \ref{fibration section},
the characteristic orbifold theory is developed in Section \ref{characteristic section} of the present paper, taking the main result of \cite{bonahon-siebenmann} as a starting point. The orbifold analogue of the results of \cite{ast} is then carried out in Section \ref{darts section}. 

After the largely foundational material in Sections \ref{prelim section}---\ref{darts section}, we proceed to material directly related to the proof of Theorem B, which occupies Sections \ref{tori section}---\ref{irr-M section}. The actual arguments needed for the proof involve a slight refinement of the approach hinted at above. If $\oldOmega$ is a closed, orientable, hyperbolic $3$-orbifold, we use a certain $2$-suborbifold $\oldTheta$ of $\oldOmega$, whose components are incompressible, and such that the components of $|\oldTheta|$ are incompressible tori in $|\oldOmega|$. (If $|\oldOmega|$ is a hyperbolic manifold or a small Seifert fibered space we take $\oldTheta=\emptyset$.) We then use the formalism of Section \ref{darts section} to bound $\vol\oldOmega$ below by a certain invariant of the orbifold $\oldOmega'$ obtained by splitting $\oldOmega$ along $\oldTheta$; this invariant is additive over components. For an arbitrary component $\oldPsi$ of $\oldOmega'$, the value of the invariant on $\oldPsi$ is bounded below by a constant multiple of the absolute value of the (orbifold) Euler characteristic of the complement  in $\oldPsi$ of the relative characteristic suborbifold of $\oldPsi$. For the case of a component $\oldPsi$ whose interior admits a hyperbolic structure, an alternative approach to estimating the invariant is to bound it below by $\vol(\inter\oldPsi)$; these two ways of producing lower bounds turn out both to be useful, and to complement each other. These lower bounds for the invariant are in turn related to mod-$2$ homology
The details of these arguments
turn out to be quite involved, and Sections \ref{tori section}---\ref{irr-M section} are the heart of the paper. A major source of difficulty is that if $\oldTheta$ is an {\it arbitrary} $2$-suborbifold  of $\oldOmega$ whose components are incompressible, and such that the components of $|\oldTheta|$ are incompressible tori in $|\oldOmega|$, then $\oldTheta$ does not necessarily yield a strictly positive lower bound for volume by the arguments that we have described, and complicated combinatorial arguments are needed to replace a given such suborbifold by one that is better adapted to the purpose.

The proof of Propositions A from the Smith inequalities, and the deduction of Proposition C from Propositions A and B, are given in Section \ref{A and C}.


The arguments in \cite{second} for the general situation in which $\oldOmega$ may be reducible involve the same basic approach, except that in place of an incompressible suborbifold $\oldTheta$ such that the components of $|\oldTheta|$ are incompressible tori, one uses an incompressible $2$-suborbifold such that the underlying surfaces of its components are essential spheres. It then becomes necessary to estimate the invariants defined in Section \ref{darts section} for the components of the orbifold $\oldOmega'$ obtained by splitting $\oldOmega$ along such a suborbifold. The boundary components of $|\oldOmega|$ are spheres, and each component of the manifold obtained from $|\oldOmega|$ by attaching $3$-balls along its boundary components is an irreducible $3$-manifold. Proposition \ref{new get lost} of the present paper, which is essentially a generalization of Proposition B, gives estimates for the values of the invariants of Section \ref{darts section} on the components of $\oldOmega$. These provide a starting point for the very involved combinatorial arguments given in \cite{second}.


I am grateful to Ian Agol, Francis Bonahon, Ted Chinburg, Marc Culler, Dave Futer, Tom Goodwillie, Tracy Hall, Christian Lange, Chris Leininger, Ben Linowitz, Shawn Rafalski, Matt Stover, and  John Voight for valuable discussions and encouragement. More specifically, Agol told me how to prove
  Proposition \ref{hepcat}; Rafalski read some of my first efforts to write about orbifolds and corrected a number of errors; and Hall and Goodwillie helped lead me to the proof of Proposition\ref {fibration-category}. Some of the work presented here was done during visits to the Technion and the IAS, and I am also grateful to these institutions
for their hospitality. Some of the work was also supported by NSF grant DMS-1207720.

\section{Conventions and Preliminaries}\label{prelim section}

\abstractcomment{\tiny A lot is tied up with the properties of the characteristic
  suborbifold. }

\Number
The set of all non-negative integers will be denoted $\NN$.
The cardinality of a finite set $S$ will be denoted $\card S$.
If $A$ and $B$ are subsets of a set, $A\setminus B$ will denote the
set of all elements of $A$ which do not belong to $B$. In the special
case where $B\subset A$, we will often use the alternative notation $A-B$.
A disjoint union of sets $A$ and $B$ will be denoted $A\discup B$: that is, writing $X=A\discup B$ means that $A\cap B=\emptyset$ and that $X=A\cup B$.
\EndNumber

\Number
If $X$ is a topological space, we will
denote by $\calc(X)$ the set of all connected components of $X$.  I
will set $\compnum(X)=\card \calc(X)$.
\EndNumber

\Number\label{injectify}
A map $f:X\to Y$ between path-connected spaces will
be termed $\pi_1$-injective if $f_\sharp:\pi_1(X)\to
\pi_1(Y)$ is injective. (Here and elsewhere, base points are
suppressed from the notation in cases where the choice of a base point
does not affect the truth of a statement.)

In general, a (continuous) map $f:X\to Y$ between arbitrary spaces $X$ will be
termed $\pi_1$-injective if each path component $C$ of $X$, the map
$f|C$ is a
$\pi_1$-injective map from $C$ to the path component of $Y$ containing
$f(C)$.

A subset $A$ of a space $X$ is termed $\pi_1$-injective if the
inclusion map from $A$ to $X$ is $\pi_1$-injective.
\EndNumber

\Number\label{in and out soda} If $X$ is a topological space having
the homotopy type of a finite CW complex, the Euler characterisistic
of $X$ will be denoted by $\chibar(X)$, and we will set
\Equation\label{in and out equation}
\chibar(X)=-\chi(X).
\EndEquation
 We will denote by $h(X)$ the dimension of the
singular homology $H_1(X;\ZZ_2)$ as a $\ZZ_2$-vector space. 
\EndNumber

\Number\label{search engine}
As the overview of our methods given in the introduction indicates, our argument depend heavily on the interaction between properties of an orientable $3$-orbifold and properties of its underlying $3$-manifold. Subsections \ref{nbhd stuff}---\ref{fhs-prop} are devoted to conventions and results concerning $3$-manifolds that will be used in this paper. While some of the concepts involved will be generalized to orbifolds later, emphasizing the manifold case here will help pave the way for using the very rich literature on $3$-manifolds that is available.
\EndNumber

\Number\label{nbhd stuff}
Let $M$ be a (topological, PL or smooth) compact manifold. A (respectively topological, PL or smooth) properly embedded, codimension-$1$ submanifold $\cals$ of $M$ will be termed {\it two-sided} if it has a neighborhood $N$ in $M$ such that the pair $(N,\cals)$ is (topologically, piecewise-linearly or smoothly) homeomorphic to $(\cals\times[-1,1],\cals\times\{0\})$.
If $\cals$ is two-sided,  we will denote by $M\cut\cals$ the compact space obtained by splitting
$M$ along $\cals$. If $M$ is a topological or PL manifold, or if $M$ is smooth and $\cals$ is closed, then $M\cut\cals$ inherits the structure of a topological, PL or smooth manifold respectively. 
We will denote
by $\rho_\cals$ the natural surjection from $M\cut\cals$ to $M$. For
each component $S$ of $\cals$,  the set $\rho_\cals^{-1}(S)$
is the
union of two components of $\tcals$. We will call these
components the {\it sides} of $S$. 
We will denote by $\grock_\cals$ the unique
involution $\grock$ of $\tcals$ which interchanges the sides of every component of $\cals$ and satisfies
$(\rho|\tcals)\circ\grock=\rho|\tcals$. 

We will regard $M\cut\cals$ as a completion of
$M-\cals$; in particular, $M-\cals$ is identified with $M\cut\cals-\rho_\cals^{-1}(\cals)$ (which is  the interior of
$M\cut\cals$ in the case where $M$ is closed). If $X$ is a union of components of $M-\cals$, its closure in
$M\cut\cals$ will be denoted $\hatX$. Note that $X\mapsto\hatX$ is a
bijection between the components of $X-\cals$ and those of
$X\cut\cals$. 
\EndNumber

\Number\label{esso}
The standard $1$-sphere $S^1\subset\CC$ inherits a PL structure from $\RR$ via the covering map $t\mapsto e^{2\pi it}$, since the deck transformations $t\mapsto t+n$ are PL. This PL structure on $S^1$ gives rise to a PL structure on $D^2$ by coning. When we work in the PL category, $S^1$ and $D^2$ will always be understood to have these PL structures. Note that the self-homeomorphism of $S^1$ or $D^2$ defined by an arbitrary element of ${\rm SO}(2)$ is piecewise linear.
\EndNumber

\Number\label{just manifolds}
From this point on, all statements and arguments about manifolds are to be interpreted in the PL category except where another category is specified.
\EndNumber

\Number\label{great day}
In large part we will follow the conventions of \cite{hempel}
regarding $3$-manifolds. We will depart slightly from these
conventions in our use of the term ``irreducible'': we define a
$3$-manifold $M$ to be {\it irreducible} if $M$ is connected, every
(tame) $2$-sphere in $M$ bounds a $3$-ball in $M$, and $M$ contains no
homeomorphic copy of $\RR P^2\times[-1,1]$. Thus $M$ is ``connected
and $P^2$-irreducible'' in the sense of \cite{hempel}. The reason for
using the term ``irreducible'' in this stronger sense is that, unlike
the more classical definition, it coincides in the manifold case with
the definition of an irreducible orbifold to be given in Subsection \ref{oops}.

We will use the word ``incompressible'' only in the context of closed surfaces. A closed surface $F$ in a $3$-manifold $M$ will be termed {\it incompressible} if $F$ is two-sided and $\pi_1$-injective in $M$, and does not bound a ball in $M$. This is also consistent with the definition to be given below for orbifolds.

As in \cite{hempel}, an irreducible $3$-manifold $M$ will be termed
{\it boundary-irreducible} if $\partial M$ is
 $\pi_1$-injective in $M$; by Dehn's lemma and the loop theorem, this is equivalent to the condition that
for every properly embedded disk
$D\subset M$ there is a disk $E\subset\partial M$ such that $\partial
E=\partial D$.  We will say that an orientable $3$-manifold $M$ is {\it acylindrical} if it is
irreducible and boundary-irreducible, and every properly embedded
$\pi_1$-injective annulus in $M$ is
boundary-parallel. 


A {\it graph manifold} is a closed, irreducible, orientable $3$-manifold
$M$ which contains a $\pi_1$-injective $2$-dimensional submanifold
$\calt$ such that every component of $\calt$ is a torus and every
component of $M\cut\calt$ is a Seifert fibered space. 

When $A$ is an annulus contained in the boundary of a solid torus $J$,
We will define the {\it winding number} of $A$ in $J$ to be the order
of the cyclic group $H_1(J,A;\ZZ)$ if this cyclic group is finite, and
to be $0$ if the cyclic group is infinite.
\EndNumber

\Definition\label{P-stuff}
If $X$ is a $3$-manifold, we will denote by $\plusX$ the $3$-manifold
obtained from $X$ by attaching a ball to each component of $\partial
X$ which is a $2$-sphere. We will say that $X$ is {\it $+$-irreducible}
if $\plusX$ is irreducible. We will say that $X$ is a {\it $3$-sphere-with-holes}
 if $\plusX$ is a $3$-sphere.
\EndDefinition

\Number\label{plus-contained}
Note that if $K$ is a compact submanifold of a compact $3$-manifold
$N$, and if every component of $\partial K$ is either a component of
$\partial N$ or a surface of positive genus contained in $\inter N$,
then $K^+$ is naturally identified with a submanifold of $N^+$.
\EndNumber

The following slightly stronger version of \cite[Corollary
5.5]{Waldhausen} will be needed in Section \ref{darts section}:

\Proposition\label{stronger waldhausen}
Let $M$ be a closed, irreducible, orientable $3$-manifold. Let $T$ be
a closed, orientable $2$-manifold, no component of which is a sphere,
and let $i:T\to M$ and $j:T\to M$
be $\pi_1$-injective embeddings. Suppose that $i$ and $j$ are homotopic, and that the
components of $i(T)$ are pairwise non-parallel in $M$. Then $i$ and
$j$ are isotopic.
\EndProposition

\Proof
We argue by induction on $\compnum(T)$. If $\compnum(T)=0$ the
assertion is
  trivial, and if $\compnum(T)=1$ it is \cite[Corollary
5.5]{Waldhausen}. Now suppose that an integer $n\ge1$ is given and
that the result is true whenever
$\compnum(T)=n$. Let $M$, $T$, $i$ and $j$ be given,
satisfying the hypothesis, with
$\compnum(T)=n+1$. We may assume without loss of generality that $T\subset M$ and that $i:T\to M$ is
the inclusion map. Thus $T$ is $\pi_1$-injective. Choose a  component $V$ of $T$, and set
$T'=T-V$. Since $j:T\to M$ is homotopic to the inclusion, so is
$j':=j|T':T'\to M$. Furthermore, $T'$ is in particular
$\pi_1$-injective, and $\compnum(T')=n$. Hence the induction
hypothesis implies that 
$j|T'$ is isotopic to the inclusion map. Thus after modifying 
$j$ within its isotopy class we may assume that $j|T'$ is the
inclusion. Now $j|V$ is homotopic to the inclusion $i_V:V\to M$. Since
$j$ and the inclusion map $i$ are embeddings, both $V$ and $j_V(V)$
are disjoint from $T'$.

We will show that $i_V$ and $j_V$ are isotopic by an ambient isotopy
of $M$ which is constant on $T'$. This will immediately imply that $i$
and $j$ are isotopic, and complete the induction.

Set $W=j_V(V)$.
We may assume after an isotopy that $W$ and $V$ intersect
transverally. We may further assumed that $j$ has been chosen within
its isotopy class rel $T'$ so as
to minimize $\compnum(W\cap V)$, subject to the
condition that $W$ and $V$ intersect transversally. We now claim:
\Equation\label{nothing there}
W\cap V=\emptyset.
\EndEquation

To prove (\ref{nothing there}), suppose that $W\cap V\ne\emptyset$. Since $j_V$ and the inclusion $i_V$ are homotopic in
$M$ and $\pi_1$-injective, it follows from
\cite[Proposition 5.4]{Waldhausen} that there
exist connected subsurfaces $A\subset V$ and $B\subset W$, and a compact submanifold
$X$ of $M$, such that
$\partial A\ne\emptyset$, $\partial X=A\cup B$,
and the pair $(X,A)$ is homeomorphic to $(A\times[0,1],A\times\{0\})$.
Since $V$ and $W$
are disjoint from $T'$, we have $T'\cap\partial X=\emptyset$; hence
every component of $T'$ is either contained in $\inter X$ or disjoint
from $X$.
Note that $X$ is a handlebody since $\partial A\ne\emptyset$. If some
component $U$ of $T'$ is contained in $\inter X$, then the inclusion
homomorphism $\pi_1(U)\to\pi_1(X)$ is an injection from a
positive-genus surface group to a free group, which is
impossible. Hence $T'\cap X=\emptyset$. On the other hand, the
properties of $A$, $B$ and $X$  listed above imply that there is
an isotopy $(h_t)_{0\le t\le 1}$ of $M$, constant outside an
arbitrarily small neighborhood of $X$, such that $h_0$ is the identity
and $\compnum(h_1(W)\cap V)<\compnum(W\cap V)$. Since $T'\cap
X=\emptyset$, we may take $(h_t)$ to be constant on $T'$. This implies
that $j_V^1:=h_1\circ j$ is isotopic rel $T'$ to $j$, and that
$W_1:=j_V^1(V)$ satisifies $\compnum(W_1\cap V)<\compnum(W\cap V)$. This
contradicts the minimality of $\compnum(W\cap V)$, and (\ref{nothing
  there}) is proved.

Since $i_V$ and $j_V$ are homotopic and are $\pi_1$-injective, and $i_V(V)$ and $j_V(V)$ are disjoint by (\ref{nothing
  there}),  it
follows from \cite[Proposition 5.4]{Waldhausen}
that 
there
exists a compact submanifold
$Y$ of $M$, and a homeomorphism $\eta :V\times[0,1]\to Y$, such that
$\eta(v,0)=i_V(v)=v$ and $\eta(v,1)=j_V(v)$ for every $v\in V$. 
Since $V$ and $W$
are disjoint from $T'$, we have $T'\cap\partial Y=\emptyset$; hence
every component of $T'$ is either contained in $\inter Y$ or disjoint
from $Y$.
If some
component $U$ of $T'$ is contained in $\inter Y$, then since $U$ is
incompressible and $Y$ is homeomorphic to $V\times[0,1]$, it follows
from \cite[Proposition 3.1]{Waldhausen} that $U$
is parallel to each of the boundary components of $Y$, and in
particular to $V$. This contradicts the hypothesis that no two
components of $T=i(T)$ are parallel. Hence $T'\cap Y=\emptyset$. On the other hand, the
properties of $Y$ and $\eta $ stated above imply that there is
an isotopy $(h_t)_{0\le t\le 1}$ of $M$, constant outside an
arbitrarily small neighborhood of $Y$, such that $h_0$ is the identity
and $h_1\circ i_V=j_V$. Since $T'\cap
Y=\emptyset$, we may take $(h_t)$ to be constant on $T'$. This implies
that $j_V$ is isotopic rel $T'$ to $i_V$, as required.
\EndProof

The following result will also be needed in Section \ref{darts section}:

\Proposition\label{snuff}
Let $M$ be a closed, irreducible, orientable $3$-manifold. Let $V$ and $W$ be closed, orientable surfaces of strictly positive genus, and let $i_V:V\to M$ and $i_W:W\to M$ be $\pi_1$-injective embeddings such that $i_V(V)\cap i_W(W)=\emptyset$. Then  for any embeddings $f:V\to M$ and $g:W\to M$ homotopic to $i_V$ and $i_W$ respectively, such that $f(V)$ and $g(W)$ meet transversally, either (a) $f(V)\cap g(W)=\emptyset$, or (b) there
exist connected subsurfaces $A\subset f(V)$ and $B\subset g( W)$, and a compact submanifold
$X$ of $M$, 
such that
$\partial A\ne\emptyset$,  $\partial X=A\cup B$, and the pair $(X,A)$ is homeomorphic to
$(A\times[0,1],A\times\{0\})$.
\EndProposition

\Proof
According to \cite[Corollary 5.5]{Waldhausen}, $g$ is isotopic to $i_W$. Hence after modifying $f$ and $g$ by a single ambient isotopy, we may assume that $g=i_W$. We may also assume that $V$ and $W$ are subsurfaces of $M$ and that $i_V$ and $g=i_W$ are the inclusion maps. The hypothesis then gives $V\cap W=\emptyset$. Then $f(V)$ meets $W$ transversally, and after a small isotopy, constant on $W$, we may assume that it also meets $V$ transversally. We may also suppose that among all embeddings in its isotopy class rel $W$, having the property that $f(V)$ meets $V$ transversally, $f$ has been chosen so as to minimize $\compnum(f(V)\cap V)$. Set $V'=f(V)$. Since $V'$ and $V$ are homotopic by hypothesis, it follows from \cite[Prop 5.4]{Waldhausen} that there
exist (not necessarily proper) connected subsurfaces $A_0\subset V'$ and $C\subset V$, a compact submanifold
$X_0$ of $M$, such that $\partial X_0=A_0\cup C$,
the pair $(X_0,A_0)$ is homeomorphic to $(A_0\times[0,1],A_0\times\{0\})$, and
$A_0\cap V=\partial A_0$. We claim:

\Claim\label{fdjt}
Either $W\cap A_0\ne\emptyset$, or Alternative (a) of the conclusion of the lemma holds.
\EndClaim

To prove \ref{fdjt}, assume that $W\cap A_0=\emptyset$. In the case where $A_0=V'$, it follows that $W\cap V'=\emptyset$, which is Alternative (a). Now suppose that $A_0$ is a proper subsurface of $V'$, so that $\partial A_0\ne\emptyset$. The properties of $A$, $C$ and $X_0$ stated above then imply that $X_0$ is a handlebody, and that there is an embedding $f_1:V\to M$, isotopic to $f$ by an ambient isotopy which is constant on an arbitrarily small neighborhood of $X_0$, such that $\compnum(f(V_1)\cap V)<\compnum(V'\cap V)=\compnum(f(V)\cap V)$. The assumption $W\cap A_0=\emptyset$, together with the fact $V\cap W=\emptyset$, implies that $W$ is disjoint from $\partial X_0$. Since $W$ is a closed $\pi_1$-injective orientable surface of positive genus, it cannot be contained in the handlebody $X_0$. Hence $W\cap X_0=\emptyset$. We may therefore take the asmbient isotopy between $f$ and $f_1$ to be constant on $W$. But then the inequality $\compnum(f(V_1)\cap V)<\compnum(f(V)\cap V)$ contradicts the minimality property of $f$. Thus \ref{fdjt} is proved.

Next, we claim:
\Claim\label{folderol}
If some component of $V'\cap W$ is homotopically trivial in $M$, then Alternative (b) of the conclusion holds.
\EndClaim

To prove \ref{folderol}, note that any homotopically trivial component of $V'\cap W$ must bound a disk in $V'$, since $V'$ is incompressible. Among all disks in $V'$ bounded by components of $V'\cap W$ choose one, $A$, which is minimal with respect to inclusion. Since $W$ is also incompressible, $\partial A$ bounds a disk $B\subset W$. The minimality of $D$ implies that $A\cap B=\partial A$, so that $A\cup B$ is a $2$-sphere; by irreducibility, $A\cup B$ bounds a $3$-ball $X\subset M$. Now the pair $(X,A)$ is homeomorphic to $(A\times[0,1],A\times\{0\})$, and hence Alternative (b) holds. Thus \ref{folderol} is established.

In view of \ref{fdjt} and \ref{folderol}, we may assume that  $W\cap A_0\ne\emptyset$ and that every component of $V'\cap W$ is homotopically non-trivial in $M$. Since $W\cap A_0\ne\emptyset$, there is a component $B$ of $W\cap X_0$ with $\partial B\ne\emptyset$. We have $B\cap C\subset W\cap V=\emptyset$, and hence $\partial B\subset A_0$. In particular, each component of $\partial B$ is a component of $V'\cap W$ and is therefore homotopically non-trivial in $M$, and in particular in $W$; this implies that $B$ is $\pi_1$-injective in the incompressible surface $W$, and is therefore $\pi_1$-injective in $M$. Thus $B$ is a properly embedded connected, $\pi_1$-injective surface in $X_0$, with $\partial B\subset A_0$. Since the pair $(X_0,A_0)$ is homemorphic to $(A_0\times[0,1],A_0\times\{0\}$, it now follows from
\cite[Proposition 3.1]{Waldhausen} that $B$ is parallel in $X_0$ to a subsurface $A$ of $A_0$. This means that there is a submanifold $X$ of $X_0$ such that $\partial X=A\cup B$ and $(X,A)$ is homeomorphic to $(A\times[0,1],A\times\{0\})$. We have $\partial A=\partial B\ne\emptyset$. This gives Alternative (b) of the conclusion.

 \EndProof

\Proposition\label{final assertion}
Let $L$ be a connected, compact, $3$-dimensional
submanifold of an irreducible, orientable $3$-manifold $M$. Suppose
that $L$ is $\pi_1$-injective, that every component of $\partial
L$ is a torus, and that $L$ is not a solid torus. Then every component of
$\partial L$ is $\pi_1$-injective in $M$.
\EndProposition

\Proof Consider any component $T$ of 
$\partial L$. Since $L$ is
$\pi_1$-injective in $M$, it suffices to show that $T$ is
$\pi_1$-injective in $L$. If it is not, there is a properly embedded
disk $D\subset L$ whose boundary does not bound a disk in $T$. Let $Y$
denote a regular neighborhood of $D$ in $L$. Then $Q:=\overline{L-Y}$
has a $2$-sphere boundary component $S$, which must bound a ball $B\subset
M$ since $M$ is irreducible. We have either $B\supset Q$ or $B\cap
Q=S$. If $B\supset Q$, then $Z=B\cup Y$ is a solid torus containing
$L$ and having boundary $T$. Since $L$ is $\pi_1$-injective in $M$ and is contained in the solid torus $Z$, it has a
cyclic fundamental group. The torus $T$ is one component of $\partial 
L$. Since $\pi_1(L)$ is cyclic, $\partial L$ cannot have a second
component of positive genus. Hence $L=Z$, which contradicts the hypothesis that $L$ is not a solid torus. On the other hand, if $B\cap
Q=S$, then $T\subset B$, so that the inclusion
homomorphism $\pi_1(T)\to\pi_1(M)$ is trivial. Since $L$ is
$\pi_1$-injective in $M$, the inclusion
homomorphism $\pi_1(T)\to\pi_1(L)$ is trivial. This is impossible,
because Poincar\'e-Lefschetz duality implies that the inclusion
homomorphism $H_1(T;\QQ)\to\pi_1(L;\QQ)$ is non-trivial. 
\EndProof

We will need the following result, the proof of which we will extract from \cite{FHS}:

\Proposition\label{fhs-prop}
Let $M$ be an orientable Riemannian $3$-manifold, let $V$ and $W$ be closed orientable $2$-manifolds, and let $f_0:V\to M$ and $g:W\to M$ be smooth embeddings, each of which has least area in its homotopy class. Suppose that for every embedding $f:V\to M$, homotopic to $f_0$,
such that $f(V)$ and $g(W)$ meet transversally, either (a) $f(V)\cap g(W)=\emptyset$, or (b) there exist connected subsurfaces $A\subset f( V)$ and $B\subset g(W)$, and a compact submanifold
$X$ of $M$, 
such that
$\partial A\ne\emptyset$, $\partial X=A\cup B$, and the pair $(X,A)$ is homeomorphic to $(A\times[0,1],A\times\{0\})$.
Then $f_0(V)\cap g(W)=\emptyset$.
\EndProposition

\Proof
This is implicit in the proof of  \cite[Lemma 1.3]{FHS}. In the language of \cite{FHS}, Alternative (b) of the hypothesis of Proposition \ref{fhs-prop} is expressed by saying that
there is a product region between $f(V)$ and $g(W)$. In the special case where $f_0(V)$ and $g(W)$ meet transversally, we may apply the hypothesis of Proposition \ref{fhs-prop}, taking $f=f_0$, to deduce that either there is a product region between $f_0(V)$ and $g(W)$, or $f_0(V)\cap g(W)=\emptyset$. But according to \cite[Lemma 1.2]{FHS}, there cannot exist a product region between two subsurfaces of a Riemannian $3$-manifold which are the images of smooth embeddings of compact surfaces, each of which has least area in its homotopy class. Hence in this case we must have $f_0(V)\cap g(W)=\emptyset$, as required.

The proof of \cite[Lemma 1.2]{FHS} depends on the observation that if $f:V\to M$ and $g:W\to M$ are smooth embeddings of closed orientable surfaces in an orientable Riemannian $3$-manifold, and if $A$, $B$ and $X$  have the properties stated in Alternative (b), then we may define piecewise smooth embeddings $f':V\to M$ and $g':W\to M$ which agree with $f$ and $g$ on $\overline{V-f^{-1}(A)}$ and $\overline{W-f^{-1}(B)}$ respectively, and such that $f'|f^{-1}(A)$ and $g'|g^{-1}(B)$ are homeomorphisms of their respective domains onto $B$ and $A$, and are homotopic in $X$, rel $f^{-1}(\partial A)$ and $g^{-1}(\partial B)$, to $f|f^{-1}(A)$ and $g|g^{-1}(B)$ respectively. For the purpose of this proof, this construction of a pair of piecewise smooth embeddings $(f',g')$ from a pair of smooth embeddings $(f,g)$, involving a product region between $f(V)$ and $g(W)$, will be called a {\it swap}. If, keeping the same assumptions and notation, $f'':V\to M$ and $g'':W\to M$ are smooth embeddings which are homotopic $f'$ and $g'$ relative to annular neighborhoods of $f^{-1}(A)$ and $g^{-1}(B)$ respectively, and satisfy $f''(V)\cap g''(W)=(f(V)\cap g(W))-f(\partial A)$, we will say that $(f'',g'')$ is obtained from $(f,g)$ by a {\it smoothed swap}. If there is a product region between smooth embeddings $f$ and $g$, then there is a pair $(f'',g'')$ obtained from $(f,g)$ by a smoothed swap such that $\area f''(V)+\area g''(W)<\area f(V)+\area g(W)$; hence either $\area f''(V)<\area f(V)$ or $\area g''(W)<\area g(W)$, so that $f$ and $g$ cannot both have least area in their respective homotopy classes.

If $f_0(V)$ and $g(W)$ do not intersect transversally, it is shown  in the proof of  \cite[Lemma 1.3]{FHS} that there exists a number $\epsilon>0$ with the following property: $f_0$ may be $C^1$-approximated arbitrarily well by an embedding $f$ such that (1) $f$ is homotopic to $f_0$; (2) $f(V)$ and $g(W)$ meet transversally; (3)  if $f_0(V)\cap g(W)\ne\emptyset$ then $f(V)\cap g(W)\ne\emptyset$; and (4) if there is a product region between $f(V)$ and $g(V)$, there is a pair $(f'',g'')$, obtained from $(f,g)$ by a smoothed swap, such that 
$\area f''(V)+\area g''(W)\le\area f(V)+\area g(W)-\epsilon$. By taking $f$ to be a good enough $C^1$-approximation to $f_0$ we can guarantee that
$\area f(V)<\area f_0(V)+\epsilon$. From (1), (2) and the hypothesis, it follows that either there is a product region between $f(V)$ and $g(W)$, or $f(V)\cap g(W)=\emptyset$. 
If there is a product region between $f(V)$ and $g(W)$, then (4) 
gives a pair $(f'',g'')$, obtained from $(f,g)$ by a smoothed swap, such that 
$\area f''(V)+\area g''(W)\le\area f(V)+\area g(W)-\epsilon<\area f_0(V)+\area g(W)$. Hence either $\area f''(V)<\area f_0(V)$, or $\area g''(W)<\area g(W)$; in either case we have a contradiction to the hypothesis that $f_0$ and $g$ have least area in their homotopy classes. We must therefore have $f(V)\cap g(W)=\emptyset$, which by (3) implies 
$f_0(V)\cap g(W)=\emptyset$, as required.
\EndProof

\Number\label{orbifolds introduced}
General references for orbifolds include \cite{bmp}, \cite{chk} and \cite{kapovich}. Although smooth orbifolds are emphasized in these books, the definition of orbifold goes through without change in the topological or smooth category. (In reading the definition in the PL category, one should bear in mind that an orthogonal action of a finite group on a Euclidean space is in particular a PL action.) 

The material from here to the end of Subsection \ref{just irreducible} is meant to be interpretable in each of the three categories, except where a restriction on category is specified.

\abstractcomment{
One point to remember is that \cite{illman} applies only to very good orbifolds; this should be OK for the apps if I bear it in mind. Another point is that \cite{lange} does  not smooth pairs, only (very good) orbifolds; again I think it will OK for the apps if I bear it in mind. Do I need to think about what stratification, and the other issues in the def of general position (for example), mean in the PL context? They should be simpler, and I'm inclined to leave the issue implicit.



}

Orbifolds will be denoted by capital fraktur letters ($\frakA, \frakB, \frakC,\ldots$).

The underlying space of an orbifold $\oldPsi$ will be denoted by $|\oldPsi|$. If $\oldPsi$ is PL, and if either $\dim\oldPsi\le2$, or $\dim\oldPsi=3$ and $\oldPsi$ is orientable, then $|\oldPsi|$ inherits the structure of a PL manifold of the same dimension as $\oldPsi$.

The singular set of an orbifold $\oldPsi$ will be denoted $\fraks_\oldPsi$. We regard it as a subset of $|\oldPsi|$.

If $\oldPsi$ is an $n$-orbifold then $|\oldPsi|$ has a canonical stratification \cite[Subsection 4.5]{chen-ruan}, in which $\fraks_\oldPsi$ is the union of all strata of dimension strictly less than $n$. If $n=3$ and $\oldPsi$ is orientable, then each component of $\fraks_\oldPsi$ is either a simple closed curve or (the
underlying space of) a trivalent graph; in the former case, the given component of $\fraks_\oldPsi$ is a single stratum, and in the latter case, each edge or vertex is a stratum.

For every
$x\in|\oldPsi|$ there exist a neighborhood $\frakU$ of $x$  in
$|\oldXi|$ and a chart map (see \cite[Subsection 2.1.1]{bmp}) $\phi:\tfrakU\to\frakU$, where
$\tfrakU\subset\RR^n$ is an open ball about $0$ and $\phi(0)=x$. By
definition there is an orthogonal action of some finite group $G$ on
$\tfrakU$ such that $\phi$ induces a homeomorphism from $\tfrakU/G$
onto $\frakU$. The group $G$ is determined up to conjugacy in $O(n)$
by the point $x$, and will be denoted $G_x$. We have $G_x=\{1\}$ if
and only if $x\notin\fraks_\oldPsi$. We will refer to the order of
$G_x$ as the {\it order} of $x$. All points of a given stratum of $\oldPsi$ have the same order; for a stratum contained in $\fraks_\oldPsi$, this order will be called the order of the stratum. 

If 
$n=2$ and $\oldPsi$ is orientable, then $G_x$ is cyclic for every $x\in\fraks_\oldPsi$. 
If 
$n=2$ and $\oldPsi$ is orientable, then $G_x$ is cyclic for every point $x$ lying in a one-dimensional stratum of $\oldPsi$.

The distinction between an orbifold $\oldPsi$ and its underlying space $|\oldPsi|$ will be rigidly observed. For example, if $\oldPsi$ is path-connected and $\star\in|\oldPsi|$ is a base point, $\pi_1(\oldPsi,\star)$ will denote the orbifold fundamental group of $\oldPsi$ based at $\star$ (denoted $\pi_1^{\rm orb}(\oldPsi,\star)$ by some authors). In contrast, $\pi_1(|\oldPsi|,\star)$ of course denotes the fundamental group of the underlying space $|\oldPsi|$ based at $\star$. 
(As in the case of spaces, we, will often suppress base points from the notation for the orbifold fundamental group in statements whose truth is independent of the choice of base point.) 
Similarly, $\partial\oldPsi$ will denote the orbifold boundary of $\oldPsi$, which is itself an orbifold.
If $\oldPsi$ is an orbifold of dimension $n\le3$ such that $|\oldPsi|$ is a topological manifold, then $\partial|\oldPsi|=|\partial\oldPsi|\cup\fraks^{n-1}_\oldPsi$, where $\fraks^{n-1}_\oldPsi$ denotes the union of all $(n-1)$-dimensional components of $\fraks_\oldPsi$. In the case of an orientable orbifold $\oldPsi$ of dimension $n\le3$, we have $\fraks^{n-1}_\oldPsi=\emptyset$ and hence $\partial|\oldPsi|=|\partial\oldPsi|$.

The Seifert van Kampen theorem
for orbifolds, which is proved in \cite[Section 2.2]{bmp}, will often be used without being mentioned explicitly.

An orbifold will be said to have {\it finite type} if it is homeomorphic to $\oldUpsilon-\frakE$, where $\oldUpsilon$ is a compact orbifold with $\fraks_\oldUpsilon\subset|\inter\oldUpsilon|$, and $\frakE$ is a union of components of $\partial\oldUpsilon$. Note that in particular,
according to this definition, a finite-type orbifold has compact boundary. 

An orbifold $\oldPsi$ of finite type has a well-defined (orbifold) Euler characteristic, which we
will denote by
$\chi(\oldPsi)$. It
is not in general equal to $\chi(|\oldPsi|)$. When $\oldPsi$ has finite type we will also set
$$\chibar(\oldPsi)=-\chi(\oldPsi)$$
in analogy with (\ref{in and out equation}).

By a {\it point} of an orbifold $\oldPsi$ we will mean simply a point
of $|\oldPsi|$. By a {\it neighborhood} of a point $x\in\oldPsi$ we
mean a 
neighborhood of $x$ in $|\oldPsi|$. If
$\oldPsi'$ is a suborbifold of an orbifold $\oldPsi$, we will regard $|\oldPsi'|$
as a subspace of $|\oldPsi|$. We will say that a suborbifold $\frakV$ is a neighborhood of $x$ if $|\oldPsi'|$ is a neighborhood of $x$.

If $\oldXi$ and $\oldUpsilon$ are orbifolds of respective dimensions $m$ and $n$, we define an {\it immersion} (or a {\it submersion}) from $\oldXi$ to $\oldUpsilon$ to be a map $\frakf$ (of sets) from $|\oldXi|$ to $|\oldUpsilon|$ such that for every $x$ there exist suborbifolds $\frakU$ and $\frakV$ of $\oldXi$ and $\oldUpsilon$, which are neighborhoods of $x$ and $y$ respectively, chart maps (see \cite[Subsection 2.1.1]{bmp}) $\phi:\tfrakU\to\frakU$ and $\psi:\tfrakV\to\frakV$, where $\tfrakU\subset\RR^m$ and $\tfrakV\subset\RR^n$ are open, and an injective (or, respectively, surjective) affine map $\alpha:\RR^m\to\RR^n$ such that $\frakf(\frakU)\subset \frakV$ and $\frakf\circ\phi=\psi\circ\alpha$.
(The only ``maps'' between orbifolds that will be considered in this
paper are immersions and submersions. An orbifold homeomorphism, or
more generally an orbifold covering map, is at once an immersion and a
submersion.) An {\it embedding} of an orbifold $\oldXi$ in an orbifold $\oldUpsilon$ is defined to be a homeomorphism of $\oldXi$ onto a suborbifold of $\oldUpsilon$; any embedding is an injective immersion, from $|\oldXi|$ to $|\oldUpsilon|$, but the converse is false.

Note that an immersion or submersion $\frakf:\oldXi\to\oldUpsilon$ is
in particular a continuous map from $|\oldXi|$ to $|\oldUpsilon|$;
thus if we ignore the orbifold structure, the immersion or submersion
$\frakf$ defines a continuous map of topological spaces, which we
denote by 
 $|\frakf|:|\oldXi|\to|\oldUpsilon|$.

\nonessentialproofreadingnote{If
  I had defined general maps, not just immersions and
  submersions---and it's hard to see why I shouldn't do so---then we
  would have a genuine functor.}

By an {\it isotopy} of an orbifold $\oldPsi$ we mean a family 
$(\frakh)_{0\le t\le1}$ of self-homeomorphisms of $\oldPsi$ such that
(i) the map $(x,t)\mapsto\frakh_t(x)$ from $|\oldPsi\times[0,1]|$
to $|\oldPsi|$ is a submersion 
from $\oldPsi\times[0,1]$ to $\oldPsi$, and (ii)
$\frakh_0$
is the identity. We will use the orbifold analogues of standard
language for isotopy of manifolds; for example, two suborbifolds
$\oldXi,\oldXi'$ will be said to be (ambiently) isotopic if there is an isotopy 
$(\frakh)_{0\le t\le1}$ of $\oldPsi$ such that
$\frakh_1(\oldXi)=\oldXi$. We will say that $\oldXi$ and $\oldXi'$ are
{\it non-ambiently isotopic} if they are isotopic when regarded as
suborbifolds of the orbifold $\oldPsi'\supset\oldPsi$ obtained from
the disjoint union of $\oldPsi$ with $(\partial\oldPsi)\times[0,1]$ by
gluing the suborbifolds $\partial\oldPsi\subset\oldPsi$ and
$(\partial\oldPsi)\times\{0\}\subset (\partial\oldPsi)\times[0,1]$ via
the homeomorphism $x\to(x,0)$.

\EndNumber

\Number\label{gen pos}
Let $\oldPsi$ be an orbifold of dimension  $m\le3$, and let $\fraks^{(0)}_\oldPsi$ denote the union of all zero-dimensional strata of $\fraks_\oldPsi$ (see \ref{orbifolds introduced}). Then $|\oldPsi|-\fraks^{(0)}_\oldPsi$ is an $m$-manifold. 
A manifold $H\subset|\oldPsi|$, having   dimension strictly less than $m$, will be said to be {\it in general position with respect to  $\fraks_\oldPsi$} if $H$ is disjoint from $\fraks^{(0)}_\oldPsi$, and intersects every positive-dimensional stratum of $\fraks_\oldPsi$ transversally in the manifold $|\oldPsi|-\fraks^{(0)}_\oldPsi$. 
\EndNumber

\begin{notationremarks}\label{obd}
Suppose that $\oldPsi$ is an orbifold, that $X$ is a manifold, and
that $f:X\to|\oldPsi|$ is a map with the property that for some orbifold $\oldXi$ we have $|\oldXi|=X$, and $f=|\frakf|$ for some orbifold immersion $\frakf:\oldXi\to\oldPsi$. In this situation, the orbifold $\oldXi$ is uniquely determined by $X$, $\oldPsi$ and $f$. In situations where it is clear from the context which $\oldPsi$ and $f$ are involved, $\oldXi$ will be denoted by $\obd(X)$. 

The most common situation in which this convention will be used is the one in which $X$ is given as a submanifold of $|\oldPsi|$ for some orbifold $\oldPsi$, in which case $f$ is understood to be the inclusion map $X\to|\oldPsi|$; in this case, $\obd(X)$ is defined if and only if $X=|\oldXi|$ for some suborbifold $\oldXi$ of $\oldPsi$, and if it is, we have $\obd(X)=\oldXi$. In the case where $m:=\dim\oldPsi\le3$, a sufficient condition for $\obd(X)$ to be defined is that $X\subset|\oldPsi|$ be a manifold which has  dimension strictly less than $m$ and is in general position (see \ref{gen pos}) with respect to  $\fraks_\oldPsi$. Furthermore, if $X$ is a closed subset of $|\oldPsi|$ such that $\obd(\Fr_{|\oldPsi|}X)$ is defined, then $\obd(X)$ is defined.

Another situation in which the convention will be used is the one in which $X$ is given as a (not necessarily proper) subset of a covering space $\toldPsi$ of $\oldPsi$. In this case $f$ will be understood to be the restriction of the covering projection to $X$. In this situation, $\obd(X)$ is defined if and only if $X=|\oldXi|$ for some subborbifold $\oldXi$ of $\toldXi$. 
In practice, the reason why such a suborbifold $\oldXi$ will usually be immediate from the context, and will be left implicit.

Still another situation in which these conditions hold (and a particularly important one) is the one in
which $\oldTheta$ is a closed two-sided $2$-suborbifold of the
interior of an orientable PL $3$-orbifold $\oldOmega$, so that
$|\oldTheta|$ is a two-sided PL $2$-submanifold of $\inter|\oldOmega|$,
and we take $\oldPsi=\oldOmega$, $X=|\oldOmega|\cut{|\oldTheta|}$, and
$f=\rho_{|\oldTheta|}$. In this case the PL orbifold
$\obd(X)=\obd(|\oldOmega|\cut{|\oldTheta|})$ will be denoted by
$\oldOmega\cut\oldTheta$, and the PL immersion $\frakf$ such that
$|\frakf|=f=\rho_{|\oldTheta|}$ will be denoted by $\rho_\oldTheta$.
\end{notationremarks}

\Number\label{unions and such}
 If $\frakA$ and $\frakB$ are
suborbifolds of an orbifold $\oldPsi$ we will set
$\frakA\cup\frakB=\obd(|\frakA|\cup|\frakB|)$, provided that
$|\frakA|\cup|\frakB|$ is the underlying set of some suborbifold of
$\oldPsi$; more generally, if $(\frakA_\iota)_{\iota\in I}$ is a
family of suborbifolds of $\oldPsi$ indexed by some set $I$, we will
set $\bigcup_{\iota\in I}\frakA_\iota=\obd(\bigcup_{\iota\in
  I}|\frakA_\iota|)$ provided that the right hand side is
defined. Likewise, for orbifolds $\frakA$ and $\frakB$ of $\oldPsi$ we
will set $\frakA\cap\frakB=\obd(|\frakA|\cap|\frakB|)$ and
$\frakA\setminus\frakB=\obd(|\frakA|\setminus|\frakB|)$, provided that
the respective right hand sides are defined, and when
$\frakB\subset\frakA$ and $\frakA\setminus\frakB$ is defined we will
write $\frakA-\frakB=\frakA\setminus\frakB$. If
$\frakf:\oldXi\to\oldUpsilon$ is an immersion or submersion of
orbifolds, and $\frakA$ is a suborbifold of $\oldXi$ (or
$\oldUpsilon$) such that 
$|\frakf|^{-1}(|\frakA|)$ (or, respectively,
$|\frakf|(|\frakA|)$) is the underlying subspace of a suborbifold of
$\oldUpsilon$ (or, respectivley. $\oldXi$), we will set
$\frakf(\frakA)=\obd(|\frakf|(|\frakA|))$ (or, respectively,
$\frakf^{-1}(\frakA)=\obd(|\frakf|^{-1}(|\frakA|))$.
\nonessentialproofreadingnote{I've done a good deal more back-and-forth between suborbifolds and submanifolds than is really necessary, and I won't be able to eliminate this completely but should try to minimize it.}
\EndNumber

\Number\label{new weak}
Let $X$ be a compact  PL subset of $|\oldPsi|$, where $\oldPsi$ is an orientable (PL) orbifold of dimension at most $3$. We define a {\it strong regular neighborhood} of $X$ in $\oldPsi$ to be a suborbifold $\frakR$ of $\oldPsi$ such that $|\frakR|$ is the second-derived neighborhood of $X$ with respect to some triangulation $\calt$ of $|\oldPsi|$, compatible with its PL structure, such that both $X$ and  $\fraks_\oldPsi$ are underlying sets of subcomplexes of $\calt$. We define a {\it weak regular neighborhood} of $X$ in $\oldPsi$ to be a suborbifold $\frakR$ of $\oldPsi$ such that (i) $\frakR$ is a neighborhood of $X$ in $\oldPsi$, and (ii) if $\oldPsi'$ is the $3$-orbifold obtained from the disjoint union $\oldPsi\discup((\partial\oldPsi)\times[0,1])$ by gluing $\partial\oldPsi\subset\oldPsi$ to $(\partial\oldPsi)\times\{0\}\subset(\partial\oldPsi)\times[0,1]$ via the homeomorphism $t\mapsto(t,0)$, then a strong regular neighborhood of $\frakR\subset\oldPsi\subset\oldPsi'$ in $\oldPsi'$ is also a strong regular neighborhood of $X$ in $\oldPsi'$.
When $\oldUpsilon$ is closed, the notions of strong and weak regular neighborhoods coincide, and we will simply use the term {\it regular neighborhood} in that case.
\EndNumber

\Number\label{other injectify}We will use the analogues for orbifolds
of the conventions of \ref{injectify}. A connected suborbifold
$\frakU$ of a connected orbifold $\frakZ$ will be
termed {\it $\pi_1$-injective} if the inclusion homomorphism $\pi_1(\frakU)\to
\pi_1(\frakZ)$ is injective. 
In general, a suborbifold $\frakU$ of an orbifold $\frakZ$ will be
termed $\pi_1$-injective if each component of $\frakZ$ is
$\pi_1$-injective in the component of $\frakU$ containing
it. A closed, connected $2$-dimensional suborbifold $\frakU$ of a
$3$-orbifold $\frakZ$ will be termed {\it incompressible} if $\frakU$
is  contained in $\inter\frakZ$, is two-sided, is $\pi_1$-injective in $\frakZ$, and has non-positive Euler characteristic.
\EndNumber

\Number
For our purposes, if $\oldPsi$ is an orbifold, $H_1(\oldPsi)$ will denote the abelianization of the orbifold fundamental group $\pi_1(\oldPsi)$, and for any abelian group $A$ we will define $H_1(\oldPsi,A)$ to be $H_1(\oldPsi)\otimes A$. We will not define (or use) higher-dimensional homology for orbifolds.
\EndNumber

\DefinitionsRemarks\label{annular etc}
An orbifold, not necessarily connected, will be said to be {\it very good} if it admits a finite-sheeted (orbifold) covering space which is a manifold. It is a standard consequence of the ``Selberg lemma'' that a compact hyperolic orbifold is very good, and this fact will often be used implicitly. Note that a suborbifold of a very good orbifold is
  very good. Note also that an orbifold with finitely many components is very
good if and only if its components are very good.

As is standard, we define an $n$-orbifold will be called {\it discal} if it is (orbifold-)homeomorphic to the quotient of $D^n$ by an orthogonal action of a finite group. It follows from the Orbifold Theorem \cite{blp}, \cite{chk} that a $3$-orbifold $\oldPsi$, such that $\fraks_\oldPsi$ has no isolated points, is discal if and only if $\oldPsi$ it admits $D^3$ as a(n orbifold) 
covering by an $n$-ball; this fact will often be used without an explicit reference. A $2$-orbifold $\oldTheta$ will be termed {\it
  spherical} if $\oldTheta$ is closed and connected and $\chi(\oldTheta)>0$; if $\oldTheta$ is very good, this is equivalent to the condition that $\oldTheta$ admits $S^2$ as a covering.
A $2$-orbifold $\oldTheta$ will be termed {\it annular} or {\it toric} if $\oldTheta$ admits $S^1\times[0,1]$ or $T^2$, respectively, as a covering; this is equivalent to the condition that $\oldTheta$ is connected, that $\chi(\oldTheta)=0$,  and that $\partial\oldTheta$ is, respectively, non-empty or empty. 
A {\it \torifold} is defined to be a $3$-orbifold which is
covered by a solid torus. 
\EndDefinitionsRemarks

\Number\label{silvering}
Let $\oldPsi$ be an orbifold, and let $\oldXi$ be a codimension-$0$
suborbifold of $\partial\oldPsi$ (so that $\oldXi$ has codimension $1$ in $\oldPsi$).
The {\it double} of $\oldPsi$ along $\oldXi$, denoted $D_\oldXi\oldPsi$, is the
closed orbifold obtained from the disjoint union of two copies of
$\oldPsi$ by gluing  together the copies of $\oldXi$ contained in their boundaries via the identity homeomorphism. We will write $D\oldPsi$ for $D_{\partial\oldPsi}\oldPsi$, which is a closed orbifold. Note that this definition includes the case where
$\oldPsi$ is a manifold;  in this case $D\oldPsi$ is a closed
manifold. 

Note that for any orbifold $\oldPsi$ and any codimension-$0$
suborbifold $\oldXi$ of $\partial\oldPsi$, 
the orbifold $D_\oldXi\oldPsi$ has a canonical involution which maps each copy of 
$\oldPsi$ on to the other, via the identity homeomorphism. The quotient of
$D_\oldXi\oldPsi$ by this involution is an orbifold $\oldPsi'$ such that
$|\oldPsi'|=|\oldPsi|$ and
$\fraks_\oldPsi'=\fraks_\oldPsi\cup\oldXi$. The orbifold $\oldPsi'$
is said to be obtained from $\oldPsi$ by {\it silvering} the suborbifold $\oldXi$ of $\partial\oldPsi$. It will be denoted
$\silv_\oldXi\oldPsi$, or by $\silv\oldPsi$ in the special case where $\oldXi=\partial\oldPsi$. (Cf. \cite[p. 65, Remark]{other-bs}.) Note that
$\silv_\oldXi\oldPsi$ is always non-orientable if $\oldXi\ne\emptyset$.

 We define $[[0,1]$ to be the $1$-orbifold $\silv_{\{0\}}[0,1]$, and we define $[[0,1]]$ to be the $1$-orbifold $\silv_{\{0,1\}}[0,1]$.
\EndNumber

\Definition\label{parallel def}
Let $\oldPsi$ be an orientable $3$-orbifold, and let $\oldXi$ be a
$2$-suborbifold of $\oldPsi$. Let $\oldPi,\oldPi'\subset\oldPsi$ be
$2$-suborbifolds whose boundaries are contained in $\oldXi$. We will say
that $\oldPi$ and $\oldPi'$ are {\it parallel in the pair $(\oldPsi,\oldXi)$}
(or simply {\it parallel in $\oldPsi$} in the case where
$\oldXi=\partial\oldPsi$, or  {\it parallel} when
$\oldXi=\partial\oldPsi$ and it is understood which orbifold $\oldPsi$
is involved) if there is
an embedding $j:\oldPi\times[0,1]\to\oldPsi$ such that
$j((\partial\oldPi)\times[0,1])\subset\oldXi$, $j(\oldPi\times\{0\})=\oldPi$, and
$j(\oldPi\times\{1\})=\oldPi'$. 
\EndDefinition

In this paper, when we consider parallelism of $2$-suborbifolds in a pair $(\oldPsi,\oldXi)$, it will always be in a context in which $\oldXi$ is a suborbifold of $\partial\oldPsi$. However, in \cite{second} we will have occasion to consider parallelism in a pair $(\oldPsi,\oldXi)$ where $\oldXi$ is a properly embedded suborbifold of $\oldPsi$.

\Definition\label{just irreducible}
A $3$-orbifold $\oldPsi$ is said to be {\it irreducible} if $\oldPsi$
is connected and every two-sided spherical $2$-suborbifold of
$\oldPsi$ is the boundary of a discal $3$-suborbifold of
$\oldPsi$. As mentioned in Subsection \ref{great day}, this
generalizes our definition for the case of a $3$-manifold.
\EndDefinition

\Number\label{categorille}In the rest of this paper and in \cite{second} we adopt the convention, generalizing the convention stated in \ref{just manifolds}, will be that all statements and arguments about orbifolds are to be interpreted in the PL category except where another category is specified. Smooth orbifolds will indeed be considered at a number of points in the paper; in particular, as hyperbolic orbifolds are smooth by definition, statements involving hyperbolic orbifolds (inlcuding some of the main results of the paper) should be interpreted in the smooth category. 

According to the main theorem of \cite{illman}, every very good smooth orbifold $\oldOmega$ admits a PL structure compatible with its smooth structure, and this PL structure is unique up to PL orbifold homeomorphism. We shall denote by $\oldOmega\pl$ the orbifold $\oldOmega$ equipped with this PL structure. The main result of \cite{lange} implies that every very good PL orbifold of dimension $n\le4$ is PL homeomorphic to $\oldOmega\pl$ for some smooth $n$-orbifold $\oldOmega$. However, as there is no result in the literature guaranteeing, even for $n=3$, that every PL suborbifold of $\oldOmega\pl$ is ambiently homemorphic to a smooth suborbifold of $\oldOmega$, a little care will be required in passing between the two categories. (A more subtle question about the relationship between the categories is addressed by Proposition \ref{fibration-category}.)

The results of \cite{bonahon-siebenmann} will play an important role in this paper. While these results are proved in the smooth category in \cite{bonahon-siebenmann}, the statements and proofs go through without change in the PL category, and it is the PL versions that will be quoted in the proofs of  Propositions \ref{when vertical} and \ref{new characteristic}. Indeed, the proofs of the PL versions of the results of \cite{bonahon-siebenmann} are simpler than the proofs of the smooth versions, in that some of the arguments involve extending an isotopy of the 2-sphere to an isotopy of  the 3-ball, and the orbifold analogue of this. In the PL category this follows from Alexander's coning trick, whereas in the smooth category it requires a deep  theorem due to Cerf. I am indebted to Francis Bonahon for these observations.

\EndNumber

\Lemma\label{cover-irreducible}
Suppose that $\oldLambda$ is a connected, orientable $3$-orbifold with non-empty boundary, and that some regular covering of $\oldLambda$ is irreducible. Then $\oldLambda$ is irreducible.
\EndLemma

\Proof
Let $p:\toldLambda\to\oldLambda$ be an irreducible regular covering of $\oldLambda$.
Let $\oldPi$ be any two-sided spherical $2$-suborbifold of $\inter\oldPi$. Then every component $\toldPi$ of $p^{-1}(\oldPi)$ is spherical and hence bounds a discal $3$-suborbifold of $\toldLambda$; this discal $3$-suborbifold is unique since $\partial\toldLambda\ne\emptyset$, and will be denoted by $\tfrakB_{\toldPi}$. If $\toldPi$ and $\toldPi'$ are distinct components of $p^{-1}(\oldPi)$, there is a deck transformation $\tau$ carrying $\toldPi$ to $\toldPi'$; since $\toldPi$ and $\tau(\toldPi)=\toldPi'$ are disjoint, and $\partial\toldLambda\ne\emptyset$, we must have either (a) $\tau(\tfrakB_{\toldPi})\subset\inter\tfrakB_{\toldPi}$, (b) $\tfrakB_{\toldPi}\subset\inter \tau(\tfrakB_{\toldPi})$, or (c) $\tau(\tfrakB_{\toldPi})\cap\tfrakB_{\toldPi}=\emptyset$. If (a) or (b) holds, there are infinitely many deck transformations mapping the compact suborbifold $\tfrakB_{\toldPi}$ of $\toldLambda$ into itself, which is impossible. Hence the discal orbifolds
$\tfrakB_{\toldPi}$ and $\tfrakB_{\toldPi'}=\tau(\tfrakB_{\toldPi})$ are disjoint. Since this holds for any two distinct components $\toldPi$ and $\toldPi'$ of $p^{-1}(\oldPi)$, there is a suborbifold $\frakB$ of $\oldLambda$ such that $p^{-1}(\frakB)=\bigcup_{\toldPi\in\calc(p^{-1}(\oldPi))}\tfrakB_{\toldPi}$; we have $\partial\frakB=\oldPi$, and for any component $\toldPi$ of $p^{-1}(\oldPi)$, the orbifold $\frakB$ is homeomorphic to the quotient of $\tfrakB_{\toldPi}$ by its stabilizer in the group of deck transformations. Since $\tfrakB_{\toldPi}$ is discal, its quotient $\frakB$ by a finite group action is discal according to \ref{annular etc}.
\EndProof

\DefinitionsRemarks\label{oops}
A $3$-orbifold $\oldPsi$ will be termed {\it weakly \simple} if it is very good,  compact,
irreducible, and non-discal, and every $\pi_1$-injective two-sided toric suborbifold
of $\inter\oldPsi$ is parallel in $\oldPsi$ to a component of $\partial\oldPsi$. A $3$-orbifold $\oldPsi$ will be termed {\it strongly
  \simple} if (I) $\oldPsi$ is  compact and non-discal, (II) some finite-sheeted regular covering of $\oldPsi$ is an irreducible $3$-manifold, and (III) $\pi_1(\oldPsi)$ has no free abelian subgroup of rank $2$.

Note that if $\oldPsi$ is strongly \simple, Condition (II) in the definition implies that $\oldPsi$ is very good. According to Lemma \ref{cover-irreducible}, Condition (II) also implies that $\oldPsi$ is irreducible. Note also that Condition (III) implies that $\oldPsi$ contains no $\pi_1$-injective toric suborbifold whatever; hence a strongly \simple\ $3$-orbifold is weakly \simple.

If $\oldOmega$ is a closed, orientable hyperbolic $3$-orbifold, then $\oldOmega\pl$ is strongly \simple. 
\EndDefinitionsRemarks

\Number\label{boundary is negative}
Let $\oldPsi$ be a compact, orientable $3$-orbifold. If $\oldPsi$ is
weakly \simple, then in particular it is irreducible and non-discal;
this implies that every boundary component of $\oldPsi$ has
non-positive Euler characteristic.

If $\oldPsi$ is strongly \simple\ and boundary-irreducible, then in fact every boundary
component of $\oldPsi$ has strictly negative Euler
characteristic. This is because any boundary component of $\oldPsi$
having Euler characteristic $0$ would be toric by definition and would
be $\pi_1$-injective by boundary irreducibility; according to an observation in \ref{oops}, this would contradict
strong \simple ity.
\EndNumber

\abstractcomment{
Weak \simple ity
mostly comes up in introducing the characteristic suborbifold. In
Bonahon and Siebenmann's result the complementary pieces are weakly
\simple\ or $S^1$-fibered.}

\DefinitionsRemarks\label{wuzza turnover}
A $2$-orbifold of finite type (see \ref{orbifolds introduced} will be termed {\it negative} if each of its
components has negative Euler characteristic. Note that the empty
$2$-orbifold is negative.

We define a {\it negative turnover} to be a negative, compact, orientable $2$-orbifold $\oldTheta$ such that $|\oldTheta|$ is a $2$-sphere and  $\card \fraks_\oldTheta=3$. 
\EndDefinitionsRemarks

\Notation\label{wuzza weight}
If $\oldPsi$ is an orbifold, and $S$ is a subset of $|\oldPsi|$ such that $S\cap\fraks_\oldPsi$ is finite, the number $\card (S\cap\fraks_\oldPsi)$ will be called the {\it weight} of $S$ in $\oldPsi$, and will be denoted $\wt_\oldPsi S$, or simply by $\wt S$ when $\oldPsi$ is understood.
\EndNotation

\Number\label{2-dim case}
Let $\otheroldLambda$ be a finite-type  $2$-orbifold. Suppose that $\fraks_\otheroldLambda$ is finite (which is in particular true if $\otheroldLambda$ is orientable). Let $x_1,\ldots,x_n$ denote the distinct points of $\fraks_\otheroldLambda$, and let $p_i$ denote the order of the singular point $x_i$ for $i=1,\ldots,n$. Then we have $\chi(\otheroldLambda)=\chi(|\otheroldLambda|)-\sum_{i=1}^n(1-1/p_i)$. In particular we have $\chi(\otheroldLambda)\le\chi(|\otheroldLambda|)$. 

Another consequence of the formula for $\chi(\otheroldLambda)$ given above, which will be used many times in this paper and in \cite{second}, often without an explicit reference, is that if $\otheroldLambda$ is an orientable annular orbifold, then either $|\otheroldLambda|$ is an annulus and $\fraks_\otheroldLambda=\emptyset$, or $|\otheroldLambda|$ is a disk and $\fraks_\otheroldLambda$ consists of two points, both of order $2$.
\EndNumber

\Number\label{cobound}An immediate consequence of the description of orientable annular orbifolds given in \ref{2-dim case} is that If $\frakB$ is an orientable annular $2$-orbifold, and if
$C\subset|\frakB|-\fraks_\frakB$ is a simple closed curve such that
$\omega(C)$ is $\pi_1$-injective in $\frakB$, then there is a
weight-$0$ annulus $A\subset|\frakB|$ having $C$ as a boundary
component, and having its other boundary component contained in
$\partial|\frakB|$. 
\EndNumber


\Proposition\label{kinda dumb}
Let $\oldPsi$ be a very good $3$-orbifold, and let $\oldTheta$ be a $2$-suborbifold of $\oldPsi$ which is either contained in  $\partial\oldPsi$, or properly embedded and two-sided. Then the following conditions are equivalent:
\begin{enumerate}[(a)]\item$\oldTheta$ is $\pi_1$-injective in $\oldPsi$; 
\item For every discal $2$-suborbifold $\frakD$ of $\oldPsi$ such that $\frakD\cap\oldTheta=\partial\frak D$, there is a discal $2$-suborbifold $\frakE$ of $\oldTheta$ with $\partial\frakE=\partial\frakD$.
\item For every discal $2$-suborbifold $\frakD$ of $\oldPsi$ such that (i) $\frakD\cap\oldTheta=\partial\frak D$ and (ii) $|\frakD|$ is in general position (see \ref{gen pos}) with respect to $\fraks_\oldPsi$, there is a discal $2$-suborbifold $\frakE$ of $\oldTheta$ with $\partial\frakE=\partial\frakD$.
\end{enumerate}
\EndProposition

\Proof
We first give the proof in the case where $\oldTheta\subset\partial \oldPsi$. Since $\oldPsi$ is very good, we may fix a regular cover $p:\toldPsi\to\oldPsi$ such that $\toldPsi$ is a $3$-manifold. Let $G$ denote the group of deck transformations of this covering.

Let us show that (a) implies (b).
Suppose that (a) holds, and that $\frakD$ is a discal $2$-suborbifold of $\oldPsi$  such that $\frakD\cap\oldTheta=\partial\frak D$. Then every component of $p^{-1}(\frakD)$ is a disk whose boundary lies in $p^{-1}(\oldTheta)$. Since $\oldTheta$ is $\pi_1$-injective in $\oldPsi$, in particular $p^{-1}(\oldTheta)$ is $\pi_1$-injective in $\toldPsi$, and so each component of $p^{-1}(\frakD)$ bounds a disk in $p^{-1}(\oldTheta)$. Among all disks in $p^{-1}(\oldTheta)$ bounded by components of $p^{-1}(\oldTheta)$, choose one, $D_0$, which is minimal with respect to inclusion. Then for every $\tau\in G$ we have either $\tau(D_0)=D_0$ or $\tau(D_0)\cap D_0=\emptyset$. Hence $\frakE:=p(D_0)$ is orbifold-homeomorphic to the quotient of $D_0$ by its stabilizer in $G$, and is therefore a discal suborbifold of $\oldTheta$ with $\partial\frakE=\frakD$. This establishes (b).

It is trivial that (b) implies (c). To show that (c) implies (a), we will suppose that (a) does not hold, and produce a discal orbifold $\frakD$ violating (c).
Since (a) does not hold, 
$p^{-1}(\oldTheta)$ is not 
$\pi_1$-injective in $\toldPsi$. It therefore follows from the equivariant loop theorem 
\cite{meeks-yau} that there is a non-empty, properly embedded $G$-invariant submanifold $\cald$ of $\toldPsi$, each of whose components is a disk, such that $\partial\cald\subset p^{-1}(\oldTheta)$, and no component of $\partial\cald$ bounds a disk in $p^{-1}(\oldTheta)$. Choose a $G$-invariant regular neighborhood $\calr$ of $\cald$ such that $\cald':=\Fr_\oldPsi\calr$ has transverse intersection with $p^{-1}(\sigma)$ for every stratum $\sigma$ of $\fraks_\oldPsi$. By choosing $\calr$ to be a sufficiently small neighborhood of $\cald$ we may guarantee that $\partial\cald'$ is contained in $ p^{-1}(\oldTheta)$, and is isotopic in $ p^{-1}(\oldTheta)$ to $\partial\cald$. Hence no component of $\partial\cald'$ bounds a disk in $p^{-1}(\oldTheta)$. If $D_1$ is a component of $\cald'$ then $\frakD:=p(D_1)$ is orbifold-homeomorphic to the quotient of $D_1$ by its stabilizer in $G$. Hence $\frakD$ is a discal suborbifold of $\oldPsi$. Since $\cald'$ is properly embedded in $\toldPsi$, has boundary
contained in $ p^{-1}(\oldTheta)$, and
has transverse intersection with $p^{-1}(\sigma)$ for every stratum $\sigma$ of $\fraks_\oldPsi$, the suborbifold $\frakD$ satisfies Conditions (i) and  (ii). 
To show that $\frakD$ violates (c), it is enough to show that there is no  discal $2$-suborbifold $\frakE$ of $\oldTheta$ with $\partial\frakE=\partial\frakD$. If $\frakE$ is such a discal $2$-suborbifold, then $\frakE:=p^{-1}(\frakE)$ is a disjoint union of disks, and $\partial\tfrakE=p^{-1}(\partial\frakD)$. In particular the component $\partial D_1$ of $\partial\cald$ bounds a disk in $p^{-1}(\oldTheta)$, a contradiction. This completes the proof of the proposition in the case where $\oldTheta\subset\partial\oldPsi$.

Now consider the case in which $\oldTheta$ is  properly embedded and two-sided. Set $\oldPsi'=\oldPsi\cut\oldTheta$, and $\oldTheta'=\rho_\oldTheta^{-1}(\oldTheta)\subset\partial\oldPsi'$. It follows from the Seifert-van Kampen theorem for orbifolds (see \cite[Section 2.2]{bmp}) that $\oldTheta$ is $\pi_1$-injective in $\oldPsi$ if and only if $\oldTheta'$ is $\pi_1$-injective in $\oldPsi'$. It is  clear that each of the conditions (b) and (c) for $\oldPsi$ and $\oldTheta$ is equivalent to the same condition with $\oldPsi$ and $\oldTheta$ replaced by $\oldPsi'$ and $\oldTheta'$. Hence the assertion of the proposition in this case follows by applying the case already proved, with $\oldPsi'$ and $\oldTheta'$ playing the respective roles of $\oldPsi$ and $\oldTheta$.
\EndProof

\Corollary\label{injective hamentash}
Let $\oldPsi$ be an orientable $3$-orbifold, and let $S\subset|\oldPsi|$ be a $2$-sphere of weight $3$, in general position with respect to $\fraks_\oldPsi$. Then $\obd(S)$ is $\pi_1$-injective in $\oldPsi$.
\EndCorollary

\Proof
We may assume after a small non-ambient isotopy that $S\subset\inter|\oldPsi|$. Set $\oldTheta=\obd(S)$. 
According to Proposition \ref{kinda dumb}, it suffices to show that if $\frakD$ is a discal $2$-suborbifold of $\oldPsi$ such that $\frakD\cap\oldTheta=\partial\frak D$, and $|\frakD|$ is in general position with respect to $\fraks_\oldPsi$, then $\partial\frakD$ is the boundary of a discal $2$-suborbifold of $\oldTheta$. Since $\oldPsi$ is orientable we have $\dim\fraks_\oldPsi\le1$, and the general position condition therefore implies  that $\fraks_\frakD=\frakD\cap\fraks_\oldPsi$ consists at most of isolated points; hence $\frakD$ is orientable, so that $\partial\frakD$ is a closed $1$-submanifold of $\oldTheta$. But since $|\oldTheta|$ is a $2$-sphere and $\fraks_\oldTheta$ has cardinality $3$, every closed $1$-submanifold of $\oldTheta$ bounds a discal suborbifold of $\oldTheta$. 
\EndProof

\Notation\label{lambda thing}
If $\oldPsi$ is a compact, orientable $3$-orbifold, we will define an
integer $\lambda_\oldPsi$ by setting $\lambda_\oldPsi=2$ if every
component of $\fraks_\oldPsi$ is an arc or a simple closed curve, and
$\lambda_\oldPsi=1$ otherwise.
\EndNotation

\Proposition\label{almost obvious}
Let $\oldPsi$ be a strongly \simple, orientable  $3$-orbifold 
containing
no embedded negative turnovers (see \ref{wuzza turnover}). Suppose that $\lambda_\oldPsi=2$. Then $|\oldPsi|$ contains no weight-$3$
sphere which is in general position with respect to $\fraks_\oldPsi$.
\EndProposition

\Proof
Suppose that $S\subset|\oldPsi|$ is a weight-$3$ sphere in general position with respect to $\fraks_\oldPsi$. By
\ref{injective hamentash}, $\omega(S)$ is $\pi_1$-injective in $\oldPsi$. If $\chi(\obd(S))=0$ then $\obd(S)$ is a $\pi_1$-injective two-sided toric suborbifold of $\oldPsi$; as observed in \ref{oops}, this contradicts the strong \simple ity of $\oldPsi$. If $\chi(\obd(S))<0$ then $\obd(S)$ is an embedded negative turnover, a contradiction to the hypothesis. There remains the possibility that $\chi(\obd(S)>0$, so that $\obd(S)$ is spherical. Since a strongly \simple\ $3$-orbifold is irreducible (see \ref{oops}), $\obd(S)$ is the boundary of a discal $3$-suborbifold $\frakB$ of $\oldPsi$. Since $\lambda_\oldPsi=2$, every component of $|\frakB|\cap\fraks_\oldPsi$ is a properly embedded arc or simple closed curve in $|\frakB|$. If $m$ denotes the number of arc components of $|\frakB|\cap\fraks_\oldPsi$, we have $\wt S=2m$, a contradiction since $\wt S=3$.
\EndProof

\begin{definitionremark}
Generalizing the definition given for manifolds in \ref{great day}, we define a $3$-orbifold $\oldPsi$ to be {\it
boundary-irreducible} if $\partial\oldPsi$ is $\pi_1$-injective in $\oldPsi$. It follows from Proposition \ref{kinda dumb} that $\oldPsi$ is boundary-irreducible if and only if the boundary of every properly
embedded discal $2$-orbifold $\oldDelta$ in $\oldPsi$, such that $|\oldDelta|$ is in general position with respect to 
$\fraks_\oldPsi$, is the boundary of a discal $2$-orbifold in $\oldPsi$. It also follows from Proposition \ref{kinda dumb} that this characterization of boundary-irreducibility remains valid if the general position condition on $|\oldDelta|$ is omitted.

\nonessentialproofreadingnote{In the presence of
irreducibility, this implies that every  properly embedded discal suborbifold is boundary-parallel.  In the one place where that fact has come up, I've just given the argument (it's
in the proof of Prop \ref{when vertical}.)
A related point is that I had been talking as if every properly embedded discal suborbifold of an orientable $3$-orbifold were two-sided. This is of course wrong, because of the example in which $\partial|\oldDelta|$ is made up of an arc in $\partial|\oldPsi|$ and an arc in $\fraks_\oldPsi$. The corresponding statement for annular suborbifolds is also false, because a M\"obius band (which is a manifold and in particular an orbifold) can be properly embedded in an orientable $3$-manifold. That misconception was the reason why I didn't have the term `two-sided'' in the statement just made. At the moment I think this note can probably be ignored...}
\end{definitionremark}

\Definition
An orbifold $\oldPsi$ will be termed {\it componentwise irreducible}, or {\it componentwise boundary-irreducible}, if each of its components is, respectively, irreducible or boundary-irreducible. We will say that $\oldPsi$ is {\it componentwise strongly \simple} if $\oldPsi$ is compact and each of its components is strongly \simple.
\EndDefinition

\Lemma\label{oops lemma} 
Let $\oldPsi$ be a compact $3$-orbifold, and let $\oldTheta$ be a properly embedded $2$-suborbifold of $\oldPsi$ whose components are two-sided, have non-positive Euler characteristic, and are $\pi_1$-injective. 
\begin{itemize}
\item If $\oldPsi$ is componentwise irreducible, or componentwise strongly \simple, then $\oldPsi\cut\oldTheta$ is, respectively, componentwise irreducible or componentwise strongly \simple.
\item If $\oldTheta$ is closed and $\oldPsi$ is componentwise boundary-irreducible, then
$\oldPsi\cut\oldTheta$ is componentwise
boundary-irreducible.
\end{itemize}
\EndLemma

\Proof
We may assume that $\oldTheta$ is connected, since the general case will follow from the connected case by induction on the number of components of $\oldTheta$. 

If $\oldPsi_0$ denotes the component of $\oldPsi$ containing $\oldTheta$, every component of $\oldPsi\cut\oldTheta$ is either a component of $(\oldPsi_0)\cut\oldTheta$, or a component of $\oldPsi$ distinct from $\oldPsi_0$; hence we may assume that $\oldPsi$ is connected.

First suppose that $\oldPsi$ is  irreducible. To prove that $\oldPsi\cut\oldTheta$ is componentwise irreducible, it suffices to prove that every two-sided spherical $2$-suborbifold $\frakV$ of $\oldPsi-\oldTheta$ bounds a discal $3$-suborbifold of $\oldPsi-\oldTheta$. Since $\oldPsi$ is  irreducible, $\frakV$ bounds a discal $3$-suborbifold $\frakB$ of $\oldPsi$. Since $\frakV\cap\oldTheta=\emptyset$, we must have either $\frakB\subset\oldPsi-\oldTheta$ or $\frakB\supset\oldTheta$. If the latter alternative holds, the $\pi_1$-injectivity of $\oldTheta$ implies that $\pi_1(\oldTheta)$ is isomorphic to a subgroup of $\pi_1(\frakB)$; this is impossible, because the hypothesis $\chi(\oldTheta)\le0$ implies that $\pi_1(\oldTheta)$ is infinite,
whereas the discality of $\frakB$ implies that $\pi_1(\frakB)$ is finite. Thus $\oldPsi\cut\oldTheta$ is indeed componentwise irreducible.

Now suppose that $\oldPsi$ is  strongly \simple. 
Since $\oldTheta$ is $\pi_1$-injective in $\oldPsi$, the immersion $\rho_\oldTheta$ is $\pi_1$-injective; since the  strong \simple ity of $\oldPsi$ implies that $\pi_1(\oldPsi)$ has no rank-$2$ free abelian subgroup (see Definition \ref{oops}), it follows that no component of $\oldPsi\cut\oldTheta$ has a fundamental group with a rank-$2$ free abelian subgroup, as required by Condition (III) of Definition \ref{oops}. 

Next note that the  strong \simple ity of $\oldPsi$ gives a finite-sheeted regular covering $p:\toldTheta\to\oldTheta$ such that $\oldPsi$ is an irreducible $3$-manifold. Note that $\toldTheta:=p^{-1}(\oldTheta)$ is properly embedded and $\pi_1$-injective in $\toldPsi$. Applying the assertion proved above, with $\toldPsi$ and $\toldTheta$ playing the respective roles of $\oldPsi$ and $\oldTheta$, we deduce that $\toldPsi\cut{\toldTheta}$ is componentwise irreducible. Every component of $\oldPsi\cut\oldTheta$ admits some component of $\toldPsi\cut{\toldTheta}$ as a regular covering, and thus satisfies Condition (II) of Definition \ref{oops}. 

To show that $\oldPsi\cut\oldTheta$ is componentwise strongly \simple, it remains only to prove that it has no discal component. If $\frakB$ is any component of $\oldPsi\cut\oldTheta$, then $\partial\frakB$ has a $\pi_1$-injective $2$-suborbifold $\frakV$ homeomorphic to some component of $\oldTheta$. The hypothesis gives $\chi(\frakV)\le0$, so that $\pi_1(\frakV)$ is infinite. Hence $\pi_1(\frakB)$ is infinite, and $\frakB$ cannot be discal.


Finally, suppose that $\oldTheta$ is closed and that $\oldPsi$ is  boundary-irreducible. To prove that $\oldPsi\cut\oldTheta$ is componentwise boundary-irreducible, we must show that $\partial(\oldPsi\cut\oldTheta)$ is $\pi_1$-injective in $\oldPsi\cut\oldTheta$. But since $\oldTheta$ is closed, each
component of $\partial(\oldPsi\cut\oldTheta)$ is either a component of $\partial\oldPsi$, which by the  boundary-irreducibility of $\oldPsi$ is $\pi_1$-injective in $\oldPsi$ and hence in $\partial(\oldPsi\cut\oldTheta)$, or a component of $\rho_\oldTheta^{-1}(\oldTheta)$. A component of the latter type is also $\pi_1$-injective since $\oldTheta$ is $\pi_1$-injective in $\oldPsi$.
\EndProof



\Definitions\label{acylindrical def}
Let $\oldPsi$ be an orientable, componentwise irreducible $3$-orbifold, and let $\oldXi$ be a
$\pi_1$-injective $2$-suborbifold of $\partial\oldPsi$. 
An orientable annular $2$-orbifold $\oldPi\subset\oldPsi$ 
will be called {\it
  essential in the pair $(\oldPsi,\oldXi)$} (or simply   {\it
  essential (in $\oldPsi$)} in the case where $\oldPsi$ is componentwise
  boundary-irreducible and $\oldXi=\partial\oldPsi$) if (1)
$\oldPi\cap\partial\oldPsi=\partial\oldPi\subset \oldXi$, (2) $\oldPi$ is
$\pi_1$-injective in $\oldPsi$, and (3) $\oldPi$ is {\it not} parallel
 in the pair $(\oldPsi,\oldXi)$ either to a suborbifold of $\oldXi$ or to a component of
$\overline{(\partial\oldPsi)-\oldXi}$. 
We define an {\it acylindrical pair} to be an ordered pair 
$(\oldPsi,\oldXi)$, where $\oldPsi$ is an orientable $3$-orbifold, $\oldXi$ is a
$\pi_1$-injective $2$-suborbifold of $\partial\oldPsi$, and
$\oldPsi$ contains no orientable annular $2$-orbifold which is    essential in
the pair
$(\oldPsi,\oldXi)$. This definition generalizes the definition of an acylindrical $3$-manifold given in \ref{great day} in the sense that an orientable, irreducible $3$-manifold $M$ is acylindrical if and only if the manifold pair $(M,\partial M)$, regarded as an orbifold pair, is acylindrical.


\EndDefinitions

\Proposition\label{what essential means}
Let $\oldPsi$ be an orientable, componentwise irreducible $3$-orbifold, let $\oldXi$ be a
$\pi_1$-injective $2$-suborbifold of $\partial\oldPsi$, and suppose that
$\oldPi$ is an   essential orientable annular $2$-orbifold 
in the pair $(\oldPsi,\oldXi)$. Let $D\subset|\oldPsi|$ be a
weight-$0$ disk such that $\beta:=D\cap|\oldPi|$ is an arc in
$\partial D$, and $D\cap\partial|\oldPsi|= \overline{(\partial
  D)-\beta}\subset\oldXi$. Then $\beta$ is the frontier in $|\oldPi|$
of a weight-$0$ disk.
\EndProposition

\Proof
Since $\oldPi$ is annular and orientable, $|\oldPi|$ is either a
weight-$0$ annulus, or a weight-$2$ disk such that the two points of
$|\oldPi|\cap\fraks_\oldPsi$ are both of order $2$. Hence if we assume
that $\beta$ is not the frontier in $|\oldPi|$
of a weight-$0$ disk, and if $\frakR$ denotes a strong regular neighborhood
of $\omega(\beta)$ in $\oldPi$, then $|\overline{\oldPi-\frakR}|$ is
either a weight-$0$ disk or a disjoint union of two weight-$1$
disks. Thus $\overline{\oldPi-\frakR}$ has either one or two
components, and each of its components is an orientable discal $2$-orbifold. Set $R=|\frakR|$ and $ J=\calc(\overline{\oldPi-\frakR})$.

Let $Z$ be a weight-$0$ ball in $|\oldPsi|$ such that (i) $D\subset
Z$, (ii)
$Z\cap|\oldPi|=R$,
(iii) $R\subset\partial Z$, (iv)
$Q:=Z\cap\partial|\oldPsi|$ is a disk contained in $|\oldXi|$, and (v)
$Q\cup R$ is a regular neighborhood of $\partial\beta$ in
$|\partial\oldPi|$ (and thus consists of two
arcs). Then $Q\cup R$ is an
annulus in $\partial Z$, and hence $\cald:=\overline{(\partial
  Z)-(Q\cup R)}$ is a disjoint union of two weight-$0$ disks.

For each $\frakV\in  J$, let 
$\oldTheta_\frakV$ denote the union of $\frakV$ with the component or
components of $\obd(\cald)$ that meet $\frakV$. Then $\oldTheta_\frakV$ is a
properly embedded $2$-orbifold in $\oldPsi$ whose boundary is
contained in $\oldXi$; since $\frakV$ is discal, and since each
component of $\cald$ is a weight-$0$ disk meeting $|\frakV|$ in an arc
or the empty set, $\oldTheta_\frakV$ is also discal. Since $\oldXi$
is $\pi_1$-injective, it follows from 
Proposition \ref{kinda dumb} (more specifically the implication (a)$\Rightarrow$(b)), applied with $\oldTheta_\frakV$ playing
the role of $\frakD$, that there is a discal $2$-suborbifold
$\frakE_\frakV$ of $\oldXi$ with
$\partial\frakE_\frakV=\partial\oldTheta_\frakV$. 
Then $\oldTheta_\frakV\cup\frakE_\frakV$
is a spherical $2$-orbifold, and since $\oldPsi$ is irreducible, $\oldTheta_\frakV\cup\frakE_\frakV$
 bounds a discal $3$-orbifold
$\oldGamma_\frakV\subset\oldPsi$. 

lf $\frakV$ and $\frakV'$ are distinct elements of $J$, we have $\oldTheta_\frakV\cap\oldTheta_{\frakV'}=\emptyset$. Hence
$\oldTheta:=\bigcup_{\frakV\in J}\oldTheta_\frakV$ is a $2$-orbifold whose components are the orbifolds $\oldTheta_\frakV$ for $\frakV\in J$. By construction we have $\oldTheta\cap\inter Z=\emptyset$, and hence $Z$ is contained in the closure of some component $U$ of $|\oldPsi-\oldTheta|$. For each $\frakV\in J$, we have $Z\cap|\oldTheta_\frakV|\supset Z\cap|\frakV|\ne\emptyset$, and hence $|\oldTheta_\frakV|\subset \overline{U}$. It follows that $\oldTheta\subset \Fr_\oldPsi\frakU$, where $\frakU=\obd(\overline{U})$. 
This in turn implies that $|\oldPi|=R\cup\overline{|\oldPi-\frakR|}\subset Z\cup|\oldTheta|\subset|\frakU|$, so that $\oldPi\subset\frakU$.

For each $\frakV\in J$, we have 
$\Fr\oldGamma_\frakV=\oldTheta_\frakV\subset\Fr\frakU$; hence either $\frakU\subset\oldGamma_\frakV$ or $\frakU\cap\oldGamma_\frakV=\oldTheta_\frakV$. If
$\frakU\subset\oldGamma_\frakV$, then in particular $\oldPi\subset\oldGamma_\frakV$, and since $\oldPi$ is $\pi_1$-injective in $\oldPsi$, it is in particular $\pi_1$-injective in $\oldGamma_\frakV$; this is impossible, since the annularity of $\oldPi$ implies that $\pi_1(\oldPi)$ is infinite, while the discality of $\oldGamma_\frakV$ implies that $\pi_1(\oldGamma_\frakV)$ is finite. Hence $\frakU\cap\oldGamma_\frakV=\oldTheta_\frakV$ for each $\frakV\in J$. It follows that $\oldGamma_\frakV$ is a component of $\overline{\oldPsi-\frakU}$ for each $\frakV\in J$. Furthermore, if $\frakV$ and $\frakV'$ are distinct elements of $J$, then since the frontiers 
$\oldTheta_\frakV$ and $\oldTheta_{\frakV'}$ of
$\oldGamma_\frakV$ and $\oldGamma_{\frakV'}$ are disjoint, $\oldGamma_\frakV$ and $\oldGamma_{\frakV'}$ are distinct components of
$\overline{\oldPsi-\frakU}$ and are therefore disjoint.
Hence
$\oldGamma:=\bigcup_{\frakV\in J}\oldGamma_\frakV$ is a $2$-orbifold whose components are the orbifolds $\oldGamma_\frakV$ for $\frakV\in J$, and $\Fr\oldGamma=\oldTheta$.

Since $\frakU\cap\oldGamma_\frakV=\oldTheta_\frakV$ for each $\frakV\in J$, we have
$\frakU\cap\oldGamma=\oldTheta$. Since $\frakZ:=\obd(Z)\subset\frakU$, it follows that 
$\frakZ\cap\oldGamma=\frakZ\cap\oldTheta=\obd(\cald)$. Thus if we set $\oldLambda=\frakZ\cup\oldGamma$, we have $\Fr_\oldPsi\oldLambda=\overline {\frakZ-\obd(\cald)}\cup\overline{\oldGamma-\obd(\cald)}=\frakR\cup\overline{\oldPi-\frakR}=\oldPi$.



Now for each 
$\frakV\in J$, since the orientable $3$-orbifold $\oldGamma_\frakV$ is discal, and the $2$-suborbifolds $\oldTheta_\frakV$ and $\overline{(\partial\oldGamma_\frakV)-\oldTheta_\frakV}=\frakE_\frakV$ of $\partial\oldGamma_\frakV$ are discal, the pair $(\oldGamma_\frakV,\frakE_\frakV)$ is homeomorphic to $(\frakE_\frakV\times[0,1],\frakE_\frakV\times\{0\})$. If we set $\frakE=\bigcup_{\frakV\in J}\frakE_\frakV=\oldGamma\cap\partial\oldPsi$, it
 follows that
$(\oldGamma,\frakE)$ is homeomorphic to $(\frakE\times[0,1],\frakE\times\{0\})$.
Since each component of $\cald$ is a disk in $|\oldTheta|$ meeting $|\partial\oldTheta|=|\partial\frakE|$ is an arc, we deduce that the triad $(\oldGamma,\frakE,\obd(\cald))$ is homeomorphic to $(\frakE\times[0,1],\frakE\times\{0\},\obd(\alpha)\times[0,1])$, where $\alpha=\cald\cap|\partial\frakE|\subset Q$. On the other hand, the triad 
$(\frakZ,\obd(Q),\obd(\cald))$ is homeomorphic to $(\obd(Q)\times[0,1],\obd(Q)\times\{0\},\obd(\alpha)\times[0,1])$. Since $\oldLambda=
\frakZ\cup\oldGamma$ and $\frakZ\cap\oldGamma=\obd(\cald)$, it now follows that
the pair $(\oldLambda,\oldPhi)$, where
$\oldPhi=\oldLambda\cap\partial\oldPsi=\oldLambda\cap\oldXi$, is homeomorphic to $(\oldPhi\times[0,1],\oldPhi\times\{0\})$. Hence $\oldPi=\Fr\oldLambda$ is
    parallel in the pair $(\oldPsi,\oldXi)$ to a suborbifold of $\oldXi$,
a contradiction to the hypothesis that $\oldPi$ is essential.
\EndProof

\Number\label{standard torifold}
Let $q\ge1$ be an integer. There is a action of the cyclic group $\langle x\,|\,x^q=1\rangle$ on the solid torus $D^2\times S^1 $ defined by $x\cdot (z,w)=(z,e^{2\pi i/q}w)$. This action is at once smooth in the standard smooth structure on $D^2\times S^1$, and PL with respect to have the standard product PL structure on $D^2\times S^1$ (see \ref{esso}). The quotient of $D^2\times S^1 $ by this action, which inherits both a smooth and a PL structure from $D^2\times S^1$, will be denoted $\frakJ_q$., Up to a homeomorphism which is at once smooth and PL, we may identify $|\frakJ_q|$ with a solid torus in such a way that $\fraks_{\frakJ_q}$ is empty if $q=1$, and is a core curve of the solid torus $|\frakJ_q|$, having order $q$, if $q>1$. 

There is also an action, which is again at once smooth and P)L, of the dihedral group $\langle x,t\,|\,x^q=1,t^2=1,txt=x^{-1}\rangle$ on $D^2\times S^1 $ will be defined by $x\cdot (z,w)=(z,e^{2\pi i/q}w)$ and $t\cdot(z,w)=(\overline{z},\overline{ w})$, where the bars denote complex conjugation.  
We will let $\frakJ'_q$ denote the quotient of $D^2\times S^1 $ by this action, which again inherits both a smooth and a PL structure from $D^2\times S^1$. We may regard $\frakJ'_q$ as the quotient of $\frakJ_q$ by the involution $(z,w)\mapsto(\overline{z},\overline{w})$, where bars denote complex conjugation in $D^2\subset\CC$ or $S^1\subset\CC$. 
Furthermore,
there is a PL homeomorphism $h:|\frakJ'_q|\to B^3$ such that (i) $h(\fraks_{\frakJ'_q})$ contains two parallel line segments $\ell_q^1,\ell_q^2$ whose endpoints are in $S^2=\partial B^3$, (ii) $s_q:=\overline {h(\fraks_{\frakJ'_q})-(\ell_q^1\cup\ell_q^2)}$ is either the empty set or a line segment having one endpoint in $\inter\ell_q^i$ for each $i\in\{1,2\}$, and (iii) for $i=1,2$, each point of $h^{-1}(\ell_q^i\setminus s_q)$ has order $2$.

We will say that an action of a finite group $G$ on $D^2\times S^1 $ is {\it standard} if either $G$ is a cyclic group of some order $q\ge1$ and the action is the first one described above, or  $G$ is a dihedral group of order $2q$ for some $q\ge1$ and the action is the second one described above. 
\EndNumber

\Lemma\label{prepre}
Let $\oldPsi$ be an irreducible, orientable $3$-orbifold, 
and let $\oldTheta\subset\inter\oldPsi$  be a two-sided toric
$2$-suborbifold of $\oldPsi$ which is not $\pi_1$-injective.
Then either (a) $\oldTheta$ bounds a $3$-dimensional
suborbifold  of $\inter\oldPsi$ which is homeomorphic to the quotient of $D^2\times S^1 $ by the standard action of a finite group, 
or (b)
$\oldTheta$ is contained in the interior of a $3$-dimensional discal suborbifold
 of
$\inter \oldPsi$.
\EndLemma

\Proof
Set $N=|\oldPsi|$ and $T=|\oldTheta|$.
Since 
 $\oldTheta$ is not
  $\pi_1$-injective,
it follows from Proposition \ref{kinda dumb} that there is a discal $2$-suborbifold $\frakD$ of $\oldPsi$, with $\partial\frakD\subset\oldTheta$, such that
$D:=|\frakD|$ is in general position with respect to $\fraks_\oldPsi$, and $\partial\frakD$ does not bound a discal suborbifold of $\oldTheta$. This means that
  $D\subset N$ is a
disk
such that (1) $D\cap T=\partial D$, 
(2) $D$ 
is
in general position with respect
to
 $\fraks_\oldPsi$
and
meets it
 in at most one point, and
(3) 
any disk in $T$  bounded by
$\gamma:=\partial D$ must meet $\fraks_\oldPsi$
 in at least two points. 
If $D\cap\fraks_\oldPsi$ consists of a single point, let $q$ denote
the order of this point in $\fraks_\oldPsi$; if $D\cap\fraks_\oldPsi=\emptyset$, set $q=1$. 

We claim:
\Claim\label{well, in that case}
Either (i) $T$ is a torus, $\wt T=0$  (see \ref{wuzza weight}), and  $\gamma$ is a non-separating curve in $T$, or (ii) $T$ is a sphere, $\wt T=4$, each point of $\fraks_\oldPsi$ has order $2$, and  $\gamma$ separates $T$ into weight-$2$ disks.
\EndClaim

To prove \ref{well, in that case}, set $n=\wt\oldTheta$ and $\fraks_\oldTheta=\{x_1,\ldots,x_n\}$, and let $p_i$ denote the order of the singular point $x_i$ for $i=1,\ldots,n$.
Since $\oldTheta$ is toric and orientable, it follows from \ref{2-dim case} that 
\Equation\label{hey, i don't make the rules}
0= \chi(\oldTheta)=\chi(T)-\sum_{i=1}^n(1-1/p_i). 
\EndEquation
In particular we have $\chi(\oldTheta)\ge0$, so that $T$ is a torus or a sphere. If $T$ is a torus, it follows from (\ref{hey, i don't make the rules}) that $n=0$; Condition (3) above then implies that $\gamma$ is a non-separating curve in $T$, and Alternative (i) of \ref{well, in that case} holds. If $T$ is a sphere, then $\gamma$ separates $T$ into two disks $\Delta_1$ and $\Delta_2$; Condition (3) implies that each $\Delta_i$ has weight at least $2$. Hence $n=\wt(\Delta_1)+\wt(\Delta_2)\ge4$. Since $\chi(T)=2$, and each term $1-1/p_i$ in (\ref{hey, i don't make the rules}) is at least $1/2$, it now follows that $n=4$ and that $1-1/p_i=1/2$, i.e. $p_i=2$, for each $i$. Hence $\wt \Delta_1=\wt \Delta_2=2$, and Alternative (ii) holds. This completes the proof of \ref{well, in that case}.


Since $D\cap T=\partial D=\gamma$, there is a ball $E$, with $D\subset E\subset N$, such that
$E\cap T$ is an annulus $A$, which is contained in $\partial E$ and is a 
regular neighborhood of $\gamma $ in $T$. 
We may choose $E$ in such a way that there exists a homeomorphism $\eta:D\times[-1,1]\to E$ such that $\eta(D\times\{0\})=D$, $\eta((\partial D)\times[-1,1])=A$, and $\eta((D\cap\fraks_\oldPsi)\times[-1,1])=E\cap\fraks_\oldPsi$. Set $D_i=\eta(D\times\{i\})\subset\partial E$ for $i=\pm1$. 
Set 
$F=\overline{T-A}\cup\overline{(\partial E)-A}$.

If Alternative (i) of \ref{well, in that case} holds, $F$ is a
  $2$-sphere.
Furthermore, in this case
  $\fraks_{\obd(F)}=F\cap\fraks_\oldPsi$ is empty if $q=1$, and consists of two points of order $q$ if $q>1$; furthermore, in the latter case, one point of $\fraks_{\obd(F)}$ is contained in $D_{1}$, and one in $D_{-1}$.

If Alternative (ii) of \ref{well, in that case} holds,
$F$ is a disjoint union of two
$2$-spheres $S_1\supset D_1$ and $S_{-1}\supset D_{-1}$.
In this case, for $i=\pm1$, the intersection $S_i\cap\fraks_\oldPsi$ consists of two points of order $2$ in $S_i-D_i$, if $q=1$; and if $q>1$ it consists of three points, two of which lie in $S_i-D_i$ and have order $2$, and one of which lies in $\inter D_i$ and has order $q$.

Hence in any event, for each component $S$ of $F$, the $2$-orbifold $\omega(S)$ is spherical. Since $\oldPsi$ is irreducible, $\omega(S)$ bounds a $3$-dimensional discal suborbifold of $\oldPsi$. For each component $S$ of $F$ we choose a discal  $3$-suborbifold $\oldUpsilon_S$ of $\oldPsi$ with $\partial\oldUpsilon_S=\omega(S)$. In particular, for each component $S$ of $F$, the manifold $|\oldUpsilon_S|$ is a $3$-ball, and the group
 $\pi_1(\oldUpsilon_S)$ is finite. 
The closures of the components of $N-S$ are $|\oldUpsilon_S|$ and
$\overline{N-|\oldUpsilon_S|}$. 

For each component $S$ of $F$, set $Y_S={(T\cup E)-S}$.
The construction gives that $Y_S=A\cup(\inter E)\cup(F-S)$; and that $F-S$ is empty if Alternative (i) of \ref{well, in that case} holds, and is a component of $F$ if Alternative (ii) of \ref{well, in that case} holds. In either case it follows that $Y_S$ is connected. Since $Y_S$ is disjoint from $S$, it must be
contained either in $|\oldUpsilon_S|$ or in
$\overline{N-|\oldUpsilon_S|}$.  Hence $T\cup E=Y_S\cup S$ is contained in either $|\oldUpsilon_S|$
or in
$\overline{N-|\oldUpsilon_S|}$. Since one of these alternatives holds for every component of $S$, we have either
$T\cup E\subset\overline{N-|\oldUpsilon_S|}$ for every component $S$ of $F$, or
$T\cup E\subset|\oldUpsilon_S|$ for some component $S$ of $F$; that is:
\Claim\label{la nouvelle monique}
Either (A) $(T\cup E)\cap|\oldUpsilon_S|=S$ for every component $S$ of $F$, or (B) 
$T\cup E\subset|\oldUpsilon_{S_0}|$ for some component ${S_0}$ of $F$.
\EndClaim

Consider the case in which Alternative (B) of \ref{la nouvelle monique}
holds. Define $\frakB$ to be a strong regular neighborhood of
$\oldUpsilon_{S_0}$ in $\oldPsi$.
Then $\frakB$ is (orbifold-)homeomorphic to $\oldUpsilon_{S_0}$ and
is therefore discal. We have $\oldTheta\subset\oldUpsilon_{S_0}\subset\inter
\frakB$. This gives Alternative (b) of the conclusion of the
proposition.

The rest of the proof will be devoted to the case in which Alternative (A) of \ref{la nouvelle monique} holds; in this case, we will show that Alternative (a) of the conclusion of the proposition holds. 

Set $Z=\bigcup_{S\in\calc(F)}|\oldUpsilon_S|$ and
$J=E\cup Z$.
It follows from Alternative (A) of \ref{la nouvelle monique} that $J$ is obtained, up to homeomorphism, from $T\cup E$ by gluing a ball to each component of $F$ along the boundary of the ball. Hence $J$ is a manifold whose boundary is $T$. 
In view of the existence of the homeomorphism $\eta:D\times[-1,1]\to E$  with the properties stated above, we can now deduce that $\frakJ:=\obd(J)$ is homeomorphic to an orbifold obtained from $\frakZ:=\obd(Z)$ by gluing together the suborbifolds $\oldDelta_1:=\omega(D_{1})$ and $\oldDelta_{-1}:=\omega(D_{-1})$ of $\partial\oldUpsilon$ by a homeomorphism. It now suffices to prove that $\frakJ:=\omega(J)$ is (orbifold-)homeomorphic to  the quotient of $D^2\times S^1 $  by the standard action of a finite group.

Consider first the subcase in which Alternative (i) of \ref{well, in that case} holds. In this subcase, we have seen that $F$ is a single $2$-sphere; and that $\fraks_{\obd(F)}$ is empty if $q=1$, and consists of two points of order $q$ if $q>1$. Furthermore, we have seen that in the latter case, one point of $\fraks_{\obd(F)}$ is contained in $D_1$, and one in $D_{-1}$. Hence $\frakZ$ is a single discal $3$-orbifold; $Z$ is a ball; and $\fraks_\frakZ$ is empty if $q=1$, and is a single unkotted arc of order $q$, with one endpoint in $D_1$ and one in $D_{-1}$, if $q>1$. It follows that $\frakJ$ is (orbifold-)homeomorphic to the orbifold $\frakJ_q$ described in \ref{standard torifold}, and is therefore homeomorphic to the quotient of $D^2\times S^1 $ by a standard action of a cyclic group of order $q$. Thus Alternative (a) of the conclusion holds.

Now consider the subcase in which Alternative (ii) of \ref{well, in that case} holds.  
In this case we have seen that $F$ is a disjoint union of two
$2$-spheres $S_1\supset D_1$ and $S_{-1}\supset D_{-1}$.
Furthermore, we have seen that if $q=1$ then for $i=\pm1$ we have $S_i\cap\fraks_\oldPsi=\{x_i,x_i'\}$ for some points $x_i,x_i'\in S_i-D_i$ having order $2$; and that if $q>1$ we have $S_i\cap\fraks_\oldPsi=\{x_i,x_i',x_i''\}$ for some points $x_i,x_i'\in S_i-D_i$ having order $2$, and some point $x_i''\in\inter D_i$ having order $q$.
Hence $\frakZ$ is a disjoint union of two discal $3$-orbifolds $\oldUpsilon_1$ and $\oldUpsilon_{-1}$. Here $Y_i:=|\oldUpsilon_i|$ is a ball for $i=\pm1$. If $q=1$ then $\fraks_{\oldUpsilon_i}$ is an unknotted arc in $Z_i$, having order $q$, and having endpoints $x_i$ and $x_i'$. If $q>1$ then $T_i:=\fraks_{\oldUpsilon_i}$ is a cone on $\{x_i,x_i',x_i''\}$; moreover, $T_i$ is unknotted in the sense that it is contained in a properly embedded disk in $Z_i$. The interiors of the arcs joining $x_i$, $x_i'$ and $x_i''$ to the cone point have orders $2$, $2$ and $q$ respectively. 

It follows that $J=|\frakJ|$ is (orbifold-)homeomorphic to the orbifold $\frakJ_q'$ described in \ref{standard torifold}, and is therefore homeomorphic to the quotient of $D^2\times S^1 $  by a standard action of a dihedral group of order $2q$. Thus Alternative (a) of the conclusion holds.
\EndProof\


\Proposition\label{three-way equivalence}
Let $\oldLambda$ be an orientable $3$-orbifold. The following conditions are mutually equivalent:
\begin{enumerate}
\item $\oldLambda$ is a \torifold;
\item $\oldLambda$ has a toric boundary component and is  strongly \simple;
\item $\oldLambda$ is homeomorphic to the quotient of $D^2\times S^1 $ by the standard action of a finite group.
\end{enumerate}
\EndProposition

\Proof
It is trivial that (3) implies (1).

If (1) holds then $\oldLambda$ is finitely covered by a solid torus $J$. Some finite-sheeted covering of $J$ is a regular covering of $\oldLambda$, and is a solid torus since $J$ is one; hence we may assume without loss of generality that $J$ is a regular covering of $\oldLambda$. Since $\partial J$ is a single torus, and the inclusion homomorphism $\pi_1(\partial J)\to\pi_1(J)$ has infinite kernel, $\partial\oldLambda$ is a single toric orbifold. To show that (2) holds, it now suffices to prove that $\oldLambda$ is strongly \simple. Since
the solid torus $J$ is an irreducible $3$-manifold, Condition (II) of Definition \ref{oops} holds with $\oldLambda$ playing the role of $\oldPsi$. Since $J$ is compact, so is $\oldLambda$; and since $\pi_1(J)$ is infinite cyclic, $\pi_1(\oldLambda)$ is infinite, so that $\oldLambda$ is non-discal, and $\pi_1(\oldLambda)$ has no rank-$2$ free abelian subgroup. Hence Conditions (I) and (III) of Definition \ref{oops} hold as well. This
completes the proof that (1) implies (2).

Now suppose that (2) holds, and let $\oldTheta$ be a toric boundary component of $\oldLambda$. Let $\oldLambda'$ denote the $3$-orbifold obtained from the disjoint union  $\oldLambda\discup((\partial\oldLambda)\times[0,1])$ by gluing $\partial\oldLambda\subset\oldLambda$ to $(\partial\oldLambda)\times\{0\}\subset(\partial\oldLambda)\times[0,1]$ via the homeomorphism $x\mapsto(x,1)$. Then $\oldLambda'$ is homeomorphic to $\oldLambda$, and is therefore strongly \simple. In particular $\oldLambda'$ is irreducible, and the toric $2$-orbifold $\oldTheta$ is not $\pi_1$-injective. Thus  the hypotheses of Lemma \ref{prepre} hold with $\oldLambda'$ 
playing the role of $\oldPsi$, and with $\oldTheta$ chosen as above. Hence either (a) $\partial\oldLambda$ bounds a $3$-dimensional
suborbifold   $\oldLambda''$  of $\inter\oldLambda'$ which is homeomorphic to the quotient of  $D^2\times S^1 $  by the standard action of a finite group, 
or (b)
$\oldTheta$ is contained in the interior of a $3$-dimensional discal suborbifold $\frakB$
 of
$\inter N$. If (a) holds, we must have either $\oldLambda''=\oldLambda$ or $\oldLambda''=\oldTheta\times[0,1]$; the latter alternative is impossible, because the solid toric $3$-orbifold $\oldLambda''$ has connected boundary. Hence $\oldLambda''=\oldLambda$, so that (3) holds. If (b) holds, then $\pi_1(\frakB)$ is finite. But the construction of $\oldLambda'$ implies that $\oldLambda$ is  $\pi_1$-injective in $\oldLambda'$, and hence $\pi_1(\oldLambda)$ is finite. But if $\toldLambda$ is a finite-sheeted manifold covering of $\oldLambda$ (which exists by Condition (II) of Definition \ref{oops}), then $\partial\toldLambda$ has a torus component since $\oldTheta$ is toric; it follows that $H_1(\toldLambda;\QQ)\ne0$, which contradicts the finiteness of $\pi_1(\oldLambda)$. Hence (a) cannot occur, and we have shown that (2) implies (3).
\EndProof

\Proposition\label{preoccupani}
Let $\oldPsi$ be a componentwise strongly \simple, orientable $3$-orbifold, and let $\oldTheta\subset\inter\oldPsi$  be a toric
$2$-suborbifold. Then either (a) $\oldTheta$ bounds a \torifold\ contained in $\inter\oldPsi$, or (b)
$\oldTheta$ is contained in the interior of a discal $3$-suborbifold of
$\inter \oldPsi$.\abstractcomment{The rest of this statement has been removed, and can be found in irreducible1.tex.}
\EndProposition

\Proof
We may assume without loss of generality that $\oldPsi$ is connected, and is therefore strongly \simple. In particular $\oldLambda'$ is irreducible, and the toric $2$-orbifold $\oldTheta$ is not $\pi_1$-injective. The conclusion is now an immediate consequence of Lemma \ref{prepre} and (the trivial part of) Proposition \ref{three-way equivalence}. \EndProof

\Corollary\label{preoccucorollary}
Let $\oldPsi$ be a componentwise strongly \simple, orientable $3$-orbifold, and let $\oldTheta\subset\oldPsi$  be a toric
$2$-suborbifold. Suppose that the image of the inclusion homomorphism $\pi_1(\oldTheta)\to\pi_1(\inter \oldPsi)$ is infinite. Then $\oldTheta$ bounds a \torifold\ contained in $\oldPsi$.
\EndCorollary

\Proof
After modifying $\oldTheta$ by a small non-ambient orbifold isotopy, we may assume $\oldTheta\subset\inter\oldPsi$. Then $\oldPsi$ and $\oldTheta$ satisfy the hypothesis of Proposition \ref{preoccupani}; but since  of the inclusion homomorphism $\pi_1(\oldTheta)\to\pi_1(\oldPsi)$ has infinite image, Alternative (b) of the conclusion of Proposition \ref{preoccupani} cannot hold. Hence Alternative (a) must hold.

This follows upon applying Proposition \ref{preoccupani} to a toric suborbifold which is contained in $\inter\oldPsi$ and is obtained from $\oldTheta$ by a small non-ambient orbifold isotopy. 
\EndProof

\tinymissingref{\tiny In this version I am reworking Proposition \ref{at least a
    sixth}, which I think was screwed up. For the old version, see
  eclipse5.tex. I fixed the statement and app. in eclipse6.tex. Here I
will try to add a lemma to make the proof go smoothly.}


\Proposition\label{silver irreducible}
Let $\oldPsi$ be a compact, orientable $3$-orbifold, and let $\oldXi$ be a
$\pi_1$-injective $2$-suborbifold of $\partial\oldPsi$. Then 
$\silv_{\oldXi}\oldPsi$ (see \ref{silvering}) is irreducible if and only if (a)
$\oldPsi$ is irreducible and (b) $\oldXi$ is $\pi_1$-injective in
$\oldPsi$.
\EndProposition

\Proof
First suppose that (a) and (b) hold. If $\oldGamma$
is a two-sided spherical $2$-suborbifold of
$\silv_{\oldXi}\oldPsi$, then since $\oldGamma$ is two-sided and closed, we may
write $\oldGamma=\silv\frakD$ for some two-sided properly embedded $2$-suborbifold
$\frakD$ of $\oldPsi$ with $\partial\frakD
\subset\oldXi$. We have
$\chi(\frakD)=\chi(\oldGamma)>0$, so that $\frakD$ is either spherical
or discal. If $\frakD$ is spherical then $\frakD\subset\inter\oldPsi$;
 since $\oldPsi$ is irreducible, $\frakD$ bounds a discal
 $3$-suborbifold in $\inter\oldPsi$, which may be identified with a
 discal suborbifold of $\silv_{\oldXi}\oldPsi$ bounded by
 $\oldGamma$. Now suppose that $\frakD$ is discal. Since $\oldXi$ is $\pi_1$-injective in $\oldPsi$, it follows from Proposition \ref{kinda dumb} (more specifically the implication (a)$\Rightarrow$(b))
that there is a discal suborbifold $\frakE$ of
$\oldXi$ with $\partial\frakE=\partial\frakD$. Thus $\frakE\cup\frakD$ is a spherical
suborbifold of $\oldPsi$, non-ambiently isotopic  to a two-sided suborbifold of $\inter\oldPsi$. Since $\oldPsi$ is irreducible, $\frakE\cup\frakD$
bounds a discal $3$-suborbifold $\frakH$ of $\oldPsi$. Then
$\frakK:=\silv_\frakE\frakH$ is a suborbifold of
$\silv_{\oldXi}\oldPsi$. Since $\frakH$, $\frakE$ and $\frakD=\overline{(\partial\frakH)-\frakE}$ are discal and orientable, the pair $(\frakH,\frakE)$ is homeomorphic to $(\frakE\times[0,1],\frakE\times\{0\})$ (where $|\frak|$ is a disk of weight at most $1$), and hence $\frakK$ is discal. We have $\partial\frakK=\oldGamma$.
Thus in either case $\oldGamma$ is the
boundary of a discal $3$-suborbifold of $\silv_{\oldXi}\oldPsi$; this proves that $\silv_{\oldXi}\oldPsi$ is irreducible.

Conversely, suppose that $\silv_{\oldXi}\oldPsi$ is irreducible. If
$\oldPi$ is a two-sided spherical $2$-suborbifold of $\inter\oldPsi$, then
$\oldPi$ may be identified with a two-sided spherical $2$-suborbifold of
$\silv_{\oldXi}\oldPsi$, which bounds a discal $3$-suborbifold
$\frakK$ of $\silv_{\oldXi}\oldPsi$ by irreducibility; we may identify
$\frakK$ with a discal $3$-suborbifold  of $\oldPsi$ bounded
by $\oldPi$. This establishes (a). To prove (b), according to Proposition \ref{kinda dumb}, it suffices to show that if $\frakD$ is a discal $2$-orbifold, properly embedded in $\oldPsi$, such that  $C:=\partial \frakD\subset\oldXi$, and $|\frakD$ is in general position with respect to $\fraks_\oldPsi$, then $C$ bounds a discal suborbifold of $\oldXi$. If we are given such a suborbifold $\frakD$, the general position property of $|\frakD|$ implies $\dim\fraks_\frakD\le\dim\fraks_\oldPsi-1$; but $\dim\fraks_\oldPsi\le1$ by the orientability of $\oldPsi$, and hence $\fraks_\frakD$ is a finite set. This implies that $\frakD$ is orientable. Hence $\silv\frakD$ is a two-sided spherical $2$-suborbifold of 
$\silv_{\oldXi}\oldPsi$, which bounds a discal $3$-suborbifold
$\frakK$ of $\silv_{\oldXi}\oldPsi$ by irreducibility. We may write
$|\frakK|=|\frakH|$ for some suborbifold $\frakH$ of $\oldPsi$, and
$\Fr_\oldPsi\frakH=\frakD$. Since $\frakK$ and $\frakD$ are discal,
$|\frakK|=|\frakH|$ is a $3$-ball and $|\frakD$ is a $2$-disk. Hence $(\partial|\frakH|)-\inter|\frakD|$ is a disk, and in particular
$\frakE:=\frakH\cap\oldXi=(\partial\frakH)-\inter\frakD$ is a
connected $2$-orbifold with boundary $C$. 
But $\frakK$ is canonically identified with $\silv_\frakE\frakH$, and hence $\frakE$ is the unique two-dimensional stratum of $\frakK$. Since $\frakK$ is discal, 
there exists a(n orbifold) covering map $p:D^3\to\frakK$ such
that every component of $p^{-1}(\frakE)$ is the fixed point set of
an orientation-reversing orthogonal involution of $D^3$. 
Since a $2$-dimensional fixed point set of
an orthogonal involution of $D^3$ is always a disk, $\frakE$ is discal, and (b) is established.
\EndProof

\Proposition\label{silver acylindrical}
Let $\oldLambda$ 
be an orientable, componentwise strongly \simple\ $3$-orbifold, and let $\oldXi$ be a 
$\pi_1$-injective $2$-suborbifold of $\partial\oldLambda$. Then
an orientable annular $2$-orbifold $\oldPi\subset\oldLambda$, with $\oldPi\cap\partial\oldLambda=\partial\oldPi\subset \oldXi$,
is
  essential in the pair $(\oldLambda,\oldXi)$ if and only if the  toric orbifold $\silv\oldPi$, which is canonically identified with a suborbifold of
$\silv_\oldXi\oldLambda$,
  is incompressible in
$\silv_\oldXi\oldLambda$ and is not parallel in $\silv_\oldXi\oldLambda$ (see \ref{parallel def}) to
a boundary component of $\silv_\oldXi\oldLambda$. \nonessentialproofreadingnote{ I think strong \simple ity is needed only for the {\it second} assertion. Of course on the face of it it costs nothing to make the assumption for the whole statement, but the question arises whether leaving it out could simplify the logic in the apps.} Furthermore,
$\silv_\oldXi\oldLambda$ is weakly \simple\ if and only if the pair $(\oldLambda,\oldXi)$
is acylindrical (see \ref{acylindrical def}).
\EndProposition

\Proof 
To prove the first assertion, first suppose that $\oldPi$
is
essential in the pair $(\oldLambda,\oldXi)$. To prove that $\silv\oldPi$ is incompressible (i.e. $\pi_1$-injective) in
$\silv_\oldXi\oldLambda$, we apply Proposition \ref{kinda dumb}, letting $\silv_\oldXi\oldLambda$ and $\silv\oldPi$ play the roles of $\oldPsi$ and $\oldTheta$. To prove $\pi_1$-injectivity, it suffices to show that Condition (c) of Proposition \ref{kinda dumb} holds. Thus we consider an arbitrary  discal orbifold $\frakD\subset
\silv_\oldXi\oldLambda$ such that (i) 
$\frakD\cap\silv\oldPi=\partial\frakD$, and (ii) $|\oldTheta|$ meets $\fraks_\oldPsi$ in general position in $|\oldPsi|$. We are required to show that $\partial\frakD$ is
the boundary of some discal suborbifold of $\silv\oldPi$. Consider first the case in which
$|\partial\frakD|\subset\inter|\silv\oldPi|=\inter|\oldPi|$. Then
$\frakD$ is identified with a discal suborbifold of $\oldLambda$ whose
boundary is contained in $\oldPi$. Since $\oldPi$ is essential, it is in particular $\pi_1$-injective. It therefore follows from Proposition \ref{kinda dumb}, this time applied with $\oldLambda$ and $\oldPi$ playing the roles of $\oldPsi$ and $\oldTheta$, that $\partial\frakD=\partial\frakE$ for some discal
suborbifold of $\oldPi$. Since
$|\partial\frakD|\subset\inter|\oldPi|$, we have $\frakE\subset\inter|\oldPi|$, and $\frakE$ may therefore be identified with
a discal suborbifold of $\silv\oldPi$; this gives the required conclusion in this case. Now consider the case in which $|\partial\frakD|$ meets
$\partial|\oldPi|\subset|\oldXi|\subset\fraks_{\silv_\oldXi\oldLambda}$.
In particular we then have $|\partial\frakD|\cap\fraks_\frakD\ne\emptyset$; since
$\frakD$ is discal, it follows that $|\frakD|$ is a disk and that
$\fraks_\frakD$ is an arc $\alpha\subset\partial|\frakD|$. 
Thus  $|\partial\frakD|$ is the arc $\beta:=
(\partial|\frakD|)-\inter\alpha$, and $|\frakD|$ is a weight-$0$ disk
meeting $\partial|\oldPsi|$ in $\alpha$ and meeting the orientable annular
orbifold $\oldPi$ in $\beta$. Now Condition (ii) above (with the definition of general position given in \ref{gen pos}) implies that $|\frakD|$ meets each one-dimensional stratum of $\fraks_\oldPsi$ only in isolated points. Hence the arc $\alpha=\fraks_\frakD\subset\fraks_\oldPsi$ must be contained in the union of the closures of the two-dimensional strata of $\oldPsi$; that is, $\alpha\subset\oldXi$. Since $\oldPi$ is essential, it now follows
from Proposition \ref{what essential means} that $\beta$ is the frontier of a weight-$0$ disk in $|\oldPi|$. 
Hence $\partial\frakD$ is the boundary of a discal suborbifold of $\silv\oldPi$ in this case as well, and we have  shown that  $\silv\oldPi$ is incompressible in
$\silv_\oldXi\oldLambda$. 

Now suppose that  $\silv\oldPi$ is parallel in $\silv_\oldXi\oldLambda$ to a boundary component
 of  $\silv_\oldXi\oldLambda$. Then by definition there exists an embedding $j:(\silv\oldPi)\times[0,1]\to\silv_\oldXi\oldLambda$ such that
$j((\silv\oldPi)\times\{0\})=\silv\oldPi$, and
$j((\silv\oldPi)\times\{1\})\subset\partial(\silv_\oldXi\oldLambda)$. Note that $|\partial\oldPi|\times[0,1]$ is a topological $2$-manifold contained in
$\fraks_{(\silv\oldPi)\times[0,1]}$; hence $|j|(|\partial\oldPi|\times[0,1])$ is contained in the union of the closures of the two-dimensional strata of $\silv_\oldXi\oldLambda$, which is $|\oldXi|$. On the other hand, $(\inter|\oldPi|)\times[0,1]$ is an open subset of $|\oldPi\times[0,1]|$ whose intersection with  $\fraks_{\oldPi\times[0,1]}$ is at most one-dimensional, and hence $|j|(\inter|\oldPi|)\times[0,1])\cap|\oldXi|=\emptyset$.  It follows that $|j|^{-1}(|\oldXi|)=
|\partial\oldPi|\times[0,1]$, and therefore that there is an embedding $j':\oldPi\times[0,1]\to\oldLambda$ such that $|j'|=|j|$. We then have
$j'((\partial\oldPi)\times[0,1])\subset\oldXi$ and
$j'(\oldPi\times\{0\})=\oldPi$. Furthermore, since
$|j|(|\oldPi|\times\{1\}) \subset|\partial(\silv_\oldXi\oldLambda)|=\overline{(|\partial\oldLambda|)-|\oldXi|}$, we have 
$j'(\oldPi\times\{1\})
\subset\overline{(\partial\oldLambda)-\oldXi}$. This shows that
$\oldPi$ is parallel in the pair $(\oldLambda,\oldXi)$ to a suborbifold of $\overline{(\partial\oldLambda)-\oldXi}$, a contradiction to the essentiality of $\oldPi$. Hence 
$\silv\oldPi$ is not parallel in $\silv_\oldXi\oldLambda$ to a boundary component
 of  $\silv_\oldXi\oldLambda$.

Conversely, suppose that
$\silv\oldPi$
  is incompressible in
$\silv_\oldXi\oldLambda$ and is not parallel to
a boundary component of $\silv_\oldXi\oldLambda$. Note that the natural  immersion $¯\alpha:\oldPi\to\silv\oldPi$ is $\pi_1$-injective, since it is the composition of the $\pi_1$-injective inclusion $\oldPi\to D(\oldPi)$ with the natural covering map $D(\oldPi)\to\silv\oldPi$. By incompressibility the inclusion map $i:\silv\oldPi\to\silv_\oldXi\oldLambda$ is also $\pi_1$-injective, and hence so is $i\circ\alpha:\oldPi\to\silv_\oldXi\oldLambda$. But we have $i\circ\alpha=\beta\circ h$, where $h:\oldPi\to\oldLambda$ is the inclusion and $\beta:\oldLambda\to\silv_\oldXi\oldLambda$ is the natural immersion. Hence $h$ is $\pi_1$-injective, i.e. $\oldPi$ is $\pi_1$-injective in $\oldLambda$.

To show that $\oldPi $ is essential in the pair $(\oldLambda,\oldXi)$, it remains to show that it is not 
 parallel in $(\oldLambda,\oldXi)$ either to a suborbifold of $\oldXi$ or to a component of
$\overline{(\partial\oldLambda)-\oldXi}$. Assume to the contrary that there is
an embedding
$\iota:\oldPi\times[0,1]\to\oldLambda$ such that 
$\iota((\partial\oldPi)\times[0,1])\subset\oldXi$,
$\iota(\oldPi)\times\{0\})=\oldPi$, and
$\iota(\oldPi\times\{1\})$ is contained either in $\oldXi$ or in $\overline{(\partial\oldLambda)-\oldXi}$.
If $\iota(\oldPi\times\{1\})\subset\oldXi$, and if we set $\frakJ=\iota(\oldPi\times[0,1])$ and $\frakB=\iota(\oldPi\times\{1\})$, then $D_\frakB\frakJ\subset D_\oldXi\oldLambda$ is a \torifold\ whose boundary is $D\oldPi$; this is impossible, because the incompressibility of $\silv\oldPi$ in $\silv_\oldXi\oldLambda$ implies that $D\oldPi$ is $\pi_1$-injective in $D_\oldXi\oldLambda$. 
Now consider the case in which
$\iota(\oldPi\times\{1\})\subset\overline{(\partial\oldLambda)-\oldXi}$. Then, since $\iota((\inter\oldPi)\times[0,1))\subset\inter\oldLambda$ and 
$\iota((\partial\oldPi)\times[0,1]\subset\oldXi$, we have $\iota^{-1}(\oldXi)=(\partial\oldPi)\times[0,1]$; hence there is an embedding $\iota':(\silv\oldPi)\times[0,1]\to\silv_\oldXi\oldLambda$ such that $|\iota'|=|\iota|$. We have
$\iota'(\silv\oldPi)\times\{0\})=\silv\oldPi$; and since $|\overline{(\partial\oldLambda)-\oldXi}|=|\partial\silv_\oldXi\oldLambda|$, we have
$\iota'(\oldPi\times\{1\})\subset
\partial\silv_\oldXi\oldLambda$. This shows that $\silv\oldPi$ is parallel in $\silv_\oldXi\oldLambda$ to a suborbifold of $\partial(\silv_\oldXi\oldLambda)$, which must be a component of $\partial(\silv_\oldXi\oldLambda)$ since $\silv\oldPi$ is closed; this is a contradiction. This completes the proof that $\oldPi$ is essential in $(\oldLambda,\oldXi)$, and the first assertion of the proposition is established.



To prove the second assertion, first suppose that
$\silv_\oldXi\oldLambda$ is weakly \simple. If $\oldPi$ is an orientable annular
orbifold which is essential in the pair$(\oldLambda,\oldXi)$, then the first
assertion implies that the 
toric orbifold
 $\silv\oldPi$
  is incompressible in
$\silv_\oldXi\oldLambda$ and is not parallel to
a boundary component of $\silv_\oldXi\oldLambda$; this contradicts the
weak
\simple ity of $\silv_\oldXi\oldLambda$. Hence $(\oldLambda,\oldXi)$
is acylindrical. 

Conversely, suppose that  $(\oldLambda,\oldXi)$
is acylindrical. Since $\oldLambda$ is strongly \simple, it is irreducible (see \ref{oops}); as $\oldXi$ is $\pi_1$-injective, it then follows from Proposition \ref{silver irreducible} that $\silv_\oldXi\oldLambda$ is
irreducible. Now suppose that $\oldTheta$ is an incompressible toric
suborbifold of $\silv_\oldXi\oldLambda$.
Since $\oldTheta$ is in
particular a closed suborbifold of $\silv_\oldXi\oldLambda$, we may
write $\oldTheta=\silv\oldPi$ for some two-sided, properly embedded suborbifold
$\oldPi$ of $\oldLambda$ with $\partial\oldPi\subset\oldXi$. We have
$\chi(\oldPi)=\chi(\oldTheta)=0$, so that $\oldPi$ is either toric or
annular. If $\oldPi$ is toric, then $\oldTheta=\oldPi$, and the
incompressibility of $\oldTheta$ in $\silv_\oldXi\oldLambda$ implies that
$\oldPi$ is incompressible in $\oldLambda$, a contradiction to the
strong \simple ity of $\oldLambda$ see \ref{oops}. Hence $\oldPi$ is annular. The acylindricity of the pair
$(\oldLambda,\oldXi)$  implies that $\oldPi$ is not
essential in the pair $(\oldLambda,\oldXi)$. It therefore follows from the first
assertion of the present proposition that $\oldTheta$ is  parallel to
a boundary component of $\silv_\oldXi\oldLambda$. This shows that $\silv_\oldXi\oldLambda$ is
weakly \simple.
\EndProof

\abstractcomment{\tiny This lemma is apparently not used. For introduction and proof summary, see cravat.tex.

\Lemma\label{graph lemma}
Let $\calt$ be a finite graph which has no isolated vertices, su y ch that each component of $\calt$ is a tree. Let $\calv$ and $\cale$ denote respectively the vertex set and edge  set of $\calt$. Let $\calv_0$ be a subset of $\calv$ such that (a) each component of $\calt$ contains at most one vertex in $\calv_0$, and (b) each vertex in $\calv_0$ has valence $1$ in $\calt$. Let $\cale_0$ be a subset of $\cale$ such that each endpoint of each  edge in $\cale$ either has valence at least $3$ or belongs to $\calv_0$. Then $\card \cale_0<(\card \cale)/2$.
\EndLemma
}


\abstractcomment{\tiny What else?}

\section{Fibrations of orbifolds}\label{fibration section}

In this section we 
will  define  orbifold fibrations and establish some of their basic properties. These will be used in developing the theory of the characteristic suborbifold in Section \ref{characteristic section}.

\DefinitionsRemarks\label{fibered stuff} (Cf.
\cite{bonahon-siebenmann}, \cite{other-bs}.) The following definitions and remarks are meant to be interpretable in the topological, PL and smooth categories.

Suppose that
$\oldLambda$ is a compact $n$-orbifold (where $n$ will be $2$ or $3$
in all the applications) and that $\frakB$ is a compact
$(n-1)$-orbifold. By a(n {\it orbifold) fibration}
$q:\oldLambda\to\frakB$ (or a(n {\it orbifold) fibration} of
$\oldLambda$ with base $\frakB$) we mean an orbifold submersion
$q:\oldLambda\to\frakB$ having the following property: for every
$x\in\frakB$,
there exist a neighborhood $\frakV$ of $v$ in
$\oldGamma$ and a(n orbifold) homeomorphism $\eta:
(D^2\times
S^1)/G
\to
q^{-1}(\frakV)
$, where 
(1)  $G$ is a subgroup of $G_1\times G_2$, for some 
subgroups $G_1$ 
and 
$G_2$  of ${\rm SO}(2)$,
and 
(2) 
if $\pi:
D^2\times
S^1
\to
(D^2\times
S^1)/G$
denotes the orbit map, then
$\eta(\pi(\{x\}\times S^1))$ 
is a fiber of $q$
for each $x\in D^2$,
and
$\eta(\pi(\{0\}\times S^1))=
r^{-1}(v)$. 
(In the PL case, we recall that the actions of $G_1$ and $G_2$ on $D^2$  and $S^1$ are PL by \ref{esso}.)

We shall call $q$ an {\it $S^1$-fibration} or an {\it
  $I$-fibration}, respectively, if we can take $J=S^1$---or,
respectively, $J=[0,1]$---for every $x\in\frakB$. (An argument similar to the connectedness argument in the proof of  Proposition \ref{new what they look like} could be used to show that every orbifold 
fibration with a connected base is either an $S^1$-fibration or an
$I$-fibration. This fact will not be needed.) 

If $q:\oldLambda\to\frakB$ is an orbifold fibration, then for every $x\in|\frakB|$, the {\it fiber over
  $x$ } is the suborbifold $q^{-1}(x)$ of $\oldLambda$;
the fiber over $x$ in an  $S^1$-fibration or
an $I$-fibration is, respectively, (orbifold-)homeomorphic
to a quotient of $S^1$ by a finite group action (and
therefore to $S^1$ or $[[0,1]]$) or to a quotient of $[0,1]$ by a
finite group action (and
therefore to $[0,1]$ or $[[0,1]$).

An $I$-fibration $q:\oldLambda\to\frakB$ will be termed {\it trivial} if there exists a homeomorphism $t:\frakB\times[0,1]\to\oldLambda$ such that $q\circ t:\frakB\times[0,1]\to\frakB$ is the projection to the first factor. 

If a compact, connected $2$-orbifold $\oldLambda$ admits an orbifold
fibration whose base is a $1$-orbifold, then
$\chi(\oldLambda)=0$. Hence if $\oldLambda$ admits an $I$-fibration,
then $\oldLambda$ is an annular orbifold (see \ref{annular etc}).

If $\oldLambda$ is a compact, orientable $3$-orbifold equipped with  an $I$-fibration $q:\oldLambda\to\frakB$, where $\frakB$ is a $2$-orbifold, 
we define the {\it vertical boundary} $\partial_v\oldLambda$ to be the
suborbifold $q^{-1}(\partial\frakB)$, of $\oldLambda$, and we define the {\it horizontal boundary}
$\partial_h\oldLambda$ to be the suborbifold of $\oldLambda$ which is the union (see \ref{unions and such} of the (orbifold-)boundaries of
all the fibers. 
We have
$\partial\oldLambda=\partial_h\oldLambda\cup\partial_v\oldLambda$, and
$\partial(\partial_h\oldLambda)=\partial(\partial_v\oldLambda)=\partial_h\oldLambda\cap\partial_v\oldLambda$. If
$\frakB$ (or equivalently $\oldLambda$) is connected, then
$\partial_h\oldLambda$ has at most two components, because each
component contains an endpoint of each fiber. 

Note that $q$ restricts to an $I$-fibration of $\partial_v\oldLambda$ whose base is $\partial\frakB$; hence $\partial_v\oldLambda$ is annular. Note also that $q|\partial_h\oldLambda:\partial_h\oldLambda\to\frakB$ is a degree-$2$ orbifold covering.
\EndDefinitionsRemarks

\Proposition\label{fibration-category}Let $q:\oldLambda\to\frakB$ be a smooth $S^1$-fibration \abstractcomment{Do I need this for $I$-fibrations too? I don't think so} of a smooth, orientable, compact, very good $3$-orbifold $\oldLambda$ over a smooth, compact, very good $2$-orbifold $\frakB$. Let $|\frakB|$ be equipped with a distance function that defines its topology, and let $\epsilon>0$ be given. Then there exist PL structures on $\oldLambda$ and $\frakB$, compatible with their smooth structures, and a sequence $(q_n')_{n\ge1}$, where $q_n':\oldLambda\to\frakB$ is a PL fibration for each $n\ge1$, such that $(|q'|)_{n\ge1}$ converges uniformly to $|q|$, and  for each compact PL subset $K$ of $|q'_n|^{-1}(|\frakB|-\fraks_\frakB)$, the sequence $(q_n'|K)_{n\ge1}$, converges to $q|K$ in the $C^1$ sense.
\EndProposition

\Proof
Since $q$ is an $S^1$-fibration and $\oldLambda$ and $\frakB$ are very good, there exist finite-sheeted regular coverings $p_\oldLambda:\toldLambda\to\oldLambda$ and $p_\frakB:\tfrakB\to\frakB$ such that $\toldLambda$ and $\tfrakB$ are connected manifolds, and a locally trivial fibration $\tq:\toldLambda\to\frakB$ whose fiber is a (possibly disconnected) closed $1$-manifold, such that $p_\frakB\circ\tq=q\circ p_\oldLambda$. Since the fiber of $\tq$ is closed, $\tq$ is a boundary-preserving map. If $G^\oldLambda$ and $G^\frakB$ denote the covering groups of the coverings $p_\oldLambda$ and $p_\frakB$ respectively, $\tq$ is equivariant in the sense that there is a homomorphism $\rho:G^\oldLambda\to G^\frakB$ such that $\tq\circ g=\rho(g)\circ\tq$ for every $g\in G^\oldLambda$. It follows from the main result of \cite{illman} that $\toldLambda$ and $\tfrakB$ have triangulations, compatible with their smooth structures, and invariant under the actions of $G^\oldLambda$ and $G^\frakB$ respectively. For the rest of the proof, the manifolds $\toldLambda$ and $\tfrakB$ will be understood to be equipped with the PL structures defined by these triangulations; note that the orbifolds $\oldLambda$ and $\frakB$ inherit PL structures, compatible with their smooth structures, and that $p_\oldLambda$ and $p_\frakB$ are then PL maps. After replacing the triangulation of $\toldLambda$ by its first barycentric subdivision, we may assume:
\Claim\label{stable is fixed}
Each simplex of $\toldLambda$ is pointwise fixed by its stabilizer under $G^\oldLambda$.
\EndClaim

Let us fix a distance function $h$ on (the total space of) the tangent bundle $T\frakB$ which determines its topology.  For any subset $K$ of $\toldLambda$ which is a union of closed simplices, any piecewise smooth maps $f,g:K\to\tfrakB$, and any $\delta>0$, we will say that $f$ and $g$ are {\it $C^1$ $\delta$-close} if for every closed simplex $\sigma\subset K$, every point $x\in\sigma$ and every tangent vector $w$ to $\sigma$ at $x$, we have $h(d(f|\sigma)(w),d(g|\sigma)(w))<\delta$.

For $k=-1,0,1,2,3$, let $\toldLambda^{(k)}$ denote the union of all simplices of dimension at most $k$ in $\toldLambda$. Note that since $\toldLambda^{(k)}$ is $G^\oldLambda$-invariant, it makes sense to speak of equivariant maps from $\toldLambda$ to $\tfrakB$. By induction, for $-1\le k\le3$, we will show:

\Claim\label{by induction}
For any $\delta>0$ there is an equivariant PL map $r^{(k)}:\toldLambda^{(k)}\to\tfrakB$ such that (1) $r^{(k)}$ is $\delta$-close to $\tq|\toldLambda^{(k)}$, and (2) $r^{(k)}(\toldLambda^{(k)}\cap\partial\toldLambda)\subset\partial\tfrakB$.
\EndClaim

Since \ref{by induction} is trivial for $k=-1$, we need only show that if $k$ is given with $0\le k\le3$, and if \ref{by induction} is true with $k-1$ in place of $k$, then it is true for the given $k$. Let $\delta'$ be a positive number less than $\delta$, which for the moment will be otherwise arbitrary; we will impose a finite number of smallness conditions on $\delta'$ in the course of the proof of \ref{by induction}. Let $r^{(k-1)}:\toldLambda^{(k-1)}\to\tfrakB$ be an equivariant PL map such that $r^{(k-1)}$ is  $\delta'$-close to $\tq|\toldLambda^{(k-1)}$, and $r^{(k-1)}(\toldLambda^{(k-1)}\cap\partial\toldLambda)\subset\partial\tfrakB$. 

Fix a complete set of orbit representatives $\cald$ for the action of $G^\oldLambda$ on the set of open $k$-simplices of $\toldLambda$. For each $\Delta\in\cald$, the stabilizer $G^{\oldLambda}_\Delta$ fixes $\overline\Delta$ pointwise by \ref{stable is fixed}. Since $\tq$ and $r^{(k-1)}$ are equivariant, the sets $\tq(\overline{\Delta})$ and $r^{(k-1)}(\partial\overline{\Delta})$ are contained in $
\Fix(\rho(G^{\oldLambda}_\Delta))$. If $\Delta$ is an element of $\cald$ such that $\Delta\subset\partial\toldLambda$, then since $\tq$ is boundary-preserving, we have $\tq(\overline{\Delta})\subset\partial\tfrakB$; furthermore, in this case we have $
r^{(k-1)}(\partial\overline{\Delta})\subset
r^{(k-1)}(\toldLambda^{(k-1)}\cap\partial\toldLambda)\subset\partial\tfrakB$. Thus if we set
$E_\Delta=
\Fix(\rho(G^{\oldLambda}_\Delta))$ for each $\Delta\in\cald$ such that $\Delta\not\subset\partial\toldLambda$, 
and $E_\Delta=
(\partial\tfrakB)\cap\Fix(\rho(G^{\oldLambda}_\Delta))$ for each $\Delta\in\cald$ such that $\Delta\subset\partial\toldLambda$, then
$\tq(\overline{\Delta})$ and $r^{(k-1)}(\partial\overline{\Delta})$ are contained in $E_\Delta$ for each $\Delta\in\cald$. 
 Note that $\Fix(\rho(G^{\oldLambda}_\Delta))$ is a PL subset of $\tfrakB$ for each $\Delta$, since the action of $G^\frakB$ on $\tfrakB$ is piecewise linear; since $\partial\frakB$ is PL, it follows that  $E_\Delta$ is a PL subset of $\tfrakB$ for each $\Delta\in\cald$. Note also that 
$r^{(k-1)}|\partial\overline{\Delta}$ is  $\delta'$-close to 
$\tq|\partial\overline{\Delta}$. If we choose $\delta'$ sufficiently small, it follows that $r^{(k-1)}|\partial\overline{\Delta}$ may be extended to a PL map $\tq_\Delta:\overline\Delta\to E_\Delta\subset\tfrakB$, and that $\tq_\Delta$ may be taken to be arbitrarily $C^1$-close to $\tq|\overline{\Delta}$.

If $\Delta$ is any open  $k$-simplex of $\toldLambda$, we may choose $g\in G^\oldLambda$ so that $\Delta_0:=g^{-1}(\Delta)\in\cald$, and define  a PL map $\tq'_\Delta:\overline{\Delta}\to\tfrakB$ by
$\tq'_\Delta=\rho(g)\circ\tq'_{\Delta_0}\circ g^{-1}$. Since 
$\tq_{\Delta_0}(\overline{\Delta_0})\subset E_{\Delta_0}\subset \Fix(\rho(G^{\oldLambda}_{\Delta_0}))$, the map $\tq_\Delta$ does not depend on the choice of $g$, and this definition of $\tq_\Delta$ specializes to the earlier one when $\Delta\in\cald$. Note  that if $\Delta$ is any $k$-simplex contained in $\partial\toldLambda$, then choosing $g$ and defining $\Delta_0$ as above, we have $\Delta_0\subset\partial\toldLambda$, and hence $E_\Delta\subset
 \partial\frakB$; it follows that $\tq'_\Delta(\overline{\Delta})\subset\rho(g)(E_\Delta)\subset\partial\frakB$.
Note also that since, for $\Delta\in\cald$, we may take $\tq_\Delta$ to be arbitrarily  $C^1$-close to $\tq|\overline{\Delta}$ if $\delta'$ is small enough, we may choose $\delta'$ so that $\tq_\Delta$ is $C^1$ $\delta$-close to $\tq|\overline{\Delta}$ for each $k$-simplex $\Delta$ of $\toldLambda$.

We may extend $\tq_{k-1}$ to a PL map $\tq_k:\toldLambda^{(k)}\to\tfrakB$ by setting $\tq_k|\overline{\Delta}=\tq_{\Delta}$ for every
open $k$-simplex $\Delta$ of $\toldLambda$. This map is equivariant by construction. Since
$r^{(k-1)}(\toldLambda^{(k-1)}\cap\partial\toldLambda)\subset\partial\tfrakB$,  and
$\tq'_\Delta(\overline{\Delta})\subset\partial\frakB$ for each $k$-simplex $\Delta\subset\partial\toldLambda$, we have
$r^{(k)}(\toldLambda^{(k)}\cap\partial\toldLambda)\subset\partial\tfrakB$.
Since $\tq_\Delta$ is $C^1$ $\delta$-close to $\tq|\overline{\Delta}$ for each $k$-simplex $\Delta$ of $\toldLambda$, and 
 $r^{(k-1)}$ is $\delta$-close to $\tq|\toldLambda^{(k-1)}$, the map  $r^{(k)}$ is $\delta$-close to $\tq|\toldLambda^{(k)}$. This proves \ref{by induction}.

It follows from the case $k=3$ of \ref{by induction} that there is a sequence $(\tq'_n)_{n\ge1}$ of boundary-preserving equivariant PL maps from $\toldLambda$ to $\tfrakB$ which converges in the $C^1$ sense to $\tq$. It follows from equivariance that for each $n\ge1$ there is a unique PL map $q'_n:\oldLambda\to\frakB$ such that $p_\frakB\circ\tq'_n=q'_n\circ p_\oldLambda$. The $C^1$-convergence of $\tq'_n$ to $\tq$ implies in particular:
\Claim\label{when I first put this uniform on}
The sequence $(|q'_n|)$ converges uniformly to $|q|$, and  for each compact PL subset $K$ of $|q'_n|^{-1}(|\frakB|-\fraks_\frakB)$, the sequence $(q_n'|K)_{n\ge1}$, converges to $q|K$ in the $C^1$ sense.
\EndClaim
On the other hand, since the smooth locally trivial fibration $\tq$ is in particular a submersion, the $C^1$-convergence of $\tq'_n$ to $\tq$ implies that $\tq'_n$ is a PL submersion for sufficiently large $n$. Since $\tq'_n$ is also boundary-preserving, we deduce:
\Claim\label{veriest dunce}
For sufficiently large $n$, the map $\tq'_n:\toldLambda\to\tfrakB$ is a locally trivial fibration whose fiber is a (possibly disconnected) closed $1$-manifold.
\EndClaim

Set $B=|\frakB|-\fraks_\frakB$. Since $B$ is connected, it follows from \ref{veriest dunce} that the map $|q'_n|\big||q'_n|^{-1}(B):|q'_n|^{-1}(B)\to B$ is a locally trivial fibration for sufficiently large $n$, and that its fiber has the form $|\frakC|$ for some orbifold quotient of the fiber of $\tq'_n$; in particular, the fiber of $|q'_n|\big||q'_n|^{-1}(B)$ is a $1$-manifold (so that each of its components is homeomorphic to $S^1$ or to $[0,1]$). We claim:

\Claim\label{peripatetics}
For sufficiently large $n$, the fiber of $|q'_n|\big||q'_n|^{-1}(B):|q'_n|^{-1}(B)\to B$ is connected.
\EndClaim

To prove \ref{peripatetics}, choose a point $x\in B$, and choose closed disk neighborhoods $D_0$ and $D$ of $x$ in $B$, with $D_0\subset\inter D$. It follows from \ref{when I first put this uniform on} that for sufficient large $n$ we have $|q'_n|^{-1}(x)\subset |q|^{-1}(D_0)$, and
$|q'_n|^{-1}(D)\supset |q|^{-1}(D_0)$. But $|q|^{-1}(D_0)$ is connected since $q$ is an $S^1$-fibration. Hence $|q|^{-1}(D_0)$ is contained in a component  $K$ of
$|q'_n|^{-1}(D)$, and in particular we have $|q'_n|^{-1}(x)\subset K$. Thus $|q'_n|^{-1}(x)$ is the fiber of the locally trivial fibration $|q'_n|\big|K:K\to D$. Since the latter fibration has a connected total space and a simply connected base, it follows from the exact homotopy sequence of this fibration that its fiber $|q'_n|^{-1}(x)$ is connected. This proves \ref{peripatetics}.

Now we claim:

\Claim\label{the rest of the story}
The PL map $q'_{n}:\oldLambda\to\frakB$ is an $S^1$-fibration for sufficiently large $n$.
\EndClaim

To prove \ref{the rest of the story}, note that, according to \ref{veriest dunce} and \ref{peripatetics}, for sufficiently large $n$, we have that (a) $\tq'_{n}:\toldLambda\to\tfrakB$ is a locally trivial fibration whose fiber is a closed $1$-manifold; and (b) the locally trivial fibration $|q'_{n}|\big||q'_{n}|^{-1}(B):|q'_{n}|^{-1}(B)\to B$ has connected  fiber. We will establish \ref{the rest of the story} by showing that $q':=q'_{n}:\oldLambda\to\frakB$ is an $S^1$-fibration when (a) and (b) hold.

To this end, let any point $x\in\frakB$ be given. Since  $p_\frakB$ is a covering map and $\tq'$ is a locally trivial fibration of manifolds having a closed $1$-manifold as fiber, there is a connected neigborhood $\frakU$ of $x$ in $\frakB$ such that for every component $\tfrakU$ of $p_\frakB^{-1}(\frakU)$, the map  $p_\frakB|\tfrakU:\tfrakU\to\frakU$ is an orbifold covering map, under which the pre-image of $x$ consists of a single point; and the map $\tq'|(\tq')^{-1}(\tfrakU):(\tq')^{-1}(\tfrakU)\to\frakU$ is a  trivial fibration whose fiber is a closed $1$-manifold. Hence for any component $\tfrakW$ of $(p_\frakB\circ  \tq') ^{-1}(\frakU)$, the map $\tq'|\tfrakW$ is a  trivial fibration, whose fiber is PL homeomorphic to $S^1$ and whose base is a component of $p_\frakB^{-1}(\frakU)$.

Since $B$ is dense in $|\frakB|$, we may select a point $y\in\frakU\cap\obd(B)$. Condition (b) above implies that $(q')^{-1}(y)$ is connected.

For each component $\tfrakW$ of
$(p_\frakB\circ\tq')^{-1}(\frakU)=(q'\circ p_\oldLambda) ^{-1}(\frakU)$, our choice of $\frakU$ guarantees that $p_\frakB\circ\tq'|\tfrakW=q'\circ p_\oldLambda|\tfrakW:\tfrakW\to\frakU$ is surjective. Now if $\frakW$ is any component of $(q')^{-1}(\frakU)$, then since $p_\oldLambda$ is an orbifold covering map, there is a component $\tfrakW$ of $(q'\circ p_\oldLambda) ^{-1}(\frakU)$ such that $p_\oldLambda(\tfrakW)=\frakW$; and the surjectivity of $q'\circ p_\oldLambda|\tfrakW:\tfrakW\to\frakU$ implies that
$q'|\frakW:\frakW\to\frakU$ is surjective. 
In particular we have $\frakW\cap (q')^{-1}(y)\ne\emptyset$. But $(q')^{-1}(y)$ is connected, and $\frakW$ is a component of $(q')^{-1}(\frakU)$; hence $\frakW\supset (q')^{-1}(y)$. Since the
component $\frakW$ of $(q')^{-1}(\frakU)$ was arbitrary, this shows that $(q')^{-1}(\frakU)$ is connected.

If we choose a component $\tfrakW_0$ of $(q'\circ p_\oldLambda) ^{-1}(\frakU)$, it follows from the  connectedness of  $(q')^{-1}(\frakU)$ that $p_\oldLambda$ maps $\tfrakW_0$ onto $\frakW_0:=(q')^{-1}(\frakU)$. On the other hand, our choice of $\frakU$ guarantees that $s:=\tq'|\tfrakW_0$ is a  trivial fibration, whose fiber is PL homeomorphic to $S^1$ and whose base is a component $\tfrakU_0$ of $p_\frakB^{-1}(\frakU)$. 
If $G$ and $H$ denote the stabilizers $G^\oldLambda_{\tfrakW_0}$ and $G^\frakB_{\tfrakU_0}$, then $p_\oldLambda|\tfrakW_0:\tfrakW_0\to\frakW_0$ and $p_\frakB|\tfrakU_0:\tfrakU_0\to\frakU$ are regular coverings whose respective covering groups are $G$ and $H$. The homomorphism $\rho$ restricts to a homomorphism $\rho':G\to H$, and $s$  is equivariant in the sense that $s\circ g=\rho'(g)\circ s$ for every $g\in G$. In particular, the action of $G$ on $\tfrakW_0$ preserves the fibration $s$ in the sense that each element of $G$ maps each fiber of $s$ onto a fiber of $s$.

Our choice of $\frakU$ guarantees that the pre-image of $x$
under the covering map
$p_\frakB|\tfrakU_0:\tfrakU_0\to\frakU$ consists of a single point $\tx$. Hence $\tx$ is invariant under $H$, which implies that the fiber $F:=s^{-1}(\tx)$ is invariant under $G$.
Since the PL action of $G$ on the PL $3$-manifold $\tfrakW_0$ preserves the fibration $s$, whose fiber is a $1$-sphere, and leaves the fiber $F$ invariant, some saturated $G$-invariant PL neighborhood $N$ of $F$ in $\tfrakW_0$ may be identified via a PL homeomorphism with $D^2\times S^1$ in such a way that the restricted action of $G$ preserves the product structure; that is, for each $g\in G$ there are periodic PL homeomorphisms $g_1$ and $g_2$ of $D^2$ and $S^1$ respectively such that $g(x,y)=(g_1(x),g_2(y))$ for every $(x,y)\in D^2\times S^1=N$. We may write $N=s^{-1}(\tfrakV)$ for some neighborhood $\tfrakV$ of $\tx$ in $\tfrakU_0$. Since $N$ is $G$-invariant, $\tfrakV$ is $G$-invariant, so that $\tfrakV=\tfrakU_0\cap p_\frakB^{-1}(\frakV)$ for some neighborhood $\frakV$ of $x$ in $\frakU$. The restriction of $p_\oldLambda$ to $N=D^2\times S^1$ induces a PL homeomorphism 
$\eta:
(D^2\times
S^1)/G
\to
(q')^{-1}(\frakV)
$, and Conditions (1) and (2) of \ref{fibered stuff} are now seen to hold with $q'$ playing the role of $q$; thus $q'$ is an $S^1$-fibration, and  \ref{the rest of the story} is proved.

According to \ref{the rest of the story}, we may assume after truncating the sequence $(q'_n)$ that $q'_n$ is an $S^1$-fibration for every $n$. With \ref{when I first put this uniform on}, this gives the conclusion of the proposition.
\EndProof

Note that, according to the convention posited in \ref{categorille}, the PL category will be the default category of orbifolds for the rest of this section.

\Lemma\label{ztimesz-lemma}
Let $\oldOmega$ be a weakly \simple\ $3$-orbifold such that every component of $\partial\oldOmega$ is toric. Suppose that $\oldOmega$ admits no (piecewise linear) $S^1$-fibration, and that no component of $\fraks_\oldOmega$ is $0$-dimensional. Then for every rank-$2$ free abelian subgroup $H$ of $\pi_1(\oldOmega)$, there is a component $\frakK$ of $\partial\oldOmega$ such that $H$ is contained in a conjugate of the image of the inclusion homomorphism $\pi_1(\frakK)\to\pi_1(\oldOmega)$. Furthermore, no finite-sheeted cover of $\oldOmega$ admits an $S^1$-fibration.
\EndLemma

\Proof
The weakly \simple\ orbifold $\oldOmega$ is by definition very good. Hence, by the main result of \cite{lange}, we may write $\oldOmega=\oldPsi\pl$ for some smooth orbifold $\oldPsi$. Since $\oldOmega$ has only toric boundary components, and has no $0$-dimensional components in its singular set, the same is true of $\oldPsi$. Since $\oldOmega$ admits no PL $S ^1 $-fibration, it follows from Proposition \ref{fibration-category} that $\oldPsi$ admits no smooth $S ^1 $-fibration  over a very good $2$-orbifold.

If $\oldTheta$ is any smooth, closed $2$-suborbifold of  $\inter\oldPsi$, then by the first assertion of the main theorem of \cite{illman}, there is a PL structure on the very good orbifold $\oldPsi$, compatible with its smooth structure, for which $\oldTheta$ is a smooth suborbifold. The second assertion of the main theorem of
\cite{illman} implies that $\oldPsi$, equipped with this PL structure, is PL homeomorphic to $\oldOmega$. This shows that  any smooth, closed $2$-suborbifold of  $\inter\oldPsi$ is (topologically) ambiently homeomorphic to a PL surface in $\oldPsi$. Since $\oldOmega$ is weakly \simple\ (and in particular irreducible), it now follows that (A) every $\pi_1$-injective smooth toric suborbifold $\oldTheta$ of $\inter\oldPsi$ is topologically boundary-parallel in $\oldPsi$, in the sense that $\oldTheta$ is the frontier of a suborbifold topologically homeomorphic to $\oldTheta\times[0,1]$; and (B) every smooth spherical suborbifold $\oldTheta$ of $\inter\oldPsi$ bounds a topological discal suborbifold of $\oldPsi$. 

Since $\oldPsi$ 
does  not admit an $S^1$-fibration, has only toric boundary components, has no $0$-dimensional components in its singular set, and has the properties (A) and (B) just stated, it follows from the Orbifold Theorem \cite{blp}, \cite{chk} (in the case where $\fraks_\oldPsi\ne\emptyset$) or from Perelman's geometrization theorem \cite{bbmbp}, \cite{Cao-Zhu}, \cite{kleiner-lott}, \cite{Morgan-Tian} (in the case where $\fraks_\oldPsi=\emptyset$) that
$\inter\oldPsi$ admits  a hyperbolic metric of finite volume. Hence  for every rank-$2$ free abelian subgroup $H$ of $\pi_1(\oldPsi)$, there is a component $\frakK$ of $\partial\oldPsi$ such that $H$ is contained in a conjugate of the image of the inclusion homomorphism $\pi_1(\frakK)\to\pi_1(\oldPsi)$. The first conclusion follows. 

To prove the second assertion, assume that $\oldPsi$ admits an $S^1$-fibration over a compact $2$-orbifold $\frakB$. The
finiteness of $\vol\oldPsi$ implies that $\pi_1(\oldPsi)$ has no abelian subgroup of finite index; hence $\chi(\frakB)<0$. This implies that $\frakB$ has a finite-sheeted cover $\tfrakB$ which is an orientable $2$-manifold with $\chi(\tfrakB)\le-2$. There is a finite-sheeted cover $\toldPsi$ of $\oldPsi$ admitting a locally trivial fibration $p:\toldPsi\to\tfrakB$ with fiber $S^1$. Since $\chi(\tfrakB)\le-2$, there is a simple closed curve $C\subset\inter\tfrakB$ which is not boundary-parallel in $\tfrakB$. But then the image of the inclusion homomorphism $\pi_1(p^{-1}(C))\to\pi_1(\tfrakB)$ is a rank-$2$ free abelian subgroup of $\pi_1(\tfrakB)$ which is not carried by any component of $\partial\tfrakB$. This contradicts the first assertion.
\EndProof

\Lemma\label{affect}
Let $\oldLambda$ be a compact, orientable $3$-orbifold, let $\frakB$ be a
compact ${2}$-orbifold, and let $q:\oldLambda\to\frakB$ be an
$I$-fibration. Then $|q|:|\oldLambda|\to|\frakB|$ is a homotopy
equivalence. 
\EndLemma

\Proof
For each $x\in|\frakB|$ we have $|q|^{-1}(x)=|\frakI|$ for some fiber $\frakI$ of $q$. As observed in \ref{fibered stuff}, $\frakI$ is homeomorphic to either $[0,1]$ or $[0,1]]$, and hence $|\frakI|$ is a topological arc. In particular $|\frakI|$ is contractible and locally contractible; hence by the theorem of \cite{smale}, $|q|$ induces isomorphisms between homotopy groups in all dimensions. As $|\oldLambda|$ and $|\frakB|$ are triangulable, it follows from Whitehead's Theorem \cite[Theorem 4.5]{hatcherbook} that $|q|$ is a homotopy equivalence.

\EndProof

\Definition\label{S-pair def}
Let $\oldLambda$ be a compact, orientable $3$-orbifold, and let $\oldXi$ be a suborbifold of $\partial\oldLambda$. We shall say that an orbifold fibration of $\oldLambda$ is {\it compatible with $\oldXi$} if either (i) the fibration is an $I$-fibration in which $\oldXi=\partial_h\oldLambda $, or (ii) the fibration is an $S^1$-fibration in which $\oldXi$ is saturated.
We define an {\it \spair} to be an ordered pair $(\oldLambda,\oldXi)$ 
where $\oldLambda$ is a compact, orientable $3$-orbifold, $\oldXi\subset\partial\oldLambda$ is a $2$-orbifold, and $\oldLambda$ admits an orbifold fibration which is compatible with $\oldXi$.
We will say that an \spair\ is {\it\pagelike} if Alternative (i) of the definition of compatibility holds, and is {\it\bindinglike} if Alternative (ii) holds; these are not mutually exclusive conditions. We will say that a \pagelike\ \spair\ is {\it\untwisted} if the $I$-fibration in Condition (i) can be taken to be a trivial $I$-fibration (and otherwise that it is twisted). 
\EndDefinition

\begin{remarksdefinitions}\label{Do I need it?}
Let $\oldPsi$ be a compact, orientable $3$-orbifold which is strongly \simple\ and
boundary-irreducible. If  $\oldLambda$ is a $3$-suborbifold of $\oldPsi$ such that
$(\oldLambda,\oldLambda\cap\partial\oldPsi)$ is an \spair, then by
definition $\oldLambda$ admits an orbifold
fibration compatible with $\oldLambda\cap\partial\oldPsi$. The definitions imply that, with respect to such a fibration, the suborbifold 
$\Fr_\oldPsi\oldLambda$ of $\oldLambda$ is saturated; and if the \spair\ $(\oldLambda,\oldLambda\cap\partial\oldPsi)$ is \pagelike\ we may take the fibration to be an $I$-fibration with $\partial_h\oldLambda=\oldLambda\cap\partial\oldPsi$, in which case $\Fr_\oldPsi\oldLambda=\partial_v\oldLambda$. In any event, the saturation of $\Fr_\oldPsi\oldLambda$ implies that each of its components is an annular or toric orbifold, and each of its components must be annular if $(\oldLambda,\oldLambda\cap\partial\oldPsi)$ is \pagelike\ (see \ref{fibered stuff}). 

We define an {\it \Ssuborbifold} of $\oldPsi$ 
to be a $3$-suborbifold $\oldLambda$ of $\oldPsi$ such that
the
components of $\Fr_\oldPsi\oldLambda$ are annular suborbifolds of $\oldPsi$
which are essential in $\oldPsi$, and
$(\oldLambda,\oldLambda\cap\partial\oldPsi)$ is an \spair. The \Ssuborbifold\
$\oldLambda$ will be called {\it\pagelike}, {\it\bindinglike} if
$(\oldLambda,\oldLambda\cap\partial\oldPsi)$ is a \pagelike\ or \bindinglike\ 
\spair, respectively.  Similarly, a \pagelike \Ssuborbifold\ $\oldLambda$ will be termed
{\it\untwisted} or {\it\twisted}   if the \pagelike\ \spair\ 
$(\oldLambda,\oldLambda\cap\partial\oldPsi)$ is 
\untwisted\ or \twisted,  respectively. 

A $3$-suborbifold $\oldLambda$ of $\oldPsi$
will be called an {\it \Asuborbifold} of $\oldPsi$ if the
components of $\Fr_\oldPsi\oldLambda$ are essential annular suborbifolds of $\oldPsi$,
and
$(\oldLambda,\oldLambda\cap\partial\oldPsi)$ is an acylindrical pair (see 
\ref{acylindrical def}).
\end{remarksdefinitions}

\Proposition\label{what i need?}
Let $\oldLambda$ be a \torifold, and let $\frakA$ be a $\pi_1$-injective $2$-suborbifold of $\partial\oldLambda$, each component of which is annular. Then $(\oldLambda,\frakA)$ is a \bindinglike\ \spair. Furthermore, if $\frakA$ is connected and the inclusion homomorphism $\pi_1(\frakA)\to\pi_1(\oldLambda)$ is an isomorphism, then $\oldLambda$ admits a trivial $I$-fibration under which $\frakA$ is a component of $\partial_h\oldLambda$.
\EndProposition

\Proof
For any integer $q\ge1$, let $\frakJ_q$ be defined as in \ref{standard torifold}. Recsall that $\frakJ_q$ has both a natural smooth structure and a natural PL structure; these are mutually compatible. Recall also that, up to a homeomorphism which is at once smooth and PL, $|\frakJ_q|$ may be identified with
$D^2\times S^1$, in such a way that $\fraks_{\frakJ_q}$ is the curve $\{0\}\times S^1$ and has order $q$ (in the sense of \ref{orbifolds introduced}) if $q>1$, and $\fraks_{\frakJ_1}=\emptyset$ if $q=1$. 
For any integer $r>1$, let $\frakD_r$ denote the $2$-orbifold such that $|\frakD_r|=D^2$, and such that $\fraks_{\frakD_r}=\{0\}$, and $0$ has order $q$. Let $\frakD_1$ denote the $2$-orbifold such that $|\frakD_1|=D^2$ and $\fraks_{\frakD_1}=\emptyset$.
For any $q\ge1$, and for any relatively prime integers $m$ and $n$ with $m\ne0$, we may define a smooth $S^1$-fibration $p_{m,n,q}:\frakJ_q\to\frakD_{qm}$ by $p_{m,n,q}(z,w)=z^mw^n$.

For the purpose of this proof, for any integer $m>0$, we define an {\it $m$-admissible system of arcs} to be a set $J\subset S^1$ which is a finite union of closed arcs, and has the property that the map $z\mapsto z^m$ from $J$ to $S^1$ is injective. 
For any $q\ge1$, for any relatively prime integers $m$ and $n$ with $|n|<m$, and for any $m$-admissible system of arcs $J$, set
$\frakA_{m,n,q,J}=\{(ut^n,t^{-m}):t\in S^1,u\in J\}\subset S^1\times S^1=\partial \frakJ_q$. 
Then 
$\frakA_{m,n,q,J}$ is  saturated in the fibration $p_{m,n,q}$, each component of $|\frakA_{m,n,q,J}|$ is a smooth annulus disjoint from $\fraks_{\frakJ_q}$, and $|\frakA_{m,n,q,J}|$ has the same number of components as $J$ (because the $m$-admissibility of $J$ implies that the map $(t,u)\mapsto(ut^n,t^{-m})$ is one-to-one on $S^1\times J$).

Suppose that $\oldLambda$ and $\frakA$ satisfy the hypotheses of the proposition. (According to our conventions, $\oldLambda$ and its suborbifold $\frakA$ are PL.)
According to Proposition \ref{three-way equivalence}, $\oldLambda$ may be identified with the quotient of $D^2\times S^1$ by a standard action of a finite group. 
In view of the description given in \ref{standard torifold} of the quotient of $D^2\times S^1$ by a standard action, 
it follows that up to PL homeomorphism we have either (i) $\oldLambda=\frakJ_q$ for some $q\ge1$, or (ii) there is a $q\ge1$ such that $\oldLambda$ is the quotient of $\frakJ_q$ by the involution $(z,w)\mapsto(\overline{z},\overline{w})$, where bars denote complex conjugation in $D^2\subset\CC$ or $S^1\subset\CC$. If  (i) holds, the hypothesis that the components of $\frakA\subset\partial\oldLambda=S^1\times S^1$ are annular and $\pi_1$-injective in $\oldLambda$ implies that $\frakA$ is piecewise smoothly isotopic to $\frakA_{m,n,q,J}$ for some relatively prime integers $m$ and $n$ 
and some $m$-admissible system of arcs $J$. 
Now since $\partial \frakJ_q$ is disjoint from $\fraks_{\frakJ_q}$, it follows from Proposition \ref{fibration-category} that there is a PL fibration $p_{m,n,q}\PL$ of $\frakJ_q$ such that $p_{m,n,q}\PL|\partial\frakJ_q$ is arbitrarily close in the $C^1$ sense to $p_{m,n,q}$. In particular we can choose $p_{m,n,q}\PL$ so that $\frakA_{m,n,q,J}$, which is  saturated in the fibration $p_{m,n,q}$, is piecewise smoothly isotopic to a PL submanifold $\frakA_{m,n,q,J}\PL$ which is saturated in $p_{m,n,q}\PL$. 
Thus $\frakA$ is piecewise smoothly, and hence piecewise linearly, isotopic to $\frakA_{m,n,q,J}\PL$, and is therefore
saturated in some PL $S^1$-fibration of $\frakJ_q=\oldLambda$. This means that 
$(\oldLambda,\frakA)$ is a \bindinglike\ \spair. This proves the first conclusion of the proposition in the case where (i) holds.

To prove the second conclusion in this case, note that if $\frakA_{m,n,q,J}$ is connected then $J$ is a single arc. Note also that $\pi_1(\frakJ_q)$ is isomorphic to $\ZZ\times(\ZZ/q\ZZ)$. Hence $\pi_1(\frakA_{m,n,q,J})$ and $\pi_1(\frakJ_q)$ cannot be isomorphic unless $q=1$. If we do have $q=1$, then $\frakJ_q$ is the manifold $D^2\times S^1$, and the image of the inclusion homomorphism
$\pi_1(\frakA_{m,n,q,J})\to\pi_1(\frakJ_q)=\pi_1(D^2\times S^1)$ has index $m$ in $\pi_1(D^2\times S^1)$. Hence if $\frakA_{m,n,q,J}$ is connected and  the inclusion homomorphism
$\pi_1(\frakA_{m,n,q,J})\to\pi_1(\frakJ_q)$ is an isomorphism, we have $q=m=1$, so that $\frakJ_q$ is a solid torus and $\frakA_{m,n,q,J}$ is an annulus in $\partial\frakJ_q$ having winding number $1$ in the solid torus $\frakJ_q$; hence the same is true of the PL annulus $\frakA$. It follows that $\frakJ_q$ admits a trivial $I$-fibration under which $\frakA$ is a component of $\partial_h\frakJ_q$. This completes the proof in the case where (i) holds.

We now turn to the case in which (ii) holds. We have $\wt|\partial\oldLambda|=4$ in this case, and every point of $\fraks_{\partial\oldLambda}$ has order $2$. 
As $\oldLambda$ is the quotient of $\frakJ_q$ by the involution $\tau_q:(z,w)\mapsto(\overline{z},\overline{w})$, the quotient map $\sigma_q:\frakJ_q\to\oldLambda$ is a degree-two (orbifold) covering map, $\tau_q$ is its non-trivial deck transformation, and $\oldLambda$ inherits a smooth structure from $\frakJ_q$, compatible with its PL structure. We claim:
\Claim\label{behind the annulus}
If $A$ is any smooth  annulus in $|\partial\oldLambda|\setminus\fraks_\oldLambda$ such that $\omega(A)$ is $\pi_1$-injective in $\oldLambda$, then some self-diffeomorphism of $\oldLambda$ carries $A$ onto an annulus $A'$ such that
$\sigma_q^{-1}(\obd(A'))=\frakA_{m,n,q,J}$ for some relatively prime integers
$m,n$
 with $|n|<m$
 and some admissible system of arcs $J$ having exactly
two components, which are  interchanged by complex
conjugation. Furthermore, the image of the
inclusion homomorphism $\pi_1(\omega(A))\to\pi_1(\oldLambda)$ has index
$2mq$ in $\pi_1(\oldLambda)$. 
\EndClaim

To prove \ref{behind the annulus}, let $E:\RR^2\to S^1\times S^1=\partial\frakJ_q$ denote the map defined by $E(x,y)=(e^{2\pi ix},e^{2\pi iy})$. Then $\sigma_q\circ E:\RR^2\to\partial\oldLambda$ is a smooth orbifold covering map whose group of deck transformations is the group $\Gamma$ generated by the integer translations of $\RR^2$ and the involution $(x,y)\mapsto(-x,-y)$. Thus $\partial\oldLambda-\fraks_{\partial\oldLambda}$ may be identified with the four-punctured sphere $(\RR^2-\Lambda)/\Gamma$, where $\Lambda=(1/2)\ZZ^2$ is the set of all fixed points of non-trivial elements of $\Gamma$. It is a standard fact in two-dimensional topology (cf. \cite[p. 243, Figure 7]{jones-reid}) that every simple closed curve $C$ in
$(\RR^2-\Lambda)/\Gamma$ which separates $(\RR^2-\Lambda)/\Gamma$ into two two-punctured disks is isotopic to the image in $(\RR^2-\Lambda)/\Gamma$ of a line of rational or infinite slope in $\RR^2$ disjoint from $\Lambda$; the slope of this line, an element of $\QQ\cup\{\infty\}$ which is uniquely determined by $C$, will be referred to as the {\it slope} of $C$. Now if $A$ satisfies the hypotheses of
 \ref{behind the annulus}, the $\pi_1$-injectivity of $A$ implies that a core curve of $A$ separates two points of $\fraks_{\partial\oldLambda}$ from the other two, and therefore has a well-defined slope. The $\pi_1$-injectivity of $A$ also implies that the slope of a core curve is
non-zero; we will write as $m/n_0$, where $m$ and $n_0$ are relatively prime integers and $m>0$. Then $A$ is isotopic in $\partial\oldLambda-\fraks_{\partial\oldLambda}=(\RR^2-\Lambda)/\Gamma$ to an annulus $A_0'$ such that
$\sigma_q^{-1}(\obd(A_0'))=
\frakA_{m,n_0,q,J}$ for some admissible system of arcs $J$ having exactly
two components, which are  interchanged by complex
conjugation. 

Now write $n_0=dm+n$, where $d,n\in\ZZ$ and $0\le n<m$. The self-diffeomorphism $\eta:(z,w)\mapsto(zw^{-d},w)$ of $S^1\times S^1=\partial\frakJ_q$ maps $\frakA_{m,n_0,q,J}$ onto $\frakA_{m,n,q,J}$ and commutes with $\tau_q$. Hence $\eta$ induces a self-homeomorphim of $\oldLambda$ which maps $A_0'$ onto an annulus $A'$ such that
$\sigma_q^{-1}(\obd(A'))=
\frakA_{m,n,q,J}$. This establishes the first assertion of \ref{behind the annulus}.

To prove the second assertion of \ref{behind the annulus},
choose a component $J_0$ of $J$, and note that according to the
definition of $\frakA_{m,n,q,J_0}$, the image of the inclusion
homomorphism $\iota:\pi_1(|\frakA_{m,n,q,J_0}|)\to\pi_1(|\frakJ_q|)$ has
index $m$ in $\pi_1(|\frakJ_q|)$. We may identify
$\pi_1(\frakJ_q)$ with $\pi_1(|\frakJ_q|)\times\ZZ_q$ in such
a way that $\iota$ is the composition of the inclusion homomorphism
$\kappa:\pi_1(\frakA_{m,n,q,J_0})\to\pi_1(\frakJ_q)$ with the
projection to the first factor. Since
$\pi_1(\frakA_{m,n,q,J_0})$ is cyclic, it follows that the image of $\kappa$ has index
$mq$ in $\pi_1(\frakJ_q)$. But $\frakA_{m,n,q,J_0}$ is a lift of
the annular orbifold $\obd(A')$
to the two-sheeted cover $\frakJ_q$ of $\oldLambda$,
and hence the image of the inclusion homomorphism
$\obd(A')\to\pi_1(\oldLambda)$ has index $2mq$ in
$\pi_1(\oldLambda)$. 
Since some self-homeomorphism of $\oldLambda$ carries
$\obd(A)$ onto 
$\obd(A')$, 
this completes the proof of \ref{behind the annulus}.


Let us now prove the first conclusion of the proposition in the subcase in which $|\frakA|$ is a single  weight-$0$ annulus. In this subcase, choose a smooth annulus $\frakA\smooth$ which is ambiently, piecewise-smoothly isotopic to $\frakA$. In view of \ref{behind the annulus}, some self-diffeomorphism of $\oldLambda$ maps $\frakA\smooth$ onto an annulus $\frakA'\smooth$ such that
$\sigma_q^{-1}(\frakA'\smooth)=\frakA_{m,n,q,J}$ for some relatively prime integers $m,n$ with $|n|<m$ and
some $m$-admissible system of arcs $J$ having exactly two components.
The definitions of $\tau_q$ and $p_{m,n,q}$ imply that the fibration $p_{m,n,q}$ is invariant under $\tau_q$, in the sense that $\tau_q$ maps each fiber onto a fiber. Hence $p_{m,n,q}$ induces a fibration of $\oldLambda$; and since $\sigma_q^{-1}(\frakA'\smooth)=\frakA_{m,n,q,J}$ is saturated in the fibration $p:=p_{m,n,q}$, the annulus $\frakA'\smooth$ is saturated in the induced fibration $\overline{p}$ of $\oldLambda$.

Now let $V$ be a compact PL neighborhood submanifold of $\oldLambda$ which is a neighborhood of $\frakA'\smooth$ and is disjoint from $\fraks_{\frakJ_q}$.  It follows from Proposition \ref{fibration-category} that there is a PL fibration $\overp\pl$ of $\frakJ_q$ such that $|\overp\pl|$ is arbitrarily close in the uniform sense to $|\overp|$, and $\overp\pl|V$ is arbitrarily close in the $C^1$ sense to $\overp|V$. In particular we can choose $\overp\PL$ so that $\frakA'\smooth$, which is contained in $\inter V$ and is saturated in the fibration $\overp$, is piecewise smoothly isotopic to a PL submanifold $\frakA'$ which is  contained in $\inter V$ and is saturated in $\overp$. Thus $\frakA$ is piecewise smoothly, and hence piecewise linearly, isotopic to $\frakA'$, and is therefore
saturated in some PL $S^1$-fibration of $\frakJ_q=\oldLambda$. 
Hence the pair $(\oldLambda,\frakA)$ is a \bindinglike\ \spair, and the first conclusion is established in this subcase.

We next prove  the first conclusion in the more general subcase in which every component $|\frakA|$ is a single weight-$0$ annulus. Since each component of $|\frakA|$ separates two points of $\fraks_{\partial\oldLambda}$ from the other two, the core curves of the components of $|\frakA|$ are all parallel, and hence there is a single $\pi_1$-injective weight-$0$ annulus $B$ containing $|\frakA|$. By the subcase proved above, letting $\frakB:=\omega(B)$ play the role of $\frakA$, the suborbifold $\frakB$ is saturated in some $S^1$-fibration $p$ of $\oldLambda$. The restriction of $p$ to $\frakB$ is an $S^1$-fibration of the annulus $\frakB$, and since the components of $\frakA$ are $\pi_1$-injective in $\oldLambda$ and hence in $\frakB$, the suborbifold $\frakA$ is isotopic in $\frakB$ to a suborbifold which is saturated in the restricted fibration of $\frakB$, and is therefore saturated in the fibration $p$ of $\oldLambda$. This shows that $(\oldLambda,\frakA)$ is a \bindinglike\ \spair, and the first conclusion is established in this more general subcase.

We can now prove the first conclusion whenever (ii) holds. Let $\frakC$ denote a strong regular neighborhood of $\partial\frakA$ relative to $(\partial\oldLambda)-\inter\frakA$. Then $\frakC$ is $\pi_1$-injective in $\oldLambda$, and each component of $|\frakC|$ is a weight-$0$ annulus. Applying the subcase just proved, with $\frakC$ playing the role of $\frakA$, we obtain an $S^1$-fibration of $\partial\oldLambda$ in which $\frakC$ is saturated. Since $\frakA$ is a union of components of $(\partial\oldLambda)-\inter\frakC$, it follows that $\frakA$ is saturated in this fibration, and the proof of the first conclusion is complete.


To prove the second conclusion in the case where (ii) holds, suppose that
 $\frakA$ is connected and that the inclusion homomorphism $\pi_1(\frakA)\to\pi_1(\oldLambda)$ is an isomorphism. 
In particular $\pi_1(\frakA)$ is isomorphic to $\pi_1(\oldLambda)$ and is therefore not cyclic; thus $|\frakA|$ is not a weight-$0$ annulus. Hence $|\frakA|$ must be a weight-$2$ disk  and the points of $\fraks_\frakA$ must have order $2$.
In particular, $\partial\frakA$ is connected. Hence, if we again let $\frakC$ denote a regular neighborhood of $\partial\frakA$ in $\partial\oldLambda-\inter\frakA$, then $|\frakC|$ is a weight-$0$ annulus. The $\pi_1$-injectivity in $\oldLambda$ of the annular orbifold $\frakA$ implies that $\frakC$ is also $\pi_1$-injective. Now let $C$ be a smooth annulus which is ambiently, piecewise-smoothly isotopic to $|\frakC|$ in $|\oldLambda|-\fraks_\oldLambda$.
In view of \ref{behind the annulus}, 
we may assume without loss of generality that 
$\sigma_q^{-1}(\obd(C))=\frakA_{m,n,q,J}$ for some relatively prime
integers $m,n$ with $|n|<m$ and some $m$-admissible system of arcs $J$ having exactly two components, which are interchanged by complex conjugation.

On the other hand, since $|\frakA|$ is a weight-$2$ disk and the points of $\fraks_\frakA$ have order $2$,  the image of the inclusion homomorphism $\pi_1(\partial\obd(C))\to\pi_1(\frakA)$ has index $2$ in $\pi_1(\frakA)$. Since
 the inclusion homomorphism $\pi_1(\frakA)\to\pi_1(\oldLambda)$ is an 
 isomorphism, the image of the inclusion homomorphism
 $\pi_1(\partial\frakA)\to\pi_1(\oldLambda)$ has index $2$ in
 $\pi_1(\oldLambda)$. Now since, by the final assertion of \ref{behind the annulus},
the image of the 
inclusion homomorphism $\pi_1(\frakC)\to\pi_1(\oldLambda)$ has index
$2mq$ in $\pi_1(\oldLambda)$, we have $2mq=2$ and hence $m=q=1$. Since
$q=1$, the orbifold $\frakJ_q$ is simply the manifold $D^2\times
S^1$, and since
$|n|<m=1$ we have $n=0$. Thus we have
$\sigma_1^{-1}(\obd(C))=\frakA_{1,0,1,J}$, which by definition means that
$\sigma_1^{-1}(\obd(C))=J\times S^1$. Since the pairs $(\oldLambda,\frakC)$ and $(\oldLambda,\obd(C))$ are piecewise smoothly, and hence piecewise linearly, homeomorphic, we may identify $\oldLambda$
with $(D^2\times S^1)/\langle\tau_1\rangle$, via a PL homeomorphism, in such a
way that $\frakC=(J\times S^1)/\langle\tau_1\rangle$.

Let $h:D^2\to [0,1]\times[-1,1]$ be a PL homeomorphism such that
$h(\overline{z})=\rho(h(z))$ for each $z\in D^2\subset\CC$, where the
bar denotes complex conjugation and $\rho$ is the reflection about the
$x$-axis in $[0,1]\times[-1,1]\subset\RR^2$. Since the components of
$J$ are interchanged by complex conjugation, we may choose $h$ in such
a way that  $h(J)=[0,1]\times\{-1,1\}$. If we define a homeomorphism
$H:D^2\times S^1\to [0,1]\times[-1,1]\times S^1$ by $H(z,w)=(h(x),w)$,
we then have $H\circ\tau_1=\tau^*\circ H$, where
$\tau^*(x,y,w)=(x,-y,\overline{w})$. Hence $H$ induces a homeomorphism
$\barH:\oldLambda\to([0,1]\times[-1,1]\times
S^1)/\langle\tau^*\rangle$. We have
$\barH(\frakC)=[0,1]\times\{-1,1\}\times
S^1)/\langle\tau^*\rangle$. Since $\frakA$ is a component of
$(\partial\oldLambda)-\inter\frakC$, we may assume, after possibly
replacing $H$ by its postcomposition with $(x,y,w)\to(1-x,y,w)$, that $\barH(\frakA)=(\{0\}\times[-1,1]\times
S^1)/\langle\tau^*\rangle$.  This shows that the pair
$(\oldLambda,\frakA)$ is PL homeomorphic to the pair
$(\frakA\times[0,1],\frakA\times\{0\})$, and the second conclusion is
proved in this case. 
\EndProof

\Corollary\label{covering annular} 
Let $p:\toldPsi\to \oldPsi$ be a covering map
of orientable $3$-orbifolds that are 
componentwise strongly \simple. Let $\oldXi$ be a $\pi_1$-injective suborbifold of $\partial\oldPsi$. Then for any essential orientable annular $2$-orbifold $\oldPi$ in the pair $(\oldPsi,\oldXi)$, every component of $p^{-1}(\oldPi)$ is
essential in the pair $(\toldPsi,p^{-1}(\oldXi))$. 
\EndCorollary

\Proof
Set $\toldXi:=p^{-1}(\toldXi)$.

Since $\oldPi$ is essential, it is in particular $\pi_1$-injective;
and since $p$ is an orbifold covering map,
$p_\sharp:\pi_1(\toldPsi)\to\pi_1(\oldPsi)$ is injective. It follows
that $p^{-1}(\oldPi)$ is $\pi_1$-injective. It remains to show that no
component of $p^{-1}(\oldPi)$ is parallel
 in the pair $(\toldPsi,\toldXi)$ either to a suborbifold of $\toldXi$ or to a component of
$\overline{(\partial\toldPsi)-\toldXi}$. 

Assume that this is false. 
Then there is a $3$-suborbifold $\toldLambda$ of $\toldPsi$ such that (i) $\Fr_{\toldPsi} \toldLambda$ is a single component of $p^{-1}(\oldPi)$, and (ii) $\toldLambda$ admits a trivial $I$-fibration under which $\Fr_{\toldPsi}\toldLambda$ is a component of $\partial_h\toldLambda$, the other component of $\partial_h\toldLambda$ is contained either in  $\toldXi$ or in 
$\overline{(\partial\toldPsi)-\toldXi}$, and
$\partial_v\toldLambda\subset\toldXi$.

We may assume that, among all $3$-suborbifolds of $\toldLambda$ for which (i) and (ii) hold, $\toldLambda$ is minimal with respect to inclusion. 
Set $\toldPi=\Fr_{\toldPsi}\toldLambda$, so that $\toldPi$ is a union of  components of $p^{-1}(\oldPi)$. Thus either (a) there is a component $\toldPi'$ of $p^{-1}(\oldPi)$ contained in $\toldLambda$ and disjoint from $\toldPi$, or (b)
$\toldLambda\cap p^{-1}(\oldPi)=\toldPi$. In each case we will obtain a contradiction.

First suppose that (a) holds.
It follows from (ii) that the inclusion homomorphism $H_1(|\toldPi|;\QQ)\to H_1(|\toldLambda|;\QQ)$ is an isomorphism. Since $|\toldPi'|$ is disjoint from $|\toldPi|$, the $2$-manifold $|\toldPi'|$ separates $|\toldLambda|$. Hence there are
$3$-suborbifolds $\toldLambda'$ and $\toldLambda''$ of $\toldLambda$ such that $\toldLambda'\cap\toldLambda''=\toldPi'$, $\toldLambda'\cup\toldLambda''=\toldLambda$, and $\toldPi\subset\toldLambda''$. 
In particular Condition (i) holds with $\toldLambda'$ in place of $\toldLambda$. 
We will show that Condition (ii) also holds with $\toldLambda'$ in place of $\toldLambda$.  This will contradict the minimality of $\toldLambda$.

As a component of $p^{-1}(\oldPi)$, the orbifold $\toldPi'$ is $\pi_1$-injective in $\toldPsi$. It therefore follows from Lemma \ref{oops lemma} that $\toldLambda'$ is strongly \simple. If $\oldTheta$ denotes the component of $\partial\toldLambda'$ containing $\toldPi'$, then $\oldTheta$ is the union of the annular orbifold $\toldPi'$ with one or more components of $\toldLambda'\cap\partial\toldPsi$. These components of $\toldLambda'\cap\partial\toldPsi$ are contained in the annular orbifold $\toldLambda\cap\partial\toldPsi$, and are bounded by components of $\partial\toldPi'$, which are $\pi_1$-injective in $\toldPsi$ since $\toldPi'$ is; hence, by \ref{cobound}, the components of $\toldLambda'\cap\partial\toldPsi$ contained in $\oldTheta$ are annular. It follows that $\oldTheta$ is toric. We now deduce from Proposition \ref{three-way equivalence} that $\toldLambda'$ is a \torifold.

Since $\toldPi'$ is $\pi_1$-injective in $\toldPsi$, we
may write $\pi_1(\toldPsi)$ as a free product with amalgamation $\pi_1(\toldLambda')\star_{\pi_1(\toldPi')}\pi_1(\toldLambda'')$,
using a base point in $\toldPi'$. Hence if $L'$, $L''$ and $P$ denote the respective images of 
$\pi_1(\toldLambda')$, $\pi_1(\toldLambda'')$ and $\pi_1(\toldPi')$ under the inclusion homomorphisms into $\pi_1(\toldPsi)$, we have $L'\cap L''=P$. But since $\toldPi\subset\toldLambda''$, and since $\toldLambda$ satisfies (ii), we have $L''=\pi_1(\toldPsi)$. Hence $L'=P$, i.e. the inclusion $\pi_1(\toldPi')\to\pi_1(\toldLambda')$ is an isomorphism. It therefore follows from
Prop. \ref{what i need?} that $\toldLambda'$ admits a trivial $I$-fibration under which $\toldPi'=\Fr_{\toldPsi}\toldLambda'$ is a component of $\partial_h\toldLambda$. 

Let $\tfrakZ$ denote the component of $\partial_h\toldLambda$ distinct from $\toldPi$ (where $\toldLambda$ has the $I$-fibration given by (ii)). Let $\tfrakZ'$ denote the component of $\partial_h\toldLambda'$ distinct from $\toldPi'$. 

Since the boundary of $\Fr_{\toldPsi}\toldLambda'=
\toldPi'$ is contained in $ \inter\toldXi$, the orbifold $\toldLambda'\cap\overline{(\partial\toldPsi)-\toldXi}$ is a union of components of $\toldLambda\cap\overline{(\partial\toldPsi)-\toldXi}$. But Condition (ii) for $\toldLambda$ implies that $\toldLambda\cap\overline{(\partial\toldPsi)-\toldXi}$ either is empty or is equal to $\tfrakZ$ (and hence  connected). Hence either $\toldLambda'\cap\overline{(\partial\toldPsi)-\toldXi}=\emptyset$ or $\toldLambda'\cap\overline{(\partial\toldPsi)-\toldXi}=\tfrakZ$. If $\toldLambda'\cap\overline{(\partial\toldPsi)-\toldXi}=\emptyset$, then $\partial_v\toldLambda'\cup\tfrakZ'\subset\inter\toldXi$, which establishes Condition (ii) for $\toldLambda'$ in this subcase. Now consider the subcase in which $\toldLambda'\cap\overline{(\partial\toldPsi)-\toldXi}=\tfrakZ$.  We then have $\tfrakZ\subset\toldLambda'\cap\partial\toldPsi\subset\toldLambda\cap\partial\toldPsi$. But Condition (ii) for $\toldLambda$ implies that 
$\toldLambda\cap\partial\toldPsi=\tfrakZ\cup\partial_v\toldLambda$, and hence $\toldLambda\cap\partial\toldPsi$ is a regular neighborhood of $\tfrakZ$ in $\partial\oldPsi$. Since $\partial(\toldLambda'\cap\partial\toldPsi)$ is $\pi_1$-injective in $\partial\toldPsi$, it follows that $\toldLambda'\cap\partial\toldPsi$ is a regular neighborhood of $\tfrakZ$ in $\partial\oldPsi$. But since $\toldPi'=\Fr_{\toldPsi}\toldLambda'$, we have
$\toldLambda'\cap\partial\toldPsi=\tfrakZ'\cup\partial_v\toldLambda'$, and hence
 $\toldLambda'\cap\partial\toldPsi$ is a regular neighborhood of $\tfrakZ'$ in $\partial\toldPsi$. Thus $\tfrakZ$ and $\tfrakZ'$ are (orbifold-)isotopic in $\toldLambda'\cap\partial\toldPsi$, and we may therefore modify the $I$-fibration of $\toldLambda'$ so as to arrange that $\tfrakZ'=\tfrakZ$. This establishes Condition (ii) in this subcase, and so we always have the required contradiction in the case where (a) holds.

 Now suppose that (b) holds. Then $\toldLambda$ is the closure of a component of $\toldPsi-p^{-1}(\oldPi)$. Hence we may write $\toldLambda=p^{-1}(\oldLambda)$, where $\oldLambda$ is the closure of some component of $\oldPsi-\oldPi$. Thus $q:=p|\toldLambda:\toldLambda\to\oldLambda$ is a(n orbifold-) covering. Since (b) holds, we have $q^{-1}(\oldPi)=\toldPi$. Since $\oldLambda$ is covered by the \torifold\ $\toldLambda$, it follows from the definition that $\oldLambda$ is itself a \torifold. Since 
$\toldPi$ is connected, and since the inclusion homomorphism $\pi_1(\toldPi)\to\pi_1(\toldLambda)$ is an isomorphism by (ii), the inclusion homomorphism $\pi_1(\oldPi)\to\pi_1(\oldLambda)$ is also an isomorphism. 
It now follows from Proposition \ref{what i need?} that $\oldLambda$ has a trivial $I$-fibration under which $\oldPi$ is a component of $\partial_h\oldLambda$. 

Let $\frakZ$ denote the component of $\partial_h\oldLambda$ distinct from $\oldPi$. Let
$\tfrakZ$ denote the component of $\partial_h\toldLambda$ distinct from $\toldPi$,
under the $I$-fibration of $\toldLambda$ given by  Condition (ii). According to that condition we have $\partial_v\toldLambda\subset\toldXi$, and $\tfrakZ$ is contained either in $\toldXi$ or in $\overline{(\partial\toldPsi)-\toldXi}$. If 
$\tfrakZ\subset\toldXi$ then $\toldLambda\cap\partial\toldPsi\subset\toldXi$; hence
$\oldLambda\cap\partial\oldPsi\subset\oldXi$, i.e. $\partial_v\oldLambda\cup\frakZ\subset\oldXi$. Thus in this subcase
$\oldPi$ is parallel in the pair $(\oldPsi,\oldXi)$ to a suborbifold of $\oldXi$,  a contradiction to the essentiality of $\oldPi$. 

If
$\tfrakZ\subset\overline{(\partial\toldPsi)-\toldXi}$, then $\tfrakZ=\toldLambda\cap\overline{(\partial\toldPsi)-\toldXi}$. But $\tfrakZ$ is connected, and the inclusion homomorphism $\pi_1(\tfrakZ)\to\pi_1(\toldLambda)$ is an isomorphism. Hence 
$\frakV:=\oldLambda\cap\overline{(\partial\oldPsi)-\oldXi}$ is connected, and the inclusion homomorphism $\pi_1(\frakV)\to\pi_1(\oldLambda)$ is an isomorphism. But we have $\frakV\subset\oldLambda\cap\partial\oldPsi=(\partial_v\oldLambda)\cup\frakZ$. The orbifold $(\partial_v\oldLambda)\cup\frakZ$ is connected, and the inclusion homomorphism $\pi_1((\partial_v\oldLambda)\cup\frakZ)\to\pi_1(\oldLambda)$ is an isomorphism. 
The inclusion homomorphism $\pi_1(\frakV)\to\pi_1((\partial_v\oldLambda)\cup\frakZ)$ is therefore also an isomorphism, so that $(\partial_v\oldLambda)\cup\frakZ$ is a 
regular neighborhood of $\frakV$ in $\partial\oldPsi$.
As $(\partial_v\oldLambda)\cup\frakZ$ is also a
regular neighborhood of $\frakZ$ in $\partial\oldPsi$, it follows that 
$\frakV$ and $\frakZ $ are isotopic in $(\partial_v\oldLambda)\cup\frakZ$. Hence after modifying the $I$-fibration of $\oldLambda$ we may assume that 
$\frakV=\frakZ $. This means that
$\oldPi$ is parallel in the pair $(\oldPsi,\oldXi)$ to the component $\frakV$ of
$\overline{(\partial\oldPsi)-\oldXi}$, and we again have a contradiction to the essentiality of $\oldPi$. 
 \EndProof

\Lemma\label{when a tore a fold}
Let $\oldLambda$ be a suborbifold of a componentwise strongly \simple,
componentwise boundary-irreducible, orientable $3$-orbifold $\oldPsi$. Then $\oldLambda$
is a \bindinglike\ \Ssuborbifold\ of $\oldPsi$ if and only if (i)
$\oldLambda$ is a \torifold\, and (ii) every component of
$\Fr_\oldPsi\oldLambda$ is an  essential annular suborbifold of
$\oldPsi$. 
\EndLemma

\Proof
Suppose that $\oldLambda$
is a \bindinglike\ \Ssuborbifold\ of $\oldPsi$. Consider first the case in which $\partial\oldLambda\ne\emptyset$, and fix a component $\oldTheta$ of $\partial\oldLambda$. Since $\oldLambda$ is a \bindinglike\ \Ssuborbifold, $\oldTheta$ admits an $S^1$-fibration over a closed $1$-orbifold, and is therefore toric. 
Since the components of $\Fr\oldLambda$ are essential, and in particular $\pi_1$-injective, in the strongly \simple\ orbifold $\oldPsi$,
it follows from Lemma \ref{oops lemma} that $\oldLambda$ is strongly \simple. Since $\oldTheta$ is toric, $\oldLambda$ satisfies Condition (2) of Proposition \ref{three-way equivalence}, and according to the latter proposition it is a \torifold. This is Condition (i) of the statement of the present lemma, and
Condition (ii) is immediate from the definition of an \Ssuborbifold.

Now consider the case in which $\oldLambda$ is closed. Then $\oldLambda$ is a component of $\oldPsi$ and is therefore strongly \simple. According to Condition (II) of Definition \ref{oops}, some finite-sheeted covering space $\toldLambda$ of $\oldLambda$ is an irreducible $3$-manifold. But by definition the \bindinglike\ \Ssuborbifold\ $\oldLambda$ admits an $S^1$-fibration, and $\toldLambda$ inherits an $S^1$-fibration $q:\toldLambda\to\frakB$ for some closed $2$-orbifold $\frakB$. Since $\oldLambda$ is irreducible, $\frakB$ is infinite and therefore has an element of infinite order. This implies that $\pi_1(\toldLambda)$ has a rank-$2$ free abelian subgroup, a contradiction to Condition (II) of Definition \ref{oops}. Hence this case cannot occur, and we have proved that if $\oldLambda$ is a \bindinglike\ \Ssuborbifold\ then (i) and (ii) hold.

Conversely, suppose that (i) and (ii) hold. It follows from (ii) that the components of $\Fr\oldPsi$ are $\pi_1$-injective in $\oldPsi$ as well as being annular. Hence Proposition \ref{what i need?} implies that $(\oldPsi,\Fr\oldPsi)$ is a \bindinglike\ \spair. As (ii) includes the assertion that the components of $\Fr\oldPsi$ are essential in $\oldPsi$, it now follows that $\oldLambda$ is a \bindinglike\ \Ssuborbifold\ of $\oldPsi$.
\EndProof

\Proposition\label{when vertical}
Let $\oldLambda$ be an orientable
$3$-orbifold equipped with an $I$-fibration over a $2$-orbifold, with $\chi(\oldLambda)<0$,
and let $\frakA$ be an essential orientable annular suborbifold of
the pair $(\oldLambda,\partial_h\oldLambda)$. Then $\oldLambda$ is isotopic in the pair $(\oldLambda,\partial_h\oldLambda)$ to a saturated annular suborbifold.
\EndProposition

\Proof
Set $\oldXi=\partial_h\oldLambda$ and $\oldPsi=\silv_\oldXi\oldLambda$, and let $\oldTheta$ denote the toric suborbifold $\silv\frakA$ of $\inter\oldPsi$. According to Proposition \ref{silver acylindrical}, $\oldTheta$ is $\pi_1$-injective in $\oldPsi$ and has no boundary-parallel component. The $I$-fibration of $\oldLambda$ gives rise to an $S^1$-fibration of $\oldXi$ in which the fibers are obtained by silvering the fibers of $\oldLambda$. The base of the $I$-fibration of $\oldLambda$, which we will denote by $\frakB$, is also the base of the $S^1$-fibration of $\oldPsi$. It follows from the PL version of 
 \cite[Verticalization Theorem 4]{bonahon-siebenmann} (see \ref{categorille}) that $\oldTheta$ is isotopic in $\oldPsi$ to a suborbifold $\oldTheta'$ which is either ``vertical,'' in the sense of being saturated in the $S^1$-fibration of $\oldPsi$, or ``horizontal,'' in the sense of being transverse to the fibers. 
If $\oldTheta'$ is horizontal then $\oldTheta$ is a(n orbifold) covering of $\frakB$; if $d$ denotes the degree of the covering, the toricity of $\oldTheta$ gives $0=\chi(\oldTheta)=\chi(\oldTheta)=d\chi(\frakB)$, so that $\chi(\frakB)=0$. But by hypothesis we have $0>\chi(\oldLambda)=\chi(\frakB)$, a contradiction. Hence $\oldTheta'$ must be vertical, i.e. saturated. We may therefore write $\oldTheta'=\silv\frakA'$ for some saturated suborbifold $\frakA'$ of $\oldLambda$. Since $\oldTheta=\silv\frakA$ is isotopic in $\oldPsi=\silv_\oldXi\oldLambda$ to $\oldTheta'=\silv\frakA'$, the suborbifolds $\frakA$ and $\frakA'$ of $\oldLambda$ are isotopic.
\EndProof

\Corollary\label{vertical corollary}
Let $(\oldLambda,\oldXi)$ be an \spair\ such that $\oldLambda$ is strongly \simple,
and let $\frakA$ be an
 essential orientable annular suborbifold of 
the pair $(\oldLambda,\oldXi)$. Then $\oldLambda$ admits an orbifold fibration compatible with $\oldXi$ in which $\frakA$ is saturated.
\EndCorollary

\Proof
In the case where
$\chi(\oldLambda)<0$, the \spair\ $(\oldLambda,\oldXi)$ must be \pagelike, since an orbifold admitting an $S^1$-fibration has Euler characteristic $0$; thus in this case the assertion follows from Proposition \ref{when vertical}. Now suppose that
$\chi(\oldLambda)\ge0$. Since $\frakA$ is annular, the (possibly disconnected) $3$-orbifold $\oldLambda':=\oldLambda\cut\frakA$ also has non-negative Euler characteristic. By Lemma \ref{oops lemma}, $\oldLambda'$ is componentwise strongly \simple; in particular, $\oldLambda'$ is componentwise irreducible and has no discal component, so $\partial\oldLambda'$ has no spherical component. Since $\chi(\partial\oldLambda')=2\chi(\oldLambda')\ge0$, it now follows that every component of $\partial\oldLambda'$ is toric. Hence by Proposition \ref{three-way equivalence}, $\oldLambda'$ is a \torifold. If we set $\tfrakA=\rho_\frakA^{-1}(\frakA)$ and $\toldXi=\rho_\frakA^{-1}(\oldXi)$ (see \ref{nbhd stuff}), , each component of $\toldXi\discup\tfrakA$ is annular and $\pi_1$-injective in $\oldLambda'$. By  Proposition \ref{what i need?}, $(\oldLambda',\toldXi\cup\tfrakA)$ is a \bindinglike\ \spair, i.e. $\oldLambda'$ has an $S^1$-fibration in which $\tfrakA$ and $\toldXi$ is saturated. Since the $S^1$-fibration of an annular $2$-orbifold is unique up to isotopy, we may choose the $S^1$-fibration of $\oldLambda'$ in such a way that $\grock_\frakA$ (see \ref{nbhd stuff})  interchanges the induced fibrations of the components of $\tfrakA$. This implies that
$\oldLambda$ has an $S^1$-fibration in which $\frakA$ and $\oldXi$ is saturated, which gives the conclusion in this case.
\EndProof

\Proposition\label{new what they look like}
Let $\oldLambda$ be a compact, connected, orientable $3$-orbifold, and let $\oldXi$ be a compact $2$-suborbifold of $\partial\oldLambda$. Then $\silv_\oldXi\oldLambda$ admits an $S^1$-fibration over a $2$-orbifold if and only if $(\oldLambda,\oldXi)$ is an \spair.
\EndProposition

\Proof 
Set $\oldOmega=\silv_\oldXi\oldLambda$. 

First suppose that  $(\oldLambda,\oldXi)$ is an \spair. By definition this means that there is a fibration $q:\oldLambda\to\frakB$ of $\oldLambda$
over a $2$-orbifold $\frakB$ such that either (i) $q$ is an
$I$-fibration and $\oldXi=\partial_h\oldLambda$, or (ii) $\oldLambda$ is an
$S^1$-fibration and $\oldXi$ is saturated. If (ii) holds, write
$\oldXi=q^{-1}(\alpha)$ for some $1$-orbifold  $\alpha\subset\frakB$, and set
$\oldGamma=\silv_\alpha\frakB$. 
If (i) holds, set
$\oldGamma=\frakB$. 
In each case 
there is a unique fibration $r:\oldOmega\to\oldGamma$ such that $|r|=|q|$. (In Case (i), the fibers of $r$ are homeomorphic to $[[0,1]]$.)

Conversely, suppose that $r$ is an $S^1$-fibration of $\oldOmega$ over a
$2$-orbifold $\oldGamma$. Consider an arbitrary point $v$  of $\oldGamma$.  It follows from the definition of an $S^1$-fibration (see \ref{fibered stuff})
that
there exist a neighborhood $\frakV$ of $v$ in
$\oldGamma$ and a(n orbifold) homeomorphism $\eta:
(D^2\times
S^1)/G
\to
r^{-1}(\frakV)
$, where 
(1)  $G$ is a subgroup of $G_1\times G_2$, for some subgroups
$G_1$ and $G_2$ of of ${\rm SO}(2)$,
and (2) 
if $\pi:
D^2\times
S^1
\to
(D^2\times
S^1)/G$
denotes the orbit map, then
$\eta(\pi(\{x\}\times S^1))$ 
is a fiber of $r$
for each $x\in D^2$,
and
$\eta(\pi(\{0\}\times S^1))=
r^{-1}(v)$. (Since we are in the PL category, this should be interpreted in terms of the PL structures on $S^1$ and $D^2$ discussed in \ref{esso}.) Now the canonical two-sheeted covering map from $D_\oldXi\oldLambda$ ro $\oldOmega$ maps $\oldXi$ homeomorphically onto a $2$-orbifold $\oldXi'\subset\oldOmega$. \abstractcomment{I was thinking it's not a suborbifold, but now after studying defs. I think it is.} Since $\oldLambda$ is orientable, the
definition of $\oldOmega=\silv_\oldXi\oldLambda$ implies 
that $\oldXi'\subset\oldOmega$ is the
union of all closures of $2$-dimensional strata of $
\fraks_\oldOmega$.
 Hence
 the pre-image of $\oldXi'\cap r^{-1}(\frakV)$ under the covering map $P=\eta\circ\pi$ 
is
the
union of all $2$-dimensional fixed point sets of elements of $G$. Each
such $2$-dimensional fixed point set is either
equal to
$D^2\times\{-1,1\}$ or has the form $\ell\times S^1$, for some diameter $\ell$ 
of $D^2$. But $P^{-1}(\oldXi'\cap r^{-1}(\frakV))$ 
is a $2$-manifold since $\oldXi'$ is a $2$-suborbifold of $\oldOmega$, and it
follows that the union expressing $P^{-1}(\oldXi'\cap r^{-1}(\frakV))$ must
consist of at most a single term. This proves:

\Claim\label{better onion}
For every $v\in\oldGamma$ there exist a neighborhood $\frakV$ of $v$ in
$\oldGamma$ and  a regular covering map $P:D^2\times S^1\to r^{-1}(\frakV)$ such that
(a) $P$ maps $\{x\}\times S^1$ onto a fiber of $r$ for each $x\in
D^2$, and (b)
$P^{-1}(\oldXi'\cap r^{-1}(\frakV))$ either (I) is empty, (II) is equal to
$D^2\times\{-1,1\}$, or (III) has the form $\ell\times S^1$ for some diameter $\ell$ 
of $D^2$. Furthermore, if (II) or (III) holds, there is a deck
transformation of the covering $P:D^2\times S^1\to r^{-1}(\frakV)$ which
interchanges the components of $(D^2\times S^1)-P^{-1}(\oldXi'\cap r^{-1}(\frakV))$.
\EndClaim 

Let $\frakU$ denote the set of all points $v\in\oldGamma$ such that there
exist a neighborhood $\frakV$ of $v$  and  a regular covering map $P$ satisfying
Condition (a), and Alternative (I) or (III) of Condition (b), of
\ref{better onion}; 
and let $\frakU'$ denote the set of all points $v\in\oldGamma$ such that there
exist a neighborhood $\frakV$ of $v$  and  a regular covering map $P$ satisfying
Condition (a), and Alternative (II) of Condition (b). Then $\frakU\cup
\frakU'=|\oldGamma|$. It is clear that $\frakU$ and $\frakU'$ are open. 
Furthermore, the
fiber over every point of $\frakU$ is either contained in $\oldXi'$ or
disjoint from $\oldXi'$, while the fiber over every point of $\frakU'$
meets $\oldXi'$ in one or two
  points. Thus $\frakU\cap\frakU'=\emptyset$. But $\oldGamma$ is connected since
$\oldLambda$ is connected, and hence either $\oldGamma=\frakU$ or
$\oldGamma=\frakU'$.

In the case where $\oldGamma=\frakU$, the local definition of $\frakU$ implies that $\oldXi'$ is saturated
in the fibration $r:\oldOmega\to\oldGamma$, and we may therefore write
$\oldXi'=r^{-1}(\alpha)$ for some $1$-suborbifold of $\oldGamma$. The local definition of $\frakU$ also implies that $|\alpha|$ is a component of $\fraks_\oldGamma$. Hence in this case there is a
$2$-orbifold
$\frakB$ such that
the orbifold with $|\frakB|=|\oldGamma|$ and
$\fraks_{\frakB}=\overline{\fraks_\oldGamma-|\alpha|}$. In the case where $\oldGamma=\frakU'$, we set $\frakB=\oldGamma$. In either case we have $|\frakB|=|\oldGamma|$, and since $|\oldOmega|=|\oldLambda|$, we may regard $|r|$ as a map of sets from $|\oldLambda|$ to $|\frakB|$.


If $z\in\frakB$ is given, then since $|\frakB|=|\oldGamma|$, we may regard $z$ as a point of $\oldGamma$. 
Choose a neighborhood $\frakV$ of $z$ in $\oldGamma$, and a regular covering map
$P:D^2\times S^1\to r^{-1}(\frakV)$, having the properties stated in
\ref{better onion}. Then $|\frakV|=|\frakW|$ for some neighborhood $\frakW$ of $z$ in $\frakB$. Since $r^{-1}(\frakW)$ is a suborbifold of $\oldOmega$, there is
a suborbifold $\frakY$ of $\oldLambda$ such that $|\frakY|=|r^{-1}(\frakW)|$.
If alternative (II) or (III) of Condition (b) of
\ref{better onion} holds, let us choose a component of $(D^2\times
S^1)-P^{-1}(|\oldXi'|)$ and denote its closure by $Z$. If Alternative (I)
holds we set $Z=D^2\times S^1$. In view of the
final assertion of \ref{better onion}, $R:=P|Z$ is a
covering map from $Z$ to $\frakY$. 

On the other hand, if Alternative (I), (II) or (III) holds, then $Z$ is naturally identified, respectively, with $D^2\times S^1$, or $D^2\times\gamma$ where $\gamma$ is a semicircular arc, or $H\times S^1$ where $H$ is a half-disk. Furthermore, in these respective cases we have $P^{-1}(\oldXi')=\emptyset$,  $P^{-1}(\oldXi')=D^2\times\partial\gamma$, or $P^{-1}(\oldXi')=\ell\times S^1$ where $\ell\subset\partial H$ is the diameter of the half-disk $H$. It follows that if $\oldGamma=\frakU$, so that (I) or (III) holds, then $|r|\circ |R|
:Z\to |\frakW|$ may be regarded as an $S^1$-fibration of $Z$ over $\frakW $ in which $R^{-1}(\oldXi')$ 
is saturated; and that if $\oldGamma=\frakU'$, so that (II) holds, then 
$|r|\circ |R|$
may be regarded as an $I$-fibration of $Z$ over $\frakW $
 in which $\partial_hX=R^{-1}(\oldXi)$. Combining this with the fact that
$R:Z \to Y$ is a
covering map, we deduce that 
if $\oldGamma=\frakU$ then  $|r|$ defines an $S^1$-fibration of $\oldLambda$ over $\frakB $ in which $\oldXi'$ 
is saturated, and that if $\oldGamma=\frakU'$ then 
$|r|$ defines an $I$-fibration of
 $\oldLambda$ over $\frakB $ for which $\partial_h\oldLambda=\oldXi'$.
In the respective cases,
$(\oldLambda,\oldXi)$ is by definition a \bindinglike\ \spair\ or a \pagelike\ \spair.
\EndProof

\Proposition\label{butthurt}
Let 
$(\oldLambda,\oldXi)$ be an acylindrical pair which 
is not an \spair. Suppose that $\oldLambda$ is strongly \simple, 
and that every component of
$(\partial\oldLambda)-\inter\oldXi$ is an annular orbifold.
Let
$p:\toldLambda\to\oldLambda$ be
a
finite-sheeted covering of $\oldLambda$, and suppose that $\toldLambda$ is strongly \simple. Then
$(\toldLambda,p^{-1}(\oldXi))$ is acylindrical and is not an \spair. 
\EndProposition

\Proof
Since $(\oldLambda,\oldXi)$ is an acylindrical pair and $\oldLambda$ is strongly \simple, it follows
from Proposition \ref{silver acylindrical} that 
$\oldOmega:=\silv_\oldXi\oldLambda$ is weakly \simple. Since $(\oldLambda,\oldXi)$
is  not an \spair, it follows from Proposition \ref{new what they look like} that 
$\oldOmega$ does not admit an
$S^1$-fibration. Furthermore, since every component of 
$(\partial\oldLambda)-\inter\oldXi$ is annular, every component of
$\partial\oldOmega$ has Euler characteristic $0$ and is therefore
toric. Since $\oldLambda$ is orientable by the definition of an acylindrical pair, $\fraks_\oldOmega$ has no $0$-dimensional components. It now follows from Lemma \ref{ztimesz-lemma} that for every rank-$2$ free abelian subgroup $H$ of $\pi_1(\oldOmega)$, there is a component $\frakK$ of $\partial\oldOmega$ such that $H$ is contained in a conjugate of the image of the inclusion homomorphism $\pi_1(\frakK)\to\pi_1(\oldOmega)$. 
Furthermore, no finite-sheeted cover of $\oldOmega$ admits an $S^1$-fibration; hence by Proposition
\ref{new what they look like},
$(\toldLambda,\toldXi)$ is not an \spair..

If we set $\toldXi=p^{-1}(\oldXi)$, then $\toldOmega:=\silv_{\toldXi}\toldLambda$ is a finite-sheeted cover of $\oldOmega$.

To show that that $(\toldLambda,\toldXi)$ is acylindrical, suppose to the contrary that there is an essential annular $2$-orbifold $\oldPi$ in the pair $(\toldLambda,\toldXi)$. Let $q:\toldLambda'\to\toldLambda$ be a finite-sheeted covering such that $\toldLambda'$ is an orientable manifold. Set $\toldXi'=q^{-1}(\toldXi)$. Then $\toldOmega':=\silv_{\toldXi'}\toldLambda'$ is a finite-sheeted covering of $\toldOmega$, and therefore of $\oldOmega$. 
Since 
the manifold $M:=D_{\toldXi'}\toldLambda'$ is a two-sheeted cover of $\toldOmega'$, it follows that $M$ is a finite-sheeted cover of $\oldOmega$. Hence every rank-$2$ free abelian subgroup of $\pi_1(M)$ is carried up to conjugacy by some component of $\partial M$.

Now fix a component $A$ of $q^{-1}(\oldPi)$. Since $\toldLambda'$ is a manifold, $A$ is an annulus. Since $\oldPi$ is essential in $(\toldLambda,\toldXi)$, it follows from Corollary \ref{covering annular} that $A$ is essential in $(\toldLambda',\toldXi')$. In particular $A$ is $\pi_1$-injective in $\toldLambda$, and is not parallel in $(\toldLambda',\toldXi')$ to an annulus in $\toldXi'$. Hence the torus $T:=DA$ is $\pi_1$-injective in $M$. Thus the image of the inclusion $\pi_1(T)\to\pi_1(M)$ is a free abelian group of rank $2$, and is therefore carried up to conjugacy by some component of $\partial M$. This implies (via a standard application of \cite[Lemma 5.1]{Waldhausen}) that $T$ is boundary-parallel in $M$; that is, there is a submanifold $K$ of $M$, homeomorphic to $T^2\times[0,1]$, with $\Fr_M K=T$. The canonical involution $\iota$ of $M=D_{\toldXi'}\toldLambda'$, which interchanges the two copies of $\toldLambda'$ in $M$, leaves $T=DA$ invariant; if $\iota$ were to interchange $K$ and $\overline{M-K}$ then $M$ would be homeomorphic to $T^2\times[-1,1]$, a contradiction to finite volume. Hence $K$ is $\iota$-invariant. Thus for some submanifold $J$ of $\toldLambda'$ we have  $K=D_BJ$, where
$B=\toldXi'\cap J$. Since $K$ is connected, so is $J$.

The essentiality of $A$ implies that $\Fr_{\toldXi'}B=\partial A$ is $\pi_1$-injective in $\toldLambda'$, and in particular in $\toldXi$; hence $B$ is $\pi_1$-injective in $\toldXi$, hence in $\toldLambda'$, and in particular in $J$. This implies that $J$ is $\pi_1$-injective in $K$, so that $\pi_1(J)$ is free abelian of rank at most $2$. But $\rank\pi_1(J)\ne2$ since $\oldLambda$ is strongly \simple, and $\rank\pi_1(J)\ne0$ since $J$ contains the $\pi_1$-injective torus $A$. Hence $\pi_1(J)$ is infinite cyclic, so that $J$ is a solid torus (see \cite[Theorem 5.2]{hempel}). Since $A$ is a $\pi_1$-injective annulus on the torus $T$, the image of the inclusion homomorphism $\pi_1(A)\to\pi_1(K)$ is a maximal cyclic subgroup of $\pi_1(K)$, and hence the inclusion homomorphism $\pi_1(A)\to\pi_1(J)$ is an isomorphism. If we set $A'=J\cap\overline{(\partial\toldLambda')-\toldXi}$, then $T':=DA'$ is a boundary torus of $K$, and hence $A'$ is an annulus which is homotopically non-trivial in $K$ and therefore in $J$. It now follows that the triad $(J,A,A')$ homeomorphic to $(A,A\times\{1\},A\times\{0\})$, so that
$A$ is parallel in $(\toldLambda',\toldXi')$ to the component $A'$ of $\overline{(\partial\toldLambda')-\toldXi}$. This contradicts the essentiality of $A$.
\EndProof

\section{The characteristic suborbifold}\label{characteristic section}

The proof of the following result, which is the basis of our treatment of the characteristic suborbifold theory, was suggested by the discussion on p. 445 of \cite{bonahon-siebenmann}.

\Proposition\label{new characteristic} Let $\oldPsi$ be a $3$-orbifold
which is orientable, componentwise strongly \simple, and componentwise
boundary-irreducible. 
Then up to (orbifold) isotopy there
exists a unique properly embedded $2$-suborbifold $\frakQ$ of $\oldPsi$
such that (1) each component of $\frakQ$ is two-sided and annular and is essential (see
\ref{acylindrical def}) in $\oldPsi$, (2) each component
$\oldLambda$ of $\overline{\oldPsi-\frakH}$, where $\frakH$ is a
strong regular neighborhood of $\frakQ$ in $\oldPsi$, 
is either an \Ssuborbifold \ or an \Asuborbifold\ of $\oldPsi$ (see Definition \ref{Do I need it?}), and (3) Condition (2) becomes false if $\frakQ$
is replaced by the union of any proper subset of its components. 
\EndProposition

The two alternatives given in Condition (2) of the above statement are not mutually exclusive: a component $\oldLambda$ of 
$\overline{\oldPsi-\frakH}$ may be both an \Ssuborbifold\ and an \Asuborbifold.

\Proof[Proof of Proposition \ref{new characteristic}]
\nonessentialproofreadingnote{I had a note here saying ``This has been revised a lot, and I was confused
  about things like the strong \simple ity of $\oldLambda$. Be very
  careful about the proofreading.'' After the mega-proofreading, I looked at that note and noticed that the sentence that now begins the proof of \ref{new-claim} was at the beginning of the proof of the lemma, where it made no sense! What am I supposed to do about these things?}
We may assume without loss of generality that $\oldPsi$ is connected,
and therefore strongly \simple\ and boundary-irreducible.

Set $\oldOmega=\silv\oldPsi$. It follows from Proposition \ref{silver
  irreducible}, applied with $\partial\oldPsi$ playing the role of $\oldXi$, that $\oldOmega$
is irreducible; indeed, Conditions (a) and (b) of Proposition \ref{silver
  irreducible} follow respectively from the strong \simple ity (cf. \ref{oops}) and the boundary-irreducibility
of $\oldPsi$. It then
follows from the PL version of \cite[p. 444, Splitting Theorem
1]{bonahon-siebenmann} (see \ref{categorille}) that there 
exists, up to (orbifold) isotopy, a unique $2$-dimensional suborbifold
$\frakR$ of $\oldOmega$ such that (1\,$'$) each component of $\frakR$
is toric and is
incompressible in $\oldOmega$, (2\,$'$) each component 
of
$\oldOmega\cut\frakR$ 
either is weakly \simple\ or admits
an $S^1$-fibration (or both), and (3\,$'$)
Condition (2\,$'$) becomes false if $\frakR$ is replaced by the union of
any proper subset of its components. 

We claim:
\Claim\label{new-claim}
Let  $\frakQ$ be any two-sided properly embedded suborbifold of $\oldPsi$.
Then Conditions (1), (2) and (3) of the statement of the present proposition hold for $\frakQ$ if and only if Conditions (1\,$'$), (2\,$'$) and (3\,$'$) hold when we set
$\frakR=\silv\frakQ$. 
\EndClaim

To prove \ref{new-claim}, 
fix a strong regular neighborhood $\frakH$ of $\frakQ$.
Since $\overline{\oldPsi-\frakH}$ is homeomorphic to
  $\oldPsi\cut\frakQ$, it follows from Lemma \ref{oops lemma} that  $\overline{\oldPsi-\frakH}$
is componentwise strongly \simple.

Set $\frakR=\silv\frakQ$ as in the
statement of \ref{new-claim}.  Note that since 
 $\oldPsi$ is strongly \simple, no component of $\frakQ$ is an incompressible
 toric suborbifold (see \ref{oops}). Hence a necessary condition for all components of $\frakR$ to
 be incompressible and toric is that all components of $\frakQ$ be
 annular. It then follows from 
the first
  assertion of Proposition \ref{silver acylindrical},  applied with
  $\oldPsi$ and $\partial\oldPsi$
  playing the roles of $\oldLambda$ and $\oldXi$, that Condition (1\,$'$)
  for $\frakR$ is equivalent to Condition (1) for $\frakQ$. (The condition in
  Proposition \ref{silver acylindrical} that $\frakR$ not be parallel to
a boundary component of $\oldOmega$ holds vacuously since $\oldOmega$ is closed.)

The second step in the proof of \ref{new-claim} is to show that if Condition (1) holds for $\frakQ$, then
$\frakR$ satisfies Condition (2\,$'$) if and only if  $\frakQ$
satisfies Condition (2). For this purpose, first observe that
 if  $\oldLambda$ is any  component of
$\overline{\oldPsi-\frakH}$, then
$\silv_{\oldLambda\cap\partial\oldPsi}\oldLambda$ is (orbifold-)homeomorphic to 
 some component
of
$\oldOmega\cut\frakR$, and conversely that 
 every component
 of
$\oldOmega\cut\frakR$
is homeomorphic to
$\silv_{\oldLambda\cap\partial\oldPsi}\oldLambda$ for some
component
$\oldLambda$  of
$\overline{\oldPsi-\frakH}$. 
 Thus Condition (2\,$'$) is equivalent to the assertion that for each component  $\oldLambda$ of  
$\overline{\oldPsi-\frakH}$, the orbifold
$\silv_{\oldLambda\cap\partial\oldPsi}\oldLambda$ either is weakly
 \simple\ or admits
an $S^1$-fibration. But since $\overline{\oldPsi-\frakH}$ is
  componentwise strongly \simple, $\oldLambda$ is strongly \simple, and
  we may therefore apply the second assertion of Proposition \ref{silver
  acylindrical}, with $\oldLambda\cap\partial\oldPsi$ playing the role
of $\oldXi$, to deduce that  $\silv_{\oldLambda\cap\partial\oldPsi}\oldLambda$ is
weakly \simple\ if and only if
the pair
$(\oldLambda,\oldLambda\cap\partial\oldPsi)$
is acylindrical.
Likewise, Proposition
\ref{new what
    they look like} implies that
 $\silv_{\oldLambda\cap\partial\oldPsi}\oldLambda$ admits
an $S^1$-fibration if and only if
$(\oldLambda,\oldLambda\cap\partial\oldPsi)$ is an \spair.
But Condition (1) for $\frakQ$ implies that  the
components of $\Fr_\oldPsi\oldLambda$ are essential annular suborbifolds of $\oldPsi$, and hence that
$(\oldLambda,\oldLambda\cap\partial\oldPsi)$
is an acylindrical pair or an \spair\ if and only if $\oldLambda$ is an \Asuborbifold\ or \Ssuborbifold\ of
$\oldPsi$, respectively.
This shows that Condition  (2\,$'$) for $\frakR$ is equivalent to Condition (2) for  $\frakQ$, provided that Condition (1) holds for $\frakQ$. 

The third and final step in the proof of \ref{new-claim} is to show
that if Condition (1) holds for $\frakQ$, then $\frakR$ satisfies Condition (3\,$'$) if and only if
$\frakQ$ satisfies Condition (3). For this purpose,
first observe that
 if $\frakQ'$ is the union
of a proper subset of the components of $\frakQ$, then
$\frakR':=\silv\frakQ'\subset\oldOmega$ is the union
of a proper subset of the components of $\frakR$; and conversely, that 
if $\frakR'$ denotes the union
of a proper subset of the components of $\frakR$, then we may write
$\frakR'=\silv\frakQ'$, where
$\frakQ'$ is the union
of some proper subset of the components of $\frakQ$. But if Condition (1) holds for $\frakQ$, and if 
$\frakQ'$ is the union
of a proper subset of the components of $\frakQ$, then the second step of the proof of \ref{new-claim}, applied with $\frakQ'$ in place of $\frakQ$, shows that 
$\frakR':=\silv\frakQ'$ satisfies Condition (2\,$'$) if and only if $\frakQ'$ satisfies Condition (2). This shows that $\frakR$ fails to satisfy Condition (3\,$'$) if and only if  $\frakQ$ fails to satisfy Condition (3), and completes the proof of \ref{new-claim}.



Now by the existence assertion of the PL version of \cite[p. 444, Splitting Theorem
1]{bonahon-siebenmann}, we may fix 
a $2$-suborbifold
$\frakR_0$ of $\oldOmega$ such that Conditions (1\,$'$)---(3\,$'$) hold with $\frakR_0$ in place of $\frakR$.
Since $\frakR_0$ is in particular a closed suborbifold of
$\oldOmega=\silv\oldPsi$, there is a unique properly
embedded $2$-suborbifold   $\frakQ_0$ of $\oldPsi$ such that
$\silv\frakQ_0=\frakR_0$.  
By \ref{new-claim}, Conditions (1)---(3) hold with $\frakQ_0$ in place of $\frakQ$. Now let $\frakQ$ be any properly embedded  $2$-suborbifold of $\oldPsi$ such that (1)---(3) hold, and set $\frakR=\silv\frakQ$. Then by \ref{new-claim}, Conditions (1\,$'$)---(3\,$'$) hold for $\frakR$. According to
the uniqueness assertion of the PL version of \cite[p. 444, Splitting Theorem
1]{bonahon-siebenmann}, $\frakR=\silv\frakQ$ is (orbifold-)isotopic to $\frakR_0=\silv\frakQ_0$, which implies (upon composing an ambient orbifold isotopy between $\silv\frakQ$ and
 $\frakR_0$
with the natural immersion $\silv\oldPsi\to\oldPsi$) that $\frakQ$ is (orbifold-)isotopic to $\frakQ_0$. 
\EndProof

\DefinitionNotation\label{oldSigma def}
If $\oldPsi$ is any componentwise strongly \simple, componentwise boundary-irreducible, orientable $3$-orbifold, we will denote by $\frakQ(\oldPsi)$ the suborbifold  $\frakQ$ of $\oldPsi$ (defined up to isotopy) given by Proposition \ref{new characteristic}, and we will denote by $\frakH(\oldPsi) $ a strong regular neighborhood of $\frakQ(\oldPsi)$ in $\oldPsi$. 
We will denote by
$\oldSigma_1(\oldPsi)$ the union of all
components  of $\overline{\oldPsi- \frakH(\oldPsi)} $ that are \Ssuborbifold s (see \ref{Do I need it?}) of $\oldPsi$,
and by $\oldSigma_2(\oldPsi)$ the union of all components of $\frakH(\oldPsi)$ that do not meet any component of $\oldSigma_1(\oldPsi)$.
We will define the (relative)
{\it characteristic suborbifold} of $\oldPsi$, denoted
$\oldSigma(\oldPsi)$, to be
$\oldSigma_1(\oldPsi)\cup\oldSigma_2(\oldPsi)$. 
Thus $\frakH(\oldPsi) $, $\oldSigma_1(\oldPsi)$, $\oldSigma_2(\oldPsi)$ and $\oldSigma(\oldPsi)$,  as well as $\frakQ(\oldPsi) $, are well-defined up to
(orbifold) isotopy. Note that each component of $\oldSigma(\oldPsi)$ is an \Ssuborbifold\ of $\oldPsi$.
\EndDefinitionNotation


\Number\label{tuesa day}
If 
$\oldPsi$ is a componentwise strongly \simple, componentwise
boundary-irreducible, orientable $3$-orbifold, we will set
$\oldPhi(\oldPsi)=\oldSigma(\oldPsi)\cap\partial\oldPsi$ and $\cala(\oldPsi)=\Fr_{|\oldPsi|}|\oldSigma(\oldPsi)|$, so that $|\oldPhi(\oldPsi)|$ and $\cala(\oldPsi)$ are
submanifolds of $\partial|\oldSigma(\oldPsi)|$ whose union is $|\partial\oldSigma(\oldPsi)|$, and
$|\oldPhi(\oldPsi)|\cap\cala(\oldPsi)=\partial |\oldPhi(\oldPsi)|=\partial\cala(\oldPsi)$.
It follows from Proposition \ref{new characteristic} and the definition of $\oldSigma(\oldPsi)$  
that every component of $\omega(\cala(\oldPsi))$ is an orientable
annular $2$-orbifold, essential in $\oldPsi$. The essentiality of the
components of $\omega(\cala(\oldPsi))$ implies that $\partial\oldPhi(\oldPsi)=
\partial\omega(\cala(\oldPsi))$ is $\pi_1$-injective in $\oldPsi$, and
in particular in $\partial\oldPsi$. 

\nonessentialproofreadingnote{I have to say I can't think of any reason why $\cala(\oldPsi)$, rather than $\obd(\cala(\oldPsi))$, should be treated as the basic object as far as notation is concerned. The question is whether I feel up to making the changes (including thinking up notation for 
$\obd(\cala(\oldPsi))$).}

By \ref{boundary is negative}, the componentwise strong \simple ity and componentwise boundary-irreducibility of $\oldPsi$ imply that each component of
$\partial\oldPsi$ has strictly negative Euler characteristic. In view of the
$\pi_1$-injectivity of
$\partial\oldPhi(\oldPsi)$ in $\partial\oldPsi$, it then follows that
each component of $\oldPhi(\oldPsi)$ or $\overline{\partial\oldPsi-\oldPhi(\oldPsi)}$ has non-positive Euler
characteristic.

Suppose that $\oldLambda$ is a component of $\oldSigma(\oldPsi)$ or of $\overline{\oldPsi-\oldSigma(\oldPsi)}$. Since $\oldLambda$ is a compact $3$-orbifold, we have $\chi(\oldLambda)=\chi(\partial\oldLambda)/2$. On the other hand, since each component of $\Fr\oldLambda$ is a component of $\obd(\cala(\oldPsi))$ and is therefore annular, we have $\chi(\Fr\oldLambda)=0$. It follows that
$\chi(\oldLambda\cap\partial\oldPsi)=\chi(\partial\oldLambda)$, and hence that
 $\chi(\oldLambda)=\chi(\oldLambda\cap\partial\oldPsi)/2$. But $\chi(\oldLambda\cap\partial\oldPsi)\le0$ since $\oldLambda\cap\partial\oldPsi$ is a union of components of $\oldPhi(\oldPsi)$ or of $\overline{\partial\oldPsi-\oldPhi(\oldPsi)}$. Hence $\chi(\oldLambda)\le0$, or equivalently $\chi(\partial\oldLambda)\le0$.

 We will denote by $\book(\oldPsi)$ the union of
$\oldSigma(\oldPsi)$ with all components $\frakC$ of $\frakH(\oldPsi)$ such that $\frakC\cap\oldSigma(\oldPsi)=\Fr_\oldPsi\frakC$.
(The notation is meant to suggest that $\book(\oldPsi)$ is an orbifold
version of a ``book of $I$-bundles.'' See \cite{hyperhaken},
\cite{canmac}, \cite{inject}.) We will
set $\kish(\oldPsi)=\overline{\oldPsi-\book(\oldPsi)}$. Thus $\book(\oldPsi)$ and
$\kish(\oldPsi)$ are suborbifolds of $\oldPsi$.
\EndNumber


\Lemma\label{i see a book}
Let $\oldPsi$ be a $3$-orbifold
which is orientable, componentwise strongly \simple, and componentwise
boundary-irreducible. Let $\frakB$ be a compact
suborbifold of $\oldPsi$. Then the following conditions are
equivalent:
\begin{itemize}
\item $\frakB$ is (orbifold-)isotopic to  $\book(\oldPsi)$ in
  $\oldPsi$.
\item There exist a properly embedded $2$-suborbifold $\oldPi$ of
  $\oldPsi$, and a strong
regular neighborhood $\oldGamma$  of $\oldPi$ in $\oldPsi$,
such that (a) each component of $\oldPi$ is two-sided and annular and is essential
in $\oldPsi$, (b) each component
$\oldLambda$ of $\overline{\oldPsi-\oldGamma}$
is either an \Asuborbifold\ or an \Ssuborbifold\ of $\oldPsi$,
and (c)
$\frakB$ is the union of $\oldGamma$ with all components of 
$\overline{\oldPsi-\oldGamma}$ that are \Ssuborbifold s.
\end{itemize}
\EndLemma

\Proof
Let $\calp$ denote the set of all properly embedded $2$-suborbifolds
$\oldPi$ of $\oldPsi$ such that every component of $\oldPi$ is two-sided and annular
and is essential in $\oldPsi$. Let $\calp_0$ denote the set of all
elements $\oldPi$ of $\calp$ such that 
 each component
of $\overline{\oldPsi-\oldGamma}$, where $\oldGamma$ is a strong regular neighborhood of $\oldPi$ in $\oldPsi$,
is either an \Asuborbifold\ or an \Ssuborbifold.
For each $\oldPi\in\calp$, if $\oldGamma$ is a strong regular neighborhood of $\oldPi$ in $\oldPsi$, let $\frakB_\oldPi$ denote the union of $\oldGamma$ with all components of 
$\overline{\oldPsi-\oldGamma}$ that are \Ssuborbifold s of $\oldPsi$; thus $\frakB_\oldPi$ is determined up to isotopy by $\oldPi$. We claim:

\Claim\label{untold us before}
If $\oldPi_0$ and $\oldPi$ are elements of $\calp_0$ and $\calp$ respectively, such that $\oldPi_0\subset\oldPi$, then $\oldPi\in\calp_0$, and $\frakB_{\oldPi}$ is isotopic to $\frakB_{\oldPi_0}$.
\EndClaim

To prove \ref{untold us before}, it is enough to consider the case in which $\oldPi-\oldPi_0$ has exactly one component, as the general case then follows by induction on the number of components of $\oldPi-\oldPi_0$. 
In this case, let $\frakA$ denote the orientable annular orbifold $\oldPi-\oldPi_0$. A strong regular neighborhood of $\oldPi$ in $\oldPsi$ may be written as $\oldGamma=\oldGamma_0\discup\oldGamma_\frakA$, where $\oldGamma_0$ and $\oldGamma_\frakA$ are strong regular neighborhoods of $\oldPi_0$ and $\frakA$ respectively. Let $\oldLambda$ denote the component of $\overline{\oldPsi-\oldGamma_0}$ containing  $\oldGamma_\frakA$. Since $\oldPi_0\in\calp_0$, the suborbifold
$\oldLambda$ of $\oldPsi$ is either an \Asuborbifold\ or an
\Ssuborbifold. 

Consider first the subcase in which  $\frakA$ is essential in the pair
$(\oldLambda,\oldLambda\cap\partial\oldPsi)$. In
this subcase, 
it follows from the definition of an \Asuborbifold\ that
$\oldLambda$ is not an \Asuborbifold\ of $\oldPsi$, and is therefore an
\Ssuborbifold\ of $\oldPsi$. Hence 
 it follows from Lemma \ref{oops lemma} that $\oldLambda$ is strongly \simple. 
Thus Corollary \ref{vertical corollary} implies that $\oldLambda$ admits an orbifold fibration, compatible with $\oldLambda\cap\partial\oldPsi$, in which $\frakA$ is saturated. Hence every component of $\overline{\oldLambda-\oldGamma_\frakA}$ is an \Ssuborbifold\ of $\oldPsi$. Since the components of $\overline{\oldPsi-\oldGamma}$ are precisely the components of $\overline{\oldPsi-\oldGamma_0}$ distinct from $\oldLambda$ and the components of $\overline{\oldLambda-\oldGamma_\frakA}$, it follows that in this subcase we have $\oldPi\in\calp_0$ and $\frakB_\oldPi=\frakB_{\oldPi_0}$.

Now consider the subcase in which $\frakA$ is non-essential in 
the pair $(\oldLambda,\oldLambda\cap\partial\oldPsi)$.  Since $\frakA\in\calc(\oldPi)$ and $\oldPi\in\calp$, the orientable annular orbifold $\frakA$ is essential in $\oldPsi$; in particular, $\frakA$ is $\pi_1$-injective in $\oldLambda$, and is not parallel in the pair
$(\oldLambda,\oldLambda\cap\partial\oldPsi)$ to a suborbifold of 
$\oldLambda\cap\partial\oldPsi$. Since $\frakA$ is not essential in $(\oldLambda,\oldLambda\cap\partial\oldPsi)$, it must be
parallel in the pair
$(\oldLambda,\oldLambda\cap\partial\oldPsi)$ to a component of $\Fr_\oldPsi\oldLambda$.
Hence $\overline{\oldLambda-\oldGamma_\frakA}$ has some component $\frakJ$ such that $\frakJ^+:=\frakJ\cup\oldGamma_\frakA$ is a strong regular neighborhood of $\frakA$ in $\oldPsi$. 
Thus $\frakJ^+$ is a \bindinglike\ \Ssuborbifold\ of $\oldPsi$. Furthermore,  $\oldLambda':=\overline{\oldLambda-\frakJ^+}$, is isotopic to
$\oldLambda$,
and is therefore either
 an \Asuborbifold\ or an \Ssuborbifold\ of $\oldLambda$. Since the components of
$\overline{\oldPsi-\oldGamma}$ are precisely $\frakJ$, $\oldLambda'$,
and the components of $\overline{\oldPsi-\oldGamma_0}$ distinct from
$\oldLambda$, we have $\oldPi\in\calp_0$. If $\oldLambda$ is an
\Ssuborbifold\ of $\oldPsi$, we have
$\frakB_\oldPi=\frakB_{\oldPi_0}$. If $\oldLambda$ is not an
\Ssuborbifold\ of $\oldPsi$, we have
$\frakB_\oldPi=\frakB_{\oldPi_0}\cup\frakJ^+$; furthermore,
$\frakB_{\oldPi_0}\cap\frakJ^+$ is then a component of
$\Fr_\oldPsi\oldLambda$, and is the component of the horizontal
boundary of a trivial $I$-fibration of $\frakJ^+$ for which the
vertical boundary is $\frakJ^+\cap\partial\oldPsi$. Hence $\frakB_\oldPi$ is isotopic to $\frakB_{\oldPi_0}$ in this situation. This completes the proof of \ref{untold us before}.

Now set $\frakQ=\frakQ(\oldPsi)$, $\frakH=\frakH(\oldPsi)$, $\oldSigma_i=\oldSigma_i(\oldPsi)$ for $i=1,2$, and $\oldSigma=\oldSigma(\oldPsi)$ (see \ref{oldSigma def}). Note that according to \ref{oldSigma def} and Proposition \ref{new characteristic} we have $\frakQ\in\calp_0$. We claim:
\Claim\label{just as you thought}
$\frakB_{\frakQ}$ is isotopic to $\book(\oldPsi)$.
\EndClaim

To prove \ref{just as you thought}, first note that we may write
$\frakH=\oldSigma_2\discup\frakH_1\discup\frakH_2$, where for $m=1,2$
we denote by $\frakH_m$  the union of all components of $\frakH$ that
meet $\oldSigma_1$ in exactly $m$ frontier components. 
The definition of $\book(\oldPsi)$ gives that $\book(\oldPsi)= \oldSigma\cup\frakH_2=\oldSigma_1\cup\oldSigma_2\cup\frakH_2$.
On the other hand, the definitions imply that $\frakH_1$ is disjoint from $\oldSigma_2$ and $\frakH_2$, and that each component of $\frakH_1$ is the strong regular neighborhood of a two-sided, properly embedded suborbifold of $\oldPsi$ and meets $\oldSigma_1$ in exactly one frontier component. Hence $\frakB_\frakQ=\frakH\cup\oldSigma_1=\oldSigma_1\cup \oldSigma_2\cup\frakH_1\cup\frakH_2$ is isotopic to $\book(\oldPsi)$, as asserted by \ref{just as you thought}.

The assertion of the present lemma may be paraphrased as saying that a
compact suborbifold $\frakB$ of $\oldPsi$ is isotopic to
$\book(\oldPsi)$ if and only if it is isotopic to $\frakB_\frakP$ for
some $\frakP\in\calp_0$. If $\frakB$ is isotopic to $\book(\oldPsi)$,
then by \ref{just as you thought}, $\frakB$ is isotopic to
$\frakB_\frakQ$, and we have $\frakQ\in\calp_0$. Conversely, suppose
that is isotopic to $\frakB_\frakP$ for some $\frakP\in\calp_0$. Among
all elements of $\calp_0$ that are
 unions of components of $\oldPi$, choose one, $\oldPi_0$,
which has the smallest number of components, and is therefore minimal with respect to inclusion. It then follows from \ref{oldSigma def} and Proposition \ref{new characteristic} that $\oldPi_0$ is isotopic to $\frakQ$. Hence by \ref{just as you thought}, $\frakB_{\oldPi_0}$ is isotopic to $\book(\oldPsi)$. But it follows from \ref{untold us before} that $\frakB_{\oldPi}$ is isotopic to $\frakB_{\oldPi_0}$, and hence to $\book(\oldPsi)$.
\EndProof

\Proposition\label{less than nothing}
Let $\oldPsi$ be a $3$-orbifold
which is orientable, componentwise strongly \simple, and componentwise
boundary-irreducible. Then every component of $\kish\oldPsi$ is either closed or has strictly negative Euler characteristic.
\EndProposition

\Proof
Let $\frakK$ is a component of $\kish\oldPsi$ with $\partial\frakK\ne\emptyset$. By Lemma \ref{oops lemma}, $\frakK$ is componentwise strongly \simple\. In particular, no component of $\partial\frakK$ is spherical. Hence if $\chi(\frakK)\ge0$ then every component of $\partial\frakK$ is toric, and by Proposition \ref{three-way equivalence}, $\frakK$ is a \torifold. By Lemmaref \ref{when a tore a fold}, $\frakK$ is a \bindinglike\ \Ssuborbifold of $\oldPsi$. But by Lemma \ref{i see a book}, no component of $\kish\oldPsi=\overline{\oldPsi-\book\oldPsi}$ is an \Ssuborbifold\ of $\oldPsi$.
\EndProof

\Proposition\label{Comment A}If $p:\toldPsi\to \oldPsi$ is a covering map
of orientable $3$-orbifolds that are 
componentwise strongly \simple\ and componentwise boundary-irreducible, then $p^{-1}(\book(\oldPsi))$ is orbifold-isotopic to
$\book(\toldPsi)$. 
\EndProposition

\Proof
We may assume without loss of generality that $\oldPsi$ and $\toldPsi$
are connected (and are therefore \simple\ and boundary-irreducible).
Set $\frakB=\book(\oldPsi)$. Then it follows from Lemma \ref{i see
  a book} that there exist a properly embedded $2$-suborbifold $\oldPi$ of
  $\oldPsi$, and a strong
regular neighborhood $\oldGamma$  of $\oldPi$ in $\oldPsi$,
such that Conditions (a), (b) and (c) of that lemma hold. Since, by Condition (a),
 each component of $\oldPi$ is two-sided and annular and is essential in
 $\oldPsi$, it follows from Corollary \ref{covering annular} that each component of the properly embedded $2$-orbifold
 $\toldPi:=p^{-1}(\oldPi)$ is two-sided and annular and is essential in
 $\toldPsi$. 

Let
$\oldLambda$ be any component of $\overline{\oldPsi-\oldGamma}$, and
set $\toldLambda:=p^{-1}(\oldLambda)$. Then $\Fr_{\toldPsi}\toldLambda$ is a union
of components of $\toldPi$, which are essential orientable annular
orbifolds in $\toldPsi$. Since $\toldPsi$ is componentwise
strongly \simple, it follows from Lemma \ref{oops lemma} that
$\toldLambda$ is strongly \simple. On the other hand, by Condition (b), either $\oldLambda$ is an \Ssuborbifold\ of
$\oldPsi$, which immediately implies that $\toldLambda$ is a \Ssuborbifold\ of
$\toldPsi$;
or $\oldLambda$ is an \Asuborbifold\ but not an \Ssuborbifold\ of
 $\oldPsi$. In the latter case,
the pair
$(\oldLambda,\oldLambda\cap\partial\oldPsi)$ is acylindrical and is
not an \spair, and each component of $\overline{(\partial\oldLambda)\setminus\partial\oldPsi}=\Fr_\oldPsi{\oldLambda}
$ is annular. 
 In this case, in view of the strong \simple ity of $\toldLambda$, it follows from Proposition \ref{butthurt} that
$(\toldLambda,\toldLambda\cap\partial\toldPsi)$ is acylindrical and is not an \spair. Thus $\toldLambda$ is either an \Ssuborbifold\ or an \Asuborbifold\ of $\toldPsi$ 
for every component $\oldLambda$ of $\overline{\oldPsi-\oldGamma}$; and $\toldLambda$ is an \Ssuborbifold\ if and only if $\oldLambda$ is an \Ssuborbifold. 
Hence 
if we denote by $\toldGamma$ the strong regular neighborhood $p^{-1}(\oldGamma)$ of
$\toldPi$, then each component of $\overline{\toldPsi-\toldGamma}$ is
either an \Asuborbifold\ or an \Ssuborbifold\ of $\toldPsi$. Furthermore, since Condition (c) of Lemma \ref{i see
  a book} asserts that $\frakB$ is the union of $\oldGamma$ with all components of 
$\overline{\oldPsi-\oldGamma}$ that are \Ssuborbifold s\ of $\oldPsi$,
it now follows that $\tfrakB:=p^{-1}(\frakB)$ is the union of $\toldGamma$ with all components of 
$\overline{\toldPsi-\toldGamma}$ that are \Ssuborbifold s of
$\toldPsi$. Thus Conditions (a), (b) and (c) of Lemma \ref{i see a
  book} hold with $\toldPsi$, $\toldPi$, 
$\toldGamma$ and $\tfrakB$ in place of
$\oldPsi$, $\oldPi$, 
$\oldGamma$ and $\frakB$. It therefore follows from Lemma \ref{i see
  a book} that $\tfrakB$ is isotopic to $\book(\toldPsi)$.
\EndProof

\Number\label{manifolds are different}We observe that the treatment of the characteristic suborbifold in this section covers only the case in which the given orbifold is strongly \simple. In the proofs of Lemma \ref{graphology} and Proposition \ref{crust chastened}, where we consider the characteristic submanifolds of $3$-manifolds which, regarded as orbifolds, are not strongly \simple, we will refer directly to the ``classical'' sources
\cite{Jo} and\cite{js}.
\EndNumber

\section{DARTS volume for orbifolds}\label{darts section}

\Definition\label{hexe}
The {\it Gromov volume} of a closed, orientable topological $3$-manifold $M$ is
defined to be $\vtet\|M\|$, where $\vtet$ denotes the volume of a
regular ideal tetrahedron in $\HH^3$, and $\|M\|$ denotes the Gromov
norm \cite{gromov}, \cite[Chapter C]{bp}. We will denote the Gromov volume of $M$ by $\volG(M)$. 
\EndDefinition 

The relevance of $\volG(M)$ to the present paper arises from Gromov's
result (see
\cite[Theorem C.4.2]{bp}) that for any closed
hyperbolic $3$-manifold $M$ we have $\volG M=\vol M$.

\begin{remarksdefinition}\label{sem ting}
According to \cite[Section 0.2]{gromov}, if $\tM$ is a $d$-sheeted cover
of a
compact, orientable topological $3$-manifold $M$, where $d$ is a positive integer,
we have 
\Equation\label{i'm silly}\volG(\tM)=d\volG(M).
\EndEquation

Now suppose that $\oldOmega$ is a closed, connected, very good topological orbifold. By definition, $\oldOmega$ admits a finite-sheeted covering which is a
manifold. 

Suppose that $\toldOmega_1$ and $\toldOmega_2$ are finite-sheeted manifold
coverings of $\oldOmega$. Then $\toldOmega_1$ and $\toldOmega_2$ admit a common
covering $\toldOmega_3$. If $d_i$ denote the degree of the covering
$\toldOmega_i$ of $\oldOmega$, then for $i=1,2$ we have $d_i|d_3$, and
$\toldOmega_3$ is a $d_3/d_i$-fold covering of $\toldOmega_i$. By (\ref{i'm
  silly}) we have
$(d_3/d_1)\volG(\toldOmega_1)=\volG(\toldOmega_3)=(d_3/d_2)\volG(\toldOmega_2)$,
so that
\Equation\label{smoking fistule}
\volG(\toldOmega_1) /d_1=\volG(\toldOmega_2) /d_2.
\EndEquation

If $\oldOmega$ is a closed, connected, very good topological orbifold, then by
(\ref{smoking fistule})   there is a well-defined invariant $\volG(\oldOmega)$ given by
setting $\volG(\oldOmega)=\volG(\toldOmega)/d$, where $\toldOmega$ is any
finite-sheeted cover of $\oldOmega$ which is a manifold, and $d$ denotes
the degree of the covering $\toldOmega$.

More generally, if $\oldOmega$ is a possibly disconnected, very good, closed,
orientable $3$-orbifold, let us define
$$\volG(\oldOmega)=\sum_{C\in\calc(\oldOmega)}\volG(C).$$
\end{remarksdefinition}

\Notation\label{t-defs}
This subsection will contain definitions of a number of invariants which,
like the invariant $\volG$ discussed above, will turn out to be closely related to 
hyperbolic volume. It is convenient to collect the definitions of all  these
invariants here; their connection with hyperbolic volume will be
established in the course of the section. Many of these invariants will be used in this paper, but there are a few whose importance will become apparent only in \cite{second}. 


For any very good, compact, orientable (but possibly disconnected) topological $3$-orbifold $\oldPsi$, we will set
$$\volorb(\oldPsi)=\frac12\volG(D\oldPsi),$$
where $D\oldPsi$ is defined as in \ref{silvering} (and is very good since $\oldPsi$ is very good and orientable). The acronym ASTA stands for Agol, Storm, Thurston and Atkinson; see \cite{ast} and \cite{atkinson}.
Note that if $\oldPsi$ is a disjoint union of suborbifolds $\oldPsi_1$ and $\oldPsi_2$, we have
$2\volorb(\oldPsi)=\volG(D\oldPsi)=\volG(\oldPsi_1)+\volG(\oldPsi_2)=2\volorb(\oldPsi_1)+2\volorb(\oldPsi_2)$, and hence $\volorb(\oldPsi)=\volorb(\oldPsi_1)+\volorb(\oldPsi_2)$. 

Note that, according to the convention posited in \ref{categorille}, the PL category will be the default category of orbifolds for the rest of this section.

For any very good, compact, orientable (but possibly disconnected) $3$-orbifold $\oldPsi$, we will set
$$\smock_0(\oldPsi)=\sup_\oldPi \volorb(\oldPsi\cut\oldPi),$$ 
where $\oldPi$ ranges
over all (possibly empty) $2$-dimensional suborbifolds of $\inter\oldPsi$
whose components are all incompressible suborbifolds of
$\oldPsi$. {\it A priori}, $\smock_0(\oldPsi)$ is an element of the extended real
number system. 

Note that the definition of $\smock_0(\oldPsi)$ implies that
$\smock_0(\oldPsi)\ge \volorb(\oldPsi\cut\emptyset)$, i.e.
\Equation\label{more kitsch}
\smock_0(\oldPsi)\ge \volorb(\oldPsi).
\EndEquation

We will also set
$$\smock(\oldPsi)= \frac{\smock_0(\oldPsi)}{0.305}.
$$ 

If $\oldPsi_1,\ldots,\oldPsi_m$ are the components of $\oldPsi$, then every
$2$-dimensional suborbifold of $\inter\oldPsi$,
whose components are all incompressible, has the form $\oldPi=\oldPi_1\cup\cdots\cup\oldPi_m$, where $\oldPi_i$ is a $2$-dimensional suborbifold of $\inter\oldPsi$
whose components are all incompressible; furthermore, we have $\volorb(\oldPsi\cut\oldPi)=\volorb((\oldPsi_1)\cut{\oldPi_1}) +\cdots+\volorb((\oldPsi_m)\cut{\oldPi_m})$. It follows that $\smock_0$ is additive over components, and hence so is $\smock$.

If $\oldPsi$ is a compact, orientable $3$-orbifold such that every boundary component of $N:=|\oldPsi|$ is a sphere, so that
$\plusN$ is closed, we will set
$$
\theta(\oldPsi)=\frac{\volG(\plusN)}{0.305}.
$$

As a hint about why the curious-looking number 
  $0.305$ appears in these
definitions, we mention that the inequality $\voct/12>0.305$, where, as in \ref{voct def}, $\voct$ denotes the volume of a regular ideal octahedron in $\HH^3$,
is used in the proof of Corollary \ref{bloody hell} below (and is close to being an equality since $\voct=3.6638\ldots$), while
Theorem 1.1 of \cite{rankfour}, which implies that
  any closed, orientable hyperbolic $3$-manifold $M$ of volume at most
  $1.22=4\cdot0.305$ satisfies $h(M)\le3$, is used in the proof of Proposition \ref{new get
    lost}. (See \ref{in and out 
    soda} above for the definition of $h(M)$.) 

If a compact, orientable $3$-orbifold $\oldPsi$ is componentwise strongly \simple\ and boundary-irreducible,
we will set
$$\sigma(\oldPsi)=
12\chibar(\kish(\oldPsi))
$$
where $\kish(\oldPsi)$ is defined by \ref{tuesa day}. (The appearance of the number $12$ in this definition is for convenience; it guarantees that certain lower bounds for $\sigma(\oldPsi)$ that appear in the course of the arguments in this paper and in \cite{second} are integers rather than fractions.)
We will also set
$$\delta(\oldPsi)=\sup_\oldPi\sigma(\oldPsi\cut\oldPi),
$$ 
where $\oldPi$ ranges
over all (possibly empty) closed, $2$-dimensional suborbifolds of $\inter\oldPsi$
whose components are all incompressible suborbifolds of
$\oldPsi$. In view of Lemma \ref{oops lemma}, the incompressibility of the components of $\oldPi$ guarantees that $\oldPsi\cut\oldPi$
is itself componentwise strongly \simple\ and boundary-irreducible, so that
$\sigma(\oldPsi\cut\oldPi)$ 
is defined.
Like $\smock_0(\oldPsi)$, the invariant $\delta(\oldPsi)$ is {\it a priori} an element of the extended real
number system.

\abstractcomment{\tiny I don't know whether there is any reason to require that
  their underlying surfaces are incompressible in $|\oldOmega|$. }

Note that the definition of $\delta(\oldPsi)$ implies that 
$\delta(\oldPsi)\ge\sigma(\oldPsi\cut\emptyset)$, i.e.
\Equation\label{kitsch}
\delta(\oldPsi)\ge\sigma(\oldPsi).
\EndEquation

Now suppose that $\oldOmega$ is a closed, orientable, hyperbolic
$3$-orbifold, let $\cals $ be a (possibly disconnected) $2$-dimensional submanifold of $M:=|\oldOmega|$, 
in general position with respect
to
$\fraks_\oldOmega$, and let $X$ be a union of components of
$M- \cals $. We will set
$t_\oldOmega(X)=\smock(\obd(\hatX))$.
(See \ref{nbhd stuff} for the definition of $\hatX$.)
If the components of $\obd(\cals)$ are  incompressible, then since $\oldOmega$ is strongly \simple\ (see \ref{oops}), it follows from Lemma \ref{oops lemma} that the
components of $\obd(\hatX)$ are strongly \simple\ and boundary-irreducible. In this case we will set 
$s_\oldOmega(X)=\sigma(\obd(\hatX))$
and 
$y_\oldOmega(X)=\delta(\obd(\hatX))$.
\abstractcomment{\tiny I've \%-ed out $\delta'$ and $y'_\oldOmega$ because they don't seem to be used.}
If every component of $\partial \hatX$ is a sphere, we will set
$q_\oldOmega(X)=\theta(\obd(\hatX))$.
\EndNotation

The following lemma, which has a simple proof and fits naturally into
the present discussion, is not quoted in the present paper but will be
used in \cite{second}.

\Lemma\label{frobisher}
Let  $\oldPsi$ be any very good, compact, orientable (possibly disconnected) $3$-orbifold, and let $\oldTheta$ be a $2$-dimensional suborbifold of $\inter\oldPsi$
whose components are all incompressible suborbifolds of
$\oldPsi$. Then 
$\smock_0(\oldPsi)\ge\smock_0(\oldPsi\cut\oldTheta)$,
and hence
$\smock(\oldPsi)\ge\smock(\oldPsi\cut\oldTheta)$.
\EndLemma

\Proof
Let $\oldPi$ be any (possibly empty) $2$-dimensional suborbifold of $\inter(\oldPsi\cut\oldTheta)$
whose components are all incompressible suborbifolds of
$\oldPsi\cut\oldTheta$. Then $\oldTheta\cup\oldPi$ is identified with a $2$-dimensional suborbifold of $\inter\oldPsi$
whose components are all incompressible suborbifolds of
$\oldPsi$. Hence the definition of
$\smock_0(\oldPsi)$ implies that $\volorb(\oldPsi\cut{\oldTheta\cup\oldPi})\le \smock_0(\oldPsi)$. Since $\oldPsi\cut{\oldTheta\cup\oldPi}$ is homeomorphic to $(\oldPsi\cut\oldTheta)\cut\oldPi$, it follows that 
$\volorb((\oldPsi\cut\oldTheta)\cut\oldPi)\le \smock_0(\oldPsi)$. Since the latter inequality holds for every
$2$-dimensional suborbifold $\oldPi$ of $\inter(\oldPsi\cut\oldTheta)$
whose components are all incompressible suborbifolds of
$\oldPsi\cut\oldTheta$, the definition of $\smock_0(\oldPsi\cut\oldTheta)$ gives $\smock_0(\oldPsi\cut\oldTheta) \le\smock_0(\oldPsi)$. The inequality $\smock(\oldPsi\cut\oldTheta) \le\smock(\oldPsi)$ then follows at once in view of the definition of the invariant $\zeta$.
\EndProof

\Theorem\label{darts theorem}
Let $\oldOmega$ be a closed, orientable, hyperbolic $3$-orbifold, and let
$\oldTheta$ be a  two-sided $2$-dimensional suborbifold of
$\oldOmega\pl$ whose components are incompressible. Then
$$\vol(\oldOmega)\ge\volorb((\oldOmega\pl)\cut\oldTheta).$$
\EndTheorem

\Proof
We first give the proof under the following additional assumptions:
\begin{enumerate}
\item $\oldOmega$ is a closed, orientable, hyperbolic $3$-manifold (so that $\oldTheta$ is a $2$-manifold whose components are incompressible in
$\oldOmega\pl$);
\item no component of $\oldTheta$ is the boundary of a twisted $I$-bundle (over a closed, non-orientable surface) contained in the manifold $\oldOmega\pl$; and
\item no two components of $\oldTheta$ are parallel in $\oldOmega\pl$.
\end{enumerate}
For this part of the proof, to emphasize that $\oldOmega$ is assumed to be a manifold, I will denote it by $M$. Because $\volorb$ is a purely topological invariant, and because of the equivalence between the PL and smooth categories for $3$-manifolds,  it suffices to prove that if $T\subset M$ is a smooth submanifold such that Conditions (1)---(3) hold with $M$ and $T$ in place of $\oldOmega\pl$ and $\oldTheta$, then 
$\vol(M)\ge\volorb(M\cut T)$. The proof of this will use the observations that Propositions \ref{stronger waldhausen} and \ref{snuff} above, and \cite[Corollary 5.5]{Waldhausen}, which are theorems in the PL category of $3$-manifolds, remaion true without change in the smooth category; this also follows from the equivalence of categories mentioned above.
For each component $V$ of $T$, it follows from
  \cite[Theorem 3.1]{schoen-yau} that the inclusion map $i_V:V\to M$ is
  homotopic to a immersion  $j_V$ which has least area in its homotopy class. By  \cite[Theorem
  5.1]{FHS}, for each component $V$ of $T$, either $j_V$ is an
  embedding,  or there is a a one-sided, closed, connected surface $F_V\subset M$ such
  that $j_V$ is a two-sheeted covering map from $V$ onto $F_V$. If the
  latter alternative holds, then $j_V$, and hence $i_V$, is homotopic
  to a homeomorphism $j'_V$ of $V$ onto the boundary of a tubular
  neighborhood $N_V$ of $F_V$, which is a twisted $I$-bundle. Since
  $i_V$ and $j_V'$ are homotopic embeddings of $V$ into $M$, it
  follows from the smooth version of \cite[Corollary 5.5]{Waldhausen} that they are
  isotopic; this implies that $V$ is the boundary of a twisted
  $I$-bundle in $M$, a contradiction to Assumption (2). Hence for
  every component $V$ of $T$, the immersion $j_V$ is an embedding.

Now define an immersion $j:T\to M$ by setting $j|V=j_V$ for each component $V$ of $T$.  We claim:
\Claim\label{two by two}
The immersion $j$ is an embedding.
\EndClaim
We have already shown that $j_V$ is an embedding for each component
$V$ of $T$. Hence, to prove \ref{two by two}, it suffices to prove
that for any distinct components $V$ and $W$ of $T$ we have  $j_V(V)\cap j_W(W)=\emptyset$. 

Since  the inclusions $i_V$ and $i_W$ are $\pi_1$-injective and have disjoint images, 
it follows from the smooth version of Proposition \ref{snuff}
that for any embeddings $f:V\to M$ and $g:W\to M$ homotopic to $i_V$ and $i_W$ respectively, such that $f(V)$ and $g(W)$ meet transversally, either (a) $f(V)\cap g(W)=\emptyset$, or (b) there
exist connected subsurfaces $A\subset f(V)$ and $B\subset g( W)$, and a compact submanifold
$X$ of $M$, 
such that
$\partial A\ne\emptyset$, $\partial X=A\cup B$, and the pair $(X,A)$ is homeomorphic to
$ (A\times[0,1],A\times\{0\})$.
Since  $j_V$ and $j_W$ are respectively homotopic to $i_V$ and $i_W$, and have least area within their homotopy classes, this shows that the hypothesis of Proposition \ref{fhs-prop} holds with $f_0=j_V$ and $g=j_W$. Thus Proposition \ref{fhs-prop} asserts that $j_V(V)\cap j_W(W)=\emptyset$, and \ref{two by two} is proved.

Now since $j_V$ is homotopic to $i_V$ for each component $V$ of $T$,
the map $j$ is homotopic to $i$. Since $i$ is an embedding, and $j$ is
an embedding by \ref{two by two}, and since Assumption (3) holds, it
follows from the smooth version of Proposition \ref{stronger waldhausen} that $i$ and $j$
are isotopic. In particular, the (possibly disconnected) manifold
$M\cut T$ is homeomorphic to $X:=M\cut{j(T)}$. 

Since each of the embeddings $j_V$ is a least area immersion, the
surface $j(T)\subset M$ is a minimal surface. The manifold $X$ inherits a metric from the hyperbolic manifold
$M$, and with this metric, $X$ is bounded by a minimal
surface. According to Theorem 7.2 of \cite{ast}, we have
$\vol(X)\ge\volorb(X)$. 
Since it is clear that $\vol(M)=\vol(X)$, and since we have seen that 
$M\cut T$ is homeomorphic to $X$, we have
$\vol(M) \ge\volorb(M\cut T)$, as required.
theorem. Thus the proof is complete under Assumptions (1)--(3).


Next, we prove the theorem under the assumptions that (1) and (2)
hold, but (3) may not. (This argument, and the rest of the proof, will be done in the PL category.) Since (1) is still assumed to hold, we continue
to write $M$ in place of $\oldOmega$, and we write $T$ in place of $\oldTheta$ (so that $T$ now denotes a PL $2$-submanifold of $M$). Under Assumptions (1) and (2), the proof will
proceed by induction on $\compnum(T)$. The assertion is trivial if
$T=\emptyset$. Now suppose that $n\ge0$ is given and that the
assertion is true in the case where $\compnum(T)=n$. Suppose that $M$
and $T$ satisfy the hypotheses of the theorem and Assumptions (1) and
(2), and that $\compnum(T)=n+1$. If (3) also holds, the conclusion has
already been established. 

Now suppose that (3) does not hold,
i.e. that $T$ has two parallel components. Then there exists a
submanifold of $M$ which is homeomorphic to the product of a closed, connected
surface with $[0,1]$, and whose boundary components are components of
$T$. Among all such submanifolds of $T$ choose one, say $Y$, which is
minimal with respect to inclusion. If $\inter Y$ contains a component
$W$ of $T$, then since $W$ is incompressible in $M$ and hence in $Y$,
it follows from \cite[Proposition 3.1]{Waldhausen} that $W$ is parallel in $Y$
to each of the components of $\partial Y$; this contradicts the
minimality of $Y$. Hence $\inter Y$ contains no component of $T$, and
therefore $Y$ is the closure of a component of $M-T$.

Fix a component $V$ of $\partial Y$.
Set $T'=T-V$. Then the abstract disjoint union $Y\discup M\cut {T'}$ is homeomorphic to
$M\cut{T}$. 
By \ref{t-defs} it follows that 
$\volorb(M\cut {T})=\volorb(Y)+\volorb(M\cut{T'})$. Since
$Y$ is homeomorphic to $V\times[0,1]$, the manifold $DY$ is
homeomorphic to $V\times S^1$. It then follows from \cite[Proposition C.3.4]{bp}) that $\volG(DY)=0$, so that $\volorb(Y)=\volG(DY)/2=0$. Thus we have 
 $\volorb(M\cut{{T'}})=\volorb(M\cut{{T}})$. Since $\compnum(T)=n$,
 the induction hypothesis gives $\vol
 M\ge\volorb(M\cut{T'})=\volorb(M\cut{T})$, and the induction is
 complete. This establishes the result under Assumptions (1) and (2).

We now turn to the proof of the theorem in the general case. Since $\oldOmega$ is hyperbolic, it is very good. Fix a finite-sheeted orbifold covering $p:N\to\oldOmega$ where $N$ is a connected $3$-manifold. Set $\calq=p^{-1}(\oldTheta)$. We claim: 
\Claim\label{even before that}
There is a two-sheeted covering $q:N'\to N$ such that $N'-q^{-1}(\calq)$ has no component whose closure  is a twisted $I$-bundle.
\EndClaim

To prove \ref{even before that}, let $\calx$ denote the set of all twisted $I$-bundles that are  closures of components of $N-\calq$. First consider the case in which the elements of $\calx$ are pairwise disjoint. In this case, let $Y$ denote the union of all elements of $\calx$, and set $Z=\overline{N-Y}$. Note that $Z$ is connected, since $N$ is connected and each element of $\calx$ has connected boundary. Since each component of $Y$ is a twisted $I$-bundle, there exist a compact (possibly disconnected) $2$-manifold $B$, and a $2$-sheeted covering map $q_Y:B\times[0,1]\to Y$, such that $q_Y$ maps $B\times\{i\}$ homeomorphically onto $\partial Y$ for $i=0,1$. Since $\partial Y=\partial Z$, we may regard $q_Y|(B\times\{i\})$ as a homeomorphism $r_i:B\times\{i\}\to\partial Z$. Let $Z_0$ and $Z_1$ be homeomorphic copies of $Z$, equipped with homeomorphisms $h_i:Z\to Z_i$. Then we may define a homeomorphism $\alpha$ from $\partial(B\times[0,1])=B\times\{0,1\}$ to $\partial(Z_1\discup Z_2)$ by setting $\alpha|(B\times\{i\})=h_i\circ r_i$ for $i=0,1$. In this case we define $N'$ to be the closed $3$-manifold obtained from the disjoint union $(B\times[0,1])\discup Z_1\discup Z_2$ by gluing $\partial(B\times[0,1])$ to $\partial(Z_1\discup Z_2)$ via the homeomorphism $\alpha$. Now we have a well-defined two-sheeted covering map $q:N'\to N$ given by setting $q|(B\times[0,1])=q_Y$ and $q|Z_i=h_i^{-1}$ for $i=0,1$. To prove the property of this covering stated in \ref{even before that}, note that the closure of each component of $N'-q^{-1}(\calq)$ either is equal to $Z$ or is a component of $B\times[0,1]$. A component of the latter type have disconnected boundary and is therefore not a twisted $I$-bundle. That $Z$ is not a twisted $I$-bundle follows from the definitions of $\calx$ and $Z$.

To complete the proof of \ref{even before that}, it remains to consider the case in which there are two elements of $\calx$ have non-empty intersection. If two such elements are denoted $\overline{X_1}$ and $\overline{X_2}$, where $X_1$ and $X_2$ are components of $N-\calq$, then since each of the twisted $I$-bundles $\overline{X_i}$ and $\overline{X_2}$ has connected boundary, and since $N$ is connected, we have $N=\overline{X_1}\cup\overline{X_2}$, and $\calx=\{\overline{X_1},\overline{X_2}\}$. We have $\overline{X_1}\cap\overline{X_2}=\partial \overline{X_1}=\partial\overline{ X_2}=\calq$.
For $j=1,2$, fix a
compact, connected $2$-manifold $B_j$, and a $2$-sheeted covering map $q_j:B_j\times[0,1]\to \overline{X_j}$, such that $q_j$ restricts to a homeomorphism 
$u_{ij}$ of $B_j\times\{i\}$ onto $\calq=\partial\overline{X_j} $ for $j=1,2$ and for $i=0,1$.
In this case we define $N'$ to be the closed $3$-manifold obtained from the disjoint union $(B_1\times[0,1])\discup( B_2\times[0,1])$ by gluing $B_1\times\{i\}$ to $B_2\times\{i\}$ via the homeomorphism 
$u_{i2}^{-1}\circ u_{i1}$, for $i=0$ and for $i=1$. 
Now we have a well-defined two-sheeted covering map $q:N'\to N$ given by setting $q|(B_j\times[0,1])=q_j$ for $j=1,2$.
To prove the property of this covering stated in \ref{even before that}, note that the closures of the components of $N'-q^{-1}(\calq)$ are $B_1\times[0,1]$ and $B_2\times[0,1]$; these have disconnected boundary, and hence neither of them is a twisted $I$-bundle. 
Thus \ref{even before that} is proved.

Let $q:N'\to N$ be a two-sheeted cover with the properties stated in \ref{even before that}, and set $\calq'=q^{-1}(\calq)$. We claim

\Claim\label{new before general}
No component of $\calq'$ is the boundary of a twisted $I$-bundle in $N'$.
\EndClaim

To prove \ref{new before general}, suppose that some component  of $\calq'$ bounds a twisted $I$-bundle $R\subset N'$. We may assume that
$R$ is
minimal with respect to inclusion among all among all twisted $I$-bundles that are bounded by  components of $\calq'$. 
According to the property of $N'$ stated in \ref{even before that},
$R$ cannot be the closure of a component of $N'-\calq'$. Hence 
$\inter R$ contains some component $Q'$ of $\calq'$. Since $H_2(R,\ZZ)=0$, there is a compact submanifold $R_1$ of $\inter R$ with $\partial R_1=Q'$. Now fix a two-sheeted covering $t:\tR\to R$, where $\tR$ is a trivial $I$-bundle over a surface. 
Since the components of
$\calq$ are incompressible,
the components of
$t^{-1}(\calq')=(pt)^{-1}(\calq)$ are also incompressible. It therefore follows from \cite[Prop. 3.1]{Waldhausen}, that each component of $\partial (t^{-1}(R_1))=t^{-1}(Q')\subset t^{-1}(\calq')$ is parallel to the boundary components of $R$. A second application of  \cite[Prop. 3.1]{Waldhausen} then shows that the two-sheeted covering $t^{-1}(R_1)$ of $R_1$ is a trivial $I$-bundle over a surface. Hence by  \cite[Prop. 3.1]{Waldhausen}, $R_1$ is a twisted $I$-bundle over a surface. But this contradicts the minimality of $R$,
and so
\ref{new before general} is proved.

Now since $N$ is a manifold, $N'$ is also a manifold, and hence Condition (1) above holds with $N'$ and $\calq'$
in place
of $\oldOmega$ and $\oldTheta$. According to \ref{new before general},
Condition (2) also holds with $N'$ and $\calq'$
in place
of $\oldOmega$ and $\oldTheta$. 
It therefore follows from the case
of the theorem already proved that
\Equation\label{tea quation}
\vol(N')\ge\volorb((N')\cut{\calq'}).
\EndEquation


Now  if $d$ denotes the degree of the covering
$q\circ p:N'\to\oldOmega$, we have $\vol(N')=d\cdot\vol(\oldOmega)$. On
the other hand, since $\calq'=(q\circ p)^{-1}(\oldTheta)$,
the manifold 
$(N')\cut{\calq'}$ is a $d$-fold covering of
$\oldOmega\cut\oldTheta$, and
$D((N')\cut{\calq'})$ is therefore a $d$-fold covering of
$D(\oldOmega\cut\oldTheta)$.
Thus
the
definitions of $\volorb((N')\cut{\calq'})$,
$\volG(D(\oldOmega\cut\oldTheta))$, and $\volorb(\oldOmega\cut\oldTheta)$ give
$$\volorb((N')\cut{\calq'})=
\frac12\volG(D ((N')\cut{\calq'}))=
\frac d2\volG(D (\oldOmega\cut\oldTheta))=
d\cdot\volorb(\oldOmega\cut\oldTheta).$$ 
Hence the conclusion of the present theorem follows from (\ref{tea
  quation}) upon dividing both sides by $d$.
\EndProof

\Corollary\label{smockollary}
If $\oldOmega$ is a
 closed, orientable, hyperbolic $3$-orbifold, we have $\vol\oldOmega=\smock_0(\oldOmega\pl)$.
\EndCorollary

\Proof 
Since by definition we have $\smock_0(\oldOmega\pl)=\sup_\oldTheta \volorb((\oldOmega\pl)\cut\oldTheta)$, where $\oldTheta$ ranges
over all $2$-suborbifolds of $\oldOmega$
whose components are all incompressible,
the inequality $\vol\oldOmega\ge\smock_0(\oldOmega\pl)$ follows
from Theorem \ref{darts theorem}.
To prove the inequality
$\vol\oldOmega\le\smock_0(\oldOmega\pl)$, note that since $\oldOmega$ is hyperbolic, it admits a $d$-sheeted manifold
cover $\toldOmega$  for some integer $d>0$. According to
\cite[Theorem C.4.2]{bp} we have   $\volG \toldOmega=\vol \toldOmega$. In
view of Definition \ref{sem ting}, we then have
$\vol\oldOmega=(\vol\toldOmega)/d=(\volG\toldOmega)/d=\volG\oldOmega$. Since
$\oldOmega$ is closed, $D\oldOmega$ is a disjoint union of two copies of
$\oldOmega$, so that the definition of the invariant $\volorb$ gives $\volorb(\oldOmega\pl)=\volG(\oldOmega)=\vol\oldOmega$. But
(\ref{more kitsch}), applied with $\oldOmega\pl$ playing the role
of $\oldPsi$, gives
$\volorb(\oldOmega\pl)\le\smock_0(\oldOmega\pl)$, and the conclusion follows.
\EndProof

\Number\label{voct def}
The volume of a regular ideal octahedron in $\HH^3$ will be denoted
$\voct$. We have $\voct=3.6638\ldots$.
\EndNumber

\Theorem\label{i want my v-8}
If $\oldPsi$ is an orientable $3$-orbifold which is
componentwise strongly \simple\ and componentwise boundary-irreducible,
we have
$$\volorb(\oldPsi)\ge\voct\chibar(\kish(\oldPsi)).$$
\EndTheorem

\Proof
Since both sides are additive over components, we may assume that
$\oldPsi$ is connected, and therefore strongly \simple\ and boundary-irreducible.

First consider the case in which $\oldPsi=M$ is a strongly 
\simple, boundary-irreducible, orientable 
$3$-manifold. In this case, because of the equivalence between the PL and smooth categories for $3$-manifolds, $M$ has a smooth structure compatible with its PL structure, and up to topological isotopy, the smooth $2$-submanifolds of $M$ are the same as its PL $2$-submanifolds. This will make it unnecessary to distinguish between smooth and PL annuli in the following discussion.

Since the components of $\omega(\cala(\oldPsi))$ are orientable annular orbifolds by \ref{tuesa day}, and since $\omega(\cala(\oldPsi))=\cala(M)$ is a $2$-manifold in this case,
the components of $\Fr_M(\calb)\subset \cala(M)$ are annuli. Let
$E_1,\ldots,E_n$ denote the components of $\kish M=\overline{M-\book M}$. Then
$\cala_i:=\Fr_ME_i$ is a union of components of $\cala(M)$ for
 $i=1,\ldots, n$. According to the manifold case of Lemma \ref{i see a book}, the pair $(E_i,E_i\cap\partial M)$ is acylindrical for $i=1,\ldots,n$. By the manifold case of Proposition \ref{less than nothing}, we have $\chi(E_i)<0$ for $i=1,\ldots,n$. It is a standard consequence of Thurston's geometrization theorem that if $(E,F)$ is an acylindrical manifold pair with $E$ compact snd $\chi(E)<0$, and if every component of $\overline{(\partial E)-F}$ is an annulus, then $(\inter E)\cup(\inter F)\subset E$ has a finite-volume hyperbolic metric with totally geodesic boundary. Thus for $i=1,\ldots,n$ the manifold
 $E_i^-:=E_i-\cala_i$ has such a
metric. According to the three-dimensional case of Theorem 4.2
of \cite{miyamoto}, we have
$\vol E_i^-\ge\voct\chibar(E_i)$. Hence  $D(E_i^-)=\inter D_{E_i'\cap\partial M}E_i\subset D_{E_i'\cap\partial M}\subset DM$ is a hyperbolic manifold
of finite volume, and 
\Equation\label{john donne}
\vol D(E_i^-)=2\vol (E_i^-)\ge2\voct\chibar(E_i)).
\EndEquation
But for each $i$, since the components of the $\cala_i$ are essential annuli, and since $DM$ is canonically identified with a two-sheeted covering of $\silv M$, it follows from the manifold case of Proposition \ref{silver acylindrical} that the components of $D\cala_i=\Fr_{DM} D(E_i^-)$ are incompressible tori in $DM$. Hence
\cite[Theorem 1]{soma}, together with the fact (see 
\cite[Section 6.5]{thurstonnotes}) that the volume of a finite-volume hyperbolic $3$-manifold is equal to the relative Gromov volume of its compact core, implies that
  $\volG(DM)\ge\vol E_1'+\cdots+\vol E_n'$ (an inequality which could be shown to be an equality). With
(\ref{john donne}), this gives
$$\volG(DM)\ge2\voct(\chibar(E_1)+\cdots+\chibar(E_n))=2\voct\chibar(\kish(M)),$$
so that
$$\volorb(M)=\frac12\volG(DM)\ge\voct\chibar(\kish(M)).$$
This completes the proof in the case where $\oldPsi$ is a manifold.

Now suppose that $\oldPsi$ is an arbitrary strongly \simple, boundary-irreducible,
orientable $3$-orbifold. According to Condition (II) of Definition \ref{oops}, $\oldPsi$ admits a finite-sheeted covering $p:\toldPsi\to\oldPsi$
such that $\toldPsi$ is an irreducible $3$-manifold. Since, by Condition (III) of Definition \ref{oops}, $\pi_1(\oldPsi)$ has no rank-$2$ free abelian subgroup,  $\pi_1(\toldPsi)$ also has no rank-$2$ free abelian subgroup. Furthermore, since Condition (III) of Definition \ref{oops} implies that $\oldPsi$ is not discal, $\toldPsi$ is also non-discal; and $\toldPsi$ is a degree-$1$ regular irreducible $3$-manifold covering of itself. Hence $\toldPsi$ is strongly \simple.
Since $\oldPsi$ is boundary-irreducible, its covering $\toldPsi$ is also boundary-irreducible. Hence by the case of the theorem for a  manifold, which was given
above, we have 
\Equation\label{summer vacation}
\volorb(\toldPsi)\ge\voct\chibar(\kish(\toldPsi)).
\EndEquation

Let $d$ denote the degree of the
covering. Doubling $p$, we obtain a $d$-fold covering map $D\toldPsi\to
D\oldPsi$. Hence by Definition \ref{sem ting}, we have
$\volG(D\toldPsi)=d\volG(D\oldPsi)$. Dividing both sides of the latter
equality by $2$, and applying the definition of $\volorb$, we obtain 
\Equation\label{geese}
\volorb(\toldPsi)=d\cdot\volorb(\oldPsi).
\EndEquation
On the other hand, it follows from Proposition \ref{Comment A} that $\kish(\toldPsi)$ is (orbifold-)isotopic to
$p^{-1}(\kish(\oldPsi))$; thus $\kish(\toldPsi)$ is homeomorphic to a $d$-fold
covering of the (possibly disconnected) orbifold $\kish(\oldPsi)$, so
that 
\Equation\label{more geese}
\chibar(\kish(\toldPsi))=d\cdot \chibar(\kish(\oldPsi)).
\EndEquation
Combining (\ref{summer vacation}) with (\ref{geese}) and (\ref{more geese}), we obtain
$$d\cdot\volorb(\oldPsi)=\volorb(\toldPsi)
\ge\voct\chibar(\kish(\toldPsi))=d\cdot\voct\chibar(\kish(\oldPsi)),$$
which gives the conclusion. 
\EndProof

\Corollary\label{bloody hell}
If $\oldPsi$ is a compact, orientable $3$-orbifold which is
componentwise strongly \simple\ and componentwise boundary-irreducible,
we have
$$\smock(\oldPsi)\ge\delta(\oldPsi).$$
\EndCorollary

\Proof
Let $\oldPi$ be any (possibly empty) closed, $2$-dimensional suborbifold of $\inter\oldPsi$
whose components are all incompressible suborbifolds of
$\oldPsi$. According to Lemma \ref{oops lemma}, the components of
$\oldPsi\cut\oldPi$ are strongly \simple\ and boundary-irreducible, and so 
Theorem \ref{i want my v-8} gives
$\volorb(\oldPsi\cut\oldPi)\ge\voct\chibar(\kish(\oldPsi\cut\oldPi))$. Hence,
using that $\voct=3.6638\ldots>3.66$, and using the definitions of
$\smock_0(\oldPsi)$ and $\smock(\oldPsi)$, we obtain
$$\sigma(\oldPsi\cut\oldPi))=12\chibar(\kish(\oldPsi\cut\oldPi))\le\frac{12}{\voct} \volorb(\oldPsi\cut\oldPi)
\le\frac{12}{\voct} \smock_0(\oldPsi)<\frac {\smock_0(\oldPsi)}{0.305}
=\smock(\oldPsi).
$$
Since the inequality $\sigma(\oldPsi\cut\oldPi))<\smock(\oldPsi)$ holds for every choice of $\oldPi$, we have
$\delta(\oldPsi)=\sup_\oldPi \sigma(\oldPsi\cut\oldPi))\le\smock(\oldPsi)$.
\EndProof

The following corollary will be used in \cite{second}.

\Corollary\label{lollapalooka}
Let $\oldOmega$ be a closed,
orientable hyperbolic $3$-orbifold. Set $M=|\oldOmega|$.
Let $\cals $ be a (possibly disconnected) $2$-dimensional submanifold of $M:=|\oldOmega|$, 
in general position with respect
to
$\fraks_\oldOmega$, such that the components of $\obd(\cals)$ are  incompressible, and let $X$ be a union of components of
$M- \cals $. Then
$$s_{{\oldOmega}}(X)\le y_{{\oldOmega}}(X)\le t_\oldOmega(X).$$ 
\EndCorollary

\Proof
According to the definitions, this
means that $\sigma(\omega(\hatX))\le\delta(\omega(\hatX))\le
\zeta(\omega(\hatX))$.  
The inequality $\sigma(\omega(\hatX)\le\delta(\omega(\hatX)$ follows
from (\ref{kitsch}), and the inequality
$\delta(\omega(\hatX)\le
\zeta(\omega(\hatX))$
follows from Corollary \ref{bloody hell}. (The conditions that $\hatX$ is componentwise strongly \simple\ and componentwise boundary-irreducible, which are needed for the application of Corollary \ref{bloody hell}, follow from Lemma \ref{oops lemma}, in view of the strong \simple ity of $\oldOmega$ and the incompressibility of the components of $\obd(\cals)$. 
\EndProof

\Proposition\label{hepcat} Let $\oldPsi$ be a very good, compact, orientable
$3$-orbifold, and set $N=|\oldPsi|$. 
\begin{itemize} 
\item If $\oldPsi$ (or equivalently $N$) is closed, then $\volG(N)\le\volG(\oldPsi)$.
\item If every component of
$\partial N$ is a sphere (so that $\plusN$ is a  closed $3$-manifold), then
$\volG(\plusN)\le\volorb(\oldPsi)$.
\end{itemize}
\EndProposition

\Proof
It follows from the definitions that the two sides of each of the asserted
inequalities are additive over components. Hence we may assume that
$\oldPsi$, or equivalently $N$, is connected.

To prove the first assertion, fix a covering $p:M\to\oldPsi  $ such that $M$ is a manifold and $d:=\deg p<\infty$. Then the definition of $\volG\oldPsi$ gives $\volG(\oldPsi)=\volG(M)/d$.
On the other hand, since $\oldPsi$ is orientable, $|p|:M\to N$ is a branched covering of degree $d$, and hence by \cite[Section 0.2]{gromov} (or \cite[Remark C.3.3]{bp}) we have $\volG M\ge d\cdot\volG N$. Thus 
$\volG(N)\le\volG(M)/d =\volG(\oldPsi)$, as required.

To prove the second assertion, first consider the case in which $N$ is closed. We have $\plusN=N$ in this case. Furthermore, $D\oldPsi$ is a disjoint union of two copies of
$\oldPsi$, and hence $\volorb(\oldPsi)=\volG(\oldPsi)$. In view of the first assertion, we now have
$\volG(\plusN)=\volG(N)\le\volG(\oldPsi)=\volorb(\oldPsi)$, as required.

Now suppose that $N$ has $n\ge1$
boundary components. Since the components of $\partial N$ are spheres,
$DN$ is homeomorphic to the connected sum of two copies of $\plusN$ and
$n-1$ copies of $S^2\times S^1$. Since Gromov volume is additive under
connected sum (see \cite[Theorem 1]{soma} or \cite[Section 3.5]{gromov}), and since $\volG(S^2\times S^1)=0$ by \cite[Proposition C.3.4]{bp}), it follows
that $\volG(DN)=2\volG(\plusN)$. Since
$|D\oldPsi |=DN$ we may write this as 
\Equation\label{deep-sea fission}
\volG(|D\oldPsi |)=2\volG(\plusN).
\EndEquation
Now the closed orbifold $D\oldPsi$ is very good since $\oldPsi$ is good, and we may apply the first assertion with $D\oldPsi$ in place of $\oldPsi$ to deduce that  $\volG(|D\oldPsi|)\le\volG(D\oldPsi)$. Combining this with the definition of
$\volorb(\oldPsi)$ and 
(\ref{deep-sea fission}) we find
$$2\volorb(\oldPsi)=\volG(D\oldPsi)
\ge \volG(|D\oldPsi |)=2\volG(\plusN),$$
from which the conclusion follows.
\EndProof

\Corollary\label{bloody hep}
If $\oldPsi$ is a very good, compact, orientable $3$-orbifold such that every
component of $\partial|\oldPsi |$ is a sphere,
we have
$$\smock(\oldPsi)\ge\theta(\oldPsi).$$
\EndCorollary

\Proof
Set $N=|\oldPsi|$.
By (\ref{more kitsch}) 
 and Proposition \ref{hepcat} we have 
$\smock_0(\oldPsi)\ge \volorb(\oldPsi) \ge \volG(\plusN)$.
By definition we have $\smock(\oldPsi)= \smock_0(\oldPsi)/0.305$, and hence 
$\smock(\oldPsi)\ge
\volG(\plusN)/0.305
=\theta(\oldPsi)$. 
\EndProof

The following corollary will be used in \cite{second}.

\Corollary\label{lollaheplooka}
Let $\oldOmega$ be a closed,
orientable hyperbolic $3$-orbifold. Set $M=|\oldOmega|$.
Let $\cals $ be a (possibly disconnected) $2$-dimensional submanifold of $M:=|\oldOmega|$, 
in general position with respect to
$\fraks_\oldOmega$, such that the components of $\obd(\cals)$ are  incompressible, and let $X$ be a union of components of
$M- \cals $.
Then
$$q_{{\oldOmega}}(X)\le t_\oldOmega(X).$$ 
\EndCorollary

\Proof
Using Corollary \ref{bloody hep} and the definitions, we find that
$q_{{\oldOmega}}(X)=
\theta(\omega(\hatX))\le
\smock(\omega(\hatX)
)=t_\oldOmega(X)$.
\EndProof

\section{Tori in $3$-manifolds}\label{tori section}

This section is exclusively concerned with $3$-manifolds. We have found it convenient to adopt classical terminology here for Seifert fibered spaces \cite{hempel}. Although these could be regarded as $3$-manifolds equipped with $S^1$-fibrations over orbifolds, we regard them here as $3$-manifolds equipped with maps to $2$-manifolds, with the local structure described in \cite[Chapter 12]{hempel}.

\Lemma\label{torus goes to cylinder}
Let $M$ be a compact, orientable, irreducible $3$-manifold. Suppose
that every component of $\partial M$ is a torus, that $M$ admits no
Seifert fibration, and that every incompressible torus in $\inter M$
is boundary-parallel in $M$. Then $M$ is acylindrical. 
\EndLemma

\Proof
We may assume that $\partial M\ne\emptyset$, as otherwise the conclusion is vacuously true.
Suppose that $A$ is a $\pi_1$-injective, properly embedded annulus in $M$. We must show that $A$ is boundary-parallel. Let $N$ denote a regular neighborhood of the union of $A$ with the component or components of $\partial M$ that meet $A$. Since each component of $\partial M$ is a torus, $N$ admits a Seifert fibration. 

Let $T$ be any component of $\Fr_M N$. Since $N$ admits a Seifert fibration, $T$ is a torus. By the manifold case of Lemma \ref{prepre}, it now follows that either 
(a) $T$ bounds a solid torus in $\inter M$,
or (b)
$T$ is contained in the interior of a ball in
$\inter M$, or (c) $T$ is $\pi_1$-injective in $M$.
But since
$A\subset M$ is a properly embedded, $\pi_1$-injective annulus, $T$ must contain an annulus which is $\pi_1$-injective in $M$; in particular, the inclusion homomorphism $\pi_1(T)\to\pi_1(M)$ is non-trivial. Hence (a) cannot hold. If (c) holds then $T$ is boundary-parallel according to the hypothesis of the present lemma. This shows that every 
 component $T$ of $\Fr_M N$ is the frontier of a submanifold $K_T$ of $M$ which either is a solid torus or is homeomorphic to $T^2\times[0,1]$. We fix such a $K_T$ for every  component $T$ of $\Fr_M N$.

For each  component $T$ of $\Fr_M N$, we must have either $K_T\supset N$ or $K_T\cap N=T$. Consider the case in which
$K_{T_0}\supset N$ for some component $T_0$ of $\Fr_MN$. The construction of $N$ implies that $N$ contains at least one boundary component of $M$, say $T_1$. Then $T_0$ and $T_1$ are distinct components of $\partial K_{T_0}$. Hence $K_{T_0}$ cannot be a solid torus, and must therefore be homeomorphic to $T^2\times[0,1]$. In particular, $T_0$ and $T_1$ are the only components of $\partial K_{T_0}$. Since $A\subset N\subset K_{T_0}$ has its boundary contained in $\partial N$, we must have $\partial A$ contained in $T_1$. It now follows from \cite[Proposition 3.1]{Waldhausen} that $A$ is boundary-parallel in $K_{T_0}$, and hence in $M$; this is the required conclusion.

There remains the case in which
 $K_T\cap N=T$ for each  component $T$ of $\Fr_M N$. If we set $\calt=\calc(\Fr_MN)$, we have $M=N\cup\bigcup_{T\in\calt}K_T$. Fix a Seifert fibration of $N$, and let $\calt'\subset\calt$ denote the set of all components $T$ of $\Fr_MN$ such that either (i) $K_T$ is homeomorphic to $T^2\times[0,1]$ or (ii) $K_T$ is a solid torus whose meridian curve is isotopic in $T$ to a fiber of $N$. Then the Seifert fibration of $N$ extends to a Seifert fibration of $M':=N\cup\bigcup_{T\in\calt'}K_T$. If $\calt'=\calt$ then $M=M'$, so that $M$ is Seifert fibered, a contradiction to the hypothesis. If $\calt'\ne\calt$ then $M$ is obtained from the Seifert fibered space $M'$ by attaching one or more solid tori along boundary components of $M'$ in such a way that each attached solid torus has a meridian curve which is attached along a fiber of $M'$. Since $\partial M\ne\emptyset$, it follows that $M$ is either  solid torus (which admits a Seifert fibration) or a reducible $3$-manifold. In either case we have a contradiction to the hypothesis.
\EndProof

\Lemma\label{and four if by zazmobile}
Let $M$ be a closed, connected, orientable $3$-manifold equipped with a Seifert
fibration. If $h(M)\ge4$, then $M$ contains a saturated incompressible torus.
\EndLemma

\Proof
Fix a Seifert fibration $p:M\to G$, where $G$ is a compact, connected $2$-manifold.
let $n\ge0$ denote the number of singular fibers, and let
$x_1,\ldots,x_n\in G$ denote the images of the singular fibers under $p$. If $G$
is not a sphere or a projective plane, then $G$ contains a simple
closed curve $C$ which does not bound a disk in $G$, and $p^{-1}(C)$ is a saturated
incompressible torus. 

Now suppose that $G$ is a sphere or a projective plane. Fix disjoint disks
$D_1,\ldots,D_n\subset G$ such that $x_i\in D_i$ for each $i$. 
For $i=1,\ldots,n$, set $V_i=p^{-1}(D_i)$.
Set
$G'=G-\bigcup_{i=1}^n\inter D_i$, and
$M'=p^{-1}(G')=\overline{M-(V_1\cup\cdots\cup V_n)}$. The map
$p':=p|M':M'\to G'$ is a fibration with fiber $S^1$, and hence $h(M')\le
h(G)+1
$. Since the $V_i$ are solid tori, the inclusion homomorphism $H_1(M';\ZZ_2)\to
H_1(M;\ZZ_2)$ is surjective, we have $h(M)\le h(M')\le h(G')+1$.

In the case where $G$ is a sphere we have $h(G')=\max(0, n-1)$, and hence $h(M)\le\max( n,1)$. The hypothesis then implies that $n\ge4$. 
Hence there is a simple
closed curve $C\subset G-\{x_1,\ldots,x_n\}$ such that each of the
disks bounded by $C$ contains at least two of the $x_i$. This implies
that $p^{-1}(C)$ is a saturated incompressible torus.

In the case where $G$ is a projective plane we have $h(G')=\max(1, n)$, and hence $h(M)\le\max(2, n+1)$. The hypothesis then implies that $n\ge3$. 
In this case, define $C\subset G$ to be a simple closed curve bounding a disk $E\subset G$ with $x_1,\ldots,x_n\in\inter E$. Since $n\ge3>2$, and since $G-\inter E$ is a M\"obius band,
$p^{-1}(C)$ is a saturated incompressible torus.
\EndProof

\Lemma\label{before graphology}
Let $\Sigma$ be a compact, orientable (but possibly disconnected)
$3$-manifold, and let $A\subset\partial\Sigma$ be a compact
$2$-manifold. Suppose that $\Sigma$ admits a Seifert fibration in
which $A$ is saturated. Let
$F$ be a closed (possibly disconnected) $2$-dimensional submanifold of
$\inter\Sigma$ whose components are two-sided and incompressible. Then
either $\Sigma$ admits a Seifert fibration in which $F\cup A$ is saturated,
or some component of $F$ is a fiber in some fibration of a component
of $\Sigma$ over $S^1$.
\EndLemma

\Proof
First consider the case in which $F$ and $\Sigma$ are both
connected. In this case, \cite[Lemma II.7.3]{js} asserts that either
(a) $F\cup A$ is saturated in some Seifert fibration of $\Sigma$, or (b)
$\Sigma$ has a finite-sheeted covering space homeomorphic to
$F\times[0,1]$. The proof shows that (b) may be replaced by the
stronger condition (b') $F$ is a fiber in some fibration of $\Sigma$
over $S^1$ for which the monodromy homeomorphism is of finite
order. In particular, the conclusion of the lemma holds in this case.

If $F$ is connected but $\Sigma$ is not, the assertion of the lemma follows upon
applying the case discussed above, with the role of $\Sigma$ played by
the component of $\Sigma$ containing $F$.

To prove the result in general, we use induction on $\compnum(F)$. If
$\compnum(F)=1$, i.e. $F$ is connected, we are in the case discussed
above, and if $\compnum(F)=0$, i.e. $F=\emptyset$, the result is
trivial. Now suppose that $\compnum(F)=n>1$ and that the result is
true for surfaces with $n-1$ components. Assume that no component of
$F$ is a fiber in a fibration of a component of $\Sigma$ over $S^1$. Fix a component $F_0$ of $F$. Since the lemma has been
proved in the case of a connected surface, there is a Seifert
fibration of $p:\Sigma\to E$, for some compact $2$-manifold $E$, in which $F_0$ is saturated. Hence $F$ is a
torus. Let $N$ denote a
regular neighborhood of $F_0$ which is saturated in the Seifert
fibration $p$. We may identify $N$
homeomorphically with $F_0\times[-1,1]$ in such a way that $F_0$ is
identified with $F_0\times\{0\}$. For $\epsilon=\pm1$, fix an annulus
$C_\epsilon\subset F_0\times\{\epsilon\}\subset\partial N$ which is
saturated in the Seifert fibration $p$. Note that the annuli $C_1$ and
$C_{-1}$ are homotopic in $N$. Set
$\Sigma'=\overline{\Sigma-N}$, and set $A'=A\cup C_1\cup
C_{-1}\subset\partial\Sigma'$, so that $A'$ is saturated with respect to
the Seifert fibration $p|\Sigma'$. Now set $F'=F-F_0$, so that
$\compnum(F')=n-1$. By the induction hypothesis, either (i) some
component $F_1'$
of $F'$ is a fiber in a fibration over $S^1$ of a component $\Sigma_1'$ of $\Sigma'$
containing $F_1'$, or (ii) $F'\cup A'$ is saturated in some Seifert fibration
of $\Sigma'$. If (i) holds, then since 
no component of
$F$ is a fiber in a fibration of a component of $\Sigma$ over $S^1$,
we must have $\Sigma_1'\subset\Sigma_0$, where $\Sigma_0$ denotes the component of $\Sigma$
containing $F_0$. This implies that
$\Sigma_1'\subset\overline\Sigma_0-N$, so that $\Sigma_1'$ is not
closed; this is impossible, because $F_1'$ is a closed surface, and
therefore cannot be the fiber in a fibration of a non-closed
$3$-manifold over $S^1$. Hence (ii) must hold. Fix a Seifert fibration
$q:\Sigma'\to B$, for some compact $2$-manifold $B$, in which $F'$ and $A'$ are saturated. In particular the
annuli $C_1$ and $C_{-1}$ are saturated in the fibration $q$. Since
$C_1$ and $C_{-1}$ are homotopic in $N=F_0\times[-1,1]$, the Seifert
fibration $q|(F_0\times\{-1,1\})$ extends to a Seifert fibration of
$N$ in which $F_0=F_0\times\{0\}$ is saturated. It follows that $q$
extends to a Seifert fibration $\overline q$ of $\Sigma$ in which $F_0$ is
saturated. Since $F'$ and $A\subset A'$ are saturated in the Seifert fibration $q$, it now follows
that $F=F'\cup F_0$ and $A$ are saturated in the Seifert fibration
$\overline q$, and the induction is complete.
\EndProof

\Lemma\label{graphology}
Let $M$ be  a closed graph manifold (see \ref{great day}). Let $X$ be a compact (possibly
empty) $2$-submanifold of $M$ whose boundary components are all
incompressible  tori in $M$. Suppose that $h(M)\ge \max(4,h(X))$, and let $m$ be an integer such that
$\max(2,h(X))\le m\le h(M)$. Then
there is a compact, connected submanifold $L$ of $ M$ such that
$X\subset \inter L$,
each component of
$\partial L$ is an incompressible  torus in $M$, and $h(L)=m$.
\EndLemma

\Proof
Let $\Sigma$ denote the characteristic submanifold of $M$. (By definition this means that $\Sigma=M$ if $M$ is Seifert-fibered, and that if $M$ is a non-Seifert-fibered graph manifold, and hence a Haken manifold, then $\Sigma$ is the characteristic submanifold as defined in \cite{js}; cf. \ref{manifolds are different}). Since $M$
is a graph manifold, each component  of $\overline{M-\Sigma}$ is
homeomorphic to $[0,1]\times S^1\times  S^1$. Since the
  components of $\partial X$ are incompressible tori,  we may assume after an isotopy that $\partial X\subset\inter\Sigma$. Since each
  component of
$\Sigma$ is a Seifert fibered space, $\Sigma$ itself admits a Seifert
fibration. Applying Lemma \ref{before graphology}, taking
$F=\partial X$ and $A=\emptyset$, we deduce that
either (i) some
component $F_0$ of $\partial X$ is a fiber in a fibration over $S^1$
of some component $\Sigma_0$  of $\Sigma$, or (ii) $\partial X$ is saturated in some Seifert fibration
of $\Sigma$. If (i) holds, then since the components of $\partial X$
are tori, $\Sigma_0$ is closed and $h(\Sigma_0)\le3$. Since the graph
manifold $M$ is by definition irreducible and therefore connected, it
follows that $M=\Sigma_0$ and hence that $h(M)\le3$, a contradiction
to the hypothesis. Hence (ii) holds. Let us fix a
Seifert fibration
  $p:\Sigma\to B$, for some compact $2$-manifold $B$,
 in which $\partial X$ is saturated.

The images under $p$ of singular fibers of $\Sigma$ will be referred to as {\it singular points} of $B$.

Let $\calK$ denote the set of all
compact, connected, non-empty, $3$-dimensional submanifolds $K$ of $ M$ such that
(1) $X\subset \inter K$ and
(2) each component of
$\partial K$ is an incompressible saturated torus in $\inter\Sigma$.

Set $\mu=\max(2,h(X))$. We claim that:
\Claim\label{oswald zdestiny}
There is an element $K$ of $\calk$ with $h(K)\le\mu$.
\EndClaim

Indeed, if $X\ne\emptyset$, then since $\partial X$ is saturated, $X$ admits a
regular neighborhood $K$ with  saturated boundary; it is clear that $K\in\calk$, and
$h(K)=h(X)\le\mu$ by the definition of $\mu$.
Hence \ref{oswald zdestiny} is true in this case. If
  $X=\emptyset$ and $\Sigma\ne M$, then
  $\partial\Sigma\ne\emptyset$. In this case, let us choose a
  saturated boundary-parallel torus
  $T\subset\inter\Sigma$, and let $K$ denote a saturated regular neighborhood
  of $T$ in $\inter\Sigma$. We clearly have $K\in\calk$, and we have  $h(K)=2\le \mu$. If
  $X=\emptyset$ and $\Sigma=M$, then since
  $h(M)\ge4$, Lemma \ref{and four if by zazmobile} gives a saturated incompressible torus in $M$. 
  If $K$ denotes a saturated neighborhood of such a torus, then  again we
  clearly have $K\in\calk$, and again we have  $h(K)=2\le \mu$. This proves \ref{oswald zdestiny}.

Now for each $K\in\calk$, since $\partial K$ is saturated in the
Seifert fibration of $\Sigma$, the $3$-dimensional submanifold
$\overline{\Sigma\setminus K}$ of $\Sigma$ is also saturated. Hence
we may write $\overline{\Sigma\setminus K}=p^{-1}(R_K)$ for
a unique $2$-submanifold $R_K$ of $B$. 
Each component of
$\partial R_K$ is contained in either $\partial B$ or $p(\partial K)$.
Let $\alpha_1(K)$ denote the number of components
of $\partial B$ contained in $R_K$, let $\alpha_2(K)$ denote the number of singular points lying in $R_K$,
and let $\gamma(K)$ denote the sum of the squares of the first betti
numbers of the components of $R_K$. Set 
$\nu(K)=(\alpha_1(K)+\alpha_2(K),\gamma(K))\in\NN^2$. The
set $\NN^2$ of pairs of non-negative integers will be endowed with
the lexicographical order. We claim:

\Claim\label{sooey pig}
Suppose that $K$ is an element of $\calk$ with non-empty boundary, and
that $T$ is a component of
$\partial K$. Then $T=p^{-1}(\ell)$ for some simple closed curve
$\ell\subset\inter B$. 
Furthermore, we have $\ell\subset\partial R_K$, and if $R_0$ denotes the component of $R_K$ whose boundary contains
$\ell$, then there exist
a subsurface $\Delta$ of
$R_0$, and a regular neighborhood $K'$ of
$K\cup p^{-1}(\Delta)$ in $M$, such that
\begin{itemize}
\item 
$K'\in\calk$,
\item
$h(K')\le h(K)+1$, and
\item $\nu(K')<\nu(K)$.
\end{itemize}
\EndClaim


To prove \ref{sooey pig}, first note that the definition of $\calk$
implies that $T$ is a saturated torus in $\inter\Sigma$, and therefore has
the form $p^{-1}(\ell)$ for some simple closed curve
$\ell\subset\inter B$. Since $T$ is a component of $\partial K$, the
curve $\ell$ is a component of $\partial R_K$. Let $R_0$ denote the component of $R_K$ whose boundary contains
$\ell$. 

If $R_0$ were a disk containing at most one singular point, then
$p^{-1}(R_0)$ would be a solid torus with boundary $T$; this would
contradict the incompressibility of $T$ in $M$. Hence one of the
following cases must
occur: (i) $R_0$ is a disk containing exactly
two singular points; (ii) $R_0$ is a disk containing at
least three singular points, or $\chi(R_0)\le0$ and $R_0$ contains at
least one singular point; (iii) $R_0$ contains no singular points
and $\chi(R_0)<0$; (iv) $R_0$ is an annulus or M\"obius band containing no
singular points, and $\partial R_0\subset p(\partial K)$; or (v) $R_0$ is an annulus containing no
singular points, and the component of $\partial R_0$ distinct from
$\ell$ is contained in $\partial B$.

In Cases (i), (iv), and (v) we will take $\Delta=R_0$. In Case (ii) we
will take $\Delta$ to be a disk contained in $R_0$, meeting
$\partial R_0$ in an arc contained in $\ell$, and containing exactly
one singular point in its interior and none on its boundary. In Case
(iii) we will take $\Delta$ to be a regular neighborhood relative to $R_0$ of
a non-boundary-parallel  arc in $R_0$ which has both endpoints in
$\ell$. 

In each case, $\Delta$ is connected, and $\Delta\cap p(K)$ is a non-empty $1$-dimensional
submanifold of $p(\partial K)$; hence $K_0':=K\cup p^{-1}(\Delta)$ is a connected
$3$-manifold. Furthermore, since $K\in\calk$, each component of
$\partial K$ is an incompressible saturated torus in
$\inter\Sigma$. Hence each boundary component of $K_0'$ is a saturated torus in
$\Sigma$. 

It is also clear from the construction of $\Delta$ in each of the five cases that
no boundary component of $\overline{R_0-\Delta}$ bounds a disk in
$R_0$ containing at most one singular point. This implies that every
boundary component of $J:=p^{-1}(\overline{R_0-\Delta})$ is
$\pi_1$-injective in $Z:=p^{-1}({R_0})$. The components of
$\partial Z$ are incompressible in $M$, because each of them is a
component of either $\partial K$ or $\partial\Sigma$. Hence $Z$ is $\pi_1$-injective
in $M$; it follows that the 
components of $\partial J$ are incompressible in
$M$. But every boundary component of $K_0'$
is a boundary component of either
$K$ or $J$. Hence all boundary
components of $K_0'$ are incompressible in $M$.

We will take $K'$ to
be a small regular neigborhood of $K_0'$ in $M$. Since the components of
$\partial K_0'$ are saturated tori in $\inter\Sigma$ and are incompressible in
$M$, we may choose  $K'$ so that its boundary components are also saturated tori in
$\inter\Sigma$ and incompressible in $M$. Since in addition  $\inter K'\supset\inter K\supset X$, we have $K'\in\calk$. 

We must show that $h(K')\le h(K)+1$, or equivalently that $h(K_0')\le h(K)+1$. 
For this purpose it suffices to show that
$\dim H_1(K_0',K;\ZZ_2)\le1$, or equivalently that $\dim H_1(L,Y;\ZZ_2)
\le1$, where
$L=p^{-1}(\Delta)$ and $Y=p^{-1}(\Delta)\cap K=
p^{-1}(\Delta\cap\partial R_0)$. This is readily verified in each of
the cases (i)---(v). Indeed, in Case (i), $L $ is a
Seifert fibered space over the disk with two singular fibers and
$Y=\partial L$. 
In Case (ii), $L $ is a solid torus and
$Y$ is a homotopically non-trivial annulus in
$\partial L$. 
In Cases (iii) and (v), the pair $(L,Y) $ is
respectively homeomorphic to
$([0,1]\times [0,1]\times S^1,[0,1]\times\{0,1\}\times S^1)$ and
$(S^1\times S^1\times[0,1] ,S^1\times S^1\times\{0\})$,
while in Case (iv), $L$ is homeomorphic to an $S^1$-bundle over an annulus or M\"obius band, and $Y=\partial L$.



To prove \ref{sooey pig}, it remains to show that $\nu(K')<\nu(K)$. For this purpose, note that,
as a consequence of the definitions and the constructions, $\alpha_1(K)-\alpha_1(K')$ is equal
 to the
  number of components of $\partial B$ 
  contained in $\Delta$, while
$\alpha_2(K)-\alpha_2(K')$ is equal to 
 the number of points of singular points lying in $\Delta$.
This gives $\alpha_1(K')\le\alpha_1(K)$ and $\alpha_2(K')\le\alpha_2(K)$ for
  $i=1,2$, and at least one of these inequalities will be strict if
$\Delta$   contains either a
component of $\partial B$ or a singular point of $B$. It follows from the construction of $\Delta$
 that this is true in Cases (i), (ii), and (v). Thus in these cases we
 have $\alpha_1(K')+\alpha_2(K')<\alpha_1(K)+\alpha_2(K)$, and in
 particular $\nu(K')<\nu(K)$. In Cases (iii) and (iv) we have
$\alpha_1(K')+\alpha_2(K')\le\alpha_1(K)+\alpha_2(K)$, and we need to
prove that $\gamma(K')<\gamma(K)$. Since $R_{K'}$ is homeomorphic to $R_{K'_0}$, it suffices
to
prove that $\gamma(K'_0)<\gamma(K)$. 

In Case (iii), $R_{K'_0}$ is obtained from $R_K$ by removing a regular
neighborhood of a properly embedded arc in the component $R_0$ of
$R_K$; this arc is not boundary-parallel in $R_0$, and its endpoints
are in the same component of $\partial R_0$. If the arc does not
separate $R_0$, then $R_{K'_0}$ is obtained from $R_K$ by replacing the
component $R_0$ by a new component whose first betti number is one less
than that of $R_0$; hence $\gamma(K'_0)<\gamma(K)$ in this subcase. If the arc 
separates $R_0$, then $R_{K'_0}$ is obtained from $R_K$ by replacing the
component $R_0$ by two new components, whose first betti numbers are
strictly positive and add up to the first betti number of $R$. Hence
$\gamma(K'_0)<\gamma(K)$ in this subcase as well. In Case (iv), $R_{K'_0}$
is obtained from $R_K$ by discarding a component which is an annulus or a M\"obius band. Hence
$\gamma(K'_0)=\gamma(K)-1$ in this case. This completes the proof of
\ref{sooey pig}.

Now let $\calk^*$ denote the subset of $\calk$ consisting of all
elements $K$ such that $h(K)\le m$. Since $m\ge\mu$ by hypothesis, it
follows from \ref{oswald zdestiny} that $\calk^*\ne\emptyset$. Since $\NN^2$
is well ordered, it now follows  that there is an element $L\in \calk^*$ such that
$\nu(L)\le\nu(K)$ for every $K\in\calk^*$. By the definition of $\calk^*$
we have $h(L)\le m$. Suppose that $h(L)<m$. Since $h(M)\ge m$ by hypothesis, we then have $h(L)<h(M)$, so
that $L\ne M$ and therefore $\partial L\ne\emptyset$. Hence
 by \ref{sooey pig}, applied with $L$ playing the role of $K$, there is an
element $L'$ of $\calk$ such that $\nu(L')<\nu(L)$ and $h(L')\le
h(L)+1\le m$. This means that $h(L')\in\calk^*$, a
contradiction. This shows that $h(L)=m$. In view of the definition of $\calk$,
it follows that $L$ has the properties asserted in the lemma.
\EndProof

\Lemma\label{even easier} Let $Y$ be a compact, connected, orientable
$3$-manifold.  Let $\calt\subset\inter Y$ be a compact $2$-manifold
whose components are all of strictly positive genus, and suppose that
$Y-\calt$ has exactly two components, $B$ and $C$. Then $h(Y)\le
h(B)+h(C)-1$.
\EndLemma

\Proof
All homology groups considered in this proof will be understood to
have $\ZZ_2$ coefficients. Let $m$ denote the number of components of
$\calt$. Consider the exact sequence
\Equation\label{yellow pig}
H_1(\calt)\to H_1(\overline B)\oplus
H_1(\overline C)\to H_1(Y)\to H_0(\calt)\to H_0(\overline B)\oplus
H_0(\overline C)\to H_0(Y)\to 0,
\EndEquation
which is a fragment of the Mayer-Vietoris sequence.
Let $V$ denote the image of the map $H_1(\calt)\to H_1(\overline B)\oplus
H_1(\overline C)$ in (\ref{yellow pig}). 
The sequence (\ref{yellow pig}) gives rise to an exact sequence
\Equation\label{purple pig}
0\to V\to H_1(\overline B)\oplus
H_1(\overline C)\to H_1(Y)\to H_0(\calt)\to H_0(\overline B)\oplus
H_0(\overline C)\to H_0(Y)\to 0.
\EndEquation
The exactness of (\ref{purple pig}) implies that
$$\begin{aligned}
0&=\dim V-\dim(H_1(\overline B)\oplus
H_1(\overline C))+\dim H_1(Y)-\dim H_0(\calt)+\dim (H_0(\overline B)\oplus
H_0(\overline C))-\dim H_0(Y)\\
&=\dim V-(h(B)+h(C))+h(Y)-m+2-1,
\end{aligned}$$
so that
\Equation\label{knute rockne}
h(Y)=(h(B)+h(C))+(m-\dim V)-1.
\EndEquation
Now let $i:H_1(\calt)\to H_1(\bar B)$ and $j:H_1(\calt)\to H_1(\bar
C)$ denote the inclusion homomorphisms. The intersection pairing on
$H_1(\calt)$ is non-singular because $\calt$ is a closed surface; but
since $\calt\subset\partial B$, this pairing is trivial on $\ker
i$. Hence $\dim\ker i\le(\dim H_1(\calt ))/2$, and therefore $\dim \image
i\ge (\dim H_1(\calt ))/2$.  Since each component of $\calt$ has strictly
positive genus, we have $\dim H_1(\calt )\ge 2m$, and so $\dim \image
i\ge m$. On the other hand, the map $H_1(\calt)\to H_1(\overline B)\oplus
H_1(\overline C)$ in (\ref{yellow pig}) is defined by
$x\mapsto(i(x),j(x))$, and so $\dim V\ge\dim\image i$. It follows that
\Equation\label{punkin}
\dim V\ge m.
\EndEquation
The conclusion of the lemma follows from (\ref{knute rockne}) and (\ref{punkin}).
\EndProof

\Lemma\label{another goddam torus lemma}
Let $U$ and $V$ be compact, connected $3$-dimensional submanifolds of an
orientable $3$-manifold. Suppose that $\inter
U\cap\inter V=\emptyset$, and that every component of $U\cap V$ is a
torus. Then $h(U\cup V)\ge h(U)-1$.
\EndLemma

\Proof
All homology groups in this proof will be understood to have
coefficients in $\ZZ_2$. Set $\calt= U\cap V$ and $Z=U\cup V$. We may
assume that $\calt\ne\emptyset$, so that $Z$ is connected. Then $\calt$ is a union of $m$ boundary tori of
$V$ for some integer $m\ge1$. The
intersection pairing on $H_1(\calt;\ZZ_2)$ is nondegenerate, and the
kernel $K$ of the inclusion homomorphism $H_1(\calt)\to
H_1(V)$ is self-orthogonal with respect to this pairing. Hence
$\dim K\le (\dim H_1(\calt;\ZZ_2))/2=m$. The exactness of the
Mayer-Vietoris fragment
$$H_1(\calt)\longrightarrow H_1(U)\oplus H_1(V)\longrightarrow
H_1(Z)$$
implies that the kernel $L$ of the inclusion homomorphism $H_1(U)\to
H_1(Z)$ is the image of $K$ under the inclusion homomorphism $H_1(\calt)\to
H_1(U)$. Hence $\dim L\le m$. Now since $U$ and $Z$ are
connected, the homology exact
sequence of the pair $(Z,U)$ gives rise to an
exact sequence
$$0\longrightarrow L \longrightarrow H_1(U) \longrightarrow H_1(Z)
\longrightarrow H_1(Z,U) \to0,$$
which implies that 
\Equation\label{glass my ass}
h(U)-h(Z)=\dim L-\dim H_1(Z,U)\le m-\dim H_1(Z,U).
\EndEquation
On the other hand, by excision we have $H_1(Z,U)\cong H_1(V,\calt)$,
and the exact homology sequence
$$H_1(V,\calt)\longrightarrow H_0(\calt)\longrightarrow H_0(V)$$
shows that $\dim H_1(V,\calt)\ge\dim H_0(\calt)-\dim H_0(V)=m-1$. Thus
$\dim H_1(Z,U)\ge m-1$, which combined with (\ref{glass my ass}) gives $h(U)-h(Z)\le1$.
\EndProof

\Lemma\label{just right}
Let $m\ge2$ be an integer, and
let $M$ be  a closed graph manifold with $h(M)>\max(3,m(m-1))$. Then there is
a compact submanifold $P$ of $M$ such that
\begin{itemize}
\item each component of $\partial P$ is an incompressible  torus in
  $M$,
\item $P$ and $M-P$ are connected, and
\item $\min(h(P),h(M-P))\ge m$.
\end{itemize}
\EndLemma

\Proof
Note that the hypothesis implies that 
$h(M)\ge 4$ and that
$2\le m\le h(M)$. Hence by Lemma \ref{graphology}, applied with $X=\emptyset$, there is a compact, connected submanifold $L$ of $ M$ such that
(a) each component of
$\partial L$ is an incompressible  torus in $M$, and (b) $h(L)=m$. It
follows from the special case $T=\emptyset$ of \cite[Theorem V.2.1]{js}
that we may
  take $L$ to be ``maximal'' among all compact, connected submanifolds
  satisfying (a) and (b) in the sense that, if $L'$ is any
  compact, connected submanifold of $\inter M$, such that
  $L\subset\inter L'$,
each component of
$\partial L'$ is an incompressible  torus in $M$, and $h(L')=m$, then
$L'$ is a regular neighborhood of $L$. Since the hypothesis implies
that $h(M)>m$, we have $L\ne M$, so that $\partial L\ne\emptyset$.

It is a standard consequence of Poincar\'e-Lefschetz duality that
the total genus of the boundary of a compact, connected, orientable
$3$-manifold $L$ is at most $h(L)$. Since the boundary components of
the manifold $L$ we have chosen are tori, it follows that $\partial L$
has at most $m$ components. Thus if $k=\compnum(\partial L)$ and
$r=\compnum(M-L)$, we have $1\le r\le k\le m$. If 
$P_1,\ldots,P_r$ are the components of $\overline{M-L}$, we may
index the components of $\partial L$ as $S_1,\ldots,S_k$, in such a
way that $S_i\subset\partial P_i$ for $i=1,\ldots,r$. 

For
$p=0,\ldots,r$, let $L_p$ denote the connected $3$-manifold obtained from the
disjoint union $L\discup P_1\discup\cdots\discup P_p$ by gluing
$S_i\subset L$ to $S_i\subset P_i$, for $i=1,\ldots,p$, via the
identity map. Then $h(L_0)=h(L)=m$. For $1\le p\le r$, it follows from
Lemma \ref{even easier} that $h(L_p)\le h(L_{p-1})+h(P_p)-1$. Hence 
$$h(L_r)\le h(L_0)+\sum_{i=1}^r h(P_i)-r= m-r+\sum_{i=1}^r h(P_i).$$
Up to homeomorphism, the manifold $M$ may be obtained from $L_r$ by
gluing the boundary components $S_{r+1},\ldots,S_k$ of $L\subset L_r$
to boundary components of $P_1\discup\cdots\discup P_r\subset
L_r$. Hence
$$h(M)\le h(L_r)+k-r\le2 (m-r)+\sum_{i=1}^r h(P_i).$$
Since $h(M)>m(m-1)$ by the hypothesis of the lemma, we obtain
$m^2-3m<-2r+\sum_{i=1}^r h(P_i)=\sum_{i=1}^r (h(P_i)-2)$, so that
$$\sum_{i=1}^r (h(P_i)-2)> m(m-3)\ge r(m-3).$$
Hence $h(P_{i_0})-2> m-3$ for some
$i_0\in\{1,\ldots,r\}$. If we set $P=P_{i_0}$ it follows that $h(P)\ge
m$. Furthermore, $P=P_{i_0}$ is by definition connected, and
$\overline{M-P}=L\cup\bigcup_{i\ne i_0}P_i$ is connected because each of
the connected submanifolds $P_i$ meets the connected submanifold
$L$. It remains only to show that $h(M-P)\ge m$.

Assume that $h(\overline{M-P})<m$. We again apply Lemma \ref{graphology},
this time taking $X=\overline{M-P}$. This gives a compact, connected submanifold $L'$ of $ M$ such that
$\overline{M-P}\subset \inter L'$,
each component of
$\partial L'$ is an incompressible  torus in $M$, and $h(L')=m$. Since
$L\subset \overline{M-P}$, we in
particular have $L\subset\inter L'$. In view of our choice of $L$,
it follows that $L'$ is a regular neighborhood of $L$. Since
$L\subset \overline{M-P}\subset \inter L'$, and since the boundary components
of $\overline{M-P}$ are incompressible, it follows that
$L'$ is also a regular neighborhood of $\overline{M-P}$. 
But this is
impossible because $h(L')=m>h(\overline{M-P})$.
\EndProof

\section{Finding useful solid tori}\label{seek 'n' find}

We indicated in the introduction that if  $\oldTheta$ is a $2$-suborbifold of a hyperbolic orbifold $\oldOmega$, such that the components of $|\oldTheta|$ are incompressible tori in $|\oldOmega|$, it may be necessary to modify $\oldTheta$ to obtain a suborbifold that is useful for finding lower bounds for the invariants introduced in Section \ref{darts section}, and hence for $\vol\oldOmega$. These modifications will be detailed in Section \ref{underlying tori section}. 

The main result of this section, Lemma \ref{uneeda}, is a technical
result on $3$-orbifolds that will be needed for the constructions given in Section \ref{underlying tori section}.
The following result, Lemma \ref{uneeda this first}, is
required for the proof of Lemma \ref{uneeda}.

\Lemma\label{uneeda this first}
Let $p:\oldUpsilon\to\frakB$ be a covering map of compact, connected $2$-orbifolds. Suppose that $\oldUpsilon$ is orientable, that $\chi(|\oldUpsilon|)=0$, and that $|p|: |\oldUpsilon| \to |\frakB|$ (see \ref{orbifolds introduced})  is $\pi_1$-injective. Then $\chi(|\frakB|)=0$.
 Furthermore, the index in $\pi_1(|\frakB|)$ of the image of $|p|_\sharp:\pi_1(|\oldUpsilon|)\to\pi_1(|\frakB|)$ is at most the degree of the orbifold covering map $p$. Finally, if $|\oldUpsilon|$ is an annulus, if every finite subgroup of $\pi_1(\frakB)$ is cyclic, and if the degree of $p$ is at least $2$, then $\card\fraks_\oldUpsilon\ne1$. 
\EndLemma


\Proof[Proof of Lemma \ref{uneeda this first}]
Let $d$ denote the degree of the orbifold covering map $p$.

In the case where $|\oldUpsilon|$ is a torus, the $\pi_1$-injectivity of $|p|$ implies that $\pi_1(|\frakB|)$ has a subgroup isomorphic to $\ZZ\times\ZZ$, and therefore that the $2$-manifold $|\frakB|$ is a torus or Klein bottle. Hence $\chi(|\frakB|)=0$. Since $|\frakB|$ is closed, we have $\dim\fraks_\frakB\le0$, so that $|p|$ is a branched covering map; in particular $|p|$ is a degree-$d$ map of closed manifolds, and hence $|\pi_1(|\frakB|):|p|_\sharp(|\pi_1(|\oldUpsilon|)|\le d$.

The rest of the proof will be devoted to the case in which
$|\oldUpsilon|$ is an annulus. Since $\frakB$ is connected, the
covering map $p$ is surjective, and hence restricts to a surjective
covering map $r:\partial\oldUpsilon\to\partial\frakB$. Since
$\partial\frakB$ is a closed $1$-orbifold, each of its components
is homeomorphic to $S^1$ or $[[0,1]]$. If a component $\oldGamma$ of
$\partial\frakB$ is homeomorphic to $[[
0,1]]$, then $r$ maps some
component $\toldGamma$ of $\partial\oldUpsilon$ to $\oldGamma$. Since
$\frakB$ is orientable, $|\toldGamma|$ is a component of
$\partial|\oldUpsilon|$, and is therefore $\pi_1$-injective in the
annulus $|\oldUpsilon|$. The $\pi_1$-injectivity of $|p|$ then
implies that $|p|\big||\toldGamma|:|\toldGamma|\to|\oldGamma|$ is
$\pi_1$-injective; this is impossible, since $|\toldGamma|$ is a
$1$-sphere and $|\oldGamma|$ is an arc. Hence every component of
$\partial\frakB$ is (orbifold)-homeomorphic to $S^1$; equivalently, $|\partial\frakB|$ is a closed $1$-manifold.

Let $\fraks^1_\frakB$ denote the union of all $1$-dimensional components of $\fraks_\frakB$. According to \ref{orbifolds introduced},  we have $\partial|\frakB|=|\partial\frakB|\cup \fraks^1_\frakB$. In the present situation, since
$|\partial\frakB|$ is a closed topological $1$-manifold,
we may write $\partial\frakB$ as a disjoint union $|\partial\frakB|\discup \fraks^1_\frakB$ of closed topological $1$-manifolds. 

Since $p$ is an
orbifold covering and $\oldUpsilon$ is orientable, we have
$|p|(\partial|\oldUpsilon|)=|p(\partial\oldUpsilon)|=|\partial\frakB|$. Since 
$|\partial\frakB|\cap\fraks^1_\frakB=\emptyset$, it follows that
$C:= |p|^{-1}(\fraks^1_\frakB)\subset\inter|\oldUpsilon|$. It follows that for each point $x$ of $C$ there exist a finite (cyclic or dihedral) subgroup $G$ of $O(2)$ containing a reflection, an orientation-preserving (cyclic and possibly trivial) subgroup $G'$ of $G$, orbifold neighborhoods $\frakU'$ and $\frakU$ of $x$ and $p(x)$ in $\oldUpsilon$ and $\frakB$ respectively, a neighborhood $V$ of $0$ in $\RR^2$, and orbifold homeomorphisms $h:V/G\to \frakU$ and $h':V/G'\to \frakU'$, such that $p(\frakU')=\frakU$, and $h^{-1}\circ p\circ h'$ is the natural quotient map from $V/G'$ to $V/G$. In particular, $\frakU'$ may be chosen so that $C\cap|\frakU'|$ is homeomorphic to the intersection of $D^2$ with a finite union of lines through the origin in $\RR^2$.  Since this holds for every $x\in C$, the set $C$ is homeomorphic to a finite graph in which every vertex has even, positive valence. In particular,  $C$ is a finite graph, without endpoints or isolated vertices, contained in the interior of the annulus $|\oldUpsilon|$; this implies:
\Claim\label{either or} Either (a) every component of $C$ is a homotopically non-trivial curve in $|\oldUpsilon|$, or (b) some component of $\inter|\oldUpsilon|-C$ is an open disk.
\EndClaim
Let us set $W=|\frakB|-\fraks^1_\frakB$ and $\tW=|p|^{-1}(W)=|\oldUpsilon|-C$. The definition of $\fraks^1_\frakB$ implies that $\omega(W)$ has only isolated singular points, and hence that $|p|\big|(\tW)$ is a branched covering map of degree $d$ from the (possibly disconnected) $2$-manifold $|\oldUpsilon|-C$ to $W$. Since $\fraks^1_\frakB\subset\partial|\frakB|$, we have $W\supset \inter|\frakB|$. Hence: 
\Claim\label{Is this it?}
For every component $X$ of $\tW$, the restriction of $p$ to $\inter X$ is a branched covering map of degree at most $d$ from $\inter X$ to $\inter|\frakB|$ (and is surjective since the $2$-manifold $|\frakB|$ is connected). 
\EndClaim
If Alternative (b) of \ref{either or} holds, so that $\tW$ has a component $X$ such that $\inter X$ is an open disk, then \ref{Is this it?} implies that
$|p|$ restricts to a (surjective) branched covering map from $\inter{X}$ to $\inter|\frakB|$, and it follows that $\inter|\frakB|$ is an open disk. This is impossible, since $|\oldUpsilon|$ is an annulus and $|p|\big||\oldUpsilon|:|\oldUpsilon|\to|\frakB|$ is $\pi_1$-injective. Hence Alternative (a) of \ref{either or} holds, i.e.
\Claim\label{It's either, not or}
Every component of $C$ is a homotopically non-trivial curve in the interior of the annulus $|\oldUpsilon|$.
\EndClaim

It follows from \ref{It's either, not or} that every component of $\tW$ is the interior of an annulus in $|\oldUpsilon|$. If we fix a component $X_0$ of $\tW$, then by \ref{Is this it?}, $|p|$ restricts to a branched covering map of degree at most $d$ from the open annulus $\inter{X_0}$ to $\inter|\frakB|$, and it follows that $|\frakB|$ is a annulus, a M\"obius band or a disk; but again, the $\pi_1$-injectivity of $|p|$ implies that $|\frakB|$ is not a disk. Thus we have shown that $|\frakB|$ is an annulus or a M\"obius band, and so $\chi(|\frakB|)=0$. Furthermore, since the branched covering map $|p|\big|\inter{X_0}\to|\frakB|$ has degree at most $d$, the image of $(|p|\big|{X_0})_\sharp$ 
has index at most $d$ in $\pi_1(|\frakB|)$. In particular, the image of 
$|p|_\sharp$ has index at most $d$ in $\pi_1(|\frakB|)$.

To prove the final assertion, suppose that every finite subgroup of $\pi_1(\frakB)$ is cyclic. Under this additional hypothesis, we claim:
\Equation\label{with additional}
 \fraks_\oldUpsilon\cap C=\emptyset.
\EndEquation
To prove (\ref{with additional}), suppose that $x$ is a point of $ \fraks_\oldUpsilon\cap C$. Set $y=p(x)$.
Since $x\in C$ we have $y\in \fraks^1_\frakB$, so that $G_y$ contains a reflection (where the notation $G_y$ is defined by \ref{orbifolds introduced}). But since $x\in\fraks_\oldUpsilon$ and $\oldUpsilon$ is orientable, $G_y$ must also contain a rotation. Hence $G_y$ is non-cyclic. But \ref{2-dim case} gives $\chi(\frakB)\le\chi(|\frakB|)=0$, which implies that $\frakB$ is a very good orbifold, and hence that $G_y$ is isomorphic to a subgroup of $\pi_1(\frakB)$. This contradicts the condition that every finite subgroup of $\pi_1(\frakB)$ is cyclic, and (\ref{with additional}) is proved. 

Now assume further that $d\ge2$ and that $\card\fraks_\oldUpsilon=1$. In particular $\fraks_\oldUpsilon\ne\emptyset$. Since $\oldUpsilon$ is orientable we have $\fraks_\oldUpsilon\subset\inter|\oldUpsilon|=|\inter\oldUpsilon|$, which with (\ref{with additional}) gives $\emptyset\ne\fraks_\oldUpsilon\subset|\inter\oldUpsilon|-C$. 
Since $p$ is a covering map we have $p(\inter\oldUpsilon)\subset\inter\frakB$ and $p(\fraks_\oldUpsilon)\subset\fraks_\frakB$. Recalling that $C=p^{-1}(\fraks^1_\frakB)$, we deduce that $\emptyset\ne p(\fraks_\oldUpsilon)\subset\fraks_\frakB\cap( |\inter\frakB|-\fraks^1_\frakB)=\fraks_\frakB\cap\inter|\frakB|$. Thus
\Equation\label{oykhie-doykh}
\fraks_\frakB\cap\inter|\frakB|\ne\emptyset.
\EndEquation
We have observed that $W\supset\inter|\frakB|$, so that $\inter|\frakB|=\inter W$, and that $\omega(W)$ has only isolated singular points.
Set $n=\card\fraks_{\omega(\inter W)}$; by (\ref{oykhie-doykh}) we have $n\ge1$. If $\fraks_{\omega(\inter W)}=\{y_1,\ldots,y_n\}$, and if $k_i$ denotes the order of $y_i$, then by \ref{2-dim case} we have $\chibar(\omega(\inter W))=\chibar(\inter W)+\sum_{i=1}^n(1-1/k_i)=\sum_{i=1}^n(1-1/k_i)\ge n/2\ge1/2$. Since  $p|\obd(\inter \tW):\obd(\inter\tW)\to\obd(\inter W)$ is a degree-$d$ (orbifold) covering, and $d\ge2$, we have $\chibar(\obd(\inter\tW))=d\cdot\chibar(\obd(\inter W))\ge1$. But we have assumed that $\card\fraks_\oldUpsilon=1$, and by (\ref{with additional}) we have  $\fraks_\oldUpsilon\subset\inter\tW$; thus
$\card\fraks_{\obd(\inter\tW)}=1$. If $\ell$ denotes the order of the unique point of $\fraks_{\inter\obd
(\tW)}$, then \ref{2-dim case} gives $\chibar(\obd(\inter\tW))=\chibar(\inter \tW)+(1-1/\ell)$.  Since each component of $\tW$ is a half-open annulus by \ref{It's either, not or}, we obtain $\chibar(\obd(\inter \tW))=1-1/\ell<1$, a contradiction. This establishes the final assertion.
\EndProof

\Lemma\label{uneeda}
Let $\frakK
$ be a strongly \simple, boundary-irreducible, orientable
$3$-orbifold. Set $K=|\frakK|$. Suppose that
$K$ is boundary-irreducible and $+$-irreducible (see
Definition \ref{P-stuff}), and that $K^+$ is not homeomorphic to
$T^2\times[0,1]$ or to a twisted $I$-bundle over a Klein bottle (where $K^+$ is defined as in \ref{P-stuff}). Suppose that $X$ is 
a $\pi_1$-injective, connected subsurface of $\partial K^+$, with $X\cap\fraks_\frakK\ne\emptyset$ and  $\chi(X)=0$,  
and that there exist
disks $G_1,\ldots,G_n\subset\inter X$, where $n\ge0$, such that
$\overline{X-(G_1\cup\ldots\cup G_n)}$ is a component of $|\oldPhi(\frakK)|$.
Then $X$ is an annulus, and there is a
$\pi_1$-injective solid torus $J\subset K^+$, 
with $\partial J\subset K\subset K^+$, such that one of the following alternatives holds:
\begin{itemize}
\item We have $(\partial J)\cap(\partial K)=X\discup X'$ for some annulus $X'\subset\partial K^+ $; furthermore, each of the annuli $X$ and $X'$ has winding
  number $1$ in $J$ (see \ref{great day}) and has non-empty intersection with $\fraks_\frakK$, and $
(\partial  J)\cap\fraks_\frakK\subset \inter X\cup \inter X'$.
\item 
We have $(\partial J)\cap(\partial K)=X$; furthermore, $X$ has winding
  number $1$ or $2$ in $J$, and we have
  $\wt (\partial J) \ge\lambda_\frakK$ (see \ref{lambda thing})  and $(\partial J)\cap\fraks_\frakK\subset\inter X $. 
\end{itemize}
\EndLemma

\Proof
Since $X$ is $\pi_1$-injective in $\partial K^+$, which is a union of components of $\partial K$, and since $K$ is boundary-irreducible, $X$ is $\pi_1$-injective in $K$.

For $i=1,\ldots,n$, the simple closed curve $\partial\obd( G_i)$ is a
component of $\partial\oldPhi(\frakK)$ and is therefore
$\pi_1$-injective in $\frakK$ by \ref{tuesa day}. Hence $\obd(G_i)$ is not discal. Since $G_i$ is a disk, it follows that $\wt(G_i)\ge2$. We have $\wt(X)\ge\sum_{i=1}^n\wt(G_i)$, and hence
\Equation\label{read this one first}
\wt(X)\ge2n.
\EndEquation

Let  $X_0$ denote the component $\overline{X-(G_1\cup\ldots\cup G_n)}$ of $|\oldPhi(\frakK)|$.
Since $X$ is $\pi_1$-injective in $K$, since the surface $X$ has Euler characteristic $0$ and is therefore non-simply connected, 
and since the inclusion homomorphism $\pi_1(X_0)\to\pi_1(X)$ is surjective, the inclusion homomorphism $\pi_1(X_0)\to\pi_1(K)$ is non-trivial.

Since $\chi(X)=0$ and $X\cap\fraks_\frakK\ne\emptyset$, it follows from \ref{2-dim case} that 
 $\chi(\omega(X))<0$. Hence if $n=0$, so that $X_0=X$, we have
 $\chi(\omega(X_0))=\chi(\omega(X))<0$. If $n>0$ then \ref{2-dim case} gives
 $\chi(\omega(X_0))\le\chi(X_0)=\chi(X)-n=-n<0$. Thus in any event we have
 $\chi(\obd(X_0))<0$. If $\oldLambda_0$ denotes the component of
$\oldSigma(\frakK)$ containing $\obd(X_0)$, then by \ref{tuesa day} we have 
$\chi(\oldLambda_0)=\chi(\oldLambda_0\cap\partial\frakK)/2$; every component of  $\oldLambda_0\cap\partial\frakK$ is a component of $\oldPhi(\frakK)$ and therefore has non-positive Euler characteristic by \ref{tuesa day}. Hence $\chi(\oldLambda_0)\le\chi(\obd(X_0))/2<0$. In particular $\oldLambda_0$ is not a \torifold, so that by Lemma \ref{when a tore a fold}, $\oldLambda$ is not a \bindinglike\ \Ssuborbifold\ of $\frakK$. It follows that $\oldLambda_0$ is a \pagelike\ \Ssuborbifold\ of
 $\frakK$. 

We may therefore fix an $I$-fibration $q_0:\oldLambda_0\to\frakB_0$, where $\frakB_0$ is some compact, connected $2$-orbifold, such that
$\calf_0:=|\partial_h\oldLambda_0|=|\oldLambda_0|\cap\partial K$. Note that $X_0$ is a component of $\calf_0$. 
 In the case where
the \pagelike\
 \Ssuborbifold\ $\oldLambda_0$ is twisted, 
$q_0$ is a non-trivial fibration, so that 
 $\calf_0$ is connected; hence we have $\calf_0=X_0$ in this case.
 In the case where
the \pagelike\
 \Ssuborbifold\ $\oldLambda_0$ is untwisted, 
$q_0$ is a trivial fibration, so that
  $\calf_0$ has exactly two components; in this case the component $\calf_0$ distinct from
$X_0$ will be denoted $X_0'$. 
In this case the inclusion map from $X'_0$ to $K$ is homotopic in $K$ to a homemorphism of $X_0'$ onto $X_0$; since the inclusion homomorphism
 $\pi_1(X_0)\to\pi_1(K)$ is non-trivial, it then follows that the
 inclusion homomorphism $\pi_1(X_0')\to\pi_1(K)$ is also non-trivial in this case.

Hence in all cases:
\Claim\label{all cases cool}
For each component $F$ of $\calf_0$, the inclusion homomorphism $\pi_1(F)\to\pi_1(K)$ is non-trivial.
\EndClaim

Let us write $|\Fr \oldLambda_0|=\cala_0\discup\cala_1$, where
$\cala_0$ (respectively $\cala_1$) denotes the union of all components
$A$ of $|\Fr \oldLambda_0|$ such that the inclusion homomorphism
$\pi_1(A)\to\pi_1(K)$ is trivial (respectively,
non-trivial). Since
each component of $\Fr \oldLambda_0\subset\cala(\frakK)$ is an orientable annular orbifold (see 
\ref{tuesa day}), each component of $\cala_i$ is an annulus or a
disk for $i=0,1$, and any annulus component of $\cala_i$ has weight $0$. It is obvious from the definition that $\cala_1$ has no disk components, and hence:
\Claim\label{lame claim}
Each component of $\cala_1$ is a weight-$0$ annulus. 
\EndClaim

If $A$ is any component of $\cala_1$, then $A$ is a component of $|\partial_v\oldLambda_0|$, and hence has a boundary component contained in each component of $|\partial_h\oldLambda_0|=\calf_0$. In particular some component $C$ of $\partial A$ is  contained in $X_0$, and is therefore a component of $\partial X_0$. On the other hand, since $A$ is an annulus by \ref{lame claim}, the definition of $\cala_1$ implies that  $C$ is homotopically non-trivial in $K$. This implies that $C$ cannot have the form $\partial G_i$ with $1\le i\le n$; hence $C$ is a component of $\partial X$. Since $C$ is  homotopically non-trivial in $K$, it is in particular homotopically non-trivial in the $2$-manifold $X$, and is therefore $\pi_1$-injective in $X$. Since by hypothesis $X$ is $\pi_1$-injective in $\partial K$, and $K$ is boudary-irreducible, it follows that $C$ is $\pi_1$-injective in $K$. As $C$ is a component of $\partial A$, the annulus $A$ is $\pi_1$-injective in $K$. This shows:
\Claim\label{omnibust}
 $\cala_1$ is $\pi_1$-injective in $K$.
\EndClaim

Since $\partial G_1,\ldots,\partial G_n$ are components of $\partial
X_0\subset\partial\calf_0$ that bound disks in $\partial K$, they are
components of $\partial\cala_0$. On the other hand, since $X$ is a
compact, orientable $2$-manifold with $\chi(X)=0$, any component of
$\partial X$ carries $\pi_1(X)$; since $X$ is non-simply connected, and is $\pi_1$-injective in $K$,
it follows that any component of $\partial
X$ is a component of $\partial\cala_1$. 
Thus:

\Claim\label{mobsters}
We have $\partial X_0\cap\cala_0=\partial G_1\cup\cdots\cup\partial G_n$, and
$\partial X_0\cap\cala_1=\partial X$.
\EndClaim

If $C$ is any component of $\partial\cala_0$, then $C$ is a simple
closed curve in $\partial\calf_0$ which is homotopically trivial in
$K$. The boundary-irreducibility of $K$ then implies that $C$ bounds a
disk in $\partial K$. A disk bounded by $C$ cannot contain a
component $F$ of $\calf_0$, since the inclusion homomorphism
$\pi_1(F)\to\pi_1(K)$ is non-trivial by \ref{all cases cool}. Hence the disk $D_C$ bounded by $C$ is unique, and satisfies $D_C\cap\calf_0=C$. Thus $D_C$ is a component of $(\partial K)-\inter\calf_0$.

According to \ref{tuesa day},
$\omega(\cala_0)\subset\omega(\cala(\frakK))$ is an orientable annular
orbifold, essential in $\frakK$, and hence every component
$C$ of $\omega(\partial\cala_0)$ is $\pi_1$-injective in
$\frakK$. This shows:
\Claim\label{hit it}
For every component  $C$ of $\partial\cala_0$ we have $D_C\cap\fraks_\frakK\ne\emptyset$.
\EndClaim

For each component $A$ of $\cala_0$, we set $S_A=A\cup\bigcup _{C\in\calc(\partial A)}D_C$. If $A$ is an annulus, then $\partial A$ has two components
$C$ and $C'$; the disks $D_C$ and $D_{C'}$ are components of $(\partial K)-\inter\calf_0$ with distinct boundaries, and are therefore disjoint.
If $A$ is a disk, then of course $\partial A$ has only one component. Hence in any case,
$S_A$ is a $2$-sphere.
Since $K^+$ is irreducible, and since it is immediate from the hypothesis that $\partial K^+\ne\emptyset$, the sphere $S_A$ bounds a unique ball $E_A\subset K^+$ for each component $A$ of $\cala_0$, and we have $\Fr_{K^+}E_A=A$. 


Since $X_0\subset|\oldLambda_0|$, and since the inclusion homomorphism $\pi_1(X_0)\to\pi_1(K)$ is non-trivial while the inclusion homomorphism $\pi_1(K)\to\pi_1(K^+)$ is an isomorphism,  the inclusion homomorphism $\pi_1(|\oldLambda_0|)\to\pi_1(K^+)$ is non-trivial. Hence if $A$ is a component of $\cala_0$, the ball $E_A$  cannot contain $|\oldLambda_0|$. We must therefore have 
\Equation\label{lamb}
E_A\cap|\oldLambda_0|=A.
\EndEquation
 Thus $E_A$ is a component of $\overline{K^+-|\oldLambda_0|}$.

If $A$ and $A'$ are distinct components of $\cala_0$, then $E_A$ and
$E_{A'}$ are components of $\overline{K^+-|\oldLambda_0|}$, and are
distinct because their frontiers $A$ and $A'$ are distinct. Hence:
\Claim\label{piggy-wig}
 $(E_A)_{A\in\calc(\cala_0)}$ is a disjoint family.
\EndClaim
 We set $J=|\oldLambda_0|\cup\bigcup_{A\in\calc(\cala_0)}E_A$. It follows from (\ref{lamb}), \ref{piggy-wig} and the definitions that
\Equation\label{where's the beef?}
\Fr\nolimits_{K^+}J=\cala_1.
\EndEquation

It follows from (\ref{where's the beef?}) that
$\partial J=\cala_1\cup(J\cap\partial{K^+})$. Since $\cala_1$ is properly embedded in $K$, it now follows that
\Equation\label{sundried}
(\partial J )\cap(\partial K)\subset\partial K^+
\EndEquation
and that
\Equation\label{two for a dime}
\partial J\subset K.
\EndEquation

It follows from \ref{omnibust} that  $\cala_1$ is $\pi_1$-injective in $K^+$, which according to the hypothesis is an irreducible $3$-manifold. With the manifold case of Lemma \ref{oops lemma} (which is also a standard fact in $3$-manifold theory), this implies:

\Claim\label{omni or bust}
$J$ is irreducible.
\EndClaim

We set $\calf=J\cap\partial K^+$. It follows from the definitions of
$J$ and of the $E_A$ that
\Equation\label{no end}
\calf=\calf_0\cup\bigcup_{C\in\calc(\partial\cala_0)}D_C.
\EndEquation

Now $X_0$ is a component of $\calf_0$, and according to \ref{mobsters},  $G_1,\ldots,G_n$ are the only disks that have the form $D_C$ for $C\in\calc(\partial\cala_0)$ and meet $X_0$. Since $X=X_0\cup\bigcup_{i=1}^nG_i$, it follows from (\ref{no end}) that
\Claim\label{change lobsters}
$X$ is a component of $\calf$.
\EndClaim

The next step involves two $3$-orbifolds $\frakZ$ and $\frakZ'$, equipped
with $I$-fibrations $r:\frakZ\to\oldDelta$
and $r':\frakZ'\to\oldDelta'$ over $2$-orbifolds, , which are
constructed as follows. We
define $\frakZ$ to be the manifold $D^2\times[0,1]$, set
$\oldDelta=D^2$, and define $r$ to the projection to the first
factor. We define $\frakZ'$ to be the quotient of $\frakZ$ by the
involution $(z,t)\mapsto(\overline{z},1-t)$, where the bar denotes
complex conjugation in $D^2$. Since the involution maps fibers to
fibers, $\frakZ'$ inherits an orbifold fibration from $\frakZ$; its
base, which we denote by $\oldDelta'$, is the quotient of $D^2$ by the
involution $z\mapsto\overline{z}$. Note that $|\oldDelta'|$ is a disk,
and that $\partial|\oldDelta'|$ is the union of the two arcs
$|\partial\oldDelta'|$ and $\fraks_\oldDelta$, which meet only in
their endpoints. Note also that $\partial_v\frakZ$ is an annulus,
while $\partial_v\frakZ'$ is an orbifold having two singular points of
order $2$, and a disk as underlying surface. 

For each component $A$ of $\cala_0$,
define $\frakZ_A$ to be a homeomorphic copy of $\frakZ$ if $A$ is an
annulus, and define  $\frakZ_A$ to be a homeomorphic copy of $\frakZ'$
if $A$ is a disk and $\card\fraks_{\omega(A)}=2$. Thus each $\frakZ_A$ has a fibration $r_A:\frakZ_A\to\oldDelta_A$, where
$\oldDelta_A$ is a $2$-orbifold homeomorphic to either $\oldDelta$ or
$\oldDelta'$, such that under suitable homeomorphic identifications of $\frakZ_A$ and $\oldDelta_A$ with $\frakZ$ and $\oldDelta$ or with $\frakZ'$ and $\oldDelta'$, the fibration $r_A$ is identified with $r$ or $r'$ respectively. 

The fibration of any given orientable annular $2$-orbifold is unique up to fiber-preserving (orbifold) homeomorphism. Hence for each component $A$ of $\cala_0$, since $\Fr (\frakZ_A)$ and $\omega(A)$ are homeomorphic by construction, there exist homeomorphisms
$\eta_A:\partial_v\frakZ_A\to\omega(A)$ and
 $h_A:\partial\oldDelta_A\to q_0(\omega(A)) $
such that 
$q_0\circ\eta_A=h_A\circ (r_A|\partial_v\frakZ_A)$. (It is worth bearing in mind that the homeomorphic $1$-manifolds  $|\partial\oldDelta_A|$ and $|\omega(q_0(A))| $ are arcs when $A$ is a disk, and are $1$-spheres when $A$ is an annulus.) Let $\oldLambda$ denote the $3$-orbifold obtained from the disjoint union of $\oldLambda_0$ with the $\frakZ_A$, where $A$ ranges over the components of $\cala_0$, by gluing $\partial_v\frakZ_A$ to $\omega(A)\subset\oldLambda_0$ via $\eta_A$ for every $A$. Let $\frakB$ denote the $2$-orbifold obtained from the disjoint union of $\frakB_0$ with the $\oldDelta_A$, where $A$ ranges over the components of $\cala_0$, by gluing $\partial\oldDelta_A$ to $\omega(q_0(A))\subset\frakB_0$ via $h_A$ for every $A$. Then there is a well-defined fibration $q:\oldLambda\to\frakB$ which restricts to $q_0$ on $\oldLambda_0$, and to $r_A$ on $\frakZ_A$ for each $A$.

For each component $A$ of $\cala_0$, both $|\frakZ_A|$ and $E_A$ are $3$-balls, and $|\partial_v\frakZ_A|\subset\partial |\frakZ_A|$ and $A\subset\partial E_A$ are either both annuli or both disks. Hence, the homeomorphism 
$|\eta_A|:|\partial_v\frakZ_A|\to A$ 
may be extended to a homeomorphism from $|\frakZ_A|$ to $E_A$; we fix such an extension $t_A:|\frakZ_A|\to E_A$ for each $A$, and define a map $T:|\oldLambda|\to J$ to be the identity on $|\oldLambda_0| $ and to restrict to $t_A$ on $\frakZ_A$ for each $A$. Since
$(E_A)_{A\in\calc(\cala_0)}$ is a disjoint family by \ref{piggy-wig}, and since
$E_A\cap|\oldLambda_0|=A$ for each $A$ by (\ref{lamb}), the map $T:|\oldLambda|\to J$ is a homeomorphism. 

It follows from (\ref{no end}) and the definition of $T$ that
\Equation\label{still no end}
T(|\partial_h\oldLambda|)=\calf.
\EndEquation


From \ref{change lobsters} and (\ref{still no end}) it follows that
$T^{-1}(X))$ is a component of $|\partial_h\oldLambda|$, so that
\Claim\label{ypsilanti}
$\oldUpsilon:=\omega(T^{-1}(X))$ is a component of $\partial_h\oldLambda$. 
\EndClaim

According to \ref{fibered stuff},
$q|\partial_h\oldLambda:\partial_h\oldLambda\to\frakB$ is a degree-$2$
orbifold covering. In particular, by \ref{ypsilanti},
$q|\oldUpsilon:\oldUpsilon\to\frakB$ is an orbifold covering of degree
at most $2$. On the other hand, the homeomorphism $T:|\oldLambda|\to
J$ maps $|\oldUpsilon|$ onto the surface $X$.  By hypothesis $X$ is
$\pi_1$-injective in $ K$ 
and hence in
$K^+$, and in particular in $J$. 
The hypothesis also gives
$\chi(X)=0$. Hence $|\oldUpsilon|$ is $\pi_1$-injective in
$|\oldLambda|$, and $\chi(|\oldUpsilon|)=0$. According to Lemma \ref{affect}, $|q|:|\oldLambda|\to|\frakB|$ is a
homotopy equivalence, and hence
$|q\big|\oldUpsilon|=|q|\big||\oldUpsilon|$ is $\pi_1$-injective. It
therefore follows from Lemma \ref{uneeda this first},  applied
with $q|\oldUpsilon$ playing the role of $p$, that $\chi(|\frakB|)=0$. Hence  $|\frakB|$  is a torus, a Klein bottle, an annulus
or a M\"obius band. 
Since $|q|:|\oldLambda|\to|\frakB|$ is a
homotopy equivalence, and $T:|\oldLambda|\to
J$ is a homeomorphism, $J$ is homotopy-equivalent to $|\frakB|$. Note also that $J\subset K=|\frakK|$, and since $\frakK$ is an orientable $3$-orbifold, $J$ is an orientable $3$-manifold.

First suppose that $|\frakB|$ is a torus or Klein
bottle.  Thus $J$ is homotopy-equivalent to a torus or Klein
bottle. Since the compact, orientable $3$-manifold $J$ is irreducible by \ref{omni or bust}, it follows from \cite[Theorem 10.6]{hempel} that $J$ is 
homeomorphic to $T^2\times[0,1]$ or to a twisted $I$-bundle over a
Klein bottle.  On the other hand, since $|\frakB|$ is closed, $\frakB$ must also be closed and $\fraks_\frakB$ must be $0$-dimensional. 
Since
$\partial_h\oldLambda$ is a covering space of $\frakB$, it follows
that $|\partial_h\oldLambda|$ is a branched cover of $|\frakB|$, and
is therefore closed. Since the components of $\partial_v\oldLambda$
are annular orbifolds, and are therefore not closed, it follows that
$\partial_v \oldLambda=\emptyset$,
i.e. $\partial|\oldLambda|=|\partial_h\oldLambda|$. Since (\ref{still no
  end}), with the definition of $\calf$, implies that 
$T(|\partial_h\oldLambda|)\subset\partial {K^+}$, it now follows that $\partial J=T(\partial|\oldLambda|)\subset\partial {K^+}$, which with the connectedness of ${K^+}$ implies that $J={K^+}$. But then ${K^+}$ is homeomorphic to $T^2\times[0,1]$ or to a twisted $I$-bundle over a Klein bottle, a contradiction to the hypothesis. 

Hence $|\frakB|$ is either an annulus or a M\"obius band. Since   $J$ is homotopy-equivalent to $|\frakB|$, it follows that $\pi_1(J)$ is cyclic. Since the compact, orientable $3$-manifold $J$ is irreducible by \ref{omni or bust}, it then follows from \cite[Theorem 5.2]{hempel} that
\Claim\label{i'll say it's solid}
$J$ is a solid torus.
\EndClaim

Since we have seen that $X$ is $\pi_1$-injective in $K $, it is $\pi_1$-injective in $J$. Furthermore, the hypothesis gives $\chi(X)=0$. In view of \ref{i'll say it's solid}, it follows that
\Claim\label{and an annulus too!}
$X$ is an annulus.
\EndClaim

Since $\pi_1(X)$ and $\pi_1(J)$ are infinite cyclic by \ref{i'll say it's solid} and \ref{and an annulus too!}, and $X$ is $\pi_1$-injective in $K$, it follows that

\Claim\label{stubbins}
$J$ is $\pi_1$-injective in $K$.
\EndClaim

Note that the properties of $X$ and $J$ stated in the lemma before the alternatives given in the two bullet points are covered by (\ref{two for a dime}), \ref{i'll say it's solid}, \ref{and an annulus too!} and \ref{stubbins}.

We have seen that
$q|\partial_h\oldLambda:\partial_h\oldLambda\to\frakB$ is a degree-$2$ orbifold covering. Furthermore, by \ref{ypsilanti}, $\oldUpsilon$ is one component of $\partial_h\oldLambda$. Hence either
(a) $\partial_h\oldLambda=\oldUpsilon$, and $q|\oldUpsilon$ is a $2$-sheeted orbifold covering map onto $\frakB$; or
(b) $\partial_h\oldLambda$ has two components, $\oldUpsilon$ and a second component $\oldUpsilon'$, and $q|\oldUpsilon$ and $q|\oldUpsilon'$ are orbifold homeomorphisms onto $\frakB$.

Suppose that (a) holds. In this case we will show that the solid torus $J$ and the annulus $X$ satisfy the second alternative in the conclusion of the lemma. 
We have $(\partial J)\cap(\partial K)=J\cap\partial K^+=\calf=T(|\partial_h\oldLambda|)$, by (\ref{sundried}) and (\ref{still no end}), and $\partial_h\oldLambda=\oldUpsilon$ by (a). 
Hence $(\partial J)\cap(\partial K) =T(|\oldUpsilon|)$. But by definition (see \ref{ypsilanti}) we have $\oldUpsilon=\omega(T^{-1}(X))$, so that $X=T(|\oldUpsilon|)$. Thus $(\partial J)\cap(\partial K)=X$. 

Next note that $\partial J=(J\cap\partial K^+)\cup\Fr_{K^+}J=X\cup\cala_1$, in view of \ref{where's the beef?}. Since 
$\cala_1\cap\fraks_\frakK=\emptyset$ by \ref {lame claim}, it follows that 
$\partial J\cap\fraks_\frakK\subset \inter X$. 

To show that $\wt (\partial J) \ge\lambda_\frakK$, or equivalently that $\wt X\ge\lambda_\frakK$, we distinguish several subcases. By hypothesis we have $X\cap\fraks_\frakK\ne\emptyset$, i.e. $\wt X\ge1$; hence the assertion is true if $\lambda_\frakK=1$. If $\lambda_\frakK=2$ and $n\ge1$, then by \ref{read this one first}
we have $\wt(X)\ge2n\ge2=\lambda_\frakK$. There remains the subcase in which $\lambda_\frakK=2$ and $n=0$. In this subcase, we must show $\wt X\ge2$, and since $\wt X\ge1$, it is enough to show that $\wt X\ne1$. Since $n=0$, it follows from \ref{mobsters} that 
$\partial X_0=\partial X$, and since $X_0\subset X$ it then follows that $X=X_0$. It also follows from \ref{mobsters} that $\partial X_0\cap\cala_0=\emptyset$; since 
each component of $\cala_0=\partial_v\oldLambda_0$ is saturated in $\oldLambda_0$, and therefore meets
$X_0=\partial_h\oldLambda_0$, it follows that $\cala_0=\emptyset$, so that $\oldLambda=\oldLambda_0$ and $\frakB=\frakB_0$.    Since $T$ is defined to be the identity on $|\oldLambda_0|$, it is the identity on its entire domain $|\oldLambda|$. Hence $\obd(X)= \obd(T(X))=\oldUpsilon$. 
To show that $\wt X\ne1$, i.e. that $\card\fraks_\oldUpsilon\ne1$, we will apply the last sentence of  Lemma \ref{uneeda this first}, with $q|\oldUpsilon$ playing the role of $p$.  We have already seen that the general hypotheses of Lemma \ref{uneeda this first} hold with this choice of $p$. By \ref{and an annulus too!}, $|\oldUpsilon|=X$ is an annulus. Since (a) holds, $p=q|\oldUpsilon$ has degree $2$.  Since $\lambda_\frakK=2$, every cyclic subgroup of $\pi_1(\frakK)$ is cyclic. But $\pi_1(\frakB)$ is isomorphic to $\pi_1(\oldLambda)$ since $\oldLambda$ has an $I$-fibration over $\frakB$, and $\oldLambda$ is $\pi_1$-injective in $\frakK$ since it is an \Ssuborbifold\ of $\frakK$. Hence every finite subgroup of $\pi_1(\frakB)$ is cyclic. It therefore indeed follows from
 the last sentence of  Lemma \ref{uneeda this first}
that $\card\fraks_\oldUpsilon\ne1$, and the proof that 
$\wt (\partial J) \ge\lambda_\frakK$ is complete in the case where (a) holds.


To complete the proof of the second alternative of the conclusion in
the case where (a) holds, it remains to show that the winding number of $X$
in $J$ (which is non-zero by the $\pi_1$-injectivity of $X$ is at most $2$. Since the homeomorphism $T:|\oldLambda|\to
J$ maps $|\oldUpsilon|$ onto $X$, this is equivalent to showing that, if $P$ denotes
the image of the inclusion homomorphism
$\pi_1(|\oldUpsilon|)\to\pi_1(|\oldLambda|)$, then the index $|\pi_1(|\oldLambda|):P|$ is at most $2$. 
Since
$|q|:|\oldLambda|\to|\frakB|$ is a homotopy equivalence by Lemma \ref{affect}, we have
$|\pi_1(|\oldLambda|):P|=|\pi_1(|\frakB|):|p|_\sharp(\pi_1(|\oldUpsilon|))|$,
where $p$ denotes the degree-$2$ orbifold covering
$q|\oldUpsilon:\oldUpsilon\to\frakB$. 
According to 
Lemma \ref{uneeda this first}, the index $|\pi_1(|\frakB|):|p|_\sharp(\pi_1(|\oldUpsilon|))|$ is bounded above by the degree of the orbifold covering $p$, which is equal to $2$. This completes the proof of the conclusion in the case where Alternative (a) holds.
Now suppose that (b) holds. In this case we will show that the solid
torus $J$ and the annulus $X$ satisfy the first alternative in the conclusion of the lemma. By definition (see \ref{ypsilanti}) we have $\oldUpsilon:=\omega(T^{-1}(X))$, so that $X=T(|\oldUpsilon|)$.  Set $X'=T(|\oldUpsilon'|)$.  Since $q|\oldUpsilon$ and $q|\oldUpsilon'$ are orbifold homeomorphisms onto $\frakB$, the orbifolds $\oldUpsilon$ and $\oldUpsilon'$ are homeomorphic. In particular, $|\oldUpsilon|$ and $|\oldUpsilon'|$ are homeomorphic, and hence so are $X$ and $X'$. Thus $X'$ is an annulus.

Since $\oldUpsilon$ and $\oldUpsilon'$ are the components of $\partial_h\oldLambda$, it follows from (\ref{still no end}) that $X$ and $X'$ are the components of $\calf$.
Since $X\cup X'=\calf=J\cap\partial {K^+}$, and since 
$\Fr_{K^+}J=\cala_1$ by \ref {where's the beef?}, each component of $\overline{(\partial J)- (X\cup X')}$ is a component of $\cala_1$, and hence by \ref{lame claim}
is a weight-$0$ annulus. It follows that
$\overline{(\partial J)- (X\cup X')}$ is a union of two weight-$0$ annuli. Hence $(\partial  J)\cap\fraks_\frakK\subset \inter X\cup \inter X'$, and the annuli $X$ and $X'$ have the same winding number in $J$. To determine this common winding number, note that since $q|\oldUpsilon:\oldUpsilon\to\frakB$ is an orbifold homeomorphism, and since
$|q|:|\oldLambda|\to|\frakB|$ is a homotopy equivalence by the first assertion of Lemma \ref{affect}, the inclusion $|\oldUpsilon|\to|\oldLambda|$ is a homotopy equivalence; since  $X=T(|\oldUpsilon|)$, the inclusion $X\to J$ is a homotopy equivalence, and hence the common winding number of $X$ and $X'$ in $J$ is $1$.

To show that  the first alternative in the conclusion of the lemma holds in this case, it remains only to show that each of the annuli $X$ and $X'$ has non-empty intersection with $\fraks_\frakK$. It follows directly from the hypothesis that $X\cap\fraks_\frakK\ne\emptyset$. To show that $X'\cap\fraks_\frakK\ne\emptyset$, we first note that by \ref{no end},
$X'$ contains a component of $\calf_0$; since $X'\ne X$, this component of $\calf_0$ is distinct from $X_0$, and is therefore equal to $X_0'$ in the notation introduced above. By \ref{no end}, $X'$ is the union of $X_0'$ with all disks of the form $D_C$, where $C$ ranges over the components of $\partial\cala_0$ contained in $X_0'$. Consider first the subcase in which $\cala_0\ne\emptyset$, and choose a component $A$ of $\cala_0$. Since $A$ is in particular a component of $|\Fr \oldLambda_0|$, it is saturated in the fibration of $\oldLambda_0$, and therefore meets
 every component of $ |\partial_h\oldLambda_0|=\calf_0$ (see \ref{fibered stuff}). Hence some component
 $C_0$ of $\partial A\subset\partial\cala_0$ is contained in $X_0'$. We therefore have $D_{C_0}\subset X'$. But by \ref{hit it} we have
 $D_{C_0}\cap\fraks_\frakK\ne\emptyset$, and hence $X'\cap\fraks_\frakK\ne\emptyset$ in this subcase.
 
 Finally, consider the subcase in which $\cala_0=\emptyset$. We then have $X=X_0$, $X'=X_0'$, $\oldLambda=\oldLambda_0$ and $\frakB=\frakB_0$.    Since $T$ is defined to be the identity on $|\oldLambda_0|$, it is the identity on its entire domain $|\oldLambda|$. Hence $\obd(X)= \obd(T(X))=\oldUpsilon$ and $\obd(X')= \obd(T(X'))=\oldUpsilon'$.  It follows that $X\cap\fraks_\frakK=\fraks_\oldUpsilon$ and that $X'\cap\fraks_\frakK=\fraks_{\oldUpsilon'}$. We have observed that $\oldUpsilon$ and $\oldUpsilon'$ are homeomorphic orbifolds, and hence $\card(X\cap\fraks_\frakK)=\card(\fraks_\oldUpsilon)=\card(\fraks_{\oldUpsilon'})=\card(X'\cap\fraks_\frakK)$. Since we have $X\cap\fraks_\frakK\ne\emptyset$, it follows that $X'\cap\fraks_\frakK\ne\emptyset$ in this subcase as well.
 \EndProof

\section{Tori in the underlying space of a
  $3$-orbifold}\label{underlying tori section}

In this section we give the construction referred to in the introduction, and at the beginning of Section \ref{seek 'n' find} for modifying a $2$-suborbifold of a hyperbolic orbifold $\oldOmega$, such that the components of $|\oldTheta|$ are incompressible tori in $|\oldOmega|$, so as to obtain a suborbifold that is useful for finding lower bounds for the invariants introduced in Section \ref{darts section}.

\begin{notationremarks}\label{wait star}
If $S$ is a subset of an orbifold $\oldPsi$ such that
$S\cap\fraks_\oldPsi$ is finite, we will define a quantity
$\wt^*_\oldPsi S$ by setting $\wt^*_\oldPsi S=\wt_\oldPsi S$ if
$\wt_\oldPsi S$ is even or $\lambda_\oldPsi=1$, and $\wt^*_\oldPsi
S=\wt_\oldPsi S+1 $ if $\wt_\oldPsi S$ is odd and
$\lambda_\oldPsi=2$. (See \ref{wuzza weight} and \ref{lambda thing} for the definitions of $\wt_\oldPsi\in\NN$ and $\lambda_\oldPsi\in\{1,2\}$.) Note that $\wt^*_\oldPsi S$ is always divisible by $\lambda_\oldPsi$. Note also that if $S$ and $S'$ are subsets of $\oldPsi$ such that $S\cap\fraks_\oldPsi$ and $S'\cap\fraks_\oldPsi$ are both finite, and if $\wt_\oldPsi S\le\wt_\oldPsi S'$, then $\wt^*_\oldPsi S\le\wt^*_\oldPsi S'$. Hence
 if $\wt^*_\oldPsi S'<\wt^*_\oldPsi S$ then $\wt_\oldPsi S'<\wt_\oldPsi S$. Moreover, 
 if $\wt_\oldPsi S=\wt_\oldPsi S'$ then $\wt_\oldPsi^* S=\wt_\oldPsi^* S'$.

We will write $\wt^* S$ for $\wt^*_\oldPsi S$ when the orbifold $\oldPsi$ is understood.
\end{notationremarks}

\Lemma\label{pre-modification} Let $\oldPsi$ be a
compact, orientable, strongly \simple, boundary-irreducible 
$3$-orbifold containing no negative
turnovers. 
Set $N=|\oldPsi |$. Suppose that each component of
$\partial N$ is a sphere, and that $N$ is $+$-irreducible (see
Definition \ref{P-stuff}). Let $K$
 be a non-empty, proper, compact, connected, $3$-dimensional
 submanifold of $N$. 
Assume that $\calt:=\Fr_NK$ is contained in $\inter N$ and is 
in general position with respect to
$\fraks_\oldPsi$,  and  that its components are all incompressible tori in $N$ (so that $K^+$ is naturally identified with a submanifold of $N^+$ by \ref{plus-contained}).
Assume that $K^+$ is not homeomorphic to $T^2\times[0,1]$ or to a twisted $I$-bundle over a Klein bottle.
Suppose that either (i) some component of $\obd(\calt)$ is
compressible in $\oldPsi$, or (ii) each component of $\obd(\calt)$ is
incompressible in $\oldPsi$ (so that $\obd(K)$ is boundary-irreducible and strongly \simple\ by Lemma \ref{oops lemma}, and hence $\kish(\obd(K))$ is defined in
view of \ref{tuesa day}) and
$$\chibar(\kish(\obd(K)))<\min\bigg(1,\frac14\wt^*\nolimits_\oldPsi(\calt)\bigg).$$
Then at least one of the following conditions holds:
\begin{enumerate}
\item There exist a disk $D\subset N$ with
$\partial D=D\cap \calt $, such that $D$ 
in general position with respect to
$\fraks_\oldPsi$,
and 
a disk $E\subset \calt $ such that $\partial E=\partial D$ and 
$\wt_\oldPsi( E)>\wt_\oldPsi(D)$; furthermore, if 
$\lambda_\oldPsi=2$,
then $\max(\wt_\oldPsi(
E)-\wt_\oldPsi(D),\wt_\oldPsi D)\ge2$. 
\item There is 
a solid torus $J\subset K^+$, $\pi_1$-injective in $N^+$,
 with $\partial J\subset
K\cap\inter N
\subset K\subset K^+$,
  such that 
  $(\partial J)\cap(\partial K)$ is a union of two disjoint annuli $X$ and $X'$ contained in $\calt$, each having winding
  number $1$ in $J$ (see \ref{great day}) and each having non-empty
  intersection with $\fraks_\oldPsi$, and 
  $(\partial J)\cap\fraks_\oldPsi\subset \inter X\cup \inter X'$.
\item 
There is 
a solid torus $J\subset K^+$, $\pi_1$-injective in $N+$, 
 with $\partial J\subset K\cap\inter N\subset K\subset K^+$,   such that 
  $(\partial J)\cap(\partial K)$ is an annulus $X\subset\calt$, having winding
  number $1$ or $2$ in $J$, and 
we have
  $\wt (\partial J) \ge\lambda_\oldPsi$ and $(\partial J)\cap\fraks_\oldPsi\subset \inter X $.
\end{enumerate}
\EndLemma

\Proof[Proof of Lemma \ref{pre-modification}]
Let us  consider the case of the lemma in which Alternative (i) of the
hypothesis holds. 
Since some component of $\obd(\calt )$ is compressible in $\oldPsi$, it
follows from Proposition \ref{kinda dumb} that there is a discal $2$-suborbifold $\frakD$ of $\oldPsi$ such that $\frakD\cap\obd(\calt)=\partial\frak D$ and  $|\frakD|$ is in general position with respect to $\fraks_\oldPsi$, but such that there is no discal $2$-suborbifold $\frakE$ of $\obd(\calt)$ with $\partial\frakE=\partial\frakD$.
Equivalently,
there is a disk $D\subset N$, in general position with respect to $\fraks_\oldPsi$,  with
$\gamma:=\partial D=D\cap \calt $, such that $\wt_\oldPsi D\le1$ and such that $\gamma$ does not bound any disk of weight at most $1$ in $ \calt $. Since by hypothesis the components of $\calt $ are incompressible,
there is a disk $E\subset \calt $ with $\partial E=\gamma$. Hence
$\wt E\ge2$. In particular we have
$\wt E>\wt D$. Furthermore, if $\wt (E)-\wt(D)=1$, then we must have
$\wt E=2$ and $\wt D=1$; this implies that $D\cup E$ is a $2$-sphere
of weight $3$ in $|\oldPsi|$, 
in general position with respect to $\fraks_\oldPsi$. If $\lambda_\oldPsi=2$, then since $\oldPsi$ 
contains
no embedded negative turnovers, it follows from Proposition \ref{almost obvious}
that $|\oldPsi|$ contains no weight-$3$
sphere in general position with respect to $\fraks_{\oldPsi}$, and hence  $\wt (E)-\wt(D)\ge2$. 
This shows that Alternative (1) of
the conclusion of the lemma holds in this case.

For the rest of the proof, we will assume that Alternative (ii) of the hypothesis
holds (so that in particular every component of $\obd(\calt)$ is
incompressible in $\oldPsi$). The strategy of the proof under this
assumption is to try to find a subsurface $X\subset K$ satisfying the
hypothesis of Lemma \ref{uneeda}, with $\frakK=\obd(K)$; when such an
annulus can be found, Alternative (2) or (3) of the
conclusion of the present lemma will be seen to hold. When the quest
for such a subsurface fails, Alternative (1) of the conclusion will turn out to hold.

Set 
$\frakU=\overline{\obd(K)-\oldSigma(\obd(K))}$. According to the definition of $\kish(\obd(K))$ (see \ref{tuesa day}), $\frakU$ is the disjoint union of $\kish(\obd(K))$ with certain components of $\frakH(\obd(K))$. Since the latter components have Euler characteristic $0$, we have $\chi(\frakU)=\chi(\kish(\obd(K)))$. Hence Alternative (ii) of the hypothesis gives
\Equation\label{bought}
\chibar(\frakU)<\min\bigg(1,\frac14\wt^*\nolimits(\calt)\bigg).
\EndEquation

We claim:
\Claim\label{a half or better}
If $\oldXi$ is any component of $\frakU\cap\obd(\partial K)$ such that
$\wt|\oldXi|$ is odd, we have $\chibar(\oldXi)\ge1/2$. 
\EndClaim

To prove \ref{a half or better}, set $q=\wt|\oldXi|$, write
$\fraks_\oldXi=\{z_,\ldots,z_q\}$, and let $e_i$ denote the order of
$z_i$ for $i=1,\ldots,q$. Then by \ref{2-dim case} we have
$\chibar(\oldXi)=\chibar(|\oldXi|)+\sum_{i=1}^q(1-1/e_i)\ge\chibar(|\oldXi|)+q/2$. If
$|\oldXi|$ is not a disk or a sphere, we have $\chibar(|\oldXi|)\ge0$,
and $q\ge1$ since $q$ is odd; hence $\chibar(\oldXi)\ge1/2$. If
$|\oldXi|$ is a disk, we have $q\ne1$ since $\partial\oldXi\subset\partial\oldPhi(\omega(K))$
is $\pi_1$-injective in $\partial(\obd(K))$ by
\ref{tuesa day}; since $q$ is odd we have $q\ge3$, and hence
$\chibar(\oldXi)\ge-1+3/2=1/2$. Now suppose that $|\oldXi|$ is a sphere. Since
$\oldXi\subset\partial\obd(K)=\partial\oldPsi\cup\calt$, and the
components of $\calt$ are tori, $\oldXi$ is a component of
$\partial\oldPsi$. In this case, if $q=1$ then $\oldXi$ is a bad $2$-orbifold; this is impossible since the strongly \simple\
$3$-orbifold $\oldPsi$ is in particular very good (see \ref{oops}). Hence $q\ne1$. Next note that since
$\oldPsi$ is strongly \simple\ and boundary-irreducible, its boundary component $\oldXi$ has negative Euler characteristic by \ref{boundary is negative}; hence if $q=3$, then $\oldXi$ is a negative
turnover, a contradiction to the hypothesis. 
Since $q$ is odd we now have $q\ge5$, and hence $\chibar(\oldXi)\ge-2+5/2=1/2$. Thus \ref{a half or better} is proved.

Now we claim:
\Claim\label{cause i may have said so}
If $Q$ is a union of components of $\partial|\frakU|$, each of which has strictly positive
genus, we have
$$\wt^*\nolimits Q\le4\chibar (\frakU).$$
\EndClaim

To prove \ref{cause i may have said so}, set $n=\wt Q$,
and write $Q\cap\fraks_\oldPsi=\{x_1,\ldots,x_n\}$. Then $x_1,\ldots,x_n$
are the singular points of the $2$-orbifold $\obd(Q)$. Let $p_i$
denote the order of the singular point $x_i$. Observing that
$p_i\ge2$ for $i=1,\ldots,n$, and that $\chibar(Q)\ge0$ since each
component of $Q$ has
positive genus, and applying \ref{2-dim case}, we find that
$\chibar(\obd(Q))=\chibar(Q)+\sum_{i=1}^n(1-1/{p_i})\ge n/2$.
Hence $n\le2\chibar(\obd(Q))$, i.e.
\Equation\label{mcfoofus}
\wt Q\le2\chibar (\obd(Q)).
\EndEquation
 Since $\obd(Q)$ is a union of components of
$\partial\frakU$, and since every  component of
$\partial\frakU$ has non-positive Euler characteristic by
\ref{tuesa day}, it follows that 
$\wt Q\le2\chibar(\partial\frakU)$. But since $\frakU$ is a compact $3$-orbifold we have
$\chibar(\partial\frakU)=2\chibar(\frakU)$, and hence
\Equation\label{almost cause i may have said so}
\wt Q\le4\chibar (\frakU).
\EndEquation
If $\lambda_\oldPsi=1$ or $\wt Q$ is even, the conclusion of
\ref{cause i may have said so} follows from (\ref{almost cause i may
  have said so}) in view of the definition of $\wt^*Q$. There remains
the case in which $\lambda_\oldPsi=2$ and $\wt Q$ is odd. In this
case, write $\partial|\frakU|=Q\discup Q'$, where $Q'$ is a union of
components of $\partial|\frakU|$. We have $\wt(|\partial\frakU|)=\wt
Q+\wt Q'$. Since $\lambda_\oldPsi=2$, each component of
$\fraks_\oldPsi\cap|\frakU|$ is an arc or a simple closed curve. Hence
$\wt\partial|\frakU|$ is equal to twice the number of arc components
of $\fraks_\oldPsi\cap|\frakU|$, and is therefore even. Since $\wt Q$
is odd, it now follows that $\wt Q'$ is odd. Each component of $Q'\cap\Fr\frakU$ is a component of $\Fr\frakU=\obd(\cala(\obd(K))$, and is therefore annular by \ref{tuesa day}; in particular, each component of $Q'\cap\Fr\frakU$ has weight $0$ or $2$. Hence $\wt(Q'\cap\partial K)$ is odd. There therefore exists a component $\oldXi$ of $\obd(Q'\cap\partial K)$ such that $\wt|\oldXi|$ is odd. According to \ref{a half or better}, we have $\chibar(\oldXi)\ge1/2$. But 
by
\ref{tuesa day},
$\partial\oldXi\subset\partial\oldPhi(\omega(K))$
is $\pi_1$-injective in
$\obd(K)$ and hence in $\obd(Q')$; since every  component of
$\partial\frakU$ has non-positive Euler characteristic by
\ref{tuesa day}, it follows that $\chibar(\obd(Q'))\ge\chibar(\oldXi)\ge1/2$. Combining this with (\ref{mcfoofus}), we find $\chibar(\partial\frakU)=\chibar(\obd(Q'))+\chibar(\obd(Q)\ge(1+\wt Q)/2$. Since $\wt^* Q=1+\wt Q$ in this case, it follows that $\wt^*Q\le2\chibar(\partial\frakU)=4\chibar(\frakU)$, as required. This completes the proof of \ref{cause i may have said so}.


Let $\calx_0$ denote the union of all components of
$|\oldSigma(\obd(K))|\cap\calt=|\oldPhi(\obd(K))|\cap\calt$ (see \ref{tuesa day}) that
 are not contained in disks in $\calt$. 
Let $\calx\subset\calt$ denote the union of $\calx_0$ with all disks
in $\calt$ whose boundaries are contained in $\partial \calx_0$. Since the
components of $\calt$ are tori, each component of $\calx$ is either an
annulus or a component of $\calt$. Furthermore, no annulus component
of $\calx$ is contained in a disk in $\calt$; that is, the annulus component
of $\calx$ are homotopically non-trivial in the torus components of
$\calt$ containing them. Hence each component of
$\caly:=\overline{\calt-\calx}$ is also  either a homotopically non-trivial
annulus in $\calt$ or a component of $\calt$. Since the
components of $\calt$ are incompressible in $N$, it now follows that the
submanifolds $\calx$ and $\caly$ are $\pi_1$-injective in $N$.

It follows from the definition of $\caly$ that every component of $(\inter
\caly)\cap|\oldSigma(\obd(K)|$ is contained in a disk in $\inter \caly$. Hence $(\inter
\caly)\cap|\oldSigma(\obd(K)|$ is contained in a (possibly empty) disjoint
union of disks
  $E_1,\ldots,E_m\subset\inter \caly$, where $\partial E_i\subset
 |\oldSigma(\obd(K)|$ whenever $1\le i\le m$. The definitions also imply
  that $\caly_0:=\overline{\caly-(E_1\cup\cdots\cup E_m)}$ is a union of components
    of $\calt\cap|\frakU|$. 
\abstractcomment{\tiny The expression $\calt\cap|\frakU|$ replaces
  $\overline{\calt\setminus|\oldPhi(\obd(K)|}$ in the version new-5point4.tex.}
%
Hence $\obd(\caly_0)$ is a union of
    components of $\obd(\calt)\cap\frakU$, which are in particular components of $\frakU\cap\partial\obd(K)$.
Since each component of $\frakU\cap\partial\obd(K)$ has non-positive Euler
    characteristic by \ref{tuesa day}, it follows that
\Equation\label{gunman}
\chibar(
\obd(\calt)\cap\frakU))\ge\chibar(\obd(\caly_0)).
\EndEquation

On the other hand, since each component of $\caly$
is an annulus or a torus, we have $\chibar(\caly_0)=m$. By
\ref{2-dim case} it follows that
\Equation\label{by crossref}
\chibar(\obd(\caly_0))\ge m.
\EndEquation
Now by \ref{tuesa day} we have
    $\chibar(\frakU)
=\chibar(\frakU\cap\omega(\partial
K))/2=(\chibar(\frakU\cap\omega(\calt))+\chibar(\frakU\cap\partial
\oldPsi))/2$. Since $\frakU\cap\partial \oldPsi$ is a union of
components of 
$\overline{\partial\obd(K)-\oldPhi(\obd(K))}$, it follows from
\ref{tuesa day} that
$\chibar(\frakU\cap\partial \oldPsi)\ge0$, and hence
$\chibar(\frakU)
\ge\chibar(\obd(\calt)\cap\frakU)/2$.
With (\ref{gunman}) and (\ref{by
      crossref}) this gives $\chibar(\frakU)\ge m/2$. Since
$    \chibar(\frakU)<1$ by (\ref{bought}), it follows that
$m<2$. Hence
\Equation\label{charmin}
m\le1.
\EndEquation


Next we claim:
\Claim \label{mabel syrble}
If a component $F$ of $|\partial\frakU|$ contains a component
of $\caly_0$, then $F$ is a torus.
\EndClaim

To prove \ref{mabel syrble}, let $Y_0$ be a component of $\caly_0$ contained
in $F$, and let $Y$ and $T$ respectively denote the components of $\caly$
and $\calt$ containing $Y_0$. Since $Y$ is a component of $\caly$, either $Y$ is an  annulus which is homotopically
non-trivial in $N$, or $Y=T$. If $Y=T$ then $Y_0$ is a genus-$1$ subsurface of
$F$. If $Y$ is a homotopically non-trivial annulus, then the
components of $\partial Y\subset Y_0\subset F$ are homotopically non-trivial
simple closed curves in $N$. Thus in any event, $F$ contains either a genus-$1$
subsurface or a simple closed curve which is homotopically non-trivial in $N$. Hence
$F$ cannot be a sphere.

If $F$ has genus at least two, then $\chibar(F)\ge2$. By \ref{2-dim case}, it follows that $\chibar(\obd(F))\ge2$. Since every
component of $ \partial\frakU$
has non-positive
Euler characteristic by \ref{tuesa day}, we have
$\chibar( \partial\frakU)
\ge\chibar(\obd(F))\ge2$, and hence $\chibar(\frakU)=\chibar( \partial\frakU)/2\ge1$.
This a
contradiction, since $\chibar(\frakU)<1$ by (\ref{bought}). Thus \ref{mabel syrble} is proved.

Now we claim:
\Claim\label{if not why not}
If $\wt^*\caly>4\chibar (\frakU)$, then Alternative (1) of the
conclusion of the lemma holds. 
\EndClaim

To prove \ref{if not why not}, 
let $Q_0$ denote the union of all
components of $|\partial\frakU|$ that contain components of
$\caly_0$. According to \ref{mabel syrble},  each component of $Q_0$ is
a torus. Since in particular every component of $Q_0$ has strictly positive genus, it follows from 
\ref{cause i may have said so} that
\Equation\label{where's the oakie?}
\wt\nolimits^*Q_0\le4\chibar(\frakU).
\EndEquation
On the other hand, since $\caly_0\subset Q_0$, we have $\wt\caly_0\le\wt Q_0$, which by \ref{wait star} implies $\wt^*\caly_0\le\wt^*Q_0$. With (\ref{where's the oakie?}), this implies $\wt^*\caly_0\le4\chibar(\frakU)$.
Now if $m=0$, we have
$\caly=\caly_0$, and hence
$\wt^*\caly\le4\chibar (\frakU)$.
Thus \ref{if not why not} is vacuously true when $m=0$. In view of
(\ref{charmin}), it remains only to consider the subcase $m=1$.

In this
subcase set $E=E_1$, so that $\caly_0=\overline{\caly-E}$. 
Let us index the components of $\caly$ as $Y^0,\ldots,Y^r$, where $r\ge0$
and $E\subset Y^0$. 
Then 
$Y^0_0:=\overline{Y^0-E}$ is a component of $\caly_0$. In particular we have
$Y^0_0\subset\calt\cap|\frakU|\subset\partial|\frakU|$.
Let $T$ denote the (torus) component of $\calt$ containing $Y^0$.
Let $F$ denote the component of
$Q_0$ containing $Y^0_0$. Then $F$ is
a torus by \ref{mabel syrble}.

Let $\Delta$ denote the component of $\overline{F-Y^0_0}$ containing
$\alpha:=\partial E\subset\partial Y^0_0$. If $Y^0=T$, then
$Y^0_0$ is a genus-$1$ surface with one boundary component, and
hence $\Delta$ is a disk. If $Y^0$ is an annulus, then $Y^0_0$ is a
connected planar surface with three boundary components; two of
these, the components of $\partial Y^0$, are homotopically non-trivial
in $N$ and hence in the torus $F$. It follows that $\alpha$, the third
boundary component of $Y^0_0$, is homotopically trivial in $F$. Thus
in any event,
$\Delta$ is a disk.

Set $Z=Y^0_0\cup\Delta$.
Then
the surfaces $Z,Y^1,\ldots, Y^r$ are contained in
$Q_0$. The surfaces $Y^1,\ldots,Y^r$ are pairwise disjoint
because they are distinct components of $\caly$. If $1\le i\le r$, we
have $Y^0\cap Y^i=\emptyset$ and hence $Z\cap
Y^i\subset\Delta\subset\inter Z$. Since
$Y^i$ is connected, and $Q_0$ is a closed $2$-manifold, we must have either $Y^i\subset\Delta$ or $Z\cap
Y^i=\emptyset$. But since $Y^i$ is a homotopically non-trivial annulus or an incompressible torus
in $N$, it cannot be contained in the disk
$\Delta$. This shows that the surfaces $Z,Y^1,\ldots, Y^r$ are
pairwise disjoint. Hence
\Equation\label{jooolie}
\wt\nolimits Z+\sum_{i=1}^r\wt\nolimits Y^i\le \wt Q_0.
\EndEquation
On the other hand, (\ref{where's the oakie?}) and the hypothesis of Claim
\ref{if not why not} give
$\wt^* Q_0\le4\chibar
(\frakU)<\wt^*\caly$, which by \ref{wait star} implies $\wt Q_0<\wt\caly$. Since $\wt\caly=\sum_{i=0}^r\wt Y^i$, it follows that
\Equation\label{not sloop}
\wt Q_0<\sum_{i=0}^r\wt Y^i.
\EndEquation
From (\ref{jooolie}) and (\ref{not sloop}) it follows that
$\wt Z< \wt Y^0$,
which may be rewritten as
$\wt Y^0_0+\wt \Delta<
\wt Y^0_0+\wt E.$
Hence 
\Equation\label{one last wisk}
\wt\Delta<
\wt E.
\EndEquation

Since the disk $\Delta$ is contained in $K$ and its boundary $\beta$ lies in
$\partial K$, we may modify $\Delta$ by a small isotopy, constant on
$\beta$, to obtain a disk $D$ which is properly embedded in $K$. In
particular we have $D\cap \calt=\beta$. We may choose the isotopy in such a
way that $D$ is in general position with respect to $\fraks_\oldPsi$ and
$\wt D=\wt\Delta$. With (\ref{one last
  wisk}), this gives $\wt D<
\wt E$. 
Now suppose that $\lambda_\oldPsi=2$ and that
$\max(\wt
(E)-\wt(D),\wt D)<2$. Since $\wt D<
\wt E$, we have either $\wt E=2$ and $\wt D=1$, or
$\wt E=1$ and $\wt D=0$; hence
$D\cup E$ is a $2$-sphere
of weight $1$ or $3$ in $|\oldPsi|$. If the sphere $D\cup E$ has weight $1$ then $\obd(D\cup E)$ is a bad $2$-suborbifold of $\oldPsi$; this is impossible since the strongly \simple\ $3$-orbifold
$\oldPsi$ is in particular very good (see \ref{oops}). If $\lambda_\oldPsi=2$, then since $\oldPsi$ 
contains
no embedded negative turnovers, it follows from Proposition \ref{almost obvious} that $|\oldPsi|$ contains no weight-$3$
sphere 
in general position with respect to $\fraks_\oldPsi$. This contradiction shows that if $\lambda_\oldPsi=2$ then $\max(\wt
(E)-\wt(D),\wt D)\ge2$. Thus Alternative (1) of the conclusion of the
lemma holds  in this situation, and Claim \ref{if not why not} is  established.


Next we claim:
\Claim\label{cornway}
Either $\calx\cap\fraks_\oldPsi\ne\emptyset$, or Alternative (1) of the
conclusion of the lemma holds. 
\EndClaim

To prove \ref{cornway}, assume that Alternative (1) of the
conclusion of the lemma does not hold.
Then by \ref{if not why not}, we have
\Equation\label{cause i said so}
\wt\nolimits^* \caly\le4\chibar (\frakU).
\EndEquation

Note that
\Equation\label{me quation}
\wt\calx+\wt \caly=\wt\calt. 
\EndEquation
If we assume that $\wt\calx=0$, then it follows from (\ref{me quation}) that 
$\wt \caly=\wt\calt$, which with \ref{wait star} implies that $\wt^* \caly=\wt^*\calt$. Combined with
(\ref{cause i said so}) this gives
$\wt^* \calt\le4\chibar (\frakU)$, which contradicts
(\ref{bought}). Hence we must have
 $\wt\calx>0$, i.e. 
$\calx\cap\fraks_\oldPsi\ne\emptyset$, and
\ref{cornway} is proved.

In view of \ref{cornway}, we will assume for the remainder of the proof
that $\calx\cap\fraks_\oldPsi\ne\emptyset$. Let us fix a component $X$ of
$\calx$ such that $X\cap\fraks_\oldPsi\ne\emptyset$. The definitions of 
$\calx_0$ and $\calx$ imply that $X$ contains a unique component $X_0$ of
$\calx_0$, and that there are disjoint disks $G_1,\ldots,G_n$ in $\inter X_0$ such that
$X_0=\overline{X-(G_1\cup \cdots\cup G_n)}$. The definition of
$\calx_0$ also implies that $X_0$ is a component of
$|\oldPhi(\obd(K))|$. As we have seen that $\calx$ is $\pi_1$-injective in $N$, in particular $X$ is $\pi_1$-injective in $K$.
Since $X\subset K$,
we have $X\cap\fraks_{\obd(K)}=X\cap\fraks_\oldPsi\ne\emptyset$. Note
also that $K$ is boundary-irreducible and $+$-irreducible, since $N$
is $+$-irreducible, $\partial N$ consists of sphere components, and the components of $\calt=\partial K^+$ are
incompressible in $N$. We have seen that each component of $\calx$ is
either a (torus) component of $\calt$ or an
annulus; in particular, we have
 $\chi(X)=0$. Furthermore, by
hypothesis,  $K^+$ is not homeomorphic to $T^2\times[0,1]$ or to a
twisted $I$-bundle over a Klein bottle. Thus
all the hypotheses of Lemma \ref{uneeda} hold with
$\frakK=\obd(K)$. Hence by Lemma 
\ref{uneeda}, $X$ is an annulus, and there is a
$\pi_1$-injective solid torus $J\subset K^+$, 
with $\partial J\subset K\subset K^+$, such that one of the following alternatives holds:
\begin{itemize}
\item We have $(\partial J)\cap(\partial K)=X\discup X'$ for some annulus $X'\subset\partial K^+ $; furthermore, each of the annuli $X$ and $X'$ has winding
  number $1$ in $J$ and has non-empty
  intersection with $\fraks_\frakK$, and 
$(\partial  J)\cap\fraks_\frakK\subset \inter X\cup \inter X'$.
\item We have $(\partial J)\cap(\partial K)=X$; furthermore, $X$ has winding
  number $1$ or $2$ in $J$, and 
we have  $\wt (\partial J) \ge\lambda_\frakK$ and $(\partial J)\cap\fraks_\frakK\subset \inter X $.
\end{itemize}

Note also that since $\frakK$ is a suborbifold of $\oldPsi$, and $J\subset K$, we have $(\partial J)\cap\fraks_\oldPsi=(\partial J)\cap\fraks_\frakK$. Hence the two alternatives above imply respectively that $(\partial  J)\cap\fraks_\oldPsi\subset \inter X\cup \inter X'$ and that  
$(\partial  J)\cap\fraks_\oldPsi\subset \inter X$. Likewise, since $\frakK$ is a suborbifold of $\oldPsi$,
it follows from the definitions that $\lambda_\frakK\ge\lambda_\oldPsi$, so that when the second alternative stated above holds, we have $\wt(\partial J)\ge\lambda_\oldPsi$. Finally, since the components of $\calt$ are incompressible in $N$, the submanifold $K^+$ of $N^+$ is $\pi_1$-injective, and hence the $\pi_1$-injectivity of $J$ in $K^+$ implies that it is $\pi_1$-injective in $N^+$. Thus one of the alternatives (2), (3) of the
conclusion holds.
\EndProof

\Lemma\label{modification} Let $\oldPsi$ be a
compact, orientable, strongly \simple, boundary-irreducible
$3$-orbifold containing no embedded negative turnovers. 
Set $N=|\oldPsi |$. Suppose that each component of
$\partial N$ is a sphere, and that $N$ is $+$-irreducible. Let $K$
 be a non-empty, proper, compact, connected, $3$-dimensional
 submanifold of $N$, and set $\calt=\Fr_N K$.  Assume that $\calt$ is contained in $N$ and is in general position with respect to
$\fraks_\oldPsi$,  and  that its components are all incompressible tori in $N$ (so that $K^+$, which by \ref{plus-contained} is naturally identified with a submanifold of $N^+$,  is irreducible and boundary-irreducible), 
 and that   either (a) $K^+$ is
acylindrical, or (b) $\overline{N-K}$ is connected and
$h(K)\ge3$. 
Suppose in addition that either (i) some component of $\obd(\calt )$ is
compressible in $\oldPsi$, or (ii) each component of $\obd(\calt )$ is
incompressible in $\oldPsi$ (so that $\obd(K)$ is boundary-irreducible and strongly \simple\ by Lemma \ref{oops lemma}, and hence $\kish(\obd(K))$ is defined in
view of \ref{tuesa day}) and
$\chibar(\kish(\obd(K)))<\lambda_\oldPsi/4$, or (ii$'$) $\wt_\oldPsi^*\calt \ge4$, each component of $\obd(\calt )$ is
incompressible in $\oldPsi$, and
$\chibar(\kish(\obd(K)))<1$. Then there is a compact,
connected $3$-manifold $K_1\subset N$ having the following properties:
\begin{itemize}
\item every component of $\calt _1:=\Fr_NK_1$ is an incompressible torus in
  $\inter N$, in general position with respect to
$\fraks_\oldPsi$;
\item $\wt_\oldPsi(\calt _1)<
\wt_\oldPsi\calt$; and if $\lambda_\oldPsi=2$, then $\max(\wt_\oldPsi(\calt )-
\wt_\oldPsi(\calt_1)
,\wt_\oldPsi(\calt_1))\ge2$; 
\item $h(K_1)\ge h(K)/2$; and
\item $h(N-K_1)\ge h(N-K)$.
\end{itemize}
In particular, $K_1$ is a proper, non-empty submanifold of $N$.

Furthermore, if (a) holds, $K_1$ may be chosen so that
$K_1^+$ and $K^+$ (which by \ref{plus-contained} are identified with
  submanifolds of $\plusN$) are isotopic in $\plusN$; and if (b) holds, $K_1$ may be chosen so that $\overline{N- K_1} $ is connected.
\EndLemma

\Proof[Proof of Lemma \ref{modification}] 
We set $P=K^+$. As was observed in the statement of the lemma, since $\calt$ is closed and has no sphere components, it follows from \ref{plus-contained} that we may
identify $P$ with a submanifold of $N^+$. It was also observed in the statement of the lemma that the $+$-irreducibility of $P$, the incompressibility of the components of $\calt$ and the hypothesis that $N$ has only sphere components, imply that $P$ is irreducible and boundary-irreducible. (This implication, which is standard in $3$-manifold theory, is also included  in the manifold case of Lemma \ref{oops lemma}.) 

Note that the hypotheses imply that $K$ has at
least one frontier component, or equivalently that it has at
least one torus boundary component. In particular,
$h(K)$ and $h(N-K)$ are both strictly positive. Hence if $K_1$ is a
submanifold with the first four properties listed, we have $h(K_1)\ge
h(K)/2>0$ and $h(N-K_1)\ge
h(N-K)>0$, and therefore $K_1\ne\emptyset$ and $N-K_1\ne\emptyset$. Thus
the first four properties listed do imply that $K_1$ is a non-empty,
proper submanifold of $N$, as asserted.

We will now turn to the proof that there is a submanifold $K_1$ with
the stated properties. 

Since $K$ is a non-empty, proper submanifold of $N$, we have $\calt\ne\emptyset$.
If $\wt\calt=0$ then the components of $\obd(\calt)$ are incompressible toric suborbifolds of the
$3$-orbifold $\oldPsi$, a contradiction since $\oldPsi$ is strongly \simple\ (see \ref{oops}). Hence
$\wt\calt>0$, which by the definition of $\wt^*\calt$ implies $\wt^*\calt>0$. Since $\wt^*\calt$ is divisible by $\lambda_\oldPsi$ according to \ref{wait star}, we have $\wt^*\calt\ge\lambda_\oldPsi$. It follows that if
Alternative (ii) of the hypothesis holds, i.e. if 
each component of $\obd(\calt )$ is
incompressible in $\oldPsi$ and
$\chibar(\kish(\obd(K)))<\lambda_\oldPsi/4\le1/2$, then
$\chibar(\kish(\obd(K)))<\min(1,{\wt^*(\calt)/4})$; this is Alternative (ii) of
the hypothesis of Lemma \ref{pre-modification}. Likewise, it is
immediate that if
Alternative (ii$'$) of the hypothesis of the present lemma holds, i.e. if 
each component of $\obd(\calt )$ is
incompressible in $\oldPsi$ and if $\wt^*\calt\ge4$ and
$\chibar(\kish(\obd(K)))<1$,
then
$\chibar(\kish(\obd(K)))<\min(1,{\wt^*(\calt)/4})$. 
Thus either of the
Alternatives (ii) or (ii$'$) of the hypothesis of the present lemma implies 
Alternative (ii) of the hypothesis of Lemma \ref{pre-modification}.
Alternative (i) of the hypothesis of the present lemma is identical to
Alternative (i) of the hypothesis of Lemma
\ref{pre-modification}. Note also that if $K^+$ 
were homeomorphic to $T^2\times[0,1]$ or to a twisted $I$-bundle over a
Klein bottle, then $K^+$ would fail to be acylindrical, and we would have
$h(K)=h(K^+)=2$, so that neither of the alternatives (a) or (b) of the
hypothesis of the present lemma would hold; thus $K^+$ is not
homeomorphic to either of these manifolds. Hence the hypotheses of the present
lemma imply those of Lemma
\ref{pre-modification}, so that under these hypotheses one of the alternatives (1)---(3) of the conclusion of Lemma
\ref{pre-modification} must hold.

Consider the case in which Alternative (1) of the conclusion of Lemma
\ref{pre-modification} holds: that is,
there is a disk $D\subset N$ with
$\partial D=D\cap \calt $, such that $D$ is in general position with respect to $\fraks_\oldPsi$ and 
there is a disk $E\subset \calt $ such that $\partial E=\partial D$ and 
$\wt E>\wt D$. Furthermore, if $\lambda_\oldPsi=2$,
we may suppose $D$ and $E$ to be chosen so that $\max(\wt(
E)-\wt(D),\wt D)\ge2$.

Since $N$ is $+$-irreducible, the sphere $D\cup
E$ bounds a ball in $\plusN$, and hence the surface $(\calt-E)\cup
D$ is isotopic to $\calt$ in $\plusN$, by an isotopy that is constant
on $\overline{\calt-E}$. Hence $\calt_1:=(\calt-E)\cup
D\subset\inter N$ bounds a submanifold $P_1$ of $\plusN$ which is isotopic to $P$ in $\plusN$, by an isotopy that is constant
on $\overline{\calt-E}$. Set $K_1=P_1\cap N$; since $\partial P_1=\calt_1\subset\inter N$, we have
$K_1^+=P_1$.

We have $\wt(\calt)-\wt(\calt_1)=\wt(E)-\wt(D)$.
Since $D$
is
in general position with respect 
to $\fraks_\oldPsi$ and
$\wt E>\wt D$,
the surface $\calt_1=\Fr_NK_1$ is 
in general position with respect to
$\fraks_\oldPsi$, and
$\wt\calt_1<\wt\calt$. Furthermore, if $\lambda_\oldPsi=2$, 
then 
$\max(\wt(\calt )-
\wt(\calt_1)
,\wt(\calt_1))\ge \max(\wt(
E)-\wt(D),\wt D)\ge2$.
Since $P_1$ is isotopic to $P$,
the components of $
\Fr_{N}K_1=
\Fr_{N^+}P_1$ are incompressible tori. 
The existence of an isotopy between $P$ and $P_1$ also implies that
  $h(K_1)=h(P_1)=h(P)=h(K)\ge h(K)/2$ and that $h(N-K_1)=h(N^+-P_1)=h(N-P)=
  h(N-K)$. 
To verify the last sentence of the lemma in this case, note that $K_1^+=P_1$ and $K^+=P$ are isotopic regardless of whether (a) or (b) holds; and that the existence of an isotopy between $P_1$ and $P$ implies that $\overline{N-K_1}=\overline{N^+-P_1}$ is homeomorphic to $\overline{N-K}=\overline{N^+-P}$, so that if (b) holds then $\overline{N-K_1}$ is connected.
Thus the proof of the lemma is complete in this case.


Now consider the case in which Alternative (2) of the conclusion of Lemma
\ref{pre-modification} holds: that is, 
there is 
a solid torus $J\subset P$, $\pi_1$-injective in $N^+$,
 with $\partial J\subset
K\cap\inter N
\subset K\subset P$, 
  such that 
  $(\partial J)\cap(\partial K)$ is a union of two disjoint annuli 
%
 $X$ and $X'$ contained in $\calt$, each having winding
  number $1$ in $J$, and we have
\Equation\label{harps}
X\cap\fraks_\oldPsi\ne\emptyset\quad\text{and}\quad X'\cap\fraks_\oldPsi\ne\emptyset,
\EndEquation
and 
  $(\partial J)\cap\fraks_\oldPsi\subset \inter X\cup \inter X'$.
Note that the set
$\Fr_{P}J=\overline{(\partial J)-(X\cup X')}$ is also a disjoint
union of two annuli, $A_1$ and $A_2$, which also have winding number
$1$ in $\calt$, and are properly embedded in $P$.  Since $(\partial 
  J)\cap\fraks_\oldPsi\subset \inter X\cup \inter X'$, we have 
\Equation\label{houris}
A_i\cap\fraks_\oldPsi
=\emptyset\quad\text{for}\quad i=1,2.
\EndEquation
Note also that since the annuli $X$, $X'$, and $A_1$ and $A_2$ have winding number
$1$ on the solid torus $J$, which is $\pi_1$-injective in $N^+$, they are themselves
$\pi_1$-injective in $N^+$.

Set $P_0=\overline{P-J}$. 
 Then $P_0$ is a compact, possibly
disconnected, $3$-dimensional submanifold of $P\subset N$. Since $\partial J\subset\inter N$,  
every component of $K\cap\partial N$ is a sphere contained in either $\inter J$ or $\inter P_0$. 
Thus if we set $K_0=P_0\cap N$, each component of $P_0$ is obtained from a component of $K_0$ by attaching balls along certain spheres in the boundary. Set
$\calt_0=\Fr_NK_0=\Fr_{N^+}P_0$. We have $\calt_0=(\calt-(X\cup X'))\cup(A_1\cup A_2)$, so that
(\ref{harps}) and (\ref{houris}) imply
\Equation\label{on travel}
\wt\calt_0\le (\wt\calt)-2.
\EndEquation

The manifold $K_1$ whose existence is asserted by the lemma will be
constructed as a suitable component of $K_0$. Before choosing a
suitable component, we will prove a number of
facts about $P_0$ and its components that will be useful in establishing the
properties of $K_1$.

Since $P=P_0\cup J$, and since
$P_0\cap J=\Fr_{P}J$ is the union of the disjoint annuli $A_1$ and $A_2$,
each of which has winding number $1$ in the solid torus $J$, we have
\Equation\label{erfer}
h(P)\le h(P_0)+1.
\EndEquation

The following property of $P_0$ will also be needed:

\Claim\label{tony baloney}
Every component of $\Fr_{N^+}P_0$ is a (possibly compressible) torus in
  $\inter {N}$, 
in general position with respect to
$\fraks_\oldPsi$.
\EndClaim

To prove \ref{tony baloney}, we need only note that since the annuli
$X$ and $X'$ are $\pi_1$-injective in the (torus) components of
$\calt$ containing them, the components of $\overline{\calt-(X\cup
  X')}$ are tori and annuli. As the closed, orientable $2$-manifold $\Fr_{N^+}P_0$ is
obtained from the disjoint union of $\overline{\calt-(X\cup
  X')}$ with $A_1$ and $A_2$  by gluing certain boundary components
in pairs, its components have Euler characteristic $0$. Since $\calt$
is by hypothesis 
in general position with respect to $\fraks_\oldPsi$, it follows from
\ref{houris} that $\Fr_{N^+}P_0$ is in general position with respect  to
$\fraks_\oldPsi$. This establishes \ref{tony
  baloney}.

We will also need:
\Claim\label{da props}
Let $L$ be any component of $P_0$. Then:
\begin{itemize}
\item $L$ is $\pi_1$-injective in ${N^+}$;
\item every component of $\Fr_{N^+}L$ is a torus in
  $\inter {N}$, in general position with respect to
$\fraks_\oldPsi$;
\item ${\wt}(\Fr_{N^+}L)\le \wt(\calt)-2$;  and
\item if $L$ is not a solid torus then every component of $\Fr_{N^+}L$
  is incompressible in ${N^+}$.
\end{itemize}
\EndClaim

To prove the first assertion of \ref{da props}, note that the
frontier of $L$ in ${P}$ consists of one or both of the
$\pi_1$-injective annuli $A_1$, $A_2$. Hence
$L$ is itself $\pi_1$-injective in ${P}$. On the other hand, since the
components of $\Fr_{N^+}{P}=\calt$ are incompressible tori, ${P}$ is in turn
$\pi_1$-injective in ${N^+}$. Hence $L$ is $\pi_1$-injective in ${N^+}$.

Next note that every component of $\Fr_{N^+}L$
is a component of $\Fr_{N^+}P_0$; hence the assertion that every component of $\Fr_{N^+}L$ is a torus in
  $\inter {N^+}$, in general position with respect to
$\fraks_\oldPsi$, is an immediate consequence of \ref{tony baloney}. 

Since $\Fr_{N^+}L\subset\Fr_{N^+}P_0=\calt_0$, the inequality
$\wt(\Fr_{N^+}L)\le \wt(\calt)-2$ is an immediate
consequence of (\ref{on travel}).

Since ${N^+}$ is irreducible, and since we have seen that $L$ is
$\pi_1$-injective in $N^+$ and that every component of $\Fr_{N^+}
L$ is a torus, the fourth assertion of \ref{da props} follows from
Proposition \ref{final assertion}, applied with $M=N^+$, and with $L$ defined as above. Thus \ref{da props} is proved.

The rest of the proof of the lemma in this case is divided into subcases. Consider the
subcase in which Alternative (a) of the hypothesis holds. 
Since $P$ is
acylindrical in this subcase, the
mutually parallel annuli $A_1$ and $A_2$
cannot be essential in $P$.
Since these annuli are $\pi_1$-injective in
${N^+}$, 
they are in particular
$\pi_1$-injective in $P$;
they must therefore be
boundary-parallel in $P$. 
It follows that some component 
of
$P_0$ is isotopic to $P$ in $N^+$. We fix such a component of $P_0$
and denote it by $P_1$. We set $K_1=P_1\cap N$.
It follows from \ref{da props},
taking $L=P_1$, that 
$\wt(\Fr_{N}K_1)=\wt(\Fr_{N^+}P_1)\le \wt(\calt)-2$.  
Since $P_1$ is isotopic to $P$,
the components of $
\Fr_{N}K_1=
\Fr_{N^+}P_1$ are incompressible tori. 
The existence of an isotopy between $P$ and $P_1$ also implies that
  $h(K_1)=h(P_1)=h(P)=h(K)\ge h(K)/2$ and that $h(N-K_1)=h(N^+-P_1)=h(N-P)=
  h(N-K)$. Thus the conclusions of the lemma (including the last sentence, which in this subcase asserts that $P=K^+$ and $P_1=K_1^+$ are isotopic) hold with this choice of $K_1$.

Next consider the
subcase in which Alternative (b) of the hypothesis holds and $K_0$ is
connected. In this subcase we take $K_1=K_0$. According to (\ref{erfer})
we have $h(K_1)=h(P_0)\ge h(P)-1= h(K)-1$; since Alternative (b) of the hypothesis
gives $h(K)\ge3$, it follows that $h(K_1)\ge h(K)/2$, as asserted in
the conclusion of the lemma. It also follows that $h(K_1)\ge2$, so that
$K_1^+$ is not a solid torus. Hence by \ref{da props}, the components of
$\Fr_{N^+} K_1$ are incompressible in ${N^+}$. 

It follows from \ref{da props} that  
$\wt(\Fr_{N^+}K_1)\le \wt(\calt)-2$. 
This in turn trivially implies that $\calt_1:=\Fr_NK_1$ has weight strictly less than $\wt\calt$, and that $\max(\wt(\calt)-
\wt_\oldPsi(\calt_1)
,\wt_\oldPsi(\calt_1))\ge
\wt(\calt)-
\wt_\oldPsi(\calt_1)\ge
2$, regardless of whether $\lambda_\oldPsi$ is $1$ or $2$; in particular, the second bullet point of the conclusion of the lemma holds in this subcase.

Since Alternative (b) holds, the manifold ${N}-K$
is connected, and hence so is ${N^+}-P$.
Since
$\overline{N^+-P_0}=\overline{N^+-P}\cup J$, and since
$\overline{N^+-P}\cap J$ is the union of the disjoint annuli $X$ and $X'$,
each of which has winding number $1$ in the solid torus $J$, the manifold 
${N^+}-P_0$ is also connected, and $h({N^+}-P_0)\ge
h({N^+}-P)$. This means that ${N}-K_1$ is connected, as required by the last sentence of the lemma when (b) holds. and that $h({N}-K_1)\ge
h({N}-K)$.
Thus all
the conclusions of the lemma are established in this subcase.

The remaining subcase of this case is the one in which Alternative (b) of the hypothesis holds and $K_0$ is
disconnected. In this subcase, $K_0$ has two components since $\Fr_PJ$
has two components $A_1$ and $A_2$. For $i=1,2$, let $K_i$ denote
the component of $K_0$ containing $A_i$. After re-indexing if
necessary we may assume that $h(K_1)\ge h(K_2)$. In this subcase,
$P_0$ has two components, and they may be indexed as $P_1$ and $P_2$,
where $P_i$ is the union of $K_i$ with certain ($3$-ball) components
of $\overline{N^+-N}$. Since the $A_i$
have winding number $1$ in $J$, the manifold $P$ is homeomorphic to
the union of $P_1$ and $P_2$ glued along an annulus. Hence $h(K)=h(P)\le
h(P_1)+h(P_2) =h(K_1)+h(K_2)\le 2h(K_1)$, i.e. $h(K_1)\ge h(K)/2$, as asserted in
the conclusion of the lemma. Since Alternative (b) of the hypothesis
gives $h(K)\ge3$, it follows that $h(K_1)\ge2$, so that
$K_1^+$ is not a solid torus. Hence by \ref{da props}, the components of
$\calt_1:=\Fr_NK_1=\partial K_1$ are incompressible in ${N^+}$. It also follows from \ref{da props} that
$\wt\Fr K_1=\wt(\Fr_{N^+}P_1)\le \wt(\calt)-2$.

To establish the conclusions of the lemma in this subcase, it remains to
prove that
$\overline{N-K_1}$ is connected and that $h(\overline{N-K_1})\ge h(\overline{N-K})$. 
Since Alternative (b) holds, 
the manifold $\overline{N-K}$
is connected, and hence so is $\overline{N^+-P}$. We have
$\overline{N^+-P_0}=\overline{N^+-P}\cup J$, where $J$ is connected, and
$\overline{N^+-P}\cap J=X\cup X'\ne\emptyset$; hence $\overline{N^+-P_0}$ is
connected, and therefore so is $\overline{N-K_0}$. Next note that
$\overline{{N}-K_1}=\overline{{N}-K_0}\cup K_2$, where $K_2$ is
connected, and has non-empty intersection with $\overline{{N}-K_0}$
because it is a component of $K_0$. Hence $\overline{{N}-K_1}$ is
connected, as required by the last sentence of the lemma when (b) holds.

To estimate 
$h(\overline{{N}-K_1})$, first note that each component of $P_0=\overline{P-J}$ contains a component of $\overline{\calt-(X\cup X')}$. Since
$P_0=\overline{P-J}$ is disconnected in this subcase, $\overline{\calt-(X\cup X')}$ is also disconnected. But each of the annuli $X$, $X'$ is $\pi_1$-injective in ${N^+}$, and therefore is non-separating in the (torus) component of $\calt$ containing it. Hence
$X$ and $X'$
must lie on the same component of $\calt$. As disjoint, homotopically
non-trivial annuli on a torus, $X$ and $X'$ are homotopic in $\calt$ and hence
in $\overline{{N^+}-P}$. We have written $
\overline{{N^+}-P_0}$ as the
union of $\overline{{N^+}-P}$ with the solid torus $J$, and the
intersection of $\overline{{N^+}-P}$ with $J$ consists of the two
annuli 
$X$ and $X'$, which have winding number $1$ in $J$. As these annuli
are homotopic in $\overline{{N^+}-P}$, we have
$h(\overline{{N^+}-P_0})=h(\overline{{N^+}-P})+1$. On the other hand, we have
$\overline{{N^+}-P_1}=\overline{{N^+}-P_0}\cup P_2$, and the components of
$ \overline{{N^+}-P_0}\cap P_2$ are components of $\Fr_{N^+}P_2$ and are
therefore tori by \ref{da props}. Applying Lemma \ref{another goddam
  torus lemma}, taking $U=\overline{{N^+}-P_0}$ (which has been seen to be
connected) and $V=P_2$, we
deduce that $h(\overline{{N^+}-P_1})\ge
h(\overline{{N^+}-P_0})-1$. Hence $h(\overline{{N^+}-P_1})\ge
h(\overline{{N^+}-P})$, or equivalently $h(\overline{{N}-K_1})\ge h(\overline{{N}-K})$. 
Thus
the conclusions of the lemma are established in this subcase.

Finally, consider the case in which Alternative (3) of the conclusion of Lemma
\ref{pre-modification} holds: that is, 
there is 
a solid torus $J\subset P$, $\pi_1$-injective in $N^+$,
 with $\partial J\subset K\cap\inter N\subset K\subset P$,   such that 
  $(\partial J)\cap(\partial K)$ is an annulus $X\subset\calt$, having winding
  number $1$ or $2$ in $J$, and we have
  $\wt (\partial J) \ge\lambda_\oldPsi$ and $(\partial J)\cap\fraks_\oldPsi\subset \inter X $.

In this case,
  $A:=\overline{(\partial J)-X}$ is a properly
 embedded,
  $\pi_1$-injective annulus in $P$, also having winding number $1$ or $2$ in $J$. If Alternative (a) of the
  hypothesis holds, then $A$ must be boundary-parallel in $P$; that
  is, some component $Z$ of $\overline{P-A}$ must be a solid torus in
  which $A$ has winding number $1$. Thus $\overline{P-Z}$ is isotopic
  to $P$.
If  $Z=\overline{P-J}$, then $P$ is homeomorphic to $J$
  and is therefore a solid torus, a contradiction to the
  incompressibility of $\calt$.
Hence we must have $Z=J$, so that $P$ is isotopic to $\overline{P-J}$. In this subcase, we set
  $K_1=(\overline{P-J})\cap N$. Setting $\calt_1=\Fr
  K_1=\partial(\overline{P-J})$, we have $\calt_1=(\calt-X)\cup A$, 
  so that $\wt\calt_1=\wt\calt-\wt X+\wt A$. Since   $\wt (\partial J)
  \ge\lambda_\oldPsi$ and $(\partial J)\cap\fraks_\oldPsi\subset \inter X $, we have
 $\wt X\ge\lambda_\oldPsi$ and $\wt A=0$. Hence
  $\wt\calt_1\le\wt\calt-\lambda_\oldPsi$; in particular, $\wt\calt_1<\wt\calt$, and if $\lambda_\oldPsi=2$ then  
$\max(\wt_\oldPsi(\calt )-
\wt_\oldPsi(\calt_1)
,\wt_\oldPsi(\calt_1))\ge
\wt_\oldPsi(\calt )-
\wt_\oldPsi(\calt_1)\ge
\lambda_\oldPsi=2$.
The existence of an isotopy between $P$ and $P_1$ also implies that
  $h(K_1)=h(P_1)=h(P)=h(K)\ge h(K)/2$ and that $h(N-K_1)=h(N^+-P_1)=h(N-P)=
  h(N-K)$. Thus the conclusions of the lemma (including the last sentence, which in this subcase asserts that $P=K^+$ and $P_1=K_1^+$ are isotopic) hold with this choice of $K_1$.

There remains the subcase in which Alternative (b) holds, i.e. 
$\overline{{N}-K}$ is connected and
$h(K)\ge3$. 
In this subcase we will set $P_1=\overline{P-J}$ and $K_1={P_1}\cap N$. Then $P_1$ is connected since $P$ and $A=\Fr_PJ$ are connected; hence $K_1$ is connected. Since $P_1\cup J=P$, and
since $P_1\cap J=A$ is connected, we have $h(P)\le
h(P_1)+h(J)=h(P_1)+1$, so that $h(K_1)=h(P_1)\ge h(P)-1=h(K)-1$. Since $h(K)\ge3$ it
follows that $h(K_1)\ge2$.

Note that $P$ is $\pi_1$-injective in ${N^+}$ since its frontier
components are incompressible, and that $P_1$ is $\pi_1$-injective in
$P$ since its frontier $A$ is $\pi_1$-injective. Hence $P_1$ is
$\pi_1$-injective in ${N^+}$, i.e. $K_1$ is
$\pi_1$-injective in ${N}$. Since $\calt_1=\Fr_{N} K_1=\Fr_{N^+} P_1=\overline{\calt
-X}\cup A$, the surfaces $\calt$ and $\calt_1$ are homeomorphic; in
particular, the components of $\calt_1$ are tori. Since
$h(P_1)\ge2$, the manifold $P_1$ is not a solid torus. It therefore
follows from Lemma \ref{final assertion}, applied with $^+$ and $P_1$ playing the respective roles of  $M$ and $L$, that the
components of $\calt_1$ are incompressible in $N^+$, and hence in $N$. This proves the first
  property of $K_1$ asserted in the conclusion of the present
  lemma. 
Next note that since   $\wt (\partial J) \ge\lambda_\oldPsi$ and $(\partial J)\cap\fraks_\oldPsi\subset X $,
we have
  $\wt X\ge\lambda_\oldPsi$ and   $\wt A=0$. Hence 
$\wt\calt_1\le\wt(\calt)-\lambda_\oldPsi=\wt(\calt)-\lambda_\oldPsi$, which implies that
$\wt\calt_1<\wt\calt$, and that
$\max(\wt(\calt )-
\wt(\calt_1)
,\wt(\calt_1))\ge \wt(\calt )-
\wt(\calt_1)\ge2$ if $\lambda_\oldPsi=2$.
 The inequality 
$h(K_1)\ge h(K)/2$ holds because $h(K_1)\ge h(K)-1$ and
$h(K)\ge3$. Since ${N}-K$ is connected by Alternative (b), ${N^+}-P$ is connected. Since, in addition, $J$ is a solid torus and $\overline{{N^+}-P}\cap J=X$ is an annulus,
the manifold $\overline{{N^+}-P_1}=\overline{{N^+}-P}\cup J$ is connected, so that $N-K_1$ is connected, as required by the last sentence of the lemma when (b) holds; and
$h(\overline{{N^+}-P_1})$ is equal either to $h(\overline{{N^+}-P})$ or to
$h(\overline{{N^+}-P})+1$. In particular, $h({N}-K_1)=h({N^+}-P_1)\ge h({N^+}-P)=h({N^+}-K)$.Thus all the conclusions are seen to hold in this final subcase.
\EndProof


\abstractcomment{\tiny
The inequality involving $\sigma$ in the hypothesis
could be weakened by replacing $\sigma$ by the corresponding thing
involving only the contribution of the tori, not the spheres. This may
mean all the results of this section remain true (and are a little
stronger) if $\sigma$ is replaced by that variant. I need to think
about whether this affects the final conclusions of the paper.

In the proof in the case where (3) holds, I could easily have gotten
$p-2$ as the upper bound for $\wt\calt_1$, since the
intersection with $J$ is an arc which must have two endpoints. The
only case where one doesn't get this is the case in which (1) holds,
but I had an idea that the non-existence of hyperbolic triangle groups might give
it there as well. This means that one may get much more info than I
thought without the assumption of no nodes.}

\Lemma\label{get a newer mop}
Let $m\ge2$ be an integer. Let $\oldPsi$ be a strongly \simple,
boundary-irreducible, orientable  $3$-orbifold containing no embedded
negative turnovers. 
Set $N=|\oldPsi |$. Suppose that each
component of $\partial N$ is a sphere, that $\plusN$ is a graph manifold, and that $h(N)\ge\max(4m-4,m^2-m+1)$.
Then
there is a compact
submanifold $K$ of $N$ such that
\begin{itemize} 
\item each component of $\Fr_N K $ is an incompressible torus in
  $\inter N $, in general position with respect to $\fraks_\oldPsi$,
\item $K$ and $N -K$ are connected, 
\item $\min(h(K),h(N-K))\ge m$, 
\item each component of $\obd(\Fr_N K )$ is incompressible in $\oldPsi$ (so that $\obd(K)$ is boundary-irreducible and strongly \simple\ by Lemma \ref{oops lemma}, and hence $\kish(\obd(K))$ is defined in
view of \ref{tuesa day}), and
\item either \begin{enumerate}[(A)]
\item $\chibar(\kish(\obd(K)))\ge1$, or 
\item  $\chibar(\kish(\obd(K)))\ge\lambda_\oldPsi/4$ and
$\wt^*
(\Fr_NK)<4$. 
\end{enumerate}
\end{itemize}
\EndLemma

\Proof
Set $M=\plusN$. 
Since $h(M)=h(N)\ge\max(4,m^2-m+1)$, the hypothesis of Lemma
\ref{just right} holds.
Hence there is a compact
submanifold $P$ of $M$ such that (I) each component of $\Fr_M P $ is an incompressible torus in
  $M$, (II) $P$ and $M-P$ are connected, and
(III) $\min(h(P),h(M-P))\ge m$. 
In particular, (III) implies that $P$ is a proper, non-empty
submanifold of $M$, so that $\Fr_MP\ne\emptyset$.
After an isotopy we may assume that (IV)
$\Fr_MP $ is contained in $\inter  N$ and is in general position with respect to $\fraks_\oldPsi$. Among all compact
submanifolds of $M$ satisfying (I)---(IV), let us choose $P$ so
as to minimize the quantity   $\wt( \Fr_MP)$. 

Set $K=P\cap N$, so that $P=K^+\subset N^+=M$.
Set 
$\calt=\Fr_MP=\Fr_NK$. Then it follows from (I)---(IV) that
each component of $\calt$ is an incompressible torus in
  $\inter N $,  in general position with respect to $\fraks_\oldPsi$; that $K$ and $N -K$ are connected; and that $\min(h(K),h(N-K))\ge m$.


In view of Condition (II), the hypotheses of Lemma \ref{even easier}
hold with $Y=M$ and with $\calt$ defined as above. Hence $h(M)\le
h(P)+h(\overline{M-P})-1$. Since $h(M)\ge 4m-4$, 
we have $h(P)+h(\overline{M-P})\ge4m-3$, so that at least one of the
integers $h(P)$ and $h(\overline{M-P})$ is at least $2m-1$. 
Since Conditions (I)---(IV), and the value of $\wt\calt$, are unaffected
when $P$ is replaced by $\overline{M-P}$, we may assume $P$ to have been chosen so that
$h(P)\ge2m-1$.
We will
show that with this choice of $P$, each component of $\obd(\calt)$
is incompressible in $\oldPsi$, and  one of the Alternatives (A) or (B)
of the statement holds. This will imply the conclusion of the  lemma.

This step is an application of Lemma \ref{modification}.  According to
the hypotheses of the present lemma, $\oldPsi$ is a compact, orientable,
strongly \simple\ $3$-orbifold containing no embedded negative turnovers, 
and each boundary component of $\partial N$
is a sphere. Since $M=\plusN$ is a graph manifold, it is by definition
irreducible, so that $N$ is $+$-irreducible. We have seen that $K\subset N$ is a compact,
connected $3$-manifold, that the components
of $\Fr_N K$ are all incompressible tori in $\inter N$, in general position with respect to
$\fraks_\oldPsi$, that $\overline{N-K}$ is connected,
and that $h(K)\ge2m-1$; since $m\ge2$, we have in particular that
$h(K)\ge3$. This gives Alternative (b) of the hypothesis of Lemma
\ref{modification}.

Now assume, with the aim of obtaining a contradiction, that
either some component of $\obd(\calt)$ is compressible in $\oldPsi$, or
that the components of $\obd(\calt)$ are incompressible in $\oldPsi$
but that both the Alternatives (A) or (B) of the conclusion of the present lemma are false. If
$\obd(\calt)$ has a compressible component in $\oldPsi$ then
Alternative (i) of Lemma \ref{modification} holds.
 If the components of $\obd(\calt)$ are incompressible in $\oldPsi$,
but (A) and (B) are
both false, then $\chibar(\kish(\obd(K)))<1$, and either
$\chibar(\kish(\obd(K)))<\lambda_\oldPsi/4$ or
$\wt^*\calt\ge4$. If 
$\chibar(\kish(\obd(K)))<\lambda_\oldPsi/4$ then Alternative (ii) of Lemma
\ref{modification} holds. If $\chibar(\kish(\obd(K)))<1$ and $\wt^*\calt\ge4$,
then 
then Alternative (ii$'$) of Lemma \ref{modification} holds. Thus in any
event, our assumption implies that one of the alternatives (i), (ii)
or (ii$'$)
of Lemma \ref{modification} holds, and hence that there is a compact,
connected $3$-manifold $K_1\subset N$ having the properties stated in
the conclusion of that lemma. According to the last sentence of Lemma
\ref{modification}, since Alternative (b) holds, we may choose
$K_1$ so that $\overline{N-K_1}$ is connected.

The conclusion of Lemma \ref{modification} gives 
$h(K_1)\ge h(K)/2$. Since $h(K)\ge2m-1$ it follows that $h(K_1)\ge
m$. The conclusion of Lemma \ref{modification} also gives  
$h(N-K_1)\ge h(N-K)$. Now since $K$ satisfies Condition (III)
above, we have $h(N-K)\ge m$; hence $h(N-K_1)\ge m$. 
Thus Condition (III) holds when $P$ is
replaced by $P_1:=K_1^+$. It is immediate from the conclusion of Lemma
\ref{modification} (including the connectedness of $\overline{N-K_1}$) that Conditions (I), (II) and
(IV) above also hold when $P$ is replaced by $P_1$. But since Lemma \ref{modification} also gives 
$\wt(\Fr_M P_1)=\wt(\Fr_N K_1)< \wt\calt$, this contradicts the minimality
of $\wt\calt$. 
\EndProof

\Proposition\label{crust chastened}
Let $\oldPsi$ be a compact,  orientable, strongly \simple, boundary-irreducible
  $3$-orbifold containing no embedded negative turnovers.
Set $N=|\oldPsi |$. Suppose that
$N$ is $+$-irreducible, and that each
component of $\partial N$ is a sphere, so that $\plusN$ is closed.
Then:
\begin{enumerate}
\item if $\plusN$ contains an incompressible torus and $h(N)\ge4$, we
  have $\delta(\oldPsi)\ge3\lambda_\oldPsi$ (see \ref{t-defs}); 
and
\item if $\plusN$ is a graph manifold, and if $h(N)\ge8$ and $\lambda_\oldPsi=2$, we have
  $\delta(\oldPsi)\ge12$.
\end{enumerate}
\EndProposition

\Proof
To prove Assertion (1), suppose, in addition to the general hypotheses of the proposition,  that
$h(N)\ge4$ and that $\plusN$ contains an incompressible torus. We claim:
\Claim\label{Henry F. Schricker}
There is a compact, connected submanifold $K$ of $N$, whose boundary
components are incompressible tori
in general position with respect to $
\fraks_\oldPsi$, such that 
every component of $\obd(\Fr_N K)$ is
incompressible in $\oldPsi$, and
$\chibar(\kish(\obd(K)))\ge\lambda_\oldPsi/4$.
\EndClaim

We will first prove \ref{Henry F. Schricker} in the case where
$\plusN$ is a graph
  manifold. In this case, the hypotheses of Lemma \ref{get a newer
    mop} hold with
$m=2$. Hence there is  a compact
submanifold $K$ of $N$ such that the conclusions of Lemma \ref{get a
  newer mop} hold (with $m=2$). 
In particular, $\Fr_NK$ is 
 in general position with respect to $
\fraks_\oldPsi$, and the components of $\obd(\Fr_N K)$ are incompressible in $\oldPsi$. Since one of the alternatives (A) or (B) of
Lemma \ref{get a newer mop} holds (and since $\lambda_\oldPsi\in\{1,2\}$), we have
$\chibar(\kish(\obd(K)))\ge\lambda_\oldPsi/4$, and the proof of (\ref{Henry
  F. Schricker}) is complete in this case.

To prove  (\ref{Henry
  F. Schricker}) in the case where $\plusN$ is not a graph
manifold, we note that since $\plusN$ contains an incompressible torus and is irreducible, it is a Haken manifold. We let  $\Sigma$ denote the characteristic submanifold of $\plusN$ in the sense of \cite{js} (cf. \ref{manifolds are different}), and we note that in this case, according to the definition of a graph manifold, there is a component $L$ of $\overline{\plusN-\Sigma}$ which is not homeomorpic to $T^2\times[0,1]$. Since $\partial L\subset\partial\Sigma$, each component of $\partial L$ is an incompressible torus in $\plusN$. 

If $L$ admits a Seifert fibration, it follows from the definition of the characteristic submanifold \cite[p. 138]{js} that the inclusion map $L\to\plusN$ is homotopic to a map of $L$ into $\Sigma$. Since $\partial L$ is incompressible, it follows that the identity map of $L$ is homotopic in $L$ to a map of $L$ into $\partial L$. This implies by \cite[Lemma 5.1]{Waldhausen} that $L$ is homeomorphic to $T^2\times[0,1]$, a contradiction. Hence $L$ does not admit a Seifert fibration.

If $V$ is an incompressible torus in $\plusN$, the definition of the characteristic submanifold implies that the inclusion $V\to\plusN$ is homotopic to a map of $V$ into $\Sigma$. Since $\plusN$ contains at least one incompressible torus, we have $\Sigma\ne\emptyset$.

If $V$ is an incompressible torus in $\inter L$, then since the inclusion map $j:V\to\plusN$ is homotopic to a map of $V$ into $\Sigma$, and since $\partial L$ is incompressible, it follows that $j$ is homotopic in $L$ to a map of $V$ into $\partial L$, which by $\pi_1$-injectivity can be taken to be a covering map from $V$ to a component $V'$ of $\partial L$. It then follows by \cite[Lemma 5.1]{Waldhausen} that $V$ and $V'$ are parallel in $L$. This shows that every incompressible torus in $\inter L$ is boundary-parallel in $L$.

Since $L$ does not admit a Seifert fibration, and every incompressible torus in $\inter L$ is boundary-parallel in $L$, it follows from Lemma \ref{torus goes to cylinder} that $L$ is acylindrical. 
After an isotopy we may assume that $\partial L\subset\inter
N$. Then $(L\cap N)^+=L$. Note also that $L$ is a proper submanifold of $\plusN$, since it is a component of $\overline{\plusN-\Sigma}$ and since $\Sigma\ne\emptyset$. In particular, there exists a compact, proper submanifold
$K$ of $N$ such that the components of $\Fr_NK$ are incompressible
tori, and $K^+$ is acylindrical. We
may choose such a $K$ so that $\Fr_NK$ is contained in $\inter N$ and  in general position with respect to
$\fraks_\oldPsi$, and so that, for every
submanifold $K'$ of $N$ such that $(K')^+$ is isotopic to $K^+$ 
in $N^+$, and $\Fr_NK'$ is  contained in $\inter N$ and in general position with respect to
$\fraks_\oldPsi$, we have $\wt(\Fr_NK')\ge \wt(\Fr_N K)$. Set $\calt=\Fr K$ and $\oldPi=\obd(\calt)$.
 
It now suffices to prove that
every component of $\oldPi$ is
incompressible in $\oldPsi$, and that $\chibar(\kish(\obd(K)))\ge\lambda_\oldPsi/4$.
Suppose to the contrary that some component of $\oldPi$ is
compressible in $\oldPsi$, or that all its components are incompressible but that $\chibar(\kish(\obd(K)))<\lambda_\oldPsi/4$.
This means that one of
the alternatives (i) or (ii) of the hypothesis of Lemma
\ref{modification} holds. Since $K^+$ is acylindrical, Alternative (a)
of \ref{modification} also holds. Hence there is a submanifold $K_1$
that satisfies the conclusions of Lemma \ref{modification}. In particular, $K_1^+$ is
isotopic to $K^+$ in $N^+$, and $\Fr K_1$ is 
contained in $\inter N$ and
 in general position with respect to $\fraks_\oldPsi$, and has weight strictly less than
$\wt\calt$. This contradiction to the minimality
of $\wt\calt$ completes the proof of \ref{Henry F. Schricker}.


Now since the components of $\oldPi:=\obd(\Fr_NK)$ are incompressible in
$\oldPsi$ by \ref{Henry F. Schricker}, the definitions of $\delta(\oldPsi)$ and $\sigma(\oldPsi\cut\oldPi)$ (see \ref{t-defs})
imply that 
\Equation\label{mr greengrass}
\delta(\oldPsi)\ge\sigma(\oldPsi\cut\oldPi)=
12\chibar(\kish(\oldPsi\cut\oldPi))
=
12 (
\chibar(\kish(\obd(K)))+\chibar(\kish(\obd(\overline{N-K})))).
\EndEquation
Since 
$\chibar(\kish(\obd(K)))\ge\lambda_\oldPsi/4$ by \ref{Henry F. Schricker}, and
$\chibar(\kish(\obd(\overline{N-K})))\ge0$ by Proposition \ref{less than nothing}, it follows from (\ref{mr greengrass}) that 
$\delta(\oldPsi)\ge3\lambda_\oldPsi$, and Assertion (1) is proved. 

Let us now turn to Assertion (2).
Suppose, in addition to the general hypotheses of the proposition,  that
$\plusN$ is a graph manifold, that $\lambda_\oldPsi=2$, and that $h(N)\ge8$. Then the hypotheses of Lemma \ref{get a newer mop} hold with
$m=3$. Hence there is  a compact
submanifold $K$ of $N$ such that the conclusions of Lemma \ref{get a
  newer mop} hold with $m=3$. Thus $K$ and $\overline{N-K}$ are
connected, and
$\min(h(K),h(N-K))\ge 3$. Furthermore, $\calt:=\Fr_NK$ is 
 in general position with respect to $\fraks_\oldPsi$; the components of $\oldPi:=\obd(\calt)$ are
 incompressible in $\oldPsi$; and either (A) $\chibar(\kish(\obd(K)))\ge1$, or (B)  $\chibar(\kish(\obd(K)))\ge\lambda_\oldPsi/4=1/2$ and
$\wt^*\calt<4$. 
By
the
definitions of $\delta(\oldPsi)$ and $\sigma(\oldPsi\cut\oldPi)$, the inequality (\ref{mr greengrass}) holds
in this context.
If (A) holds, then since $\chibar(\kish(\obd(\overline{N-K})))\ge0$,
the right hand side of (\ref{mr greengrass}) is bounded below by $12$; hence
$\delta(\oldPsi)\ge12$, so that the conclusion of (2) is true when (A) holds.

Now suppose that (B) holds. We claim:
\Equation\label{football heh heh}
\chibar(\kish(\obd(\overline{N-K})))\ge1/2.
\EndEquation

To prove (\ref{football heh heh}), suppose that 
$\chibar(\kish(\obd(\overline{N-K})))<1/2$. Since (B) holds we have $\wt^*\calt<4$; since in addition $\wt^*\calt$ is divisible by $\lambda_\oldPsi=2$ by \ref{wait star}, we have $\wt^*\calt\le2$. 
Set
$K^*=\overline{N-K}$. We will apply Lemma \ref{modification}, with
$K^*$ playing the role of $K$ in that lemma. We have observed that
$K^*$ is connected, and that the components of $\calt=\Fr_NK^*$ are
incompressible tori,  in general position with respect to $\fraks_\oldPsi$. Since
$K^*=\overline{N-K}$ is connected, and $h(K^*)\ge3$, Alternative (b) of
the hypothesis of Lemma \ref{modification} holds with $K^*$ playing the role of $K$ in that lemma. Since the components
of $\oldPi:=\obd(\Fr_NK^*)$ are incompressible in $\oldPsi$, and since we
have assumed that $\chibar(\kish(\obd(\overline{N-K})))<1/2=\lambda_\oldPsi/4 $,
Alternative (ii) of the hypothesis of Lemma \ref{modification}
holds. Hence there is a compact, connected $3$-manifold $K_1^*\subset N$
such that the conclusions of Lemma \ref{modification} hold
with $K_1^*$ and $K^*$ playing the roles of $K_1$ and $K$
respectively. In particular, $K_1^*$ is a proper, non-empty
submanifold of $N$, so that $\calt_1:=\Fr_NK_1^*\ne\emptyset$; every
component of $\Fr_N K_1^*$ is an incompressible toric suborbifold of $\inter N$;
and we have $\wt\calt_1<\wt\calt\le2$. Since $\lambda_\oldPsi=2$,
Lemma \ref{modification} also guarantees that $\max(\wt_\oldPsi(\calt )-
\wt_\oldPsi(\calt_1),\wt_\oldPsi(\calt_1))\ge2$. 
Since $\wt\calt_1<2$, it follows that $\wt_\oldPsi(\calt _1)\le
\wt_\oldPsi(\calt)-2\le0$,
i.e. $\calt_1\cap\fraks_\oldPsi=\emptyset$. Thus any
component of $\obd(\calt_1)=\calt_1\ne\emptyset$ is an
incompressible torus in $\oldPsi$; this is a contradiction since $\oldPsi$
is strongly \simple (see \ref{oops}). This
  completes the proof of (\ref{football heh heh}). 

It follows from (\ref{football heh heh}), together with the inequality
  $\chibar(\kish(\obd(K)))\ge1/2$ which is contained in Condition (B),
  that the right hand side of (\ref {mr greengrass}) is again at least $12$, so that
$\delta(\oldPsi)\ge12$ as required.
\EndProof

\Lemma\label{get the basin}
Let $\oldPsi$ be a compact,  orientable, strongly \simple, boundary-irreducible
$3$-orbifold containing no embedded negative turnovers. Suppose that
$\lambda_\oldPsi=2$. 
Set $N=|\oldPsi |$. Suppose that each
component of $\partial N$ is a sphere,
that $N$ is $+$-irreducible, and that
there
is an acylindrical manifold $Z\subset \plusN$ such that 
\begin{itemize}
\item $\Fr_N Z$ is a single torus which is
incompressible in $\plusN$, and 
\item$h(Z)\le h(N)-2$.
\end{itemize}
Then $\delta(\oldPsi)\ge12$. \abstractcomment{\tiny I hope I have fixed up the notation in the
  application to be correct and to fit with
  the statement.}
\EndLemma

\Proof
After an isotopy we may assume that $\partial Z\subset
\inter N$. Then $(Z\cap N)^+=Z$. In particular, there exists a submanifold
$W$ of $N$ such that $W^+$ is acylindrical (and in particular connected), $\Fr_NW$ is a single torus which is
contained in $\inter N$ and is incompressible in $\plusN$, and $h(W)\le h(N)-2$.

We
may choose such a $W$ so that $T:=\Fr_NW$ is in general position with respect to
$\fraks_\oldPsi$, and  so that, for every
submanifold $W'$ of $N$ such that $\Fr_NW$ is contained in $\inter N$ and in general position with respect to
$\fraks_\oldPsi$, and such that $(W')^+$ is isotopic to $W^+$ 
in $N^+$, we have $\wt(\Fr_NW')\ge \wt T$.

We will apply Lemma \ref{modification},
letting $W$ play the role of $K$. Since $W^+$ is acylindrical, Alternative (a) of Lemma
\ref{modification} holds. If one of the alternatives (i), (ii) or (ii$'$)  of Lemma
\ref{modification} also holds, then Lemma \ref{modification} gives a
submanifold $W_1$ of $N$ such that $W_1^+$ is isotopic to $W^+$ in  $\plusN$, and $\Fr_N W_1$ is contained in $\inter N$, is in general position with respect to $\fraks_\oldPsi$ and has weight strictly less than
$\wt T$, a contradiction to the
minimality of $\wt T$. Hence Alternatives  (i), (ii), and (ii$'$) of Lemma
\ref{modification} must all fail to hold. Since Alternative  (i)
fails to hold, $\obd(T)$ is
incompressible in $\oldPsi$. Since Alternative  (ii)
fails to hold, we have
\Equation\label{obadiah}
\chibar(\kish(\obd(W)))\ge1/2.
\EndEquation
Since Alternative  (ii$'$)
fails to hold, either
$\wt_\oldPsi^* T<4$ or 
$\chibar(\kish(\obd(W)))\ge1$. But by \ref{wait star},
$\wt^*T$ is divisible by $\lambda_\oldPsi=2$; the definition of $\wt^*T$ given in \ref{wait star} also implies that $\wt T\le\wt^*T$. Hence
\Equation\label{jeremiah}
\chibar(\kish(\obd(W)))\ge1\quad \text{or}\quad \wt T\le2.
\EndEquation
Now since $\obd(T)$ is incompressible in
$\oldPsi$ the definitions of $\delta(\oldPsi)$ and $\sigma(\oldPsi\cut{\obd(T)})$ (see \ref{t-defs})
imply that 
\Equation\label{isaiah}
\delta(\oldPsi)\ge\sigma(\oldPsi\cut{\obd(T)})=12\chibar(\kish(\oldPsi\cut{\obd(T)}))
=12
(\chibar(\kish(\obd(W)))+\chibar(\kish(\obd(\overline{N-W})))).
\EndEquation
If the first alternative of (\ref{jeremiah}) holds, i.e. if
$\chibar(\kish(\obd(W)))\ge1$, then it follows from (\ref{isaiah})
that $\delta(\oldPsi)\ge12$, which is the conclusion of the lemma. 

It
remains to consider the case in which the second alternative of
(\ref{jeremiah}) holds, i.e. $\wt T\le2$.
In this case, we will make a second application of Lemma \ref{modification}, this time
taking $K=\overline {N-W}$. According to Lemma \ref{even easier}, we
have $h(N)\le h(W)+h(N-W)-1$, so that $h(N-W)\ge h(N)-h(W)+1$. Since $h(W)\le h(N)-2$, it follows that $h(N-W)\ge 3$. Since $W$ is connected, Alternative (b) of Lemma
\ref{modification} is now seen to hold with $K=\overline{N-W}$. If Alternative (ii)  of Lemma
\ref{modification} also holds with this choice of $K$, then since $\lambda_\oldPsi=2$, Lemma
\ref{modification} gives a proper, non-empty, compact 
submanifold $K_1$ of $N$, whose frontier components are incompressible
tori in $N$, such that $T_1:=\Fr_N
K_1$ has weight less than
$ \wt T$. 
Since 
 $\lambda_\oldPsi=2$, Lemma
\ref{modification} also guarantees that 
$\max(\wt T-\wt T_1,\wt T_1)
\ge2$.
But $\wt T_1<\wt T\le2$, and hence $\wt T-\wt T_1\ge2$, i.e. $\wt T_1\le \wt T-2\le0$.
This means that
$T_1\cap
K_1=\emptyset$, and so any component of $\obd(\calt_1)=\calt_1$ is an
incompressible toric suborbifold of $\oldPsi$. This is a contradiction since $\oldPsi$
is strongly \simple (see \ref{oops}).
Hence Alternative  (ii) of Lemma
\ref{modification} must fail to hold with $K=\overline{N-W}$. As we
have already seen that $\obd(T)=\obd(\Fr_N(\overline{N-W}))$ is incompressible, this
means that 
\Equation\label{obadobah}
\chibar(\kish(\obd(\overline{N-W})))\ge1/2.
\EndEquation
It now follows from (\ref{isaiah}),  (\ref{obadiah}) and (\ref{obadobah}) that $\delta(\oldPsi)\ge12$.
\EndProof

\section{Homology of underlying
  manifolds}\label{irr-M section}

\Proposition\label{agol-plus}
If an orientable hyperbolic manifold $M$
has at least two cusps then $\vol M\ge\voct$ (see \ref{voct def}).
\EndProposition

\Proof
We may of course assume that $\vol M<\infty$. If $M$ has exactly two
cusps, the result follows from \cite[Theorem 3.6]{twocusps}. If $M$
has more than two cusps, it is a standard consequence of Thurston's
hyperbolic Dehn filling theorem \cite[Chapter E]{bp} that there is a hyperbolic manifold $M'$ having
exactly two cusps which can be obtained from $M$ by Dehn filling, and
that $\vol(M')<\vol(M)$. Since $\vol(M')\ge\voct$, we have
$\vol(M')>\voct$ in this case.
\EndProof

\Proposition\label{new get lost}
Let $\oldPsi$ be a compact, strongly \simple, boundary-irreducible,
orientable $3$-orbifold containing no embedded negative
turnovers. Set $N=|\oldPsi |$. Suppose that every
component of $\partial N$ is a sphere, that $N$ is
$+$-irreducible, and that $\smock_0(\oldPsi)\le3.44$. (See \ref{t-defs}.)
Then the following conclusions hold.
\begin{enumerate}
\item If $\lambda_\oldPsi=2$ then $h(N)\le7$.
\item If $\lambda_\oldPsi=2$ and $h(N)\ge 4$ then either $\theta(\oldPsi)>4$ or $\delta(\oldPsi)\ge 6$.
\item If $\lambda_\oldPsi=2$ and $h(N)\ge6$ then either $\theta(\oldPsi)>10$ or
  $\delta(\oldPsi)\ge6$.
\item If $h(N)\ge 4$ then either $\theta(\oldPsi)> 4$ or $\delta(\oldPsi)\ge3$.
\end{enumerate}
\EndProposition

\Proof
The hypothesis and the definition of $\smock(\oldPsi)$ (see \ref{t-defs}) give
\Equation\label{for future reference}
\smock(\oldPsi)= \frac{\smock_0(\oldPsi)}{0.305}
\le 
\frac{3.44}{0.305}
<11.
\EndEquation

Set $M=\plusN$, so that $M$ is a closed,
orientable $3$-manifold. The $+$-irreducibility of $N$
means that $M$ is irreducible. 
We have $h(M)=h(N)$.

By (\ref{more kitsch}) and the
hypothesis of the present proposition, we have
$\volorb(\oldPsi)\le
 \smock_0(\oldPsi)\le3.44$. Hence by Proposition \ref{hepcat}, we have
\Equation\label{slob cat}
\volG M\le \volorb(\oldPsi)\le 3.44.
\EndEquation

Because of the equivalence between the PL and smooth categories for $3$-manifolds, $M$ has a smooth structure compatible with its PL structure, and up to topological isotopy, the smooth $2$-submanifolds of $M$ are the same as its PL $2$-submanifolds. This will make it unnecessary to distinguish between smooth and PL tori in the following discussion.

Consider the case in which $M$ (regarded as a smooth manifold) is
hyperbolic. In this case, we have $\vol M=\volG M$ by
\cite[Theorem C.4.2]{bp}, and (\ref{slob cat}) gives
$$
\vol M\le \volorb(\oldPsi)\le 3.44. 
$$
According to  \cite[Theorem 1.7]{fourfree},
  any closed, orientable hyperbolic $3$-manifold $M$ of volume at most
  $3.44$ satisfies $h(M)\le7$. This establishes Conclusion (1) in this
  case.

According to \cite[Theorem 1.1]{rankfour},
  any closed, orientable hyperbolic $3$-manifold $M$ of volume at most
  $1.22$ satisfies $h(M)\le3$. Hence if $h(N)\ge4$, we have 
$$\vol
  M>1.22.$$
 But since $\volG M=\vol M$ by
\cite[Theorem C.4.2]{bp}, the definition of
 $\theta(\oldPsi)$ (see \ref{t-defs}) gives
$\theta(\oldPsi)=
\vol(M)/0.305$.
Hence
$\theta(\oldPsi)>
1.22/{0.305}
=4$.
This establishes conclusions (2) and (4) in this case.

According to \cite[Theorem 1.2]{lastplusone}, 
  any closed, orientable hyperbolic $3$-manifold $M$ of volume at most
  $3.08$ satisfies $h(M)\le5$. 
Hence if $h(N)\ge6$, we have $\vol
  M>3.08$, and therefore
$\theta(\oldPsi)=\vol(M)/{0.305}>
3.08/0.305>10$.
This establishes conclusion (3) in this case, and completes the proof
of the proposition in the case where $M$ is hyperbolic.

Now consider the case in which  $M$ is
not hyperbolic. 
Since $M$ is irreducible, it follows from Perelman's geometrization theorem
\cite{bbmbp}, \cite{Cao-Zhu}, \cite{kleiner-lott}, \cite{Morgan-Tian}, that there is a (possibly empty) $2$-manifold
$\calt\subset M$, each component of which is an incompressible torus,
such that for each component $C$ of $M-\calt$, either $\hatC$ is a Seifert
fibered space, or $C$ admits
a
hyperbolic metric of finite volume; in the latter case, we will say
more briefly that $C$ is hyperbolic. 

We claim:
\Claim\label{found a doughnut}
If $h(N)\ge4$ then $M$ contains at least one incompressible torus.
\EndClaim

Note that \ref{found a doughnut} is obvious if $\calt\ne\emptyset$. If
$\calt=\emptyset$, then since $M$ is not itself hyperbolic, $M$ must
be a Seifert fibered space. But since $h(N)\ge4$, it follows from Lemma \ref{and four if by zazmobile}
that $M$
contains an incompressible torus. Thus \ref{found a doughnut} is established.

If
$h(N)\ge4$, then by \ref{found a doughnut} and 
Assertion (1) of Proposition \ref{crust
  chastened}, we have
$\delta(\oldPsi)\ge3\lambda_\oldPsi$. 
 Thus in particular we have $\delta(\oldPsi)\ge3$; and if, in addition to assuming
$h(N)\ge4$, we assume $\lambda_\oldPsi=2$, then 
$\delta(\oldPsi)\ge6$. 
This establishes Conclusions (2), (3) and (4) in this case.

To prove Conclusion (1) in this case, let $\ch$ denote the set of all hyperbolic components of
$M-\calt$. (The set $\ch$ may be empty.) Since $M$ is not itself
hyperbolic, we have $\partial \hatC\ne\emptyset$ for any
$C\in\ch$. 

It follows from \cite[Theorem 1]{soma}, together with the fact (see \cite[Section 6.5]{thurstonnotes}) that the volume of a finite-volume hyperbolic $3$-manifold is equal to the relative Gromov volume of its compact core, that $\sum_{C\in\ch}\vol( C)=\volG M$. With (\ref{slob cat}), this gives
\Equation\label{last cat standing}
\sum_{C\in\ch}\vol( C)
\le 3.44.
\EndEquation
Since $\partial\hatC\ne\emptyset$ for each $C\in\ch$, each 
hyperbolic manifold in the collection $\ch$ has at least one cusp. Hence by
\cite[Theorem 1.1]{cao-m}, we have
\Equation\label{what do you want? four fewer years?}
\vol( C)\ge2\calv=2.029\ldots \quad\text{for each }C\in\ch,
\EndEquation
where $\calv$
is the volume of the ideal regular
tetrahedron in hyperbolic 3-space. It follows from
(\ref{last cat standing}) and (\ref{what do you want? four fewer
  years?}) that $\card \ch\le1$.

Let us now distinguish two subcases, according as $\ch=\emptyset$ or
$\card \ch=1$. In each of these subcases, we will show that Conclusion
(1) of the proposition holds.

In the subcase $\ch=\emptyset$, the manifold $M=\plusN$ is by definition a
graph manifold. 
To prove Conclusion (1) in this subcase, note that if
$\lambda_\oldPsi=2$ and $h(M)\ge8$, we may apply Assertion (2) of
Proposition \ref{crust chastened} to deduce that 
$\delta(\oldPsi)\ge12$. By Corollary \ref{bloody hell},  we therefore have
$\smock(\oldPsi)\ge12$,
%
%
a contradiction to (\ref{for future reference}). Thus (1) holds in this subcase.

There remains the subcase in which $\card \ch=1$. Let $C_0$ denote the
unique element of $\ch$. Then $ C_0$ is a (connected)
finite-volume, orientable hyperbolic $3$-manifold with at least one
cusp. According to (\ref{last cat standing}) we have
$\vol( C_0)
\le 3.44$. On the other hand, if $C_0$
had at least two cusps, we would have $\vol C_0\ge\voct=3.66\ldots$ by
Proposition \ref{agol-plus}. Hence $ C_0$ has exactly one cusp,
i.e. $\partial \hatC_0$ is a single torus.

According to \cite[Theorem 6.2]{hodad}, if $Y$ is a complete,
finite-volume, orientable hyperbolic 3-manifold having exactly one
cusp, and if
$h(Y)\ge6$, then $\vol(Y)>5.06$. Since $ C_0$ has exactly one cusp, and
since $\vol( C_0)\le3.44$ by (\ref{last cat standing}),
we have $h(C_0)\le5$. Thus if we assume that $\lambda_\oldPsi=2$ and $h(N)\ge7$, we have
$h(C_0)\le h(N)-2$. Since $C_0$ is acylindrical by Lemma \ref{torus
  goes to cylinder}, 
we may apply Lemma
\ref{get the basin}, with $Z=C_0$, to deduce that 
$\delta(\oldPsi)\ge12$. By Corollary \ref{bloody hell},  we therefore have
$\smock(\oldPsi)\ge12$,  a contradiction to (\ref{for future reference}).
Hence in this subcase $\lambda_\oldPsi=2$ implies $h(M)\le6$, and in particular
Conclusion (1)
holds in this 
subcase as well.

\abstractcomment{\tiny Check carefully for editing errors. In particular check that
  it's clear that the hypotheses of Lemma \ref{get the basin} hold for
that last app. I'm confusing about the role of the boundary spheres
that constitute the difference between $N$ and $M$. How am I making
the transition? I guess the answer is contained in the proofs of
Proposition \ref{crust chastened} and Lemma \ref{get the basin}, which are about
$\plusX$ (or is it $\plusN$? I don't know what I was talking about
here exactly), but depend on earlier lemmas about $X$ (or is it $N$?).}

\EndProof

The following result was mentioned in the introduction as Proposition B:

\Proposition\label{lost corollary}
Let $\oldOmega$ be a closed,
orientable, hyperbolic $3$-orbifold such that $\oldOmega\pl$ contains no embedded negative
turnovers. 
Suppose that $M:=|\oldOmega|$ is irreducible.
Then:
\begin{itemize}
\item If $\lambda_\oldOmega=2$ and $\vol(\oldOmega)\le3.44$ then $h(M)\le7$. 
\item If  $\lambda_\oldOmega=2$ and $\vol(\oldOmega)\le1.22$, then $h(M)\le3$.
\item If $\lambda_\oldOmega=2$ and  $\vol(\oldOmega)<1.83$, then $h(M)\le5$.
\item If $\vol(\oldOmega)<0.915$, then $h(M)\le3$. 
\end{itemize}
\EndProposition

\Proof Note that the hypothesis of each of the assertions implies that $\vol\oldOmega\le3.44$. According to Corollary \ref{smockollary} and the hypothesis, we then
have $\smock_0(\oldOmega\pl)= \vol\oldOmega\le3.44$. Furthermore, the hyperbolicity of $\oldOmega$ implies that $\oldOmega\pl$ is strongly \simple\ (see \ref{oops}). Thus the hypotheses of
Proposition \ref{new get lost} hold with $\oldOmega\pl$ and $M$ playing the
respective roles of $\oldPsi$ and $N$. Note that $\lambda_{\oldOmega\pl}=\lambda_\oldOmega$. According to Assertion (1)
of Proposition \ref{new get lost}, if $\lambda_\oldOmega=2$, we have $h(M)\le7$, which
is the first conclusion of the present proposition. 

To prove the remaining
conclusions, note that 
Corollary \ref{bloody hep} gives 
$\smock(\oldOmega\pl)\ge\theta(\oldOmega\pl)$, and that
Corollary \ref{bloody hell} gives
$\smock(\oldOmega\pl)\ge\delta(\oldOmega\pl)$. Hence $\smock(\oldOmega\pl)\ge\max(\theta(\oldOmega\pl),\delta(\oldOmega\pl))$.
In view of the definition of $\smock(\oldOmega\pl)$ (see \ref{t-defs}), this
means that
$\smock_0(\oldOmega\pl)/{0.305} \ge\max(\theta(\oldOmega\pl),\delta(\oldOmega\pl))$. Again using that
$\vol\oldOmega=\smock_0(\oldOmega\pl)$ by  Corollary \ref{smockollary}, we deduce:
\Equation\label{that's not an all}
\vol(\oldOmega) \ge0.305\max(\theta(\oldOmega\pl),\delta(\oldOmega\pl)).
\EndEquation

If  $\lambda_\oldPsi=2$ and and $\vol(\oldOmega)\le1.22$, it follows from (\ref{that's not an all})
that $\max(\theta(\oldOmega\pl),\delta(\oldOmega\pl))\le4$, which by
Assertion (2) of Proposition \ref{new get lost} implies that
$h(M)\le3$. This is the second assertion of the present proposition.

If  $\lambda_\oldPsi=2$ and and $\vol(\oldOmega)<1.83$, it follows from (\ref{that's not an all})
that $\max(\theta(\oldOmega\pl),\delta(\oldOmega\pl))<6$, which by
Assertion (3) of Proposition \ref{new get lost} implies that
$h(M)\le5$. This is the third assertion of the present proposition.

If $\vol(\oldOmega)<0.915$, it follows from (\ref{that's not an all})
that $\max(\theta(\oldOmega\pl),\delta(\oldOmega\pl))<3$, which by
Assertion (4) of Proposition \ref{new get lost} implies that
$h(M)\le3$. This is the fourth assertion of the present proposition.
\EndProof

\section{A bound for orbifold homology}\label{A and C}

The following result was stated in the introduction as Proposition A.

\Proposition\label{boogie-woogie bugle boy}
Let $\oldOmega$ be a 
closed, orientable, hyperbolic 3-orbifold such that $\fraks_\oldOmega$
is a link. Then $\oldOmega$ is covered with degree at most $2$ by some
orbifold $\toldOmega$ such that
$$\dim_{\ZZ_2}H_1(\oldOmega;\ZZ_2)\le1+ h(|{\toldOmega}|)+h(|{\oldOmega}|).$$
\EndProposition

\Proof
Set $M=|\oldOmega|$. Let $\eta:H_1(\oldOmega;\ZZ_2)\to H_1(M;\ZZ_2)$
denote the 
canonical surjection.
Let $T_1,\ldots,T_m$ denote the components of $\fraks_\oldOmega$. For
$i=1,\ldots,m$, let $C_i$ denote a simple closed curve bounding a disk
which meets $T_i$ transversally in one point, and is disjoint from
$T_j$ for each $j\ne i$. Then $C_i$ defines an element $c_i$ of
$H_1(\oldOmega;\ZZ_2)$, 
and the elements
$c_1,\ldots,c_m$ span the vector space $K:=\ker\eta$. Hence $k:=\dim K\le m$,
and after possibly
re-indexing the $T_i$, we may assume that $c_1,\ldots,c_k$ form a
basis for $K$. 

Now set $h=h_1(M)$. Since $\eta$ is surjective, there are elements
$e_1,\ldots,e_h$ of $H_1(\oldOmega;\ZZ_2)$ such that
$\eta(e_1),\ldots,\eta(e_h)$ form a basis of $H_1(M;\ZZ_2)$. Then
$c_1,\ldots,c_k,e_1\ldots,e_h$ form a basis of $H_1(\oldOmega;\ZZ_2)$. If
$k=0$ then $K=0$, so that $\dim(H_1(\oldOmega;\ZZ_2))=h(M)$. Hence in this
case, the conclusion is trivially true if we take $\toldOmega=\oldOmega$, a
one-sheeted covering. We may therefore assume that $k>0$.

Define a homomorphism $\phi:H_1(\oldOmega;\ZZ_2)\to\ZZ_2$ by setting
$\phi(c_i)=1$ for $i=1,\ldots,k$ and $\phi(e_j)=0$ for
$j=1,\ldots,h$. After possibly re-indexing the $c_i$ for $k<i\le m$,
we may assume that there is an integer $r$ with $k\le r\le m$ such
that $\phi(c_i)=1$ for $1\le i\le r$ and $\phi(c_i)=0$ for $r<i\le
m$. The  homomorphism $\phi:H_1(\oldOmega;\ZZ_2)\to\ZZ_2$ is non-zero
since $k>0$, and therefore defines a
two-sheeted covering of the orbifold $\oldOmega$, say $p:\toldOmega\to\oldOmega$. If we
regard $p$ as a map from $\tM:=|\toldOmega|$ to $M$, it is a branched covering map whose branch
locus is $T_1\cup\ldots\cup T_r$. Hence $\tT_i:=p^{-1}(T_i)$ is a
simple closed curve for $i=1,\ldots,r$, and the non-trivial deck transformation
is an involution $\tau:\tM\to\tM$ with $\Fix\tau=\tT_1\cup\cdots\cup\tT_r$.

According to the Smith inequality (see \cite[p. 126, Theorem 4.1]{bredon}), we
have
\Equation\label{a shitte hier, a shitte her}
\sum_{d\ge0}\dim H_d(\Fix(\tau);\ZZ_2)\le \sum_{d\ge0}\dim
H_d(\tM;\ZZ_2).
\EndEquation
(In general this inequality holds if $\tau$ is an order-$p$
homeomorphism for a given prime $p$, and homology with coefficients in
$\ZZ_p$ is used; in this application we have $p=2$.)

Since $\Fix(\tau)$ is a disjoint union of $r$ simple closed curves,
the left hand side of (\ref{a shitte hier, a shitte her}) is $2r$. By
Poincar\'e duality, the right hand side of (\ref{a shitte hier, a
  shitte her}) is $2+2h(\tM)$. Hence
$r\le1+h(\tM)$. 
Since
$c_1,\ldots,c_k,e_1\ldots,e_h$ form a basis of $H_1(\oldOmega;\ZZ_2)$, we
have
$$\dim H_1(\oldOmega;\ZZ_2)=h+k\le h+r\le h+1+h(\tM),$$
which gives the conclusion.
\EndProof

The following result was stated in the introduction as Proposition C.

\Proposition\label{orbifirst}Let $\oldOmega$ be a 
closed, orientable, hyperbolic 3-orbifold such that 
$\fraks_\oldOmega$ is a link, and such that $\pi_1(\oldOmega)$
contains no hyperbolic triangle group. Suppose that $|\oldOmega|$ is irreducible, and that $|\toldOmega|$ is irreducible for every two-sheeted (orbifold) covering $\toldOmega$ of $\oldOmega$. If $\vol \oldOmega\le1.72$ then
$\dim H_1(\oldOmega;\ZZ_2)\le
15$. Furthermore, if
$\vol \oldOmega\le1.22$ then
$\dim H_1(\oldOmega;\ZZ_2)\le
11$, and if $\vol \oldOmega\le0.61$ then
$\dim H_1(\oldOmega;\ZZ_2)\le
7$. 
\EndProposition

\Proof
According to Proposition
\ref{boogie-woogie bugle boy}, $\oldOmega$ is covered with degree at most $2$ by an
orbifold $\toldOmega$ such that 
\Equation\label{number me}
\dim_{\ZZ_2}H_1(\oldOmega;\ZZ_2)\le1+ h(|{\toldOmega}|)+h(|{\oldOmega}|).
\EndEquation
Since $\oldOmega$ has a link as its singular set,
so does $\toldOmega$; thus $\lambda_\oldOmega=\lambda_{\toldOmega}=2$. Since  $\pi_1(\oldOmega)$
contains no hyperbolic triangle group, it follows from Corollary \ref{injective hamentash}
that neither $\oldOmega\pl$ nor $\toldOmega\pl$ contains
any embedded negative turnover. The hypothesis implies that $|\oldOmega|$ and $|\toldOmega|$ are both irreducible. If $\vol\oldOmega\le1.72$, then $\vol\oldOmega$ and $\vol\toldOmega=2\vol\oldOmega$ are both bounded above by $3.44$; it then follows from Proposition \ref{lost corollary} that $h(|\oldOmega|)$ and $h(|\toldOmega|)$ are both bounded above by $7$. We therefore have $\dim_{\ZZ_2}H_1(\oldOmega;\ZZ_2)\le15$ by (\ref{number me}). 
If $\vol\oldOmega\le1.22$, then  Proposition \ref{lost corollary} gives 
$h(|\oldOmega|)\le3$; and since 
$\vol\toldOmega\le2.44<3.44$,  Proposition \ref{lost corollary} also gives 
$h(|\toldOmega|)\le7$, so that $\dim_{\ZZ_2}H_1(\oldOmega;\ZZ_2)\le11$ by (\ref{number me}). Likewise, if $\vol\oldOmega\le0.61$, then $\vol\oldOmega$ and $\vol\toldOmega=2\vol\oldOmega$ are both bounded above by $1.22$; it then follows from Proposition \ref{lost corollary} that $h(|\oldOmega|)$ and $h(|\toldOmega|)$ are both bounded above by $3$, so that $\dim_{\ZZ_2}H_1(\oldOmega;\ZZ_2)\le7$ by (\ref{number me}). 
\EndProof
\bibliographystyle{plain}
\bibliography{orbilink}

\end{document}